\documentclass[a4paper,12pt]{amsart}
\usepackage{import}
\usepackage[utf8]{inputenc}
\usepackage{datetime} 
\usepackage[margin=1in]{geometry}
\usepackage{amsmath, amssymb, amsthm, extarrows, float,mathrsfs, mathtools, amsfonts}
\usepackage{color}
\usepackage{tikz, tikz-cd}
\usepackage{quiver}
\usepackage{graphicx}
\usepackage{tensor}
\usepackage[colorinlistoftodos]{todonotes}
\usepackage{latexsym}
\usepackage{mdwlist}
\usepackage[backend=biber,
            isbn=false,
            doi=false,
            eprint=false,
            url=false,
            maxbibnames=9,
            giveninits=true,
            style=alphabetic,
            citestyle=alphabetic]{biblatex}
\AtEveryBibitem{\clearlist{language}}
\bibliography{References}
\renewbibmacro{in:}{}

\newcommand{\changelocaltocdepth}[1]{%
  \addtocontents{toc}{\protect\setcounter{tocdepth}{#1}}%
  \setcounter{tocdepth}{#1}%
}
\usepackage{color}
\definecolor{e-mail}{rgb}{0,.40,.80}
\definecolor{burgundy}{RGB}{128,0,32}
\definecolor{citation}{rgb}{0,.40,.80}

\usepackage[colorlinks=true,
            linkcolor=burgundy,
            citecolor=citation,
            urlcolor=e-mail]{hyperref}
\usepackage{preamble}

\title{Cohomological Donaldson-Thomas theory for local systems on the $3$-torus }
\author{\v{S}arūnas Kaubrys}
\address{School of Mathematics, University of Edinburgh, Edinburgh, UK}
\email{S.Kaubrys@ed.ac.uk}
\begin{document}
\begin{abstract}
    This paper studies the Cohomological Donaldson-Thomas theory of $G$-local systems on the topological three torus. Using an exponential map we prove cohomological integrality for $\GL_n$-local systems using the statement of cohomological integrality for the tripled Jordan quiver from \cite{davison2020cohomological}. Using this result we prove a version of cohomological integrality for $\SL_n$ and $\PGL_n$ for prime $n$. Finally, for prime $n$, we prove a Langlands duality statement for the $\SL_n$ and $\PGL_n$ cohomological Donaldson-Thomas invariants.
\end{abstract}
\maketitle
\tableofcontents
\newpage

\section{Introduction}
\changelocaltocdepth{1}
Let $M$ be a closed connected oriented $n$-manifold and $G$ a connected reductive group over $\mathbb{C}$. A $G$-\emph{local system} is a homomorphism $\pi_{1}(M) \to G$ from the fundamental group of $M$ to $G$. We denote the stack of $G$-local systems on $M$ by $\Loc_{G}(M)$. In this paper we will study the \emph{cohomological Donaldson-Thomas (DT) invariant}  $\HHf(\Loc_{G}(M), \varphi_{G}(M))$, where $M$ is $3$-dimensional and $\varphi_{G}(M)$ is the \emph{DT perverse sheaf} on $\Loc_{G}(M)$ as defined in \cite{ben2015darboux}. The construction of $\varphi_{G}(M)$ uses the derived enhancement $\LocB_{G}(M)$ of $\Loc_{G}(M)$ with its natural $(-1)$-shifted symplectic structure as in \cite{PTVV}. The goal of this paper is to compute $\HHf(\Loc_{G}(T^{3}), \varphi_{G}(T^{3}))$ for the topological example of local systems on the real $3$-torus and $G = \GL_n , \SL_n , \PGL_n$. \par  The  DT sheaf can be viewed as a categorification of the DT invariant originally defined for moduli spaces of sheaves on a $3$-Calabi-Yau variety $X$ by Thomas \cite{thomas_dt_og}. The $(-1)$-shifted symplectic structure on the derived moduli spaces of sheaves on $X$ comes from \emph{Serre duality} on $X$. Meanwhile the $(-1)$-shifted symplectic structure on derived moduli of local systems comes from \emph{Poincar\'e  duality} on a $3$-manifold $M$. This $(-1)$-shifted symplectic structure is a derived analogue of the symplectic structure due to Atiyah-Bott on character varieties of surfaces. \par
  
 The good moduli space of $X_{G}(M)$ is called the \emph{character variety}. When we take a surface $\Sigma$, the spaces $\Loc_{G}(\Sigma)$ and $X_{G}(\Sigma)$ have been extensively studied including in the context of nonabelian Hodge theory, geometric Langlands, quantization and low-dimensional topology.  The singular cohomology of certain twisted $G$ character varieties has been investigated. It has been shown that the cohomology of twisted $\SL_n$ and $\PGL_n$ character satisfy a type of \emph{Langlands duality} or \emph{topological mirror symmetry} \cite[Section 4]{hausel2011globaltopologyhitchin}. Another goal of this paper is to prove an analogue of the Langlands duality or topological mirror symmetry statement for cohomological DT invariants of the $3$-torus. 
\subsection{Background in DT theory.} 
We now turn to stating the results of the paper more precisely. We start by recalling the general definition of cohomological Donaldson-Thomas invariants. Let $\XB$ be any $(-1)$-shifted symplectic scheme. Brav-Bussi-Joyce \cite{brav2019darboux} have proved a Darboux theorem, which says that locally $\XB$ is a derived critical locus $\mathbf{crit}f$ of some function $f \colon U \to \mathbb{A}^{1}$ with $U$ smooth. Brav-Bussi-Ben-Bassat-Joyce have extended this theorem to derived Artin stacks \cite{ben2015darboux}. In particular, the classical truncation of a derived $(-1)$-shifted symplectic scheme $\XB$ has the structure of a \emph{d-critical locus}, so around every point it is the classical critical locus of some function. Classical truncations of $(-1)$-shifted symplectic stacks also have a similar $d$-critical structure. To capture the singularities of $\mathbf{crit}f$ we can construct the perverse sheaf of vanishing cycles $\varphi_{f}$ on the truncation of $\mathbf{crit}f$. It is then proven in \cite{ben2015darboux} that up to the existence of an orientation data on $\XB$ there is a \emph{global perverse DT sheaf} $\varphi_{\XB}$. 
\subsection{Moduli spaces of Local systems on manifolds} 
 We call the variety of homomorphisms $\Loc^{\ff}_{G}(M) = \{ \pi_{1}(M) \to  G \} $ the \emph{representation variety} of $G$-local systems on $M$. We can also view this as the moduli space of $G$-local systems with a trivialisation at a chosen point. The stack of local systems is given by $\Loc_{G}(M) \cong \Loc^{\ff}_{G}(M)/G$. The character variety of $M$, $X_{G}(M) = \{ \pi_{1}(M) \to G \} /\!\!/ G$, is given by taking the affine GIT quotient.  Recently, orientation data for all $3$-manifolds has been provided by Naef-Safronov using Reidemeister torsion in \cite{naef2023torsion}. Therefore, the DT sheaf is well defined on $\LocB_{G}(M)$ and we denote it by $\varphi_{G}(M)$. In this paper we study the DT sheaf on $\LocB_{G} (T^{3})$ for $G = \GL_n , \SL_n $ and $\PGL_n$. We have a very concrete description of the classical stack $\Loc_{G}(T^{3}) \cong \CC_{3}(G)/G$ where 
 $$\CC_{3}(G)= \{ (A,B,C) \in G^{3} \mid [A,B]=[A,C]=[B,C]=I\}$$ 
 is the scheme of $3$ pairwise commuting elements and $G$ acts by conjugation. The character variety is also concretely given by the affine GIT quotient $H^{3} /\!\!/ W \coloneqq \Spec \mathcal{O}(H^{3})^{W}$, where $W$ is the Weyl group of $G$. For $\GL_n$, the character variety is $\So^{n} \mathbb{G}^{3}_{m}$, the $n$-th symmetric power of $\mathbb{G}^{3}_{m}$.  \subsection{Tripled Jordan quiver and cohomological integrality}
Let $Q$ be a quiver. Moduli spaces of quiver representations give examples where one can compute cohomological DT invariants. A \emph{potential} $W$ is a linear combination of elements in the path algebra up to cyclic permutation. The potential $W$ defines a function $\tr(W) \colon \MM_{Q} \to  \mathbb{A}^{1}$ from the stack of quiver representations $\MM_{Q}$. We can consider the tripled quiver $\widetilde{Q}$ by doubling $Q$ and adding a loop for each vertex. There is a canonical potential $\widetilde{W}$ on $\widetilde{Q}$. Using the direct sum of quiver representations we can define a convolution symmetric monoidal structure $\boxdot$ on the good moduli space 
 of quiver representations. Summarising work of Davison-Meinhardt and Davison we will now discuss the following structural result about DT invariants of  tripled quivers called cohomological integrality.
\begin{thm} \cite[Theorem A]{davison2020cohomological}  \footnote{ The theorem is much more general and there are versions for general quivers with potential.} \label{dm_cohintquiver}
Let $Q$ be a quiver and $\widetilde{Q}$ the tripled quiver. Let $\pi \colon M_{\widetilde{Q}} \to X_{\widetilde{Q}}$ be the map to the good moduli space. Take the perverse sheaf of vanishing cycles $\varphi_{\tr \widetilde{W}}$. Consider the pushforward $\pi_{*} \varphi_{\tr \widetilde{W}}$. There exists a perverse sheaf $\BPS \cong \pH^{1} (  \pi_{*} \varphi_{\tr \widetilde{W}})$ such that we have the following isomorphism
\begin{align}
   \pi_{*} \varphi_{\tr \widetilde{W}} \cong \Sym_{\boxdot}(\BPS \otimes \HHf(\B \mathbb{G}_{m})[-1]). 
\end{align}
\end{thm}
Furthermore, using the results of \cite{bps_less_perverse} we can deduce that $\pi_{*} \varphi_{\tr \widetilde{W}}$ has the structure of a pure complex of mixed Hodge modules. Of particular interest for us is the tripled Jordan quiver $\widetilde{Q}_{\operatorname{Jor}}$ with the canonical potential $\widetilde{W} = x[y,z]$, denoting the loops of $\widetilde{Q}_{\operatorname{Jor}}$ by $x,y,z$. The critical locus of $M_{\widetilde{Q}_{\operatorname{Jor}}} \xrightarrow{\tr(\widetilde{W})} \mathbb{A}^{1}$ is then given by $\coprod_{n} \CC_{3}(\mathfrak{gl}_{n})/ \GL_{n}$, where we again have commuting elements up to conjugation. We view this space as the \emph{additive} version of the moduli of local systems on the $3$ torus. See Section \ref{coh_int_quiver} for more details.
\subsection{Exponential map and cohomological integrality} The goal of this paper is to prove an analogue of Theorem \ref{dm_cohintquiver} above for the stack of local systems on the three torus $T^3$. More precisely we will consider the disjoint union of the stacks of local systems of $\GL_n$ for all $n$. The main difficulty is the fact that $\LocB_{\GL_{n}}(T^{3})$ is \emph{not} a global critical locus as a $(-1)$-shifted symplectic stack. To work around this difficulty we consider an exponential map of associated analytic stacks
\begin{equation}
    \exp \colon \CC_{3}(\GL^{2}_{n} , \mathfrak{gl}_{n})/ \GL_{n} \to \Loc_{G}(T^{3})
\end{equation}
where the map sends $(A,B,x) \in \CC_{3}(\GL^{2}_{n}, \mathfrak{gl}_{n})$ to $(A,B, \exp(x)) \in \Loc_{\GL_n}(T^{3})$.  Note that here $\CC_{3}(\GL^{2}_{n} , \mathfrak{gl}_{n})/ \GL_{n}$ \emph{is} a global critical locus. 
This leads to the first main theorem.
\begin{thm}[= Theorem \ref{exp_dcrit}]
The map $\exp \colon \CC_{3} (\GL^{2}_{n}, \mathfrak{gl}^{\et}_{n})/ \GL_{n} \to \Loc_{\GL_{n}} (T^{3})$ is a map of $d$-critical loci, where we have restricted to an open locus where the exponential map is \'etale.
\end{thm}
We prove this theorem by first showing that the map on formal completions at a point preserves closed $2$-forms and then lifting the result to complex analytic maps. \par Using this theorem we can prove purity as a complex of mixed Hodge modules on $\So^{n} \mathbb{G}^{3}_{m}$ of the pushforward $\pi_{*} \varphi_{\LocB_{\GL_{n}}(T^{3})}$. Once we have purity, Saito's theory of mixed Hodge modules then gives a splitting into a direct sum of intersection cohomology sheaves.  We then have the theorem
\begin{thm}[Cohomological integrality for $\GL_n$ = Theorem \ref{sym_coh_3tor}] \label{sym_coh_3tor_intro}
Take the coproduct over all dimensions of the maps to the good moduli space 
\begin{equation*}
    \pi_{m} \colon \coprod_{n \geq 0} \Loc_{\GL_{n}}(T^{3}) \to \coprod_{n \geq 0 } \So^{n} \mathbb{G}^{3}_{m}
\end{equation*}
and $\varphi_{m} = \bigoplus_{n \geq 0}\varphi_{\GL_{n}}(T^{3})$ the DT sheaf on $\coprod_{n} \Loc_{\GL_{n}}(T^{3})$. We have an equivalence in $\operatorname{D}^{+}_{c}(\So \mathbb{G}^{3}_{m})$ 
\begin{equation} \label{intro_coh_eq}
    \JHmp \Jm \cong \Sym_{\boxdot} (\BPS_{m} \otimes \HHf(\B \mathbb{G}_{m})[-1])
\end{equation}
with 
\begin{equation*}
    \BPS_{m} \coloneqq \pH^{1} \JHmp \Jm  \cong \bigoplus_{n} \BPS_{\GL_n}
\end{equation*}
and $\BPS_{\GL_n} \cong (\Delta \colon \mathbb{G}^{3}_{m} \to \So^{n} \mathbb{G}^{3}_{m})_{*}  \mathbb{Q}_{\mathbb{G}^{3}_{m}}[3]$.
\end{thm} 
Using this Theorem we can compute the cohomological DT invariants using the isomorphism deduced from equation \eqref{intro_coh_eq}
\begin{equation}
    \bigoplus_{n} \HHf(\Loc_{\GL_n}(T^{3}), \varphi_{\GL_n}(T^{3})) \cong \Sym( \bigoplus_{n} \BPSo_{\GL_n} \otimes \HHf(\B \mathbb{G}_{m})[-1]),
\end{equation}
denoting the cohomology of $\BPS_{\GL_n}$ by $\BPSo_{\GL_n}$. We now give some details about the proof of Theorem \ref{sym_coh_3tor_intro}. The proof is by an explicit computation of the local systems appearing in the Saito decomposition of $\pi_{\GL_n *} \varphi_{\GL_n}$. The intersection complex perverse sheaves that appear in the decomposition are associated to local systems with respect to a stratification $\So^{n}_{\lambda} \mathbb{G}^{3}_{m}$ by partitions $\lambda$ of $n$ of the character variety $\So^{n} \mathbb{G}^{3}_{m}$. $\So^{n}_{\lambda} \mathbb{G}^{3}_{m}$ consists of elements of $\So^{n} \mathbb{G}^{3}_{m}$ that are allowed to repeat according to the partition $\lambda$. For each partition $\lambda$ we have a standard Levi 
 of $\GL_n$: $L_{\lambda} = \prod \GL_{\lambda_{i}}$. We compute these local systems by a reduction to the Levis along the natural map $\LocB_{L_{\lambda}}(T^{3}) \to \LocB_{\GL_n}(T^{3}) $ induced by the inclusion $L_{\lambda} \to \GL_n$. We consider the induced map on good moduli spaces $X_{L_{\lambda}}(T^3) \to X_{\GL_n}(T^3)$. Taking the preimage of $\So^{n}_{\lambda}(\mathbb{G}^{3}_{m})$ we obtain a cover by the relative Weyl group $W_{L_{\lambda}}$. Using this cover we can compute the local systems appearing in the Saito decomposition using an easier description on $X_{L_{\lambda}}(T^3)$. See Section \ref{section_coh_int} for more precise details.
 \subsection{Cohomological integrality for $\SL_n$ and $\PGL_n$} For a general reductive group $G$ it is not clear how to define a symmetric monoidal structure on the character variety $X_{G}$. Therefore, it is not clear how to define a symmetric algebra in a similar way as Theorem \ref{sym_coh_3tor_intro}. However, we can view cohomological integrality as a splitting of the pushforward of the DT sheaf into contributions from the Levi subgroups of $\GL_{n}$. This version is easier to generalise to more general reductive groups. In this paper we prove such a version of cohomological integrality for $\SL_n$ and $\PGL_n$. A similar formula was first explained to the author by Tasuki Kinjo.
\begin{thm}[Cohomological Integrality for $\SL_n$, and $\PGL_n$ for prime $n$  = Theorem \ref{coh_proof_pgsln}] \label{conj_intro}
Pick representatives $L$ in each conjugacy class of the Levi subgroups of $G= \SL_n$ and $\PGL_n$. Denote by $\LocB^{1}_{G}(T^{3})$ the connected component of the trivial local system in $\LocB_{G}(T^{3})$, $\pi_{G} \colon \LocB^{1}_{G}(T^{3}) \to X^{1}_{G}(T^{3})$ the good moduli space and $\varphi_{G}$ the restriction of the DT sheaf to $\LocB^{1}_{G}(T^{3})$. Define the map $\theta \colon X^{1}_{L}(T^{3}) \to X^{1}_{G}(T^{3})$ induced from the inclusion $L \to G$. Then we have
    \begin{equation} \label{intro_cohint_eq}
        \pi_{G,*} \varphi_{G}(T^{3}) \cong \bigoplus_{L \subseteq G} (\theta_{*} \BPS_{L} 
 \otimes \HHf(\B \ZZ(L))[- \dim \ZZ(L)])^{W_{L}}
    \end{equation}
Here $\BPS_{L} = \Delta_{*} \mathbb{Q}_{\ZZ(L)^{3}} [ 3\dim \ZZ(L) ]$ with $\Delta \colon \ZZ^{3}(L) \to X^{1}_{L}$ induced by the inclusion $\ZZ(L) \to L$ and $W_{L}$ is the relative Weyl group of $L$.
\end{thm}
Note that the Theorem does not make claims about contributions coming from the \emph{twisted} stacks of local systems, which are components in $\LocB_{G}$ that arise from non-trivial torsion elements in $\pi_{1}(G)$. For $\SL_n$ there is only the component of the trivial local system and the theorem is proven by deducing purity through the $\GL_n$ version and then again using a reduction to the Levis procedure. The $\PGL_{n}$ version follows by similar methods.  Note that the formula in Theorem \ref{conj_intro} makes sense for any connected reductive group so we can ask if such a version of cohomological integrality holds in other types.
\subsection{Geometric Langlands for 3 manifolds}
In \cite{witten_kapustin} Kapustin-Witten define  4D topological field theories $\ZZ^{\Psi}_{G}$ depending on a reductive group $G$ and parameter $\Psi$. S-duality implies an equivalence $\ZZ^{\Psi}_{G} \cong \ZZ^{\Psi^{L}}_{G^{L}}$, where $G^{L}$ is the Langlands dual of $G$ and $\Psi^{L}$ is the dual parameter to $\Psi$. The theories $\ZZ^{\Psi}_{G}$ and $\ZZ^{\Psi^{L}}_{G^{L}}$ assign to a surface the categories appearing in the Geometric Langlands Conjecture as defined in \cite{arinkin_gaits} and \cite{betti_gl}. It is expected\footnote{This was explained to the author by Pavel Safronov and is joint work of Ben-Zvi-Gunningham-Jordan-Safronov.} that for generic $\Psi$ the vector space assigned to a $3$-manifold is modelled by $\HHf(\Loc_{G}(M), \varphi_{G}(M))$. This inspires the following conjecture
\begin{conj}\label{langlands_conj}
    Let $G$ be a reductive group, $G^{L}$ its Langlands dual and $M$ a closed oriented $3$-manifold. We have an isomorphism
    $$ \HHf(\Loc_{G}(M),\varphi_{G}(M)) \cong \HHf(\Loc_{G^{L}}(M), \varphi_{G^{L}}(M)).$$
\end{conj}
Using cohomological integrality for $\SL_n$ and $\PGL_n$ and a computation of the twisted components of $\LocB_{\PGL_{n}}(T^{3})$ we get the following theorem.
\begin{thm}[Langlands duality for prime $n$ = Corollary \ref{langlands} ]\label{langlands_intro}
Let $n$ be prime. We have an isomorphism of graded vector spaces 
\begin{equation}
    \HHf(\Loc_{\SL_{n}}(T^{3}) , \varphi_{\SL_{n}}(T^{3})) \cong \HHf(\Loc_{\PGL_{n}}(T^{3}), \varphi_{\PGL_{n}}(T^{3})).
\end{equation}
\end{thm}
Let us now sketch the Langlands duality in the special case of $\SL_2$ and $\PGL_2$. First,  the statement of Theorem \ref{langlands_intro} implies that the BPS sheaves $\BPS_{G}$ for $G = \SL_n$ or $\PGL_n$ are well-defined and are constant sheaves of rank $1$ supported on the centre $\ZZ^{3}(G)$. Denote the cohomology of BPS sheaves $\BPS_{G}$ by $\BPSo_{G}$. In particular, $\BPS_{\SL_2}$ with cohomology 
$$\BPSo_{\SL_2} \cong \mathbb{Q}^{|\mu^{3}_{2}|}$$ is a skyscraper supported on the $8$ disjoint points in $X_{\SL_{2}}(T^{3})$ corresponding to $\ZZ^{3}(\SL_2) \cong (\mu_{2})^{3}$. Meanwhile $\BPS_{\PGL_2}$ is a constant sheaf supported on the trivial local system in $X_{\PGL_2}(T^{3})$ since the centre of $\PGL_2$ is trivial. Therefore it has cohomology
\begin{equation}
    \BPSo_{\PGL_2} \cong \mathbb{Q}.
\end{equation}
Looking at equation \eqref{intro_cohint_eq} we see that for both $\SL_2$ and $\PGL_2$ there will be two contributions to the direct sum. The first contribution comes the maximal torus, while the second comes from the trivial Levi and is the contribution of the BPS sheaves.
For $\PGL_2$ to compute the full cohomology of the DT sheaf, we will explicitly compute all the components of $\LocB_{\PGL_2}(T^{3})$. These turn out to be points with finite stabiliser, one for each non-trivial element of $\HH^{2}(T^{3}, \pi_{1}(\PGL_2)) \cong (\mu_{2})^{3}$. Each component contributes $\mathbb{Q}$ to the total cohomology.  Therefore we get
\begin{equation}
    \BPSo_{\SL_2} \cong \BPSo_{\PGL_2}  \oplus \mathbb{Q}^{|\mu^{3}_{2}|-1},
\end{equation}
which comes from the Langlands duality isomorphism $\ZZ(\SL_{2}) \cong (\pi_{1}(\PGL_2))^*$, where $*$ is the Pontryagin dual. Recall that Langlands duality interchanges maximal tori, so $H_{\SL_2}$ and $H_{\PGL_2}$ are isomorphic and hence the torus contribution is identified. We can summarise Theorem \ref{langlands_intro} for $n=2$ via the following diagram


\begin{center}
\begin{tikzcd}[column sep=tiny]
	{\HHf(\Loc_{\SL_2}(T^{3}), \varphi_{\SL_2}(T^{3})) } & \cong & {\BPSo_{\SL_2}  } & \oplus & {H_{\SL_2}  \text{ contribution}} \\
	{\HHf(\Loc_{\PGL_2}(T^{3}), \varphi_{\PGL_2}(T^{3}))} & \cong & {\BPSo_{\PGL_2}  \oplus \mathbb{Q}^{|\mu^{3}_{2}|-1} } & \oplus & {H_{\PGL_2}  \text{ contribution}}
	\arrow["{ \cong}", tail reversed, from=1-1, to=2-1]
	\arrow["\cong", tail reversed, from=1-3, to=2-3]
	\arrow["\cong", tail reversed, from=1-5, to=2-5]
\end{tikzcd}
\end{center}
The restriction to only prime dimensions in Theorem \ref{langlands_intro} is necessary to be able to compute the contributions of non-trivial components of  $\LocB_{\PGL_{n}}(T^{3})$, which are related to \emph{twisted character stacks} $\LocB^{\operatorname{tw}}_{\SL_{n}}(T^{3})$. When $n$ is not prime the contributions of the non-trivial components of $\LocB_{\PGL_n}(T^{3})$ are more complicated and the author plans to explore this question in future work.
\subsection{Relation to complexified instanton Floer homology and skein modules}
In \cite{floer} Abouzaid and Manolescu have defined an $\SL_{2}(\mathbb{C})$ version of instanton Floer homology for a closed, connected, oriented $3$-manifold $M$. Recall that we can express $\Loc_{G}(M) \cong \Loc^{\ff}_{G}(M) / G$ with $\Loc^{\ff}_{G}(M)$ the representation variety of $G$-local systems. Define $\varphi^{\#}_{G}(M)$ to be the pullback of $\varphi_{G}(M)$ to $\Loc^{\ff}_{G}(M)$. The \emph{framed complexified instanton Floer homology} is defined as the cohomology of $\varphi^{\#}_{\SL_2}(M)$. There is an action of $G$ on $\HHf(\Loc^{\ff}(M), \varphi^{\#}_{G}(M))$ and the $G$-equivariant cohomology $\HHf_{G}(\Loc^{\ff}_{G}(M), \varphi^{\#}_{G}(M))$\footnote{ More precisely we view $\HHf(\Loc^{\ff}_{G}(M), \varphi^{\#}_{G}(M))$ as a constructible sheaf on $\B G$. Constructible sheaves on $\B G$ can be identified with modules over the chains $\CC_{*}(G)$ on $G$. The equivariant cohomology is then given by taking invariants.} is identified with $\HHf(\Loc_{G}(M), \varphi_{G}(M))$.

\begin{conj} \cite[Conjecture D]{gunningham2023deformation} \label{pavel_sam}
    We have an isomorphism
    \begin{equation}
        \HH^{0}(\Loc^{\ff}_{G}(M),\varphi^{\#}_{G} (M)) \otimes_{\mathbb{C}} \mathbb{C}(q^{1/d}) = \operatorname{Sk}^{\operatorname{gen}}_{G}(M)
    \end{equation}
    here $\operatorname{Sk}^{\operatorname{gen}}_G(M)$ is the skein module  with generic quantum parameters associated to the $3$-manifold.
\end{conj}
Therefore, the skein module is a certain part of the full cohomology of the DT sheaf. From this point of view Conjecture \ref{langlands_conj} is a derived version of the following Conjecture.
\begin{conj}\cite[Conjecture 1.1]{jordan2023langlands} \label{jordan_langlands}
    Let $G$ be a reductive group, $G^{L}$ its Langlands dual and $M$ a closed, connected, oriented $3$-manifold. We have an isomorphism
    $$\operatorname{Sk}^{\operatorname{gen}}_{G}(M) \cong \operatorname{Sk}^{\operatorname{gen}}_{G^{L}}(M).$$
\end{conj}
We note that even assuming Conjecture \ref{pavel_sam} we cannot immediately deduce Conjecture \ref{jordan_langlands} for the $3$-torus from Theorem \ref{langlands_intro} since we can only compute the $G$-equivariant cohomology. We leave the computation of $\HHf(\Loc^{\ff}_{G}(M), \varphi^{\#}_{G}(M))$ for future work. 
\subsection{Relation to the work of Kinjo-Park-Safronov  }
While this paper was being completed a relative Cohomological Hall algebra structure on  $\pi_{m,*}  \varphi_{m}$ was defined in \cite{kinjo2024cohomological}. This provides a different strategy to prove Theorem \ref{sym_coh_3tor_intro} by constructing directly a map 
\begin{equation} \label{kps_eq}
  \Sym_{\boxdot}(\BPS_{m} \otimes \HHf(\B \mathbb{G}_{m})[-1]) \to \pi_{m,*}  \varphi_{m}  
\end{equation}
 using the Cohomological Hall algebra structure and a natural embedding 
 \begin{equation*}
     \BPS_{m} \otimes \HHf(\B \mathbb{G}_{m})[-1] \xhookrightarrow{} \pi_{m,*}  \varphi_{m}.
 \end{equation*}
 One can then pullback the map in equation \eqref{kps_eq} via the exponential and deduce that it is an isomorphism using the Cohomological integrality theorem for the tripled Jordan quiver. In this paper we do not require the existence of the Cohomological Hall algebra structure for any of the proofs. In future work the author will consider applications of the cohomological Hall algebra structure to DT invariants of $3$-manifolds. 
\subsection*{Acknowledgements}
I would like to thank my advisers Ben Davison and Pavel Safronov for numerous discussions and generous help. I would also like to thank Lucien Hennecart, Shivang Jindal, Danil Ko\v{z}evnikov, Karim R\'ega, Sebastian Schlegel Mejia and Nikola Tomi\'c for useful discussions. Special thanks are due to Damien Calaque and Tasuki Kinjo for independently suggesting to work with the exponential map. The author was supported by the Carnegie Trust for the Universities of Scotland for the duration of this research. 
\changelocaltocdepth{2}
\section{Preliminaries and shifted symplectic structures}
\subsection{Preliminaries in algebraic geometry}
In this paper we work with derived prestacks $\operatorname{dPreStk}$ over $\mathbb{C}$. Namely functors $\operatorname{cdga}^{\leq 0,\op} \to \operatorname{Spc}$. Here $\operatorname{Spc}$ is the $\infty$-category of spaces. There is a truncation functor $\tt_{0}\colon \operatorname{dPreStk} \to \operatorname{PreStk}$ to classical (higher) prestacks, which has a fully faithful right adjoint $\iota\colon \operatorname{PreStk} \to \operatorname{dPreStk}$. We can further truncate to $1$-prestacks $\operatorname{PreStk}^{\leq 1}$, where $\operatorname{PreStk}^{\leq 1}$ are functors $\operatorname{CAlg}^{\op}_{\mathbb{C}} \to \operatorname{Gpd}$. We will use the $\infty$-category of quasicoherent sheaves $\QCoh(\XB) = \lim_{\Spec R \to \XB} (\Mod R)$ and the subcategory of perfect complexes $\Perf (\XB)$. There is an internal Hom in $\operatorname{dPreStk}$ denoted by $\Map(\XB, \YB) \in \operatorname{dPreStk}$. $\Map(\XB, \YB)$ is defined by sending $R$ to the mapping space $\Hom(\XB \times \Spec R, \YB) \in \operatorname{Spc}$.
\begin{notation}
We will denote derived prestacks by bold letters such as $\XB$ and their truncation $\tt_{0} \XB \in \operatorname{PreStk}^{\leq 1} $ by unbolded letters $X$.
\end{notation}
We denote the $\infty$-category of derived stacks for the \'etale topology by $\operatorname{dStk}$ a derived stack $\XB$ is an Artin stack if it is $n$-geometric for some $n$ as in \cite[Appendix B7]{calaque2022aksz} and locally of finite type. \par
We will also need to work with complex analytic stacks. We define complex analytic stacks $\operatorname{Stk}^{\an}$ as in \cite[Definition 3.1.1]{sun_analytic}. In particular, $X$ is a complex analytic stack if it is a stack over the site of complex analytic spaces with analytic topology, there is a smooth surjective map $U \to X$ from a complex analytic space $U$ and the diagonal of $X$ satisfies a representability and finiteness condition. There is an analytification functor from finite type Artin stacks $(-)_{\an}\colon \operatorname{ArtStk}^{\leq 1} \to \operatorname{Stk}^{\an}$, see \cite[Section 3.2.2]{sun_analytic}. 
\begin{ex}
In this paper the main example of stacks we will use is the following. Let $G$ be an algebraic group acting on a finite type scheme $X$. Then $([X/G])_{\an} = [X_{\an} / G_{\an}] $. Here $[X_{\an} / G_{\an}]$ is the quotient of the groupoid $G_{\an} \times X_{\an} \rightrightarrows X_{\an}$.
\end{ex}
\begin{remark}
    There is also a theory of derived complex analytic stacks but we will not use this notion.
\end{remark}
\begin{notation}
If it is clear from context we will abuse notation and denote a stack $X$ and its analytification $X_{\an}$ by the same symbol.
\end{notation}
\subsection{Perverse sheaves and mixed Hodge modules} \label{mhm_background}
 In this paper we work with sheaves of vector spaces with coefficients in $\mathbb{Q}$, so we drop the coefficients from the notation. We recall some definitions and theorems that we use in this paper. Let $X$ be a finite type $\mathbb{C}$-scheme or complex analytic space, when considering (constructible) sheaves on these spaces we are implicitly considering them over $X_{\red}$. We can take the derived category of constructible sheaves $\DD^{b}_{c}(X)$, with the perverse $t$-structure and heart $\Perv(X)$. We will also use the bounded below constructible category $\DD^{+}_{c}(X)$ and pushforward and pullback functors. \par  We define the perverse \textbf{nearby cycles} and \textbf{vanishing cycles} functors 
 $${}^{p} \psi_{f}, {}^{p}\varphi_{f}\colon \DD^{b}_{c}(U) \to \DD^{b}_{c}(f^{-1}(0)) $$ for any function $f \colon U \to \mathbb{C}$. These functors restrict to functors $\Perv(U) \to \Perv(f^{-1}(0))$. We will always drop the $p$ from the notation and always assume we are taking the perverse vanishing cycles. \par  The category $\Perv(X)$ is artinian with simple objects given by $\textbf{intersection cohomology}$ complexes $\IC_{Y}(\mathcal{L})= j_{!*} (\mathcal{L}[\dim Y])$ with $j \colon Y \to X$ a smooth locally closed subvariety of $X$ and $\mathcal{L}$ a local system on $Y$. We will use the following characterisation of $\IC$ sheaves.
 \begin{lem}\cite[Lemma 3.3.3]{achar2021perverse} \label{ic_achar}
     Let $j \colon X \to Y$ be a locally closed embedding. $\IC_{Y}(\mathcal{L})$ can be characterised as the unique perverse sheaf on $Y$ such that
     \begin{enumerate}
         \item $\IC_{Y}(\mathcal{L})$ is supported on $\overline{Y}$
         \item $j^{*} \IC_{Y}(\mathcal{L}) \cong \mathcal{L}$
         \item $\IC_{Y}(\mathcal{L})$ has no subobjects or quotients supported on $\overline{Y} - Y$.
     \end{enumerate}
 \end{lem}
 We can also define the derived category of mixed Hodge modules $\DD^{b}(\MHM (X))$  on a reduced separated scheme $X$ with heart $\MHM(X)$. There is a forgetful functor 
 $$\operatorname{rat} \colon \MHM(X) \to \Perv(X),$$ which induces a functor $$\operatorname{rat}\colon \DD^{b}(\MHM(X)) \to \DD^{b}(\Perv(X)) \cong \DD^{b}_{c}(X).$$ Each $\mathcal{F} \in \MHM(X)$ carries an ascending weight filtration $W_{\bullet}$. $\mathcal{F}$ is pure of weight $n$ if $W_{k} \mathcal{F} = \mathcal{F}$ for $k \geq n$ and $W_{k} \mathcal{F} = 0$ for $k < n$. An object $\mathcal{F} \in \DD^{b}(\MHM (X))$ is pure if each $\mathcal{H}^{i}(\mathcal{F})$ is pure of weight $i$. \par We can also define the $\textbf{intersection cohomology}$ mixed Hodge module given a pure polarisable variation of Hodge structure $\mathcal{L}$ on $X$. \par We will need one final upgrade to monodromic mixed Hodge modules $\MMHM(X)$. This can be viewed as a certain Serre quotient of the category of $\MHM(X \times \mathbb{A}^{1})$.
\begin{defn}[Tate twists]
Define the pure mixed Hodge structure of weight $0$ and  cohomological degree $2$ as $\mathbb{L} =\HHf_{c}(\mathbb{A}^{1}) = \mathbb{Q}[-2](1)$. \\
In general we cannot construct a square root $\mathbb{L}^{1/2}$ in the category of mixed Hodge modules. However, $\mathbb{L}^{1/2}$ exists in $\MMHM(pt)$.
\end{defn}
\begin{thm}[ \cite{saito}]
The category of pure mixed Hodge modules is semisimple. Let $\mathcal{F} \in \DD^{b}(\MHM X)$ be pure. Then 
\begin{equation}
    \mathcal{F} = \bigoplus_{i \in \mathbb{Z}} \mathcal{H}^{i}(\mathcal{F})[-i]
\end{equation}
and
\begin{equation}
    \mathcal{H}^{i} \mathcal{F} = \bigoplus_{i \in J } \IC_{Z_{i}}(\mathcal{L}_{i})
\end{equation}
here $Z_{i}$ are locally closed smooth subvarieties and $\mathcal{L}_{i}$ are pure weight $i$ simple variations of Hodge structure.
\end{thm}
We will also need to use unbounded categories of mixed Hodge modules. We will use this in Section \ref{coh_int_quiver}.  We define $\DD^{+}(\MHM (X))$ as in \cite[Section 2.1.4]{bps_less_perverse}. We will also need perverse sheaves and mixed Hodge modules on stacks. All the stacks in this paper are global quotient stacks so we can work with the $G$-equivariant versions of the constructions previously mentioned. See \cite[Section 2]{davison2023purity} for more on mixed Hodge modules. In particular, let $M$ be a stack with a good moduli space $\pi\colon M \to X$. Then the pushforward $\pi_{*} \mathcal{F}$ of $\mathcal{F} \in \DD^{b}(\MHM (M))$ exists as an object in $ \DD^{+}(\MHM (X))$. 
\subsection{Differential forms and shifted symplectic structures}
In this section we will not explain the whole theory of $n$-shifted closed forms. We will only recall what we will use and refer the reader to \cite{PTVV} for more. Define for a connective cdga $R$ the de Rham complex
$\mathbf{DR}(R) = \Sym_{R}(\mathbb{L}_{R}[-1])$. Recall that a graded mixed complex is a complex of $\mathbb{C}$ vector spaces $(E,d)$ equipped with a decomposition $E = \bigoplus_{p \in \mathbb{Z}} E(p)$, where $E(p)$ is called the weight $p$ part. The differential $d$ preserves the weight grading and there is another differential $\epsilon$, which increases the weight grading and the cohomological grading. We denote the category of graded complexes by $\operatorname{dg}^{\gr}$ and graded mixed complexes by $\epsilon-\operatorname{dg}^{\gr}$ . There is a functor $(-)^{\#}\colon \epsilon-\operatorname{dg}^{\gr} \to \operatorname{dg}^{\gr}$ given by forgetting the $\epsilon$ differential. The de Rham algebra $\mathbf{DR}(R)$ is a graded mixed complex where $\epsilon$ action is given by the de Rham differential. This defines an $\infty$-functor $\mathbf{DR}(-)\colon \operatorname{cdga}^{\leq 0 } \to \epsilon-\operatorname{dg}^{\gr}$ that satisfies \'etale descent. For a general derived stack $\XB$, the graded mixed cdga $\mathbf{DR}(\XB)$ is then defined by right Kan extension along the map $\operatorname{cdga}^{\leq 0 } \to \operatorname{dStk}^{\op}$.  Concretely we have $\mathbf{DR}(\XB) = \lim_{\Spec R \to \XB} \mathbf{DR}(R)$. We can define as in \cite[Section 1.2]{PTVV} the functors of spaces of $n$-shifted $p$-forms $\mathcal{A}^{p}(-,n)$ and $n$-shifted closed $p$-forms $\mathcal{A}^{p,\cl}(-,n)$.
These functors satisfy \'etale descent hence one can define the same spaces for a derived stack $\XB$. Let $\XB$ be a derived stack that admits a cotangent complex, then there is a canonical map $\Gamma(\XB, \Sym(\mathbb{L}_{\XB}[-1])) \to \mathbf{DR}(\XB)$.
\begin{thm}\cite[Theorem 2.6]{cal_saf} \label{forms_artin}
Let $\XB$ be a derived  prestack that admits a cotangent complex. Then the above map is an equivalence of graded cdgas. 
\end{thm}
For Artin stacks locally of finite type, this theorem was already proven in \cite[Proposition 1.14]{PTVV}. Apart from Section \ref{s0_formal_sect}, we will work with Artin stacks locally of finite type.
Let us briefly recall a more concrete description of closed forms when $\XB$ is derived Artin using this theorem. Denote by $A^{p}(\XB)$ the complex of $p$-forms, which are all the elements of weight $p$ in $\mathbf{DR}(\XB)$. Define $n$-shifted $p$-forms by $ \HH^{n}(A^{p}(\XB))$.
Construct the complex 
$$A^{p,\cl}(\XB) = \prod_{i \geq 0} A^{p}(X)[-i].$$ Then define closed $n$-shifted $p$ forms to be $ \HH^{n}(A^{p,\cl}(X))$.  In particular, if $p = 0$ then $\HHf(A^{0,\cl}(\XB)) = \HHf_{\dR}(\XB)$. 
More explicitly, an $n$-shifted closed $p$-form is the data of a power series $\omega = \sum_{i \geq p} \omega_{i}$ with $\omega_i$ a $p+i$ form of degree $n+p$ that is closed under the total differential $d + d_{\dR}$.
Written explicitly in increasing weight we have 
\begin{align*}
    d \omega_{0} &= 0 \\
    d_{\dR} \omega_{0} + d \omega_{1} & = 0 \\
    d_{\dR} \omega_{i} + d \omega_{i+1} & = 0 
\end{align*}
The forgetful and de Rham differential maps can be described as 
\begin{align*}
    \pi\colon A^{p,\cl}(\XB) & \to A^{p}(\XB) \\
    (\omega_{0} , \omega_{1}, \omega_{2} , \cdots ) & \mapsto \omega_{0}  \\
     d_{\dR}\colon A^{p,\cl}(\XB) & \to A^{p+1, \cl}(\XB) \\
      \omega_{0} & \mapsto (d_{\dR}\omega_{0} , 0 , 0 , \cdots , )
\end{align*}
Two closed $p$ forms $\omega = \sum_{i \geq p } \omega_{i}$ and $\omega^{'} = \sum_{i \geq p } \omega^{'}_{i}$ are equivalent $\omega \sim \omega^{'}$ if there exists a homotopy $h = \sum h_{i}$, with $h_i$ $p+i$-forms  such that
\begin{equation}
    \omega_{i} - \omega^{'}_{i} = d h_{i} + d_{\dR} h_{i-1}. 
\end{equation}
If two classes $\omega$ and $\omega^{'}$ are homotopic then they define the same cohomology class in $\HH^{n}(A^{p,\cl}(\XB))$.
\begin{remark}
    We have that $\pi_{0} \mathcal{A}^{p,\cl}(\XB,n) \cong \HH^{n}(A^{p,\cl}(\XB)$. By abuse of notation we write $\omega \in \mathcal{A}^{p,\cl}(\XB, n)$ for $\omega \in \pi_{0} \mathcal{A}^{p,\cl}(\XB,n)$
\end{remark}
In classical algebraic geometry, a symplectic structure is a nondegenerate closed $2$-form $\omega$ on a scheme or manifold $X$. We can express the nondegeneracy condition by saying that the form gives an isomorphism $\To X \cong \To^{*}X$ between the tangent and cotangent bundles. This version of symplectic structure can be readily generalised as follows.
\begin{defn}[Shifted symplectic structure]\label{shifted_symp_def}
An $n$-shifted symplectic structure on a derived Artin stack $\XB$ is a closed $2$-form $\omega \in \mathcal{A}^{2,\cl}(\XB,n) $ along with the \textbf{non-degeneracy} condition that the induced map 
$$\mathbb{T}_{\XB} \to \mathbb{L}_{\XB}[n]$$
is a quasi-isomorphism.
\end{defn}
\begin{defn}[Symplectic maps]
Let $f\colon \XB \to \YB$ be a map of derived stacks with $n$-shifted symplectic structures $\omega_{\XB}$ and $\omega_{\YB}$. Then we say that the map is $n$-symplectic or simply symplectic if $f^{*} \omega_{\YB} \sim \omega_{\XB}$.
\end{defn}
\begin{defn}[Lagrangians and Lagrangian correspondences]
Let $f\colon \XB \to \YB$ be a map of derived stacks and let $\YB$ have a $n$-shifted symplectic structure. Then $f$ is $n$-shifted Lagrangian if 
\begin{enumerate}
    \item $f$ is isotropic : there is a homotopy $\gamma \colon f^{*} \omega_{\YB} \sim 0$
    \item the null homotopic sequence $\mathbb{T}_{\XB} \to f^{*} \mathbb{T}_{\YB} \cong f^* \mathbb{L}_{\YB}[n] \to \mathbb{L}_{\XB}[n]   $ induced by $\gamma $ is a fiber sequence.
\end{enumerate}
A correspondence $\XB \xleftarrow{} \ZB \to \YB$, with $\XB$ and $\YB$ $n$-shifted symplectic is $n$-shifted Lagrangian if the induced map $\ZB \to \XB \times \overline{\YB}$ is a $n$-shifted Lagrangian. Here the shifted symplectic structure on $\XB \times \overline{\YB}$ is given by $\pr^{*}_{\XB} \omega_{\XB} - \pr^{*}_{\YB} \omega_{\YB}$.
\end{defn}
The rest of this subsection consists of various examples that will come up in various proofs in the rest of the paper.
\begin{ex}[Shifted (co)tangent bundles and derived critical loci] \label{dcrit_co}
 Let $\XB$ be a derived stack and $\mathbb{L}_{\XB}$ the cotangent complex then define the $n$-shifted cotangent bundle $\To^{*}[n] \XB = \Tot(\mathbb{L}_{\XB} [n] )$ and $n$-shifted tangent bundle $\To[n] \XB = \Tot(\mathbb{T}_{\XB} [n] )$. Recall that a map $S \to \To^{*}[n]\XB$ is given by a map $x \colon S \to \XB$ and a section $s \in \Gamma(x^{*} \mathbb{L}_{\XB})[n]$. Take the identity map $\To^{*}[n] \XB  \xrightarrow{\id} \To^{*}[n] \XB $ this map corresponds to the projection $\pi \colon \To^{*}[n]\XB \to  \XB$ and a section $l_{\XB} \in \Gamma(\pi^{*} \mathbb{L}_{\XB}[n])$. Define the Liouville form $\lambda_{\XB}$ to be the image of $l_{\XB}$ under the map $\pi^{*} \mathbb{L}_{\XB}[n] \to \mathbb{L}_{\To^{*}[n] \XB}[n]$.  It is shown in \cite{calaque_cot} that $d_{\dR} \lambda_{\XB}$ defines an $n$-shifted symplectic structure on $\To^{*}\XB [n]$. \par
 Let $f \colon U \to \mathbb{A}^{1}$ be a function from a smooth variety $U$. The derived critical locus is defined as the intersection $\critB f \coloneqq U \times_{\To^*U} U$, where the first map to $\To^{*}U$ is the zero section and the second is induced by $\operatorname{d} \!f$. Assume that an algebraic group $G$ acts on $U$ and $f$ is $G$-invariant. Then $f$ descends to a function $f/G : U/G \to \mathbb{A}^{1}$ and the derived critical locus is defined in the same way as $\critB (f/G) \coloneqq U/G \times_{\To^*(U/G)} U/G$. There is a canonical $(-1)$-shifted symplectic structure in both cases given by derived Lagrangian intersection. See \cite[Section 2.2]{PTVV}.
\end{ex}
\begin{ex}[Quotient and Classifying stacks] \label{quotient_class_forms}
Let $X$ be a scheme with an action of a reductive group $G$. We have $\mathbf{DR}(X/G) \cong (\Sym ( \mathbb{L}_{X} [-1] )\otimes \Sym \mathfrak{g}^{*}[-2])^{G}$.  \par
Consider  $X = \pt $. We know that $\mathbb{L}_{\B G} = \mathfrak{g}^{*}[-1]$ with coadjoint action.
We have 
$$\mathbf{DR}(\B G) \cong \Gamma(\B G, 
\Sym\mathbb{L}_{\B G}[-1]) = \Gamma(\B G, \Sym \mathfrak{g}^{*}[-2]) = (\Sym \mathfrak{g}^{*}[-2])^{G}$$
This follows since we are taking derived global sections, which in the case of $G$ reductive is just $G$-invariants. Therefore, $(2)$-shifted closed forms on $ \B G$ are exactly invariant bilinear forms $\omega \in (\Sym^{2}(\mathfrak{g}^{*}))^{G}$ and shifted symplectic structures correspond to nondegenerate forms.
\end{ex}
\begin{ex}[Tangent and cotangent stacks of BG] \label{tangent_bg}
We can show that $\To^{*}[1] \B G \cong [\mathfrak{g}^{*} / G]$ and $\To[-1] \B G \cong [ \mathfrak{g} / G]$. Now the $(2)$-shifted symplectic structure on $\B G$ can be thought of as the $G$-equivariant isomorphism $\psi \colon \mathfrak{g}^{*} \to \mathfrak{g}$. This then gives us an isomorphism $\psi \colon \To^{*}[1] \B G \to \To[-1] \B G$. $\To^{*}[1] \B G$ has a canonical $1$-shifted symplectic structure given by the exact $2$-form $\gamma_{0}$ and we can define a $1$-shifted symplectic structure on  $\To[-1] \B G 
 \cong [\mathfrak{g} / G]$ using the $2$-shifted symplectic structure on $\B G$.
\end{ex}
\subsection{AKSZ construction}\label{aksz_sect}
Now let us describe how one can induce $n$-forms, closed $n$- forms and symplectic structures from $\YB$ to $\Map(\XB, \YB)$ via transgression. This is called in \cite{PTVV} the AKSZ construction. We require a compactness condition on $\XB$.
\begin{defn}
Let $\XB$ be a derived stack, then it is $\mathcal{O}$-compact if the following hold
\begin{enumerate}
    \item $\mathcal{O}_{X}$ is a compact object in $\QCoh \XB$ 
    \item for any perfect complex $\mathcal{F}$ we have that $\Hom(\mathcal{O}_{\XB},\mathcal{F}) \in \Perf \mathbb{C}$ is a perfect complex of vector spaces. 
\end{enumerate}
\end{defn}
This is the compactness condition used in \cite{PTVV}. For us it will also be useful to recall the following finiteness conditions from \cite{naef2023torsion}.
\begin{defn}[Finiteness assumption] \cite[Assumption 1.1]{naef2023torsion} 
Let $\XB$ be a prestack such that
\begin{enumerate}
    \item $\mathcal{O}_{\XB}$ is a compact object of $\QCoh ( \XB)$
    \item The functor $p^*$ has a left adjoint $p_{\#}  \colon \QCoh (\XB) \to \Vect$. 
\end{enumerate}
\end{defn}
\begin{prop}[Consequences of assumption 1.1]  \cite[Proposition 1.4,1.5]{naef2023torsion} \label{mapping_cot}
Let $\XB$ be a prestack that satisfies assumption 1.1 then we have
\begin{enumerate}
    \item every perfect complex is compact
    \item For a derived scheme $\SB$ we have that for $\pi \colon \SB \times \XB \to \SB$ a left adjoint $\pi_{\#}$ to $\pi^{*}$ exists
    \item $\XB$ is $\mathcal{O}$-compact
    \item Let $\YB$ be a prestack with perfect cotangent complex. Let us also define the maps $\pi \colon \XB \times \Map(\XB, \YB) \to \Map(\XB, \YB)$ and $\ev\colon \XB \times \Map(\XB, \YB) \to \YB$. Then we have that $\mathbb{L}_{\Map(\XB,\YB)} = \pi_{\#} \ev^{*} \mathbb{L}_{\YB}$. 
\end{enumerate}
\end{prop}

\begin{defn}\cite[Definition 2.1]{PTVV} and  \cite[Definition 1.21]{naef2023torsion} 
Let $\XB$ be an $\mathcal{O}$-compact stack, a $\mathcal{O}$-\textbf{orientation} on $\XB$  of degree $d$ is a morphism $p_{*} \mathcal{O}_{\XB} \to \mathbb{C}[-d]$ with a non-degeneracy condition. \par
Now let $\XB$ further satisfy assumption 1.8 and 1.15 in \cite{naef2023torsion}. Then a \textbf{fundamental class} of degree $d$ is a morphism $\mathbb{C} \to p_{\#} \mathcal{O}_{\XB}[-d]$. The existence of a fundamental class gives a natural isomorphism $p_{*} \to p_{\#}[-d]$. Then  \cite[Proposition 1.22]{naef2023torsion} proves that a fundamental class gives an $\mathcal{O}$-orientation.
\end{defn}
We will not specify the assumptions 1.8 and 1.11 in \cite{naef2023torsion} since they will be satisfied in the following example, which is the only example we will use.
\begin{ex}[Betti stacks] \label{betti_stack}
Our main example will be the Betti prestack $M_B$ for $M$ a $n$-dimensional closed, connected, oriented manifold. $M_{B}$ is the constant prestack given by sending any cdga $A$ to $M$ considered just as a topological space. In particular, it is shown in \cite{naef2023torsion} that $M_B$ satisfies the finiteness assumptions above. By \cite[Proposition 3.19]{naef2023torsion} the fundamental class induces a fundamental class on the stack $M_{B}$. If it is clear from context, we will denote $M$ and $M_B$ by the same symbol.
\end{ex}
\begin{thm}[AKSZ] \cite[Theorem 2.5]{PTVV} \cite[Section 3.2]{calaque2022aksz}\label{aksz}
Let $\XB$ be an $\mathcal{O}$-compact stack with an orientation of degree $d$ and $\YB$ a $n$-shifted symplectic derived artin stack, then there is a canonical $(n-d)$-shifted symplectic structure on $\Map(\XB,\YB)$. Furthermore, the construction is covariantly functorial in the target: namely, given an $n$-shifted symplectic map $\YB_{1} \to \YB_{2}$ we get an $(n-d)$-shifted symplectic map $\Map(\XB, \YB_{1}) \to \Map(\XB, \YB_{2})$.
\end{thm}

\section{D-critical loci and DT sheaves}
In this section we will define d-critical structures which we view as classical truncations of  $(-1)$-shifted symplectic stacks and schemes. We will briefly recall the construction of the global DT sheaf on a d-critical scheme or stack. The cohomology of this sheaf in the case of the moduli space of $G$-local systems on the $3$-torus is the main object of study in this paper.
\subsection{D-critical loci on schemes}
\begin{prop}\cite[Theorem 6.1]{brav2019darboux} \label{s_sheaf}
Let $X$ be a $\mathbb{C}$-scheme or a complex analytic space. For every open $R \xrightarrow{j} X$ with a closed embedding $R \xrightarrow{i} U$ into a smooth scheme $U$ we have a short exact sequence
    $$0 \to I_{R,U} \xrightarrow{i} i^{-1} \mathcal{O}_{U} \to \mathcal{O}_{R} \to 0.$$
There exists a sheaf $S_{X}$ of $\mathbb{C}$-vector spaces on $X$ such that the following hold 
\begin{enumerate}
    \item there is an exact sequence $0 \to S_{X} |_{R} \xrightarrow{i} i^{-1} \mathcal{O}_{U} / I_{R,U}^{2} \xrightarrow{d} i^{-1} \Omega_{U}/ I_{R,U} \cdot i^{-1} \Omega_{U}$
    \item there is a morphism $\beta_{X} \colon S_{X} \to \mathcal{O}_{X}$ inducing a decomposition $S_{X} = \mathbb{C}_{X} \oplus S^{0}_{X}$, with 
$S^{0}_{X} = \ker(S_{X} \xrightarrow{\beta_{X}}  \mathcal{O}_{X} \to \mathcal{O}_{X,\red})$ and $\mathbb{C}_X$ the constant sheaf.
 \item we have an exact sequence
$$0 \to \HH^{-1}(\mathbb{L}_{X}) \to S_{X} \to \mathcal{O}_{X} \to \Omega_{X} $$
and a similar one for $S^{0}_{X}$.
\end{enumerate}

\end{prop}

\begin{defn}[D-critical locus] \cite[Definition 3.1]{ben2015darboux}
A $d$-critical structure on a scheme $X$ or complex analytic space is the data of a section $s \in \HH^{0}(S^{0}_{X})$ and for each point $p \in X$ a \textbf{critical chart} $ (R, U , f , i) $ with $x \in R$ an open of $X$, $i \colon R \to U$ a closed embedding,  $U$ smooth and $U \xrightarrow{f} \mathbb{C}$ a function satisfying $s|_{R} = i^{-1}f + I^{2}_{R,U}$ and $i(R) = \crit f$. \\
Let $f \colon X \to Y$ be a morphism of schemes or complex analytic spaces with d-critical loci structures, then there is an induced map $f^{\star} \colon f^{-1} S^{0}_{Y} \to S^{0}_{X}$. We say $f$ is a morphism of d-critical loci if $f^\star s_Y = s_X$.  
\end{defn}
\begin{defn}[Canonical bundle and orientation]   \cite[Theorem 6.4]{brav2019darboux}
Let $(X,s)$ be a $d$-critical locus. Then there exists a line bundle $K_{X}$ on $X_{\red}$ with the property that for every critical chart $(R,U,f,i)$ there is an isomorphism
$$i \colon K_{X}|_{R_{\red}} \to (\omega_{U}^{\otimes 2}) |_{R_{\red}}.$$
An \textbf{orientation} for a $d$-critical locus is the data $(L,\theta)$ with $L$ a line bundle on $X$ and  an isomorphism $\theta \colon L^{\otimes 2} \to K_{X}$. An isomorphism of two orientations $\psi \colon (L_{1}, \theta_{1}) \to (L_{2}, \theta_{2})$   is given by an isomorphism $\psi \colon L_{1} \to L_{2}$ which satisfies $\theta_{2} \circ \psi^{\otimes 2}  = \theta_{1}$.
\end{defn}
Now we are ready to state the theorem about the passage from $(-1)$-shifted symplectic derived schemes to $d$-critical loci.
\begin{thm} \cite[Theorem 6.6]{brav2019darboux}\label{truncation}
Let $\XB$ be a $(-1)$-shifted symplectic derived scheme.
Then the truncation $\tt_{0}(\XB)$ has a natural structure of $d$-critical locus with $(\det \mathbb{L}_{\XB} )|_{\tt_{0}\XB_{\red}} = K_{\tt_{0}(\XB)}$.
\end{thm}
\begin{remark}
There are different $(-1)$-shifted symplectic derived schemes that give the same d-critical locus. This is explained in \cite[Example 2.17]{Joyce_dcrit}.
\end{remark}
Summing up, we have that the classical truncation of our $(-1)$-shifted symplectic scheme can be expressed locally as a critical locus, in the form of critical charts. 
We recall here some first examples of d-critical schemes 
\begin{ex}[Smooth space]
Let $U$ be smooth, then the function $U \xrightarrow{0} \mathbb{C}$ defines a global chart. Note that we could define more complicated critical charts, this is just the obvious one.
However, $(U,0)$ is the unique $d$-critical structure on $U$. This can be seen by considering the exact sequence relating $\mathbb{L}_{U}$ and definition of $S^{0}_{U}$ to see that $S^{0}_{U} = 0$. In this case, there is an orientation for $K_{U} = \omega^{\otimes 2}_{U}$, since we can take $L =  \omega_{U}$, however, we can tensor $L$ with a non-trivial $2$-torsion line bundle and also change the isomorphism, so there are many choices for orientation a priori.
\end{ex}
We now turn to more precisely describe the procedure of gluing up vanishing cycles on d-critical schemes. First let us define vanishing cycles on a critical chart. See \cite[Definition 2.12]{brav2012symmetries}
\begin{defn}
Consider $f \colon U \to \mathbb{C}$ and $X = \crit f$. Define $X_{c} = f^{-1}(c) \cap X$. Then the vanishing cycles sheaf on a chart is
$$PV_{U,f} = \bigoplus_{c \in f(X)} \varphi_{f-c}(\mathbb{Q}_{U}[\dim U])|_{X_{c}}.$$
\end{defn}
The sheaf $PV$ can be shown to be in the category $\Perv(X)$. Perverse sheaves form a stack and hence can be glued on an open cover. However, to glue the sheaves $PV$ we will need to use the orientation we defined earlier. It turns out one needs to define an extra $\mathbb{Z} / 2 \mathbb{Z}$ bundle and an associated  local system.
\begin{defn}[Orientation principal bundle]\label{orientation_loc}
Let $X$ be a $d$-critical locus with canonical bundle $K_{X}$ and orientation $(K^{1/2}_{X}, \theta)$ and $C = (U,R,f,i)$ a critical chart.
Define a principal $\mathbb{Z} / 2 \mathbb{Z}$ bundle $\pi \colon Q_{C} \to R$, which has sections that are maps $s \colon K^{1/2}_{X} |_{R_{\red}} \to i^{*}K_{U}$ with $s \otimes s = \theta$. So sections of this bundle are square roots of the isomorphism $\theta$.
\end{defn}
\begin{thm}\cite[Theorem 6.9]{brav2012symmetries}
Let $(X,s)$ be a $d$-critical locus with an orientation $K^{1/2}_{X}, \theta$. Then there is a perverse sheaf $\varphi_{X}$ on $(X,s)$ such that for a critical chart $C = (U,R,f,i)$ we have 
$$\varphi_{X} |_{R} \cong i^{*}(PV_{U,f}) \otimes Q_{C} $$
where $Q_{C}$ is the orientation principal bundle defined above. $\varphi_X$ can be upgraded to a mixed Hodge module on $X$.
\end{thm}
\begin{ex}[Global critical locus] \label{global_crit_locus}
Let $X = \crit f \xrightarrow[]{\iota} U$ for $f \colon U \to \mathbb{C}$, then we have an isomorphism $K_{X} \cong \iota^{*} K^{\otimes 2}_{U} |_{X_{\red}} $. Furthermore, we can set $K^{1/2}_{X} = \iota^{*} K_{U} |_{X_{\red}}$ and take $K^{1/2}_{X} \xrightarrow[]{id} \iota^{*} K_{U} |_{X_{\red}}$. Then the orientation isomorphism clearly has a global square root, from which we can conlude that the orientation local system is trivial. Therefore the DT sheaf $\varphi_{X}$ is $\varphi_f$ the sheaf of vanishing cycles.
\end{ex}

\subsection{DT sheaf on stacks}
A similar story holds for $(-1)$-shifted symplectic derived Artin stacks. There is a truncation to $d$-critical stacks and similarly a perverse sheaf. We will briefly recall the constructions. First, we need to define sheaves on Artin stacks. As in \cite[Section 2.7]{Joyce_dcrit} we work with the site Lis-\'et$(X)$. A sheaf $\mathcal{F}$ on $X$ will be the data of an \'etale sheaf $\mathcal(F)_{T}$ for every smooth map $f \colon T \to X$ from a scheme $T$ with some compatibility conditions.
\begin{prop}[D-critical structures for stacks]  \cite[Corollary 2.52]{Joyce_dcrit} \label{d_critc_stacks}
Let $X$ be an Artin stack or complex analytic stack.
\begin{enumerate}
    \item We have the following 
\begin{enumerate}
    \item there exists a sheaf $S_{X}$ of $\mathbb{C}$-vector spaces on $X$ such that for each smooth morphism $f \colon T \to X$ we have an isomorphism $\theta_{f} \colon f^{*} S_{X} \to S_{T}$ 
    \item there is a canonical splitting $S_{X} = \mathbb{C}_{X} \oplus S^{0}_{X} $.
\end{enumerate}
\item $X$ has a d-critical structure if there is a section $s_{X} \in S^{0}_{X}$ such that for each smooth morphism $f \colon T \to X$ we have that $f^{*} s_{X}$ defines a d-critical structure on $T$. We call $X$ a d-critical stack.
\item Let $X$ be a d-critical stack. Then there is a canonical line bundle $K_{X}$ on $X_{\red}$. An \textbf{orientation} on a d-critical stack $X$ is the data $(L, \theta)$ of a line bundle $L$ on $X_{\red}$ and an isomorphism $\theta \colon L^{\otimes 2} \to K_{X}$.
\end{enumerate}
\end{prop}
Again we have a truncation theorem
\begin{thm}\cite[Theorem 3.18]{ben2015darboux} \label{truncation_dcrit_stacks}
Let $\XB$ be a $-1$-shifted symplectic derived artin stack. Then the truncation $\tt_{0}(\XB)$ has a natural structure of $d$-critical stack and $(\det \mathbb{L}_{\XB} )|_{(\tt_{0}\XB)_{\red}} \cong K_{\tt_{0}(\XB)}$. 
\end{thm}
Then similarly one can define the DT sheaf on stacks.
\begin{thm} \cite[Theorem 4.8]{ben2015darboux}  \label{joyce_sheaf_stacks}
Let $(X,s)$ be a $d$-critical stack with an orientation $(K^{1/2}_{X}, \theta)$. Then there is a perverse sheaf $\varphi_{X}$ on $(X,s)$ such that for each smooth map $f \colon T \to X$ we have $f^{*}[d] \varphi_{X} \cong \varphi_{T} $. Here $d$ is the relative dimension of $f$ and $T$ has the induced $d$-critical structure from $X$. Furthermore, there is an upgrade of $\varphi$ to a mixed Hodge module on $\XB$.
\end{thm}
\begin{ex}[Products] \cite[Proposition 4.3]{kinjo2024cohomological} \label{dcrit_products}
Let $\XB$ and $\YB$ be $(-1)$-shifted symplectic oriented stacks. Then  the d-critical locus structure on $X \times Y$ is given by $s_{X} \oplus s_{Y}$ and we have that $\varphi_{X \times Y} \cong \varphi_{X} \boxtimes \varphi_{Y} $. Also see \cite[Remark 5.23]{joyce_conje_lino_ben}.
\end{ex}
We will need to be a bit more explicit about d-critical structures on quotient stacks. We recall \cite[Section 3.2]{ben2015darboux}, where it is explained that $d$-critical structures on a quotient stack are the same as $G$-equivariant $d$-critical structures. Let $G$ be an algebraic group acting on a scheme $X$ and denote the action map by $a \colon G \times X \to X$ and the projection map by $\pi \colon G \times X \to X$. Then a $G$-equivariant $d$-critical structure is a section $s \in \Gamma(X, S^{0}_{X})$, with the property that $\pi^{\star} s = a^{\star} s \in \Gamma(G \times X, S^{0}_{G \times X})$. Equivalently for each $g \colon X \to X$ we have that $g^{\star}s = s$, where $g$ is the map induced by the action of $g \in G$.
\begin{defn}[Equivariant orientation data]
Let $X$ be a $G$-equivariant critical locus. Then the canonical bundle $K_X$ has a canonical $G$-equivariant structure. A $G$-equivariant orientation data is the data of a $G$-equivariant line bundle $L$ together with a $G$-equivariant isomorphism $L \otimes L \to K_X$.
\end{defn}
In particular, we have $\varphi_{X/G}$ on $X/G$ as defined in Theorem \ref{joyce_sheaf_stacks}. This sheaf satisfies the equation $p^{*}[\dim G] \varphi_{X/G} = \varphi_{X}$ for $p \colon X \to X/G$. So the perverse sheaf $\varphi_{X}$ is $G$-equivariant. Furthermore, here we can upgrade to a $G$ equivariant mixed Hodge module, using the mixed Hodge module structure on $\varphi_{X}$.
Let $X = Z/G$ be a global critical locus with $Z= \crit f$ for $f \colon U \to \mathbb{A}^{1}$ and $f$ $G$-invariant. We have $K_{Z} \cong K_{U}|_{\red}$ as $G$-equivariant sheaves. Local sections of the orientation principal bundle $P$ are given by local morphisms $s \colon K^{1/2}_{Z} \to K_{Z}$. Therefore, we can act on the set of $s$ by $G$, giving $P$ a $G$-equivariant structure and thus the associated  $\mathbb{Z}/ 2 \mathbb{Z}$ local system is $G$-equivariant as well.

\section{\texorpdfstring{$S^{0}$}{S 0} sheaves for formal completions} \label{s0_formal_sect}
In this section we establish some folklore results about $S$ sheaves on formal completions as well as comparisons between the $S$ sheaves of an Artin stack $X$, its analytification $X_{\an}$ and its formal completion at a point $\widehat{X}^{x}$. The results in this section are technical in nature and will only be used in the proof of Theorem \ref{exp_dcrit}. \par
Let us start by recalling some facts about formal completions.\begin{defn}[Formal completions] \label{formal_comp_forms}
    Let $f \colon \YB \to \XB$ be a map of derived prestacks. Define the completion $\widehat{\YB}^{f} = X \times_{\YB_{\dR}} \XB_{\dR}$. Here $(\XB)_{\dR}(R) = \XB(R_{\red})$ for more on the de Rham stack see \cite[Definition 2.1.3]{cptvv}. \par  
Completion is compatible with truncation in the sense that $\tt_{0}(\widehat{\XB}^{x}) = \widehat{\tt_{0} (\XB)}^{x}$.
\end{defn}
Note that if we consider a classical stack and embed it into derived prestacks it is generally only locally almost of finite type. Similarly formal completions of finite type schemes are locally almost of finite type by \cite[Corollary 6.3.2]{dg_indschemes}.  First, let $X$ be a finite type scheme and $\widehat{X}^{x}$ its completion at a point. Then we have the following lemma, which shows that the formal completion defined in terms of the de Rham stack as in Definition \ref{formal_comp_forms} agrees with the more classical definition of completing along a closed subscheme.

\begin{lem}\cite[Lemma B.1.2]{cptvv} \label{cptvv_formal_lemma}
    Let $R$ be a noetherian classical ring and denote $X = \Spec R$. Let $X_{n} = \Spec (R/I^{n})$, with $I$ being the ideal defining $x$ in $\Spec R$, then we have an equivalence of prestacks
    \begin{equation}
        \widehat{X}^{x} \coloneqq X \times_{X_{\dR}} \pt \cong \colim_{n \in \mathbb{N}} X_{n}
    \end{equation}
\end{lem}
This allows us to identify $\QCoh (\widehat{X}^{x}) \cong \lim_{n \in \mathbb{N}} \QCoh (X_{n})$ and the functor $\eta^{*}$ induced by $\eta \colon \widehat{X}^{x} \to X$ is $\eta^{*}(M) = (M / I^{n} M)_{n \in \mathbb{N}}$. We also have the completion $\widehat{M} = \lim M/I^{n}M$. We can identify $\Gamma(\widehat{X}^{x}, \eta^{*} M) = \widehat{M}$.
When considering formal completions of algebraic stacks we can reduce to the case of quotient stacks in the following way. By \cite[Theorem 4.12]{AHR_luna} there is an \'etale map $  Y / G_{x} \to X$ for $Y = \Spec R$.  Therefore, the formal completions of $X$ and $Y/ G_{x}$ will coincide and we have $\widehat{X}^{\B G_{x}} = \widehat{Y}^{x} / G_{x}$ and $\widehat{X}^{x} = \widehat{Y}^{x} / \widehat{G}^{1}_{x}$. Here $\widehat{X}^{\B G_{x}}$ is the completion along the map $\B G_{x} \to X$.  \par
\begin{defn}[Stack of exact $2$ forms]
We can define the following functor:
\begin{align*}
    \mathcal{A}^{2,\exa}(-,-1) \colon \operatorname{cdga}^{\leq 0 }  & \to \operatorname{Spc} \\
    R & \mapsto |\cone(\mathbf{DR}(R)(0) \xrightarrow{\epsilon} \mathbf{DR}(R)(1)[1])[-1]|
\end{align*}
\end{defn}
Note that for a cdga $R$ the map $\mathbf{DR}(R)(0) \xrightarrow{\epsilon} \mathbf{DR}(R)(1)[1]$ can be written as the map 
$$ \cone(R \xrightarrow{d_{\dR}} \mathbb{L}_{R})[-1].$$  This functor satisfies \'etale descent and therefore we can define the space of exact forms for a derived stack $\XB$ by right Kan extension. By right Kan extension we can then also write $\mathcal{A}^{2,\exa}(\XB) = |\cone(\mathbf{DR}(\XB)(0) \xrightarrow{\epsilon} \mathbf{DR}(\XB)(1)[1])[-1]|$. In particular, we are interested in $\XB = \widehat{X}^{x}$ a completion of a finite type scheme at a point. In this case, we can use Theorem \ref{forms_artin} to deduce that
\begin{align*}
    \mathcal{A}^{2,\exa}(\widehat{X}^{x}) & \cong |\cone(\Gamma(\widehat{X}^{x},\mathcal{O}_{\widehat{X}^{x}} \xrightarrow{d_{\dR}} \mathbb{L}_{\widehat{X}^{x}}))[-1]| \\
    \mathcal{A}^{2,\exa}(X) & \cong |\cone(\Gamma(X,\mathcal{O}_{X} \xrightarrow{d_{\dR}} \mathbb{L}_{X}))[-1]|.
\end{align*}
The generality of Theorem \ref{forms_artin} is necessary here since we cannot use \cite[Proposition 1.14]{PTVV} since $\widehat{X}^{x}$ is not Artin and $X$ is not of finite type as a derived prestack.
\begin{prop}  \cite[Proposition 3.2]{kinjo2024cohomological}  \label{ses_closed_forms}
Let $X$ be a derived Artin stack locally almost of finite type. Then we have an exact sequence of stacks
$$\mathbb{C} \to \mathcal{A}^{2,\exa}(-,-1) \xrightarrow{d_{\dR}} \mathcal{A}^{2,\cl}(-,-1) $$
where $\mathbb{C}$ is the constant prestack that assigns $R \mapsto \mathbb{C}$. Furthermore, there is a splitting of this sequence.
\end{prop}
\begin{proof}
    For finite type stacks, it is enough to prove this affine locally for which we can use \cite[Proposition 5.6,5.7]{brav2019darboux}. In particular, if $R$ is a cdga, then $\pi_{0}\mathcal{A}^{2,\exa}( \Spec R,-1)$ consists of pairs $(f, \alpha)$ with $f$ of degree $0$ and $\alpha$ a one form of degree $-1$. The splitting is given by restricting the function $f$ to $ (\tt_{0} \Spec R)_{\red}$. For the extension to stacks locally almost of finite type see \cite{kinjo2024cohomological}.
\end{proof}

We will now give a definition of sections of the $S$-sheaf which works for formal completions.
\begin{defn}[S sheaf for formal completions] \label{new_S_sheaf_def}
Let $X$ be a classical Artin stack. If $\widehat{X}^{x}$ is the completion of an Artin stack at a point $x$, we define the vector spaces $S_{\widehat{X}^{x}} = \pi_{0}\mathcal{A}^{2,\exa}(\widehat{X}^{x},-1)$ and $S^{0}_{\widehat{X}^{x}} = \pi_{0}\mathcal{A}^{2,\cl}(\widehat{X}^{x},-1)$. 
\end{defn} 
The following proposition now ensures that the definition above is compatible with the original Definition \ref{s_sheaf}.
\begin{prop} \label{closed_forms_ssheaf}
 Let $X$ be a classical scheme or Artin stack. Then we have 
 \begin{equation}
     \pi_{0}\mathcal{A}^{2, \exa }(X,-1) \cong \Gamma(X, S_{X}),
 \end{equation}
 where the sheaf $S_{X}$ is as defined in \ref{s_sheaf}.   \par 
 Let $\eta_{X} \colon \widehat{X}^{x} \to X$ be the formal completion of a classical scheme at a point $x \colon \pt \to X$.
    Fix an open neighbourhood $R \subseteq X$ of $x$ and a closed immersion $R \xrightarrow{i} U$ with ideal $I$. This induces a map $\widehat{\mathcal{O}}^{x}_{U} \to \widehat{\mathcal{O}}^{x}_{R} = \widehat{\mathcal{O}}^{x}_{X}$. Then  $S_{\widehat{X}^{x}}$ fits into the following short exact sequence of vector spaces
    \begin{equation} \label{formal_s_sheaf_ses}
        0 \to S_{\widehat{X}^{x}} \to \widehat{\mathcal{O}}^{x}_{U} / \widehat{I}^{2} \to  \widehat{\Omega}^{x}_{U} / \widehat{I} \widehat{\Omega}^{x}_{U}.
    \end{equation}
\end{prop}
\begin{proof}
The equivalence of the two definitions of sections of the $S$ sheaf follow by \cite[Remark 2.2b]{Joyce_dcrit}. In particular, one can consider the truncation of the cotangent complex of $X$. Given an embedding $i \colon X \to U$ for $U$ smooth we have
\begin{equation} \label{cot_truncated}
    \tau_{\geq -1} \mathbb{L}_{X} = I/I^{2} \to i^{*} \Omega_{U}
\end{equation} 
We can form the following exact sequence of complexes
\begin{equation} \label{ses_S_sheaf_cone}
\begin{tikzcd}
	{I/I^{2}} & {\mathcal{O}_{U}/I^{2}} & {\mathcal{O}_{X}} \\
	{i^{*}\Omega_{U}} & {i^{*}\Omega_{U}} & 0
	\arrow[from=1-1, to=1-2]
	\arrow[from=1-1, to=2-1]
	\arrow[from=1-2, to=1-3]
	\arrow[from=1-2, to=2-2]
	\arrow[from=1-3, to=2-3]
	\arrow[from=2-1, to=2-2]
	\arrow[from=2-2, to=2-3]
\end{tikzcd}
\end{equation}
Shifting we can show that $\cone(\mathcal{O}_{X} \to \tau_{\geq -1}\mathbb{L}_{X}) \cong \mathcal{O}_{U}/I^{2} \to i^{*}\Omega_{U}$. From this we can see that \begin{equation*}
     \pi_{0}\mathcal{A}^{2,\exa}(X,-1) \cong \Gamma(X, S_{X}).
 \end{equation*}
We will now repeat the same proof for formal completions. Note for the purposes of formal completion we can work affine locally so we can assume $X = \Spec A$ is affine and we have a closed embedding $X \to U$ with $U= \Spec R$ smooth and affine. Then we have an ideal $I \subseteq R$ such that $A = R/I$. \par
    The second property will follow from a description of the cotangent complex of the formal completion of $X$.  We also get induced maps on formal completions that make the following square commute.
    \begin{equation}
\begin{tikzcd}
	X & U \\
	{\widehat{X}^{x}} & {\widehat{U}^{x}}
	\arrow["i"', from=1-1, to=1-2]
	\arrow["{\eta_{X}}", from=2-1, to=1-1]
	\arrow["{\hat{i}}", from=2-1, to=2-2]
	\arrow["{\eta_{U}}", from=2-2, to=1-2]
\end{tikzcd}
    \end{equation}
 Since we are considering $\mathcal{A}^{2,\exa}(\widehat{X}^{x}, -1)$, it is enough to consider the truncation $\tau_{\geq -1} \mathbb{L}_{\widehat{X}^{x}}$. 
The maps $\eta_{X}$ and $\eta_{U}$ are formally \'etale so we get $\tau_{\geq -1}\mathbb{L}_{\widehat{X}^{x}} =  \eta^{*}_{X} \tau_{\geq -1} \mathbb{L}_{X}$. The inverse systems $\eta^{*} I/I^{2}$ and $\eta^{*} i^{*}\Omega_{U}$ are Mittag-Leffler since all the maps in the inverse system are surjective, therefore the limit functor does not have any higher cohomology. Then using Lemma \ref{cptvv_formal_lemma} and equation \eqref{cot_truncated} we can deduce that on global sections on $\widehat{X}^{x}$ we have
\begin{equation} \label{cot_truncated_formal}
    \tau_{\geq -1} \mathbb{L}_{\widehat{X}^{x}} = \widehat{I}/\widehat{I}^{2} \to i^{*} \widehat{\Omega}_{U}.
\end{equation} 
Now we can consider the following short exact sequence of complexes, which comes from completion of the analogous exact sequence \ref{ses_S_sheaf_cone}. Using the Mittag-Leffler condition again we get
\begin{equation}
\begin{tikzcd}
	{\widehat{I}/\widehat{I}^{2}} & {\mathcal{O}_{\widehat{U}}/\widehat{I}^{2}} & {\mathcal{O}_{\widehat{X}}} \\
	{\widehat{i}\Omega_{\widehat{U}}} & {\widehat{i}\Omega_{\widehat{U}}} & 0
	\arrow[from=1-1, to=1-2]
	\arrow[from=1-1, to=2-1]
	\arrow[from=1-2, to=1-3]
	\arrow[from=1-2, to=2-2]
	\arrow[from=1-3, to=2-3]
	\arrow[from=2-1, to=2-2]
	\arrow[from=2-2, to=2-3]
\end{tikzcd}
\end{equation}
By shifting this short exact sequence it follows that we have an isomorphism $\cone(\mathcal{O}_{\widehat{X}} \to \tau_{\geq -1}\mathbb{L}_{\widehat{X}}) \cong (\widehat{\mathcal{O}}^{x}_{U}/ \widehat{I}^{2} \to \widehat{\Omega}^{x}_{U})$. Since we defined the space of exact $-1$ forms to be the cohomology of the cone we get the desired description of $S^{0}_{\widehat{X}^{x}}$. 
\end{proof}
From this point of view we can express the induced $d$-critical structure on $\tt_{0} \XB$ in Theorem \ref{truncation_dcrit_stacks} as the one induced by the map $\mathcal{A}^{2,\cl}(\XB,-1) \to \mathcal{A}^{2,\cl}(\tt_{0} \XB , -1)$. Proposition $\ref{ses_closed_forms}$ now shows that we have a decomposition $S_{\widehat{X}^{x}} = S^{0}_{\widehat{X}^{x}} \oplus \mathbb{C}$.
\begin{remark}
    Note that in the classical or analytic setting it makes sense to also define the $S$ sheaf on $X$ for the \'etale or Zariski topologies on $X$. For formal completions at a point the underlying space is just a point so we only have a vector space. 
\end{remark}
\begin{lem} \label{ss_ssheaf_inj}
    Let $X$ be a finite type scheme. There is a map $\Gamma(X,S^{0}_{X}) \to \Gamma(X_{\an},S^{0}_{X_{\an}})$ and also an injective map on stalks $S^{0}_{X,x} \to S^{0}_{X_{\an},x}$. There are injective maps $S^{0}_{X,x} \to  S^{0}_{\widehat{X}^{x}}$ and $S^{0}_{X_{\an},x} \to  S^{0}_{\widehat{X}^{x}}$. Furthermore, there is a commutative diagram
    \begin{equation} \label{stalk_alg_formal_an}
\begin{tikzcd}
	& {S^{0}_{X,x}} \\
	{S^{0}_{\widehat{X}^{x}}} & {S^{0}_{X_{\an},x}}
	\arrow[hook, from=1-2, to=2-1]
	\arrow[hook, from=1-2, to=2-2]
	\arrow[hook, from=2-2, to=2-1]
\end{tikzcd}
    \end{equation}
    \end{lem}
\begin{proof}
Denote by $h \colon X_{\an} \to X$ the inclusion map. To define the map $\Gamma(X,S^{0}_{X}) \to \Gamma(X_{\an},S^{0}_{X_{\an}})$, as in \cite[Section 3.1]{Joyce_dcrit} we can cover $X$ by opens $R$  such that $R \xhookrightarrow{} U$ is a closed embedding into $U$ smooth. We can then use the following diagram
    \begin{equation}
\begin{tikzcd}
	R & U & {\mathcal{O}_{U}/I^{2}_{U}} & {\Omega_{U}/I_{U}\Omega_{U}} \\
	{R_{\an}} & {U_{\an}} & {h_*(\mathcal{O}_{U_{\an}}/I^{2}_{U_{\an}})} & {h_{*}(\Omega_{U_{\an}}/I_{U_{\an}}\Omega_{U_{\an}})}
	\arrow["i", hook, from=1-1, to=1-2]
	\arrow[from=1-3, to=1-4]
	\arrow[from=1-3, to=2-3]
	\arrow[from=1-4, to=2-4]
	\arrow[from=2-1, to=1-1]
	\arrow["{i_{\an}}", hook, from=2-1, to=2-2]
	\arrow[from=2-2, to=1-2]
	\arrow[from=2-3, to=2-4]
\end{tikzcd}
    \end{equation}
    which induces the map of short exact sequences
    \begin{equation}
\begin{tikzcd}
	0 & {S_{X}|_{R}} & {i^{-1}\mathcal{O}_{U}/I^{2}} & {i^{-1} \Omega_{U}/Ii^{-1}\Omega_{U}} \\
	0 & {h_{*}(S_{X_{\an}}|_{R_{\an}})} & {h_{*}i^{-1}_{\an}\mathcal{O}_{U_{\an}}/I^{2}_{\an}} & {h_{*}i^{-1}_{\an} \Omega_{U_{\an}}/I_{\an}i^{-1}_{\an}\Omega_{U_{\an}}}
	\arrow[from=1-1, to=1-2]
	\arrow[from=1-2, to=1-3]
	\arrow[dashed, from=1-2, to=2-2]
	\arrow[from=1-3, to=1-4]
	\arrow[from=1-3, to=2-3]
	\arrow[from=1-4, to=2-4]
	\arrow[from=2-1, to=2-2]
	\arrow[from=2-2, to=2-3]
	\arrow[from=2-3, to=2-4]
\end{tikzcd}
    \end{equation}
    This defines a map $S^{0}_{X} \to h_{*} S^{0}_{X_{\an}}$ and thus a map $\Gamma(X, S^{0}_{X}) \to \Gamma(X_{\an}, S^{0}_{X_{\an}})$.
    If the map $i^{-1} \mathcal{O}_{U} / I^{2} \to i^{-1} \mathcal{O}_{U_{\an}} / I^{2}_{\an}$ is injective, then the map $S_{X}|_{R} \to S_{X_{\an}|_{R_{\an}}}$ is injective. Note that $I = i^{-1}I_{U}$. The  sheaf $O_{U}/I^{2}_{U}$ is coherent so the canonical map to the analytification is injective. This also means that the map on stalks $S^{0}_{X,x} \to S^{0}_{X_{\an},x}$ is injective. 

        The maps $S^{0}_{X,x} \to  S^{0}_{\widehat{X}^{x}}$ and $S^{0}_{X_{\an},x} \to  S^{0}_{\widehat{X}^{x}}$ are defined by taking colimits over analytic opens $U \subseteq X_{\an}$ or Zariski opens $U \subseteq X$ of the maps $\Gamma(U,S^{0}_{X_{\an}}) \to \Gamma(X, S^{0}_{\widehat{X}^{x}})$ or $\Gamma(U,S^{0}_{X}) \to \Gamma(X, S^{0}_{\widehat{X}^{x}})$ respectively.
    
    To prove the maps to $S^{0}_{\widehat{X}^{x}}$ are injective we can use the argument in \cite[Proposition 3.12]{ricolfi_savvas}, where it is proven in the algebraic case using the description of $S^{0}_{\widehat{X}^{x}}$ in \eqref{formal_s_sheaf_ses}. The map $S^{0}_{X,x} \to  S^{0}_{\widehat{X}^{x}}$ is then induced by the map $\mathcal{O}_{U,x}/ I^{2}_{x} \to \widehat{\mathcal{O}}_{U,x}/ \widehat{I}^{2}_{x}$.  We can repeat the argument of \cite{ricolfi_savvas} also in the complex analytic case because the map $\mathcal{O}_{U_{\an},x} \to \widehat{\mathcal{O}}_{U_{\an},x} \cong \widehat{\mathcal{O}}_{U,x}$ is still faithfully flat. This follows because $\mathcal{O}_{U_{\an},x}$ is still a noetherian ring despite $\mathcal{O}_{U_{\an}}(U_{\an})$ not being noetherian in general. \par
    To prove that the diagram \eqref{stalk_alg_formal_an} commutes we can again consider the local models of $S$ sheaves and the commutative diagram
    \begin{equation}
\begin{tikzcd}
	& {\mathcal{O}_{U,x}/I^{2}_{x}} \\
	{\widehat{\mathcal{O}}_{U,x}/\widehat{I}^{2}_{x}} & {\mathcal{O}_{U_{\an},x}/I^{2}_{an,x}}
	\arrow[from=1-2, to=2-1]
	\arrow[from=1-2, to=2-2]
	\arrow[from=2-2, to=2-1]
\end{tikzcd}
    \end{equation}
This diagram commutes because of the fact that a completion of the algebraic functions and analytic functions at a point is the same. 
\end{proof}
The constructions of the above lemma are functorial in the sense that for a map of schemes $f \colon X \to Y$ the following diagrams commute
    \begin{equation}\label{an_to_formal}
\begin{tikzcd}
	& {S^{0}_{X_{\an},x}} & {S^{0}_{Y_{\an},x}} \\
	& {S^{0}_{\widehat{X}^{x}}} & {S^{0}_{\widehat{Y}^{y}}} \\
	{S^{0}_{X,x}} &&& {S^{0}_{Y,y}}
	\arrow[hook, from=1-2, to=2-2]
	\arrow["{f^*}"', from=1-3, to=1-2]
	\arrow[hook', from=1-3, to=2-3]
	\arrow["{\widehat{f}}"', from=2-3, to=2-2]
	\arrow[hook, from=3-1, to=1-2]
	\arrow[hook, from=3-1, to=2-2]
	\arrow[hook', from=3-4, to=1-3]
	\arrow[hook, from=3-4, to=2-3]
\end{tikzcd}
\end{equation}
\begin{equation}\label{functoriality}
\begin{tikzcd}
	{S^{0}_{\widehat{X}^{x}}} & {S^{0}_{\widehat{Y}^{x}}} \\
	{\Gamma(X,S^{0}_{X})} & {\Gamma(Y,S^{0}_{Y})} \\
	{S^{0}_{X,x}} & {S^{0}_{Y,y}}
	\arrow[from=1-2, to=1-1]
	\arrow[from=2-1, to=1-1]
	\arrow[from=2-1, to=3-1]
	\arrow[from=2-2, to=1-2]
	\arrow[from=2-2, to=2-1]
	\arrow[from=2-2, to=3-2]
	\arrow[bend left = 60, hook', from=3-1, to=1-1]
	\arrow[bend right = 60, hook, from=3-2, to=1-2]
	\arrow[from=3-2, to=3-1]
\end{tikzcd}
\end{equation}
 The diagram \eqref{functoriality} commutes either for complex analytic spaces or schemes.
We now explain how to upgrade the previous lemma to stacks. For quotient stacks $X/G$ we will by abuse of notation denote $\Gamma(X/G, S^{0}_{X/G})$ by  $\Gamma(X, S^{0}_{X})^{G}$. Strictly speaking $S^{0}_{X}$ is not $G$-equivariant in the usual sense since we are working with the lisse-\'etale site. 
\begin{lem} \label{stack_S_sheafcomm}
    Let $X$ be an Artin stack, $U \to X$ be an atlas and $x \in U$ a $\mathbb{C}$ point with stabiliser $G_{x}$.  Then there is a commutative diagram
    \begin{equation}
\begin{tikzcd}
	{S^{0}_{\widehat{X}^{\B G_x}}} & {\Gamma(X,S^{0}_{X})} \\
	{S^{0}_{\widehat{U}^{x}}} & {\Gamma(U,S^{0}_{U})} \\
	& {S^{0}_{U,x}} & {}
	\arrow[hook, from=1-1, to=2-1]
	\arrow[from=1-2, to=1-1]
	\arrow[hook, from=1-2, to=2-2]
	\arrow[from=2-2, to=2-1]
	\arrow[from=2-2, to=3-2]
	\arrow[hook, from=3-2, to=2-1]
\end{tikzcd}
    \end{equation}
    The same diagram commutes for $X$ an analytic stack.
\end{lem}
\begin{proof}
We can reduce to the case of a quotient stack with $U/G$ and a point $x \in U$ with stabiliser $G$. In particular, we have $\widehat{U/G}^{\B G} = \widehat{U}^{x}/G$. Then we immediately get the injectivity of the map $S^{0}_{\widehat{U}^{x}/G} =(S^{0}_{\widehat{U}^{x}})^{G}  \to S^{0}_{\widehat{U}^{x}}$. The action on  $S^{0}_{\widehat{U}^{x}}$ is given by the natural $G$ action on the space of closed $(-1)$-shifted two forms. The commutativity of the square follows from functoriality and the triangle commutes already from the previous lemma. The same argument works for the complex analytic case.
\end{proof}
This lemma says that given two sections $s_1$ and $s_{2} \in \Gamma(X,S^{0}_X)$ we can check if they agree at a point $x \in X$ by checking if they agree on the formal completion at that point. This follows from the lemma by the commutativity of the diagram and the injectivity of the maps. Note that for comparing sections of $\Gamma(X,S^{0}_{X})$ we can work on a subspace of $\Gamma(U,S^{0}_{U})$ by \cite[Proposition 2.54]{Joyce_dcrit}. Furthermore, even though $S^{0}_{U}$ is an \' etale sheaf the global sections are the same as the associated Zariski sheaf.

\section{Cohomological integrality for tripled quivers} \label{coh_int_quiver}
In this section we recall cohomological integrality for quivers with potential, which we will later use to deduce cohomological integrality for the $3$-torus. Let $Q$ be a quiver with vertex set $Q_0$ and edge set $Q_1$. Let us take a dimension vector $v \in \mathbb{N}^{Q_{0}}$. The representation variety of the quiver with respect to the dimension vector is given by $\Rep_{v}(Q) = \prod_{i \to j \in Q_{1}} \Hom(\mathbb{C}^{v_i} , \mathbb{C}^{v_{j}})$. There is a conjugation action on this variety by the gauge group $\GL_v = \prod \GL_{v_{i}}$. Define the stack of quiver representations $\MM_{Q} = \coprod_{v \in \mathbb{N}^{Q_{0}}} \Rep_{v}(Q) / \GL_{v}$ and denote each dimension component by $\MM_{Q,v}$. \par

Define the \textbf{doubled quiver} $\overline{Q}$ to be the quiver with the same vertices as $Q$ and for each edge $a \colon i \to j \in Q_1$ an opposite added arrow $a^{*} \colon j \to i$. We have an equivalence $\To^{*} \Rep_{v} Q \cong \Rep_{v} \overline{Q}$ and we will denote elements in $\To^{*} \Rep_{v} Q$ by $(\rho(a), \rho(a^{*}))$ for $a \in Q_{1}$. \par  Define the \textbf{tripled quiver} $\widetilde{Q}$ to have the same vertices as $Q$ but for each arrow $a \colon i \to j$ we add the opposite arrow $a^*  \colon j \to i$ and a loop $\omega_i$ for each vertex $i \in Q_0$. Denote by $\mathbb{C}Q$ the path algebra of $Q$. A \textbf{potential} $W$ on $Q$ is given by a linear combination of cyclic words, where two cyclic words are considered the same if they can be cyclically permuted to each other. If $W$ is a single cyclic word, then we define 
$$\frac{\partial W}{\partial a} = \sum_{W = cac^{'}} cc^{'},$$
where $c$, $c^{'}$ are paths in $Q$. The Jacobi algebra of a quiver $Q$ with potential $W$ is defined by $\Jac(Q,W) = \mathbb{C}Q / \langle \partial W / \partial a  \mid a \in Q_{1} \rangle$. For the tripled quiver we define a potential $\widetilde{W} = \sum_{a \in Q_{1}}[a,a^{*}] \sum_{i \in Q_{0}} \omega_{i}$.
\begin{ex}[Jordan quiver]
Let $Q_{\Jor}$ be the quiver with one vertex and one loop. The main example we will consider will be the tripled Jordan quiver. The tripled potential in this case is $\omega[a,a^*]$ giving us $\Jac(\widetilde{Q}_{\Jor},\widetilde{W}) = \mathbb{C}[a,a^{*},\omega]$ the polynomial algebra in $3$ variables.
\end{ex}
From now on we work with tripled quivers with the canonical potential, so we write $\Jac = \Jac(\widetilde{Q}, \widetilde{W})$. Using the potential for the tripled quiver we can define a function 
\begin{align*}
    \tr(\widetilde{W})_{v} \colon \Rep_{v}(\widetilde{Q}) & \to \mathbb{A}^{1} \\
    \rho & \mapsto \tr(\sum_{i \in Q_{0}} \rho(\omega_{i})\sum_{a \in Q_{1}}{[\rho(a),\rho(a^*)]} ).
\end{align*}
This function is $\GL_{v}$ invariant, and so descends to $\tr(\widetilde{W})_v/ \GL_v \colon \Rep_{v}(\widetilde{Q})/ \GL_v \to \mathbb{A}^{1}$. Taking all dimension vectors at once we get an induced function $\tr(\widetilde{W})  \colon \MM_{Q} \to \mathbb{A}^{1}$. We have 
$$\MM_{\Jac} \cong \coprod_{v} \crit (\tr(\widetilde{W})_v) / \GL_v$$ where $\MM_{\Jac}$ is a component of the classical truncation of the moduli of objects as in \cite{toen2007moduli}. Each component  of $\MM_{\Jac}$ is a global critical locus. Therefore, as in Example \ref{global_crit_locus} it is automatically oriented with trivial orientation twist and the DT sheaf is $\varphi_{\MM_{\Jac,v}} = \varphi_{\tr(\widetilde{W})_{v}}$. In this case we can also upgrade this to a monodromic mixed Hodge module. Let us define the good moduli space map $\pi_{Q} \colon \MM_{\Jac} \to X_{\Jac} $ with $X_{\Jac} = \coprod X_{\Jac, v}$. We can also describe the derived enhancements of $\MM_{\Jac}$ in the following way. Recall the definition of derived critical loci in Example \ref{dcrit_co}.
\begin{lem} \label{stacky_kinjo}
    We have an isomorphism of $(-1)$-shifted symplectic stacks
    \begin{equation}
        \To^{*}[-1] \To^{*} (\Rep_{v}(Q)/ \GL_{v}) \cong \crit (\tr(\widetilde{W})_v / \GL_{v}).
    \end{equation}
\end{lem}
\begin{proof}
    Let $\YB$ be a derived stack with an action of an algebraic group $G$. There is an induced action on $\To^{*}[n] \YB$ and a moment map $\To^{*}[n]\YB/G \to \mathfrak{g}^{*}[n]/G$. Then we have the following formula from \cite[Example 2.6]{calaque_anel}
    \begin{equation}
        \To^{*}[n](\YB / G )=\To^{*}[n]\YB/G \times_{\mathfrak{g}^{*}[n]/G} \B G.
    \end{equation}
    We will work with the case $\YB = \Rep_{v} (Q)$ and $G = \GL_{v}$. There is a canonical moment map
    \begin{align*}
        \mu \colon \To^{*} \Rep_{v}(Q) & \to \mathfrak{gl}_{v}  \cong \mathfrak{gl}^{*}_{v} \\
        (\rho(a),\rho(a)^{*}) \in \To^{*} \Rep_{v}(Q) & \mapsto \sum_{a \in Q_{1}} [\rho(a),\rho(a^{*})]
    \end{align*}
     where we have used the identification $\mathfrak{gl}_{v} \cong \mathfrak{gl}^{*}_{v}$ using the trace form. Using this moment map we construct the $\GL_{v}$-invariant function $\tr(\widetilde{W})_v \colon \To^{*} \Rep_{v}(Q) \times \mathfrak{gl}_{v} \cong \Rep_v ( \widetilde{Q}) \to \mathbb{A}^{1}$.
    Now there is a $(-1)$-shifted symplectic equivalence, which follows from \cite[Lemma 2.7]{dim_red_kinjo}.
    \begin{equation} \label{tasuki_lemma}
        \textbf{crit} \tr(\widetilde{W})_v \cong \To^{*}[-1] \ZB(\mu). 
    \end{equation}
Here $\ZB(\mu) $ is the derived zero locus. There is also an equivalence of $(-1)$-shifted moment maps 
\begin{align*}
    \textbf{crit} \tr(\widetilde{W})_v & \to \mathfrak{gl}^{*}_{v}[-1] \\
    & \text{and} \\
    \To^{*}[-1] \ZB(\mu) & \to \mathfrak{gl}^{*}_{v}[-1]. 
\end{align*}
Then by using \cite[Theorem A]{calaque_anel} in the first isomorphism we can conclude that
\begin{align*}
      \mathbf{crit} (\tr(\widetilde{W})_{v} / \GL_{v}) & \cong (\textbf{crit}  \tr(\widetilde{W})_v )/\GL_{v} \times_{\mathfrak{gl}^{*}_{v}[-1]/\GL_{v}}\B \GL_{v} \\   & \cong (\To^{*}[-1] \ZB(\mu))/\GL_{v} \times_{\mathfrak{gl}^{*}_{v}[-1]/\GL_{v}}\B \GL_{v}  \\ & \cong \To^{*}[-1](\mathbf{Z}(\mu)/\GL_{v})  \\ & \cong \To^{*}[-1] \To^{*} (\Rep_{v} (Q) / \GL_{v}).
\end{align*}
\end{proof}
There is a symmetric monoidal structure $\boxdot$ on $\DD^{+}_{c}(X_{\Jac})$ by convolution along the direct sum morphism $\oplus \colon X_{\Jac} \times X_{\Jac} \to X_{\Jac}$ given by the formula 
$$\mathcal{F} \boxdot \mathcal{G} = \oplus_{*} \mathcal{F} \boxtimes \mathcal{G}.$$ 
\begin{defn}[BPS sheaves for quivers]
Define 
$$\BPS_v  \coloneqq \begin{cases}
    \varphi_{\tr(\widetilde{W})} \IC_{X_{\widetilde{Q},v}} \text{ if } X_{Q,v, \text{simple}} \neq \emptyset \\
    0 \text{ otherwise }
    
\end{cases}$$ 
Here $X_{\widetilde{Q},v, \text{simple}}$ is the space of simple quiver representations of $\widetilde{Q}$ of dimension vector $v$. Again we have an upgrade to monodromic mixed Hodge modules. Upon taking cohomology we denote $\BPSo_{v} = \HHf(X_{Q,v}, \BPS_{v} )$. Finally $\BPS = \bigoplus_{v \in \mathbb{N}^{Q_{0}}} \BPS_{v}$ and $\BPSo = \bigoplus_{v \in \mathbb{N}^{Q_{0}}} \BPSo_{v}$.
\end{defn}
 We have the following theorem of Davison-Meinhardt
\begin{thm}\cite[Theorem A]{davison2020cohomological} \label{coh_int_dm}
 We have the following decomposition
\begin{equation}\label{integrality_quivers}
   \pi_{Q*} \varphi_{\tr(\widetilde{W})} \cong \Sym_{\boxdot}(\BPS \otimes \HHf(\B \mathbb{G}_{m})_{\vir}). 
\end{equation}
Here $\HHf(\B \mathbb{G}_{m})_{\vir} = \HHf(\B \mathbb{G}_{m}) \otimes \mathbb{L}^{1/2}$ . Furthermore, $\pH^{1} ( \pi_{Q*} \varphi_{\tr(\widetilde{W})}) \cong \BPS$ as perverse sheaves.
\end{thm}
Here one needs to make sense of mixed Hodge modules on stacks as in \cite[Section 2]{bps_less_perverse}. In the special case of tripled quivers with potential we have the following theorem 
\begin{thm}\cite[Corollary 3.9]{bps_less_perverse}
The BPS sheaf $\BPS_d$ is a pure monodromic mixed Hodge module. As a consequence $\pi_{Q*} \varphi_{\tr(\widetilde{W})}$ is a pure complex of mixed Hodge modules. 
\end{thm}
The last part of the theorem holds because $\BPS_v$ is only monodromic up to tensoring with $\mathbb{L}^{1/2}$. Since we take symmetric algebra of  $\BPS_{v} \otimes \HHf(\B \mathbb{G}_{m})_{\vir} $ this half Tate twist cancels out and we get a monodromy-free mixed Hodge module. To say more about the BPS sheaves appearing in the theorem we need the following support lemma
\begin{lem}[Support lemma] \cite[Lemma 4.1]{davison2022integrality} \label{support_lemma_quivers}
Let $x \in X_{\Jac}$ that corresponds to a $\Jac$ module $\rho$ that is in the support of $\BPS$. Then the set of generalised eigenvalues of $\rho(\omega_{i})$ for $i \in Q_0$ contains only one element.
\end{lem}
In other words, this means that the BPS sheaf is supported on the locus of $X_{\Jac}$ that corresponds to the eigenvalues of $\omega_{i}$ being the same. 
\subsection{Tripled Jordan quiver} \label{jordan}
We will now specialise to the case of the Jordan quiver $Q_{\operatorname{Jor}}$. This is the additive or Lie algebra version of the moduli stack of local systems on the $3$-torus. We can explicitly describe the stack of representations of the Jacobi algebra as $\MM_{\Jac,n} \cong \CC_{3}(\mathfrak{gl}_{n})/ \GL_n$, where $\CC_{3}(\mathfrak{gl}_{n})= \{ (x,y,z) \mid [x,y]=[x,z]=[y,z]=0 \}$ is the scheme of $3$ pairwise commuting matrices in $\mathfrak{gl}_{n}$. Similarly we can describe $X_{\Jac,n} \cong  \So^{n} \mathbb{G}_{a}^{3}$. As explained in \cite[Section 5]{davison2022integrality} we can use the support lemma \ref{support_lemma_quivers} $3$ times to deduce that the BPS sheaves $\BPS_{\Jac, n}$ must be supported on the diagonal $\mathbb{G}^{3}_{a} \to  \So^{n} \mathbb{G}_{a}^{3}$ and it is proven that $\BPS_{\Jac} = \bigoplus \Delta_{*} \mathbb{Q}_{\mathbb{G}^{3}_{a}}[3]$.
\begin{notation}
From now on denote the stack of representations of the Jacobi algebra of the Jordan quiver of dimension $n$ by $\MM_{\mathfrak{gl}_n}$, the good moduli space by $\pi_{\mathfrak{gl}_n} \colon \MM_{\mathfrak{gl}_n} \to X_{\mathfrak{gl}_n}$, the DT sheaf by $\varphi_{\mathfrak{gl}_n}$ and the BPS sheaf by $\BPS_{\mathfrak{gl}_n}$. We establish this notation to have analogous notation to the one we will use for the stack of local systems on the $3$-torus in Section \ref{section_coh_int}.
\end{notation}
 Note that the purity of $\pi_{Q*} \varphi_{\tr(W)}$ implies that this complex splits into a direct sum of $\IC$ sheaves supported on smooth locally closed subvarieties of $\So^{n} \mathbb{G}^{3}_{a}$. We will from now on forget the additional structure of mixed Hodge modules and just consider the underlying perverse sheaves.
Next we want to explicitly determine these sheaves in terms of the decomposition of equation \eqref{integrality_quivers}. We will do so by determining the perverse cohomology over each dimension $n$. Pick standard Levis inside $\GL_n$, they correspond to partitions of $n$. Let us call these Levis $L_{\lambda}  = \prod_{\lambda_{i} \in \lambda} \GL_{\lambda_i}$ with Lie algebra $\mathfrak{l}_{\lambda} = \oplus_{i} \mathfrak{gl}_{\lambda_{i}}$. Then also define the stack $M_{\mathfrak{l}_{\lambda}} = \CC_{3}(\mathfrak{l}_{\lambda}) / L_{\lambda}$ and $X_{\mathfrak{l}_{\lambda}} = \prod_{i} \So^{\lambda_{i}} \mathbb{G}^{3}_{a}$.
Let us define a stratification of $X_{\mathfrak{gl}_{n}}$ by setting
\begin{align}
    X^{\lambda}_{\mathfrak{gl}_n} & = \So^{n}_{\lambda} \mathbb{G}^{3}_{a} = \{ \sum^{l}_{i} \lambda_{i} x_{i} \mid \lambda_{i} \in \lambda =(\lambda_{1}, \dots, \lambda_{l}), x_{i} \neq x_{j} \in \mathbb{G}^{3}_{a} \} \\
 \text{with }\sum^{l}_{i} \lambda_{i} x_{i} & = \{ \underbrace{x_{1}, \dots , x_{1} , }_{\lambda_{1}\text{ times}}  \underbrace{x_{2}, \dots , x_{2} , }_{\lambda_{2}\text{ times}} \dots \underbrace{x_{l} , \dots , x_{l}}_{\lambda_{l}\text{ times}} \} \in \So^{n} \mathbb{G}^{3}_{a} . \nonumber
\end{align}
This is a locally closed smooth subscheme. We define the stack $M^{\lambda}_{\mathfrak{gl}_n} = M_{\mathfrak{gl}_n} \times_{X_{\mathfrak{gl}_{n}}} X^{\lambda}_{\mathfrak{gl}_{n}}$. \par
Note that the relative Weyl group $W_{L_{\lambda}}$ naturally acts on $\HHf(\B \ZZ(L_{G, \lambda}))$. We can split $\HHf(\B \ZZ(L_{G, \lambda}))$ by cohomological degree into subspaces $V_i$ $i \geq 0$. The natural $W_{L_{\lambda}}$ action on $\HHf(\B \ZZ(L_{G, \lambda}))$ perserves cohomological degree so each $V_{i}$ is a $W_{L_{\lambda}}$ subrepresentation.
\begin{lem}[BPS sheaves for Levis in $\mathfrak{gl}_n$]\label{additive_computation_lemma}
The following properties hold for BPS sheaves on Levis:
\begin{enumerate}
    \item $\pi_{\mathfrak{l}_{\lambda}*} \varphi_{\mathfrak{l}_{\lambda}}$ has perverse cohomology bounded below. We define $$\BPS_{\mathfrak{l}_{\lambda}} = \pH^{l} \pi_{\mathfrak{l}_{\lambda}*} \varphi_{\mathfrak{l}_{\lambda}}$$  here $l = \dim \ZZ(L_{\lambda})$. We have 
    $$\BPS_{l_{\lambda}} \cong \BPS_{\lambda_{1}} \boxtimes \cdots \boxtimes \BPS_{\lambda_{l}}.$$
    Furthermore, $\BPS_{\mathfrak{l}_{\lambda}}$ is a constant sheaf of rank $1$ supported on 
 $$\operatorname{supp}(\BPS_{\mathfrak{l}_{\lambda}}) = \im ( \Delta^{\lambda} \colon \ZZ^{3}(\mathfrak{l}_{\lambda}) \xhookrightarrow{} X_{\mathfrak{l}_{\lambda}}) .$$
 Finally, the components of the Saito decomposition of $\pi_{\mathfrak{l}_{\lambda}*} \varphi_{\mathfrak{l}_{\lambda}}$ with supports given by $\ZZ^{3}(\mathfrak{l}_{ \lambda})$ are $\BPS_{\mathfrak{l}_{\lambda}} \otimes \HHf(\B \ZZ(L_{ \lambda}))[- \dim \ZZ(L_{\lambda})]$.
 \item   Consider the subspace $V_{i}$ of cohomological degree $i$ in $\HHf(\B \ZZ(L_{ \lambda}))[- \dim \ZZ(L_{\lambda})]$.  The term 
 \begin{equation} \label{equiv_bps_add}
     \BPS_{\mathfrak{l}_{ \lambda}}  \otimes V_{i} 
 \end{equation}
  has a natural action of $W_{L_{\lambda}}$, which corresponds to the finite dimensional representation $V_{i}$. Pushing forward by $\theta \colon X_{L_{G,\lambda}} \to X_{G}$ and taking invariant part we get
  \begin{equation} \label{add_levi_ic_comp}
      (\theta_{*} \BPS_{\mathfrak{l}_{ \lambda}} 
 \otimes V_{i})^{W_{L_{\lambda}}} \cong \IC_{X^{\lambda}_{\gl_n}}(\mathcal{L}^{\lambda}_{i})[-i]
  \end{equation}
  where $\mathcal{L}^{\lambda}_{i}$ is some local system.
\end{enumerate}
 \end{lem}
 \begin{proof}
    By example \ref{dcrit_products} we can write $\varphi_{\mathfrak{l}_{\lambda}} \cong \varphi_{\mathfrak{gl}_{\lambda_{1}}} \boxtimes \dots \boxtimes \varphi_{\mathfrak{gl}_{\lambda_{l}}}$. Then the description of the BPS sheaf follows from the one for each $\mathfrak{gl}_{\lambda_{i}}$. We want to compute the pushforward $\theta_{*}( \Delta_{\lambda, *}\BPS_{\mathfrak{l}_{\lambda}} [ - \dim \ZZ(\mathfrak{l}_{\lambda})] \otimes\HHf(\B \ZZ (L_{ \lambda})))$. So we are pushing forward by the map $\ZZ^{3}(\mathfrak{l}_{\lambda}) \to X_{\mathfrak{gl}_{n}}$. The image of this map is the closure of $X^{\lambda}_{\mathfrak{gl}_{n}}$.   Therefore, we have the following pullback diagram 
\begin{equation}\label{computation_diagram}
\begin{tikzcd}
	{\widetilde{\ZZ}^{3}(\mathfrak{l}_{\lambda})} & {\ZZ^{3}(\mathfrak{l}_{\lambda})} \\
	{X^{\lambda}_{\mathfrak{gl}_{n}}} & {\overline{X}^{\lambda}_{\mathfrak{gl}_{n}}}
	\arrow["\theta", from=1-2, to=2-2]
	\arrow["j"', hook, from=2-1, to=2-2]
	\arrow["{j^{'}}", from=1-1, to=1-2]
	\arrow["{\theta^{'}}"', from=1-1, to=2-1]
\end{tikzcd}
\end{equation}
We know that $\ZZ^{3}(\mathfrak{l}_{\lambda}) = \prod^{l}_{i=1}\mathbb{C}^{3}_{}$. We can compute directly from the diagram above that $\widetilde{\ZZ}^{3}(\mathfrak{l}_{\lambda}) = \{ (x_{i}, \dots , x_{l}) \in \ZZ^{3}(\mathfrak{l}_{\lambda}) \mid x_{i} \neq x_{j} \text{ for } i \neq j \}$. This implies that the map from $\widetilde{\ZZ}^{3}(\mathfrak{l}_{\lambda}) \to X^{\lambda}_{\mathfrak{gl}_{n}}$ is a $W_{L_{\lambda}}$ cover.  \par
Let us prove that $\theta_{*}(\mathcal{L}[3\dim \ZZ(\mathfrak{l}_{\lambda})]) \cong \IC_{X^{\lambda}_{\mathfrak{gl}_{n}}}(\theta^{'}_{*} j^{'*} \mathcal{L})$ for $\mathcal{L}$ a local system on $\ZZ^{3}(\mathfrak{l}_{\lambda})$. Starting with
\begin{align*}
    \IC_{X^{\lambda}_{\mathfrak{gl}_{n}}}(\theta^{'}_{*} j^{'*} \mathcal{L}) & = \im (\pH^{0} j_{!} \theta^{'}_{*} j^{'*} \mathcal{L} \to \pH^{0} j_{*} \theta^{'}_{*} j^{'*} \mathcal{L}) \\
    & \cong \im (\pH^{0} \theta_{*} j^{'}_{!} j^{'*} \mathcal{L} \to \pH^{0} \theta_{*} j^{'}_{*}  j^{'*} \mathcal{L}) \qquad \text{(by finiteness of $\theta$)} \\
    & \cong \theta_{* } \im (\pH^{0}  j^{'}_{!} j^{'*} \mathcal{L} \to \pH^{0}  j^{'}_{*}  j^{'*} \mathcal{L}) \qquad \text{(by finiteness of $\theta$)}  \\
   & \cong \theta_{*} (\mathcal{L}[3\dim \ZZ(l_{\lambda})]) \qquad \text{($\IC_{U}(\mathcal{L}_{U}) \cong \mathcal{L}[\dim X]$ for $U \subseteq_{\text{open}} X$ smooth}) .
\end{align*}
Now there is a residual $W_{L_{\lambda}}$ action on the pushforward by $\theta$. Since the $W_{L_{\lambda}}$ action on $X^{\lambda}_{\mathfrak{gl}_{n}}$ is trivial we can decompose the pushforward into a direct sum of sheaves tensored by simple 
 $W_{L_{\lambda}}$ representations. In particular, we have $\theta^{'}_{*} j^{'*} \mathcal{L}[3\dim \ZZ(\mathfrak{l}_{\lambda})] = \bigoplus \mathcal{L}_{i} \otimes \rho_{i}$. Taking invariant parts we get
\begin{align*}
    (\theta_{*}(\mathcal{L}[3\dim \ZZ(\mathfrak{l}_{\lambda})]))^{W_{L_{\lambda}}} \cong (\IC_{X^{\lambda}_{\mathfrak{gl}_{n}}}(\theta^{'}_{*} j^{'*} \mathcal{L}))^{W_{L_{\lambda}}} & \cong (\IC_{X^{\lambda}_{\mathfrak{gl}_{n}}}(\bigoplus \mathcal{L}_{i} \otimes \rho_{i}))^{W_{L_{\lambda}}} \\
    & \cong (\bigoplus \IC_{X^{\lambda}_{\mathfrak{gl}_{n}}}(\mathcal{L}_{i}) \otimes \rho_{i})^{W_{L_{\lambda}}} \\
    & \cong \IC_{X^{\lambda}_{\mathfrak{gl}_{n}}}(\mathcal{L}_{\operatorname{triv}}) \otimes \operatorname{triv} \\
    & \cong  \IC_{X^{\lambda}_{\mathfrak{gl}_{n}}}((\theta^{'}_{*}j^{'*} \mathcal{L}[ 3\dim \ZZ(\mathfrak{l}_{\lambda})])^{W_{L_{\lambda}}}).
\end{align*}
Now we can use that for a principal $G$ bundle $\pi \colon P \to X$ with $G$ discrete $(\pi_{*}(V_{P}))^{G} = V_{X}$. Here $V_{P}$ is the constant sheaf on $P$ with fibre $V$. $V$ also has the structure of a $G$ representation, so $\pi_{*} V_{P}$ has a residual $G$ action. We apply this to the $W_{L_{\lambda}} \colon 1$ cover $\theta^{'}$.
 \end{proof}
 Now we give a different formulation of cohomological integrality that will be useful in Section \ref{section_coh_int} later. The next proposition is the restriction of Theorem \ref{coh_int_dm} for the tripled Jordan quiver to a fixed dimension vector $n$.
\begin{prop}\label{coh_add_prop}
     Cohomological integrality for the tripled Jordan quiver is equivalent to the following statement for all $n$.
    \begin{equation}\label{coh_add_gln}
        \pi_{\mathfrak{gl}_{n}*} \varphi_{\mathfrak{gl}_{n}} = \bigoplus_{L_{\lambda} \subseteq G} (\theta_{*}\BPS_{\mathfrak{l}_{\lambda}} 
 \otimes \HHf(\B \ZZ(L_{\lambda}))[- \dim \ZZ(L_{\lambda})])^{W_{L_{\lambda}}}
    \end{equation}
    Here the map $\theta \colon X_{\mathfrak{l}_{\lambda}} \to X_{\mathfrak{gl}_{n}}$ is induced by the inclusion $\mathfrak{l}_{\lambda} \to \mathfrak{gl}_{n}$ and $W_{L_{\lambda}}$ is the relative Weyl group.
\end{prop}
\begin{proof}
     Note that $W_{L_{\lambda}}$ is just a product of symmetric groups. Firstly we can define $W_{L_{\lambda}}$ invariant sheaves using idempotents as in \cite{mhm_sym}. Denote $\BPS_{n} \otimes \HHf(\B \mathbb{G}_{m})[-1]$ by $V_n$ and $V = \bigoplus_{n} V_{n}$. Then $\Sym V = \bigoplus_{k} (V^{\otimes k})^{S_k} $, where the tensor product is given by the monoidal structure $\boxdot$. The $n$-th graded part of $\Sym V$ can be written as
    \begin{equation}
          (\bigoplus_{\substack{(\lambda_{1}, \dots, \lambda_{l}) \in \mathbb{Z}^{l}_{> 0} \\ \sum \lambda_{i} = n}} V_{\lambda_{1}} \otimes \dots \otimes V_{\lambda_{l}})^{S_l} = \bigoplus_{\sum \lambda_{i} = n} (V_{\lambda_{1}} \otimes \dots \otimes V_{\lambda_{l}})^{S_{\lambda}}
    \end{equation}
    Here the direct sum on the left is over all tuples $(\lambda_{1}, \dots \lambda_{l})$ of some length $l \geq 1$, while the direct sum on the right is over partitions. $S_{\lambda}$ is the subgroup of $S_l$ that preserves the partition $\lambda = (\lambda_{1}, \dots , \lambda_{l})$. Note that $S_{\lambda}$ is exactly $W_{L_{\lambda}}$. Noting that 
    $$V_{\lambda_{1}} \otimes \cdots \otimes V_{\lambda_{l}} = \theta_{*} \BPS_{l_{\lambda}} \otimes \HHf(\B \ZZ(L_{\lambda}))[- \dim \ZZ(L_{\lambda})]$$ we are done.
\end{proof}
In light of the previous two results we have the following decomposition for $\pi_{\mathfrak{gl}_{n}*} \varphi_{\mathfrak{gl}_{n}} $
\begin{equation} \label{add_double_decomp}
    \pi_{\gl_{n}*} \varphi_{\gl} \cong \bigoplus_{i \geq 0} \bigoplus_{\lambda} \IC_{X^{\lambda}_{\gl_{n}}}(\mathcal{L}^{\lambda}_{i})
\end{equation}
where the second direct sum is over partitions of $n$.

\section{Stacks of local systems}
Let $M$ be a closed, connected, oriented $n$-manifold. For a reductive group $G$ we define the derived stack of local systems to be $\LocB_{G}(M) = \Map(M_{B}, \B G)$, with $M_{B}$ the Betti stack as in Example \ref{betti_stack}.  Using the AKSZ construction of Theorem \ref{aksz} we see that $\LocB_{G}(M)$ is $(2-n)$-shifted symplectic. Furthermore, if we consider a parabolic $P \subseteq G$ with Levi factor $L$ we have the $(2-n)$-shifted Lagrangian correspondence 
\begin{equation}\label{lagrag_og}
    \LocB_{G}(M) \xleftarrow{} \LocB_{P}(M) \to \LocB_{L}(M).
\end{equation}
The maps are induced from the $2$-shifted Lagrangian correspondence $\B G \xleftarrow{} \B P \to \B L$. This is proven in \cite[Lemma 3.4]{pavel_impl}. The purpose of this section is to establish some results about the structure of the stack of local systems we will need to use as well as define orientation data. \par We have the following well known description of the cotangent complex of the stack of local systems. Let $\mathcal{L}$ be a $G$-local system on $M$ and consider the adjoint action of $G$ on $\mathfrak{g}$ then we denote by $\ad_{\mathcal{L}} \mathfrak{g}$ the $\GL(\mathfrak{g})$-local system given by the composition $\pi_{1}(M) \to G \to \GL(\mathfrak{g})$. We can repeat the same construction for the coadjoint action on $\mathfrak{g}^{*}$ or other representations of $G$.
\begin{prop}[(Co)tangent complex of stack of local systems]
Using Proposition \ref{mapping_cot} we can express $\mathbb{L}_{\LocB_{G}(M)} \cong \pi_{\#} \ev^{*} \mathbb{L}_{\B G}$.  Let $\mathcal{L} \in \LocB_{G}(M)$ be a $G$-local system on $M$ corresponding to a $\mathbb{C}$-point. Then we have
\begin{equation}
    \mathbb{T}_{\LocB_{G}(M), \mathcal{L}} \cong \CC^{*}(M, \ad_{\mathcal{L}} \mathfrak{g}[1]) \quad \text{and } \quad \mathbb{L}_{\LocB_{G}(M), \mathcal{L}} \cong \CC_{*}(M, \ad_{\mathcal{L}} \mathfrak{g}^{*}[-1]).
\end{equation}
\end{prop}
Pick a CW structure on $M$. Denote by $A$ the set of cells of $M$.
\begin{defn}[Euler structure]
    An Euler structure on $M$ is a singular $1$-chain $\zeta$ such that 
    \begin{equation}
        d \zeta = \sum_{a \in A} (-1)^{\dim a} \alpha_{a}
    \end{equation}
where $\alpha_{a} \in a$. Two Euler structures $\zeta$ and $\eta$ with $d \zeta = \sum_{a \in A} (-1)^{\dim a} \alpha_{a}$ and $d \eta = \sum_{a \in A} (-1)^{\dim a} \beta_{a}$ are equivalent if for some paths $x_{a} \colon [0,1] \to a$ from $\alpha_{a}$ to $\beta_{a}$ the $1$-cycle $\zeta - \eta  + \sum_{a \in A} (-1)^{\dim a} x_{a}$ is a boundary.
\end{defn}
It can be shown that Euler structures exist if and only if $\chi(M) = 0$. Roughly speaking an Euler structure is a choice of paths to a base point from any cell in our cell decomposition. We will use this to more explicitly describe the cotangent complex of the stack of local systems.
\begin{ex}[Local systems on the circle]\label{s1_forms}
First let us describe the cotangent complex of $\Loc_{G}(S^{1})$. Let $G$ be a linear reductive group and fix a $G$ invariant nondegenerate symmetric bilinear form $(-,-)$ on $\mathfrak{g}$. Let us more explicitly describe what closed forms look like on $\XB = \LocB_{G}(S^{1}) = G/G$. Here $G$ acts on $G$ by right conjugation. Pick a standard cell structure of a $0$-cell $p$ and $1$-cell $a$. Pick an Euler structure by picking a path from the $1$-cell given anticlockwise with respect to the standard orientation of $S^{1}$. Let us describe the differential for the tangent complex. The Euler structures allows us to write down an explicit model for $\CC^{*}(S^{1},\mathcal{L})$ for a local system as a module over $\mathbb{C} \pi_{1}(S^{1}) \cong \mathbb{C}[t^{\pm 1}]$.
\begin{equation}
    \mathcal{L}_{p} \xrightarrow{t\cdot} \mathcal{L}_{p}
\end{equation}
Then pulling back the cotangent complex $\mathbb{L}_{\LocB_{G}(S^{1})} $ along the map $G \to G/G$ we get $\CC^{*}(S^{1}, \ad_{\mathcal{L}_{S^{1}}} \mathfrak{g}[1])$ where $\mathcal{L}_{S^{1}}$ is the universal local system given by the map $G \to G/G$
\begin{align*}
    \mathcal{O}_{G} \otimes \mathfrak{g} &\xrightarrow{X^{-1}vX-v} \mathcal{O}_{G} \otimes \mathfrak{g}  
\end{align*}
Where $v \in \mathfrak{g}$ and $X$ is an element of $\mathcal{O}(G)$ thinking of it as a matrix of variables satisfying the relations cutting out $G$ inside $\GL_n$.  \par
Let us describe the $1$-shifted symplectic structure on $\LocB_{G}(S^{1})$. We follow  \cite[Section 2.3]{pavel_qhr}.
Define the Maurer-Cartan forms $\theta , \overline{\theta} \in \Omega_{G} \otimes \mathfrak{g}$ by
\begin{equation}
    \iota_{v} \theta = (x \in G \mapsto \LL_{x^{-1},*} v_{x}), \quad  \iota_{v} \overline{\theta} = (x \in G \mapsto \RR_{x^{-1},*} v_{x})
\end{equation}
for a vector field $v \in \Gamma(G, \mathbb{T}_{G})$. Here $\RR$ and $\LL$ are the left and right multiplication maps and $\LL_{x^{-1},*}$ is pushforward of vector fields. Theorem \ref{forms_artin} gives that $2$ forms on $[G/G]$ are
\begin{equation}
    \bigwedge^{2} \Omega_{G} \oplus \Omega_{G} \otimes \mathfrak{g}^{*}[-1] \oplus \mathcal{O}_{G} \otimes \Sym^{2}(\mathfrak{g}^{*})[-2]
\end{equation}
Define a two form $\omega_{0}(y) = - \frac{1}{2} (\theta + \overline{\theta} , y)$ for all $y \in \mathfrak{g}$. In \cite{pavel_qhr} it is proven that $\omega_{0}$ is $d$ closed but it is not $d_{\dR}$ closed. We can define a three form $\omega_{1} =  \frac{1}{12} (\theta, [\theta, \theta])$. Then we can show that $d_{\dR} \omega_{0} + d \omega_{1}=0$.
\end{ex}
\begin{defn}[Moduli of framed local systems] \label{framed_locsys}
Picking a point $x \in M$ we can define a map $\LocB_{G}(M) \to \B G$ then we define $\LocB^{\ff}_{G}(M) \coloneqq \LocB_{G}(M) \times_{\B G} \pt$. We call this the \textbf{ moduli of framed local systems}. This allows us to write
\begin{equation}
    \LocB_{G}(M) \cong \LocB^{\ff}_{G}(M) / G.
\end{equation}
\end{defn}
In the next example we will describe the fibers of the tangent complex of $\LocB_{G}(T^{3})$.
\begin{ex}[(Co)tangent complex for the 3 torus]\label{cotangent_t3}
Let $G \subseteq \GL_n$ be a linear algebraic group. 
Let $M = T^{3}$, which we can view as a cube with opposite faces identified. Pick the standard cell structure with one $0$-cell $p$, three $1$-cells $a_{i}$, three $2$-cells $b_{i}$ and one $3$-cell $c$. Also pick an Euler structure, we pick a corner of the cube as our base point $p$. Then we pick straight line paths from the center of each adjacent face, edge and the center of the cube. This allows us to pick an explicit model for $\CC_{*}(T^{3}, \mathcal{L})$ for some local system $\mathcal{L}$. The Euler structure allows us to use the parallel transport of $\mathcal{L}$ to identify $\Gamma(a,\mathcal{L}|_a)$ with $\mathcal{L}_{p}$ for any cell $a$. Therefore, we can use a Koszul resolution of $\mathcal{L}_{p}$ as a $\mathbb{C}[t^{\pm 1}_{1}, t^{\pm 1}_{2},t^{\pm 1}_{3}]$ module to explicitly write down the differentials that appear in terms of the actions of the three monodromy operators $t_{1},t_{2},t_{3}$
\begin{equation}
\begin{tikzcd}
	{\Gamma(p,\mathcal{L}|_{p})} & {\oplus^{3}_{i=1}\Gamma(a_i,\mathcal{L}|_{a_{i}})} & {\oplus^{3}_{i=1}\Gamma(b_i,\mathcal{L}|_{b_{i}})} & {\oplus^{3}_{i=1}\Gamma(c,\mathcal{L}|_{c})} \\
	{\mathcal{L}_{p}} & {\mathcal{L}^{3}_{p}} & {\mathcal{L}^{3}_{p}} & {\mathcal{L}_{p}}
	\arrow["{d_{0}}", from=1-1, to=1-2]
	\arrow["\cong", from=1-1, to=2-1]
	\arrow["{d_{1}}", from=1-2, to=1-3]
	\arrow["\cong", from=1-2, to=2-2]
	\arrow["{d_{2}}", from=1-3, to=1-4]
	\arrow["\cong", from=1-3, to=2-3]
	\arrow["\cong", from=1-4, to=2-4]
	\arrow["{d_{0}}", from=2-1, to=2-2]
	\arrow["{d_{1}}", from=2-2, to=2-3]
	\arrow["{d_{2}}", from=2-3, to=2-4]
\end{tikzcd}
\end{equation}
with maps 
\begin{align} \label{differentials_t3}
    d_{0}(v) & = t_{1}v_{3} - t_{2} v_{2}+t_{3}v_{1} \\
    d_{1}(v_{1}) & = t_{1} v_{2} -t_{2}v_{1}  \nonumber \\
        d_{1}(v_{2}) & = t_{1} v_{3} -t_{3}v_{1} \nonumber \\
        d_{1}(v_{3}) & = t_{2} v_{3} -t_{3}v_{2} \nonumber \\
        d_{2}(v_{1})&= t_{1}v, \quad d_{2}(v_{2})= t_{2}v, \quad d_{2}(v_{3})= t_{3}v. \nonumber
\end{align} 
Here for $v \in \mathcal{L}_{p}$ we write $v_{1}= (v,0,0) \in \mathcal{L}^{3}_{p}$ and similarly for $v_{2},v_{3}$. Now we write down the fibers of the cotangent complex at a local system $\mathcal{L}$, specified by three matrices $(X_1,X_{2},X_{3}) \in G^{3}$. The differentials will be given by $d_{i}$ so we will just specify how the $t_{i}$ act.
\begin{align*}
     \CC_{*}(T^{3}, \ad_{\mathcal{L}} \mathfrak{g}^{*}[-1]) & \cong \mathfrak{g}^{*} \xrightarrow{d_{0}}  (\mathfrak{g}^{*})^{3}  \xrightarrow{d_{1}} (\mathfrak{g}^{*})^{3}  \xrightarrow{d_{2}}  \mathfrak{g}^{*} \\
     f \in \mathfrak{g}^{*}, \, t_{i}f & \mapsto \ad^{*}_{X_i}f -f 
\end{align*}
where $\ad^*$ is the coadjoint action.
\end{ex}
\subsection{Orientation data for 3-manifolds}
In this section we recall some theorems from \cite{naef2023torsion} that define orientation data on the stack of local systems of a $3$-manifold. In particular, we will use their theory of volume forms. After this we will be able to define the DT sheaf on $\LocB_{G}(M)$. 
\begin{defn}[Torsion]  \cite[Section 1.1]{absolute_torsion} \label{torsion}
Let $C$ be a bounded cochain complex $\mathbb{C}$-vector spaces. We call the \textbf{torsion of }$C$ the following isomorphism
\begin{equation}
   \psi(C) \colon \det   C  \to \det   \HHf(C) 
\end{equation}
To define $\psi(C)$ first pick $c \in \det C$.
 Choose for each $k=0,...,m$  
non-zero elements
$c_k\in \det \, C_k $ and $h_k\in \det  \HH^k(C)$. Set 
$c=c_0\otimes c_1^{-1}\otimes c_2\otimes \dots \otimes c_m^{(-1)^m}\in
\det \, C $ and  $h=h_0\otimes h_1^{-1}\otimes h_2\otimes \dots
\otimes h_m^{(-1)^m} \in \det \HHf(C)$.
We define $\psi(C)$ by
\begin{equation}
    \psi(C) (c) = (-1)^{N(C)}\, [c:h]\,h, 
\end{equation}
where  $ [c:h]$ is a
nonzero element of $\mathbb{C}$, defined by 
\begin{equation}
    [c:h] =  \prod_{k=0}^m
 \det [d(b_{k+1})\hat h_k b_k/\hat c_k]^{(-1)^{ k+1 }}
\end{equation}
Here:
\begin{align*}
    b_k & \text{ is a sequence of vectors of }C^k \text{ whose
 image }d_{k}(b_k) \text{
is a basis of } \im  d_{k}, \\
\hat h_k & \text{ is a sequence of
cycles in } C^k  \\
& \text{ such that the wedge product of their cohomology classes 
equals } h_k, \\
\hat c_k & \text{ is a basis of }
C^k \text{  whose wedge product  
equals } c_k, \\
[d(b_{k+1})\hat h_k b_k/\hat  c_k] & \text{ is the change of basis matrix
 from } \hat  c_k \text{ to the 
basis } d(b_{k+1})\hat h_k b_k \text{ of } C^{k}.
\end{align*}
$N(C)$  is defined by
$$N(C)=\sum_{k=0}^m\alpha_k(C)\beta_k(C) \, \modo\,2, $$
where
$$\alpha_k(C) = \sum_{j=0}^k \dim  C^j\, \modo\,2,
\quad  \beta_k(C) = \sum_{j=0}^k \dim  \HH^{j}(C)\,\modo\,2.$$
\end{defn}
\begin{defn}[Volume form] \label{naef_saf_vol}
    Let $\XB$ be a derived stack with perfect cotangent complex. Define the dimension of $\XB$ to be $\chi(\mathbb{L}_{\XB})$.
    \begin{enumerate}
        \item Then a \textbf{volume form} is an isomorphism $\mathcal{O}_{\XB} \to \det \mathbb{L}_{\XB}$.
        \item Let $f \colon \XB \to \YB$ be an an \'etale map of derived stacks equipped with volume forms $\omega_{\XB}$ and $\omega_{\YB}$. Then there is an induced pullback volume form $f^{*} \omega_{\XB}$ on $\YB$, which differs from $\omega_{\YB}$ by some invertible function $g_{f}$ induced by the quasi-isomorphism $\det (f^{*} \mathbb{L_{\XB}} \cong \mathbb{L}_{\YB})$.
        \item  Let $X$ be a $d$-critical stack. Then we call an isomorphism  $\mathcal{O}_{X_{\red}} \to K_{X}$ a volume form. Similarly for $X$ a complex analytic $d$-critical stack we call an isomorphism $\mathcal{O}_{X_{\red}} \to K_{X}$ a volume form.
    \end{enumerate}
\end{defn}
We will later use Definition \ref{torsion} to explicitly compare volume forms under \'etale maps. \par
Recall that if $\XB$ is $(-1)$-shifted symplectic with $X = \tt_{0} \XB$, then we have $K_{X} \cong \det \mathbb{L}_{\XB} |_{X_{\red}}$. Therefore, a volume form on $\XB$ induces a volume form on $X$ and thus an orientation on $X$. Then we can also induce a volume form and orientation on $X_{\an}$. \par If we have an \'etale $(-1)$-shifted symplectic map $f \colon \XB \to \YB$ of stacks with volume forms such that $g_{f}$ has a square root, then the map $f$ is oriented and we have $f^{*} \varphi_{\YB} \cong \varphi_{\XB}$. 
We now state a theorem for the existence of volume forms on mapping stacks. We will combine several theorems from \cite{naef2023torsion} and only state the parts we will use later.
\begin{prop} \label{orientation_pavel}
\begin{enumerate}
    \item Let $M$ be closed, connected oriented $n$-manifold and $\YB$ a derived stack of dimension $0$ with a volume form. Then $\Map(M, \YB)$ has a volume form. \cite[Theorem 2.8, Proposition 3.19]{naef2023torsion}.
    \item Let $M$ be a closed, connected oriented $3$-manifold and $G$ an algebraic group. If the Lie algebra $\mathfrak{g}$ of $G$ has a $G$-invariant volume form, then $\LocB_{G}(M)$  has a volume form. At a fiber $\mathcal{L} \in \LocB_{G}(M)$, this volume form is given by the Reidemeister torsion as defined in \cite{absolute_torsion}. Scaling the volume form on $\mathfrak{g}$ by a scalar $\lambda \in \mathbb{C}^{*}$ scales the volume form on $\Map(M, \B G)$ by $\lambda^{\chi(M)}$. \cite[Proposition 3.22, Proposition 3.23]{naef2023torsion}.
    \item Finally, the construction is functorial in the sense that if $f \colon \XB \to \YB$ is an \'etale map of derived stacks equipped with volume forms $\omega_{\XB}$ on $\XB$ and $\omega_{\YB}$ on $\YB$ such that $f^{*}\omega_{\YB} = \omega_{\XB}$ then we have an equality $\Map(f)^{*} \omega_{\Map(M, \YB)} = \omega_{\Map(M, \XB)}$ for the map $\Map(f) \colon \Map(M, \XB) \to \Map(M, \YB)$.   \cite[Theorem 3.28]{naef2023torsion}.
\end{enumerate}
\end{prop}
\begin{defn}[Modular character]
    Let $G$ be a an algebraic group. The \textbf{modular character} $\Delta_{G}$ is the character $G \xrightarrow{\operatorname{ad}} \GL(\mathfrak{g}) \xrightarrow{\det} \mathbb{G}_{m}$. A group $G$ is called unimodular if $\Delta_{G}  = 1$ or in other words the representation $\det \mathfrak{g}$ is trivial. It can be shown that any connected reductive group is unimodular.
\end{defn}
\begin{remark}
    A volume form on the Lie algebra exists if and only if $G$ is unimodular. While an Euler structure exists if and only if $\chi(M) = 0$.
\end{remark}
 In this paper we will be interested in the groups $G= \GL_{n}, \SL_{n}, \PGL_{n}$, for which a volume form on $\mathfrak{g}$ exists. Now we can use part $2$ of  Proposition \ref{orientation_pavel} and Theorem \ref{joyce_sheaf_stacks} to define the \textbf{DT sheaf} on $\LocB_{G}(M)$. 

\subsection{Local systems on the \texorpdfstring{$3$}{TEXT}-torus and group actions} \label{loc_stack_t3}
Let us specialise to the setting where $M = T^3$. From this point we will denote $\LocB_{G}(T^3)$ by $\LocB_G$ and the DT sheaf on $\LocB_{G}$ by $\varphi_{G}$. In this section we will examine how the DT sheaf defined in the previous subsection behaves under various group actions and Lagrangian correspondences.
We start by describing the good moduli spaces of $\Loc_{G}$, which are given by the affine GIT quotient $\Spec \mathcal{O}(\Loc_G)$.
\begin{warning}
Note that for $G=\SL_n$ or $\PGL_n$ it is generally not known if the good moduli space $\Spec \mathcal{O}(\Loc_{G})$ is $H^{3} /\!\!/ W$. In general, there is a map $H^{3} /\!\!/ W \to \Spec \mathcal{O}(\Loc_{G}) $, which is a bijection on $\mathbb{C}$-points. The difference between the two spaces is essentially a question about the reducedness of the ring of invariants. Therefore, apriori there is not even a map $\Loc_{G} \to H^{3} /\!\!/ W$. However, since the DT sheaf lives on the underlying complex topological space we abuse notation by identifying $\Loc_{G}$ and $\Loc_{G} \times_{ \Spec \mathcal{O}(\Loc_{G})} H^{3} /\!\!/ W$.
\end{warning}
With the above remark in mind we have the following definition 
\begin{defn} \label{good_moduli_loc}
    Let $G$ be a connected reductive group. Set $X_{G} = H^{3}_{G} /\!\!/ W$, where $H_{G}$ is the maximal torus of $G$.  We will denote the good moduli space map by $\pi_{G} \colon \Loc_{G} \to X_{G}$.
\end{defn}
\begin{notation}
    We will only consider good moduli spaces for trivial components of $\Loc_{G}$. If $\Loc_{G}$ has components, for clarity we will sometimes write $X^{1}_{G}$ to emphasize this. Similarly we will then write $\pi_{G} \colon \Loc^{1}_{G} \to X^{1}_{G}$.
\end{notation}
  In this paper we will mainly be interested in $G = \GL_{n}, \SL_{n}$ or $\PGL_{n}$ or one of their Levi subgroups. We can think of elements of this space as triples of diagonal matrices up to permutation by the Weyl group. Let us also fix for the rest of the paper unless stated otherwise, a non-degenerate  $\GL_n$-invariant symmetric bilinear form given by $(x,y) = \tr(xy)$ for $x,y \in \mathfrak{gl}_{n}$. The same pairing restricts to a non-degenerate invariant pairing  $\mathfrak{sl}_{n}$. Let us describe the Levi subgroups of $\GL_n, \SL_n, \PGL_n$.
 \begin{defn}[Levi subgroups] \label{levi_partitions}
Let $\lambda$ be a partition of $n$ of length $l$. We pick standard Levis $L_{\GL_{n}, \lambda} $. 
 \begin{enumerate}
     \item $\mathbf{\GL}_n$: we set $L_{\GL_{n}, \lambda} = \prod_{i} \GL_{\lambda_{i}}$, with Lie algebra $\mathfrak{l}_{\mathfrak{gl}_{n}, \lambda} = \oplus^{l}_{i=1} \mathfrak{gl}_{\lambda_{i}}$.
     \item $\mathbf{\SL}_n $: we set $L_{\SL_{n}, \lambda} = L_{\GL_{n},\lambda} \cap \SL_{n}$  with Lie algebra $\mathfrak{l}_{\mathfrak{sl}_{n},\lambda} = \mathfrak{l}_{\mathfrak{gl}_{n},\lambda} \cap \mathfrak{sl}_n$.
     \item $\mathbf{\PGL}_n$: we set $L_{\PGL_{n}, \lambda} = L_{\SL_{n}, \lambda} / \mu_{n}$ or $L_{\GL_{n}, \lambda} / \mathbb{G}_{m}$, with Lie algebra $\mathfrak{l}_{\mathfrak{pgl}_{n},\lambda} = \mathfrak{l}_{\mathfrak{sl}_{n},\lambda}$.
 \end{enumerate}
 \end{defn}
 The trace pairing on $\mathfrak{gl}_{n}$ restricts to the Levi Lie algebras. Associated to the Levis we have the relative Weyl group $W_{L_{\lambda}} = \NN(L_{G,\lambda})/ L_{G,\lambda}$ with $\NN(L_{G,\lambda})$ being the normaliser in $G$. We can show that $W_{L_{\lambda}}  \cong \prod^{l_{\gamma}}_{i = 1} S_{n_{i}}$, where $n_{i}$ is the number of times that the number $\lambda_{i}$ is repeated in the partition $\lambda$. So the Weyl group of the Levi is given by a product of symmetric groups of size the number of repeated blocks of the same dimension. We will denote relative Weyl groups by $W_{L_{\lambda}}$ since for all three groups they are isomorphic.
\begin{defn}[Actions]
Let $G$ be an algebraic group and $L \subseteq G$ a Levi subgroup.
\begin{enumerate}
    \item \textbf{Relative Weyl group actions}: \\ 
    $W_{L}$ naturally acts on $\LocB_{L}$ in the following way. Firstly, $\NN_{G}(L)$ acts on $L_{G, \lambda}$ by conjugation, then $L \subseteq \NN_{G}(L)$ acts by inner automorphisms thus acts trivially on $\B L$. Therefore, the action descends to $W_{L} = N_{G}(L)/L$ and  we can construct isomorphisms $\B L \to \B L$ for each element of $W_{L}$. Finally we can construct isomorphisms $\LocB_{L_{G, \lambda}} \to \LocB_{L}$ by taking mapping stacks $\Map(T^{3},-)$. Explicitly this map is given by sending an $S$ valued point $T^{3} \times S \to \B L$ to $T^{3} \times S \to \B L \xrightarrow{g} \B L$.
    \item \textbf{Center actions}: \\ 
    We have an action of $\LocB_{\ZZ(G)}$ on $\LocB_{G}$ given by taking mapping stacks of the action $\B \ZZ(G) \times \B G \to \B G$. 
\end{enumerate}
\end{defn}
Since $\ZZ(G)$ is abelian we have that $\Loc_{\ZZ(G)} \cong \ZZ^{3}(G) \times \B \ZZ(G)$. By restriction we then get an action of $\Loc^{\ff}_{\ZZ(G)} \cong \ZZ^{3}(G)$. In the case $G = \GL_n$, the action of $\LocB_{\mathbb{G}_{m}}$ can be seen as tensoring by a rank $1$ local system. For $G = \GL_n$ we can split the normaliser short exact sequence
\begin{equation}
    1 \to L_{\GL_n, \lambda} \to \NN(L_{\GL_n, \lambda}) \to W_{L_{\lambda}} \to 1
\end{equation}
by defining the map $W_{L_{\lambda}} \to \NN(L_{\GL_n , \lambda})$ by sending a permutation to its corresponding permutation matrix. This gives an isomorphism $\NN(L_{\GL_n, \lambda}) \cong L_{\GL_{n},\lambda} \rtimes W_{L_{\lambda}}$ and so we get an action of $W_{L_{\lambda}}$ on $L_{\GL_n , \lambda}$ by conjugation. Since the corresponding Levi for $\PGL_n$ is a quotient, we also get an action of $W_{L_{\lambda}}$. When $G = \SL_n$ we can view the Levi $L_{\SL_n , \lambda}$ as a subgroup of $L_{\GL_n , \lambda}$ and the $W_{L_{\lambda}}$ action restricts since conjugation does not change the determinant. Note however, that the normaliser exact sequence for $\SL_n$ in general \emph{will not} split. 
\begin{remark}[Equivariant sheaves with respect to group stacks]
Recall that $\Loc_{\ZZ(G)} = \ZZ^{3}(G) \times \B \ZZ(G)$. Via the next theorem by restriction we will get an action of $\ZZ^{3}(G)$ and the group stack $\B \ZZ(G)$ on the DT sheaf of $\Loc_{G}$. This means that there is an action of $\HHf(\B \ZZ(G))$ on $\varphi_{G}$. In particular, for $\SL_n$ we have that $\ZZ(\SL_n) = \mu_n$ and $\HHf(\B \mu_{n}) \cong \mathbb{Q}$ so in this case we do not get any extra structure.
\end{remark}
\begin{prop} \label{center_Wact}
    Let $G$ be a connected reductive group. The $\LocB_{\ZZ(G)}$-action on $\LocB_{G}$ preserves the $(-1)$-shifted symplectic structure and is oriented. By restriction there is an action of $\ZZ^{3}(G)$ on the DT sheaf $\varphi_{G}$. \par 
    Now let $G = \GL_n , \SL_n $ or $\PGL_n$.
    The $W_{L_{\lambda}}$ action on $\LocB_{L_{G,\lambda}}$ is $(-1)$ symplectic and oriented. Finally, there is an action of $W_{L_{\lambda}}$ on $\varphi_{L_{G, \lambda}}$.
\end{prop}
\begin{proof}
Consider the map $\LocB_{G} \to \LocB_{G/ \ZZ(G)}$. Note that there is a natural map $\LocB_{G / \ZZ(G)} \to \LocB_{\B \ZZ(G)}$ induced from the map $\B(G / \ZZ(G)) \cong \B G / \B \ZZ(G) \to \B^{2} \ZZ(G)$. Namely, to every $G/ \ZZ(G)$-local system we can associate a local system of $\ZZ(G)$ gerbes. Then we get diagram where both squares are pullbacks 
\begin{equation}
\begin{tikzcd}
	{\LocB_{G}} & \pt \\
	{\LocB^{1}_{G / \ZZ(G)}} & {\B \LocB_{\ZZ(G)}} \\
	{\LocB_{G / \ZZ(G)}} & {\LocB_{\B \ZZ(G)}}
	\arrow[from=1-1, to=1-2]
	\arrow[from=1-1, to=2-1]
	\arrow[from=1-2, to=2-2]
	\arrow[from=2-1, to=2-2]
	\arrow[from=2-1, to=3-1]
	\arrow[from=2-2, to=3-2]
	\arrow[from=3-1, to=3-2]
\end{tikzcd}
\end{equation}
Here $\LocB^{1}_{G/ \ZZ(G)}$ is the component of the trivial $G/ \ZZ(G)$ local system and is given by pulling back along the trivial local system of gerbes. The outer square is a pullback diagram because of the equations  
\begin{align*}
    \Map(T^{3}, \B (G / \ZZ(G))) \times_{ \Map(T^{3}, \B^{2} \ZZ(G))} \pt  & \cong \Map(T^{3}, \B (G / \ZZ(G)) \times_{\B^{2} \ZZ(G)} \pt)  \\
     &  \cong \Map(T^{3}, \B G) \quad \text{(because $\B G  \cong \B (G / \ZZ(G) ) \times_{ \B^{2} \ZZ(G)} \pt$  )}.  
\end{align*}
By the $2$ out of $3$ property for pullbacks the upper square is a pullback as well. Therefore, the map  $\LocB_{G} \to \LocB^{1}_{G/ \ZZ(G)}$ is an effective epimorphism since it is a pullback of the effective epimorphism $\pt \to \B \LocB_{\ZZ(G)}$. Taking \v{C}ech nerve we get the action groupoid of $\LocB_{\ZZ(G)}$ on $\LocB_{G}$ using that mapping stacks commute with limits in the target and that $\B G \times_{ \B (G/ \ZZ(G))} \B G = \B \ZZ(G) \times \B G$. This allows us to identify $\LocB^{1}_{G/ \ZZ(G)}$ with $\LocB_{G}/ \LocB_{\ZZ(G)}$. This implies that $\LocB_{G} / \LocB_{\ZZ(G)}$ is $(-1)$-symplectic and oriented since $\LocB^{1}_{G/ \ZZ(G)}$ is. Therefore, there is a DT sheaf on the quotient stack $\LocB_{G} / \LocB_{\ZZ(G)}$ and so $\varphi_{G}$ is $\ZZ^{3}(G)$ equivariant by restriction. \par
Let $G = \GL_n, \SL_n$ or $ \PGL_n$. Let us consider $W_{L_{\lambda}}$ acting on $L_{G,\lambda}$, then each element $\sigma \in W_{L_{\lambda}}$ gives an automorphism $\B L_{G,\lambda} \to \B L_{G,\lambda}$. We wish to show this map preserves $2$-shifted symplectic structures. The induced map on quasicoherent sheaves is the restriction functor $\Rep L_{G, \lambda} \to \Rep L_{G, \lambda}$. Note that the $2$-shifted symplectic structure on $\B L_{G,\lambda}$ is given by an isomorphism of $L$ representations $\mathfrak{l}_{\mathfrak{g}, \lambda} \to \mathfrak{l}^{*}_{\mathfrak{g},\lambda}$. We need to show that this map is also $W_{L_{ \lambda}}$-equivariant. But this follows from the fact that permuting the blocks in a Levi does not affect the trace of the matrix. Note that $\sigma$ acts as a permutation on $\mathfrak{l}_{\mathfrak{g}, \lambda}$ so it has determinant $\pm 1$. Therefore, $\sigma$ scales the volume form on $\B L_{G,\lambda}$ by $\pm 1$. By Proposition \ref{orientation_pavel} the induced map scales the volume form on $\Loc_{L_{G,\lambda}}$ by $(\pm 1)^{\chi(T^{3}) }= 1$ since $\chi(T^{3}) = 0$. We can now conclude that the action of $W_{L_{\lambda}}$ is oriented. The oriented $(-1)$-shifted symplectic automorphisms then define an action of $W_{L_{\lambda}}$ on $\varphi_{L_{G,\lambda}}$.
\end{proof}
$\ZZ^{3}(G)$ also acts on the good moduli space and the good moduli space map is equivariant with respect to this action. The action is well defined since the action of $\ZZ^{3}(G)$ on $H^{3}_{G}$ commutes with the action of $W$ on $H^{3}_{G}$. Since the map $\pi_G$ is $\ZZ^{3}(G)$ equivariant we also have an action on $\pi_{G*} \varphi_{G}$ and also on each perverse cohomology of $\pi_{G*} \varphi_{G}$. \par
 Because $\PGL_{n}$ is no longer simply connected $\LocB_{\PGL_{n}}$ will have several connected components with $\pi_{0} \LocB_{\PGL_n} \cong \mu^{3}_{n}$. Let $\zeta = (\zeta_{1}, \zeta_{2}, \zeta_{3}) \in \mu^{3}_{n}$. Write $\omega$ for the generator of $\mu_{n}$. The component $\LocB^{\zeta}_{\PGL_{n}}$ is defined to be the component where we have $(A,B,C) \in \Loc_{\PGL_{n}}$ with $[\widetilde{A}, \widetilde{B}] = \zeta_{1}$, $[\widetilde{A}, \widetilde{C}] = \zeta_{2}$ and $[\widetilde{B}, \widetilde{C}] = \zeta_{3}$. Here $(\widetilde{A}, \widetilde{B}, \widetilde{C})$ are representatives in $\GL_n$. In general, each component $\Loc^{\zeta}_{\PGL_{n}}$ is $(-1)$ symplectic and has its own DT sheaf. When $n$ is a prime we have a particularly simple description.
\begin{lem} \label{prime_pgln}
Let $n$ be prime. Then $\Loc^{\zeta}_{\PGL_{n}} \cong \B \mu^{2}_{n}$ where $\zeta = (\zeta_{1}, \zeta_{2}, \zeta_{3}) \in \mu^{3}_{n} - (1,1,1)$.   
\end{lem}
\begin{proof}
    Recall that we have a natural map $\LocB_{\PGL_{n}} \to \LocB_{\B \mu_{n}} \cong \HH^{2}(T^{3}, \mu_{n}) \times \B \mu_{n} \times \B^{2} \mu_{n} $. The latter isomorphism follows from computing the triple loop space $LLL \B \mu_{n} \cong \Loc_{\B \mu_{n}}(T^3)$ 
    \begin{align*}
        LLL \B \mu_{n} & \cong LL (\B \mu_{n} \times \B^{2} \mu_{n}) \\
        & \cong L( \mu_{n} \times \B \mu_{n} \times \B \mu_{n} \times \B^{2} \mu_{n}  ) \\
        & \cong \mu^{3}_{n} \times (\B \mu)^{3}_{n} \times \B^{2} \mu_{n} \\
        & \cong \HH^{2}(T^{3}, \mu_{n}) \times (\B \mu_{n})^{3} \times \B^{2} \mu_{n}.
    \end{align*}
    Here we have repeatedly used that $L \B G = G \times \B G$ for $G$ an abelian group and that loop spaces of products are products of the loop spaces. In particular, the component $\LocB^{\zeta}_{\PGL_{n}}$ is the preimage of $\zeta \in \HH^{2}(T^{3}, \mu_{n}) \cong \mu^{3}_{n}$. The mapping class group of $T^{3}$, $\Mod(T^{3}) \cong \SL_{3} (\mathbb{Z})$ acts by automorphisms on $\LocB_{\PGL_{n}}$. Fix a generator $\omega \in \mu_{n}$ and an isomorphism $\mu_{n} \to \mathbb{Z}/ n \mathbb{Z}$ given by $\zeta^{m} \mapsto m$. For any $\zeta \in \HH^{2}(T^{3}, \mu_{n})$ there is a mapping class $\gamma \in \Mod(T^{3})$ that maps 
    $$\zeta = (\omega^{n_{1}} , \omega^{n_{2}} , \omega^{n_{3}}) \to  \widetilde{\zeta} = (\omega^{gcd(n_{1}, n_{2} , n_{3})}, 1 , 1 ).$$ This claim follows since we can perform the Euclidean algorithm on $(n_{1},n_{2},n_{3})$ by multiplication with matrices that add multiples of one row to another, which are in $\SL_3 (\mathbb{Z})$. Therefore, for any $\zeta$ there is an automorphism $\gamma \colon \LocB_{\PGL_{n}} \to \LocB_{\PGL_{n}}$ that identifies $\LocB^{\zeta}_{\PGL_{n}}$ with $\LocB^{\widetilde{\zeta}}_{\PGL_{n}}$. Therefore, it is enough to determine $\LocB^{\zeta}_{\PGL_{n}}$ for $\zeta$ of the form $(\omega^{m},1,1)$ for some $1 \leq m < n$. Let $(A,B,C) $ be a commuting triple in $\PGL_{n}$ and $(\widetilde{A}, \widetilde{B}, \widetilde{C})$ their representatives in $\GL_{n}$. Without loss of generality assume $[\widetilde{A}, \widetilde{B}] = \zeta_{1} \neq 1$. Then by \cite[Lemma 2.2.15]{hausel_mixed} we have that up to conjugation $\widetilde{A} = \gamma_{1} \cdot \operatorname{diag}(1, \zeta_{1} , \cdots \zeta^{n-1}_{1})$ and $\widetilde{B} = \gamma_{2} P$. Here $P$ is the permutation matrix of the cycle $(12 \dots n)$ and $\gamma_{1} , \gamma_{2} \in \mathbb{G}_{m}$. Now if $[\widetilde{A}, \widetilde{C}] = [\widetilde{B}, \widetilde{C}] = 1$, then by \cite[Lemma 2.2.6]{hausel_mixed} $\widetilde{C}$ must be central. Therefore, in $\PGL_{n}$ we have a unique such point up to conjugation. Let us now compute the stabiliser of this point. This is an element $D$ in $\PGL_{n}$, which commutes with $(A,B,C)$. Then there are a priori $\mu^{3}_{n}$ choices for commutators of the representatives $(\widetilde{A}, \widetilde{B}, \widetilde{C})$ and $\widetilde{D}$. However, since $\widetilde{C}$ is central this reduces the number of choices to $\mu^{2}_{n}$. This implies that the stabiliser is $\mu^{2}_{n}$. 
\end{proof}
We now introduce some generic loci of our moduli spaces which we will need to be able to compare DT sheaves on $\LocB_{G}$ and $\LocB_{L_{G,\lambda}}$.
\begin{defn}[Generic loci] \label{generic_loci} \,
\begin{enumerate}
    \item \textbf{Generic locus for} $\GL_n$: \\ Say that an element $x \in L_{\GL_{n},\lambda} = \prod^{l}_{i = 1} \GL_{\lambda_{i}}$ is \textbf{generic} if the eigenvalues of the $\GL_{\lambda_{i}}$ blocks are pairwise distinct. \par
Define an open subvariety $X^{g}_{L_{\GL_{n},\lambda}} \subseteq X_{L_{\GL_n ,\lambda}}$ by the condition that for any $x = ( A_{1},A_{2},A_{3}) \in X^{g}_{L_{_{\GL_{n}}\lambda}}$ at least one of the $A_{i} $ $1 \leq i \leq 3$ is generic. \par 
Define $\Loc^{g}_{L_{\GL_{n},\lambda}} = \Loc_{L_{\GL_{n},\lambda}} \times_{X_{L_{\GL_{n},\lambda}}} X^{g}_{L_{\GL_{n},\lambda}}$. This is an open substack of $\Loc_{L_{\GL_{n},\lambda}}$.
\item \textbf{Generic locus for} $\SL_n$: \\
First, define an open subvariety $X^{g}_{L_{\SL_{n}, \lambda}} = X_{L_{\SL_{n}, \lambda}} \cap X^{g}_{L_{\GL_n}}$. \par  
Define $\Loc^{g}_{L_{\SL_{n},\lambda}} = \Loc_{L_{\SL_{n},\lambda}} \times_{X_{L_{\SL_{n}, \lambda}}} X^{g}_{L_{\SL_{n}, \lambda}}$. This is an open substack of $\Loc_{L_{\SL_{n},\lambda}}$. 
\end{enumerate}
\end{defn}
The inclusion $L_{G, \lambda} \subseteq G$ induces a map $\Theta \colon \LocB^{1}_{L_{G, \lambda}} \to \LocB^{1}_{G}$.
\begin{prop}\label{lag_corresp_gen}
 Let $G= \GL_n$ or $\SL_n$. The map $\Theta \colon \LocB_{L_{G, \lambda}}  \to \LocB_{G}$ preserves symplectic forms. Restricted to the generic locus $\Theta^{g} \colon \LocB^{g}_{L_{G, \lambda}}  \to \LocB_{G} $ is \'etale and oriented. The $W_{L_{\lambda}}$ action restricts to $\LocB^{g}_{L_{G , \lambda}}$ and is also oriented.\end{prop}\begin{proof}
The fact that $\Theta$ preserves symplectic forms follows from the fact that the map $\B L_{G, \lambda} \to \B G $ preserves closed $2$-forms because the pairing on $\mathfrak{g}$ restricts to the pairing on $\mathfrak{l}_{\mathfrak{g} , \lambda}$. Using Theorem \ref{aksz} we can then deduce the statement for mapping stacks. \par 
To show that the map is \'etale we can prove that $\Theta^{g}(S^{1}) \colon \LocB^{g}_{L_{G, \lambda}}(S^{1})  \to \LocB_{G}(S^{1}) $ is \'etale and then take the mapping stack $\Map(T^2,-)$. \par To start proving that $\Theta^{g}(S^{1})$ is \'etale, decompose the Lie algebra  $\mathfrak{g}$ as $\mathfrak{u}_{-} \oplus \mathfrak{l} \oplus  \mathfrak{u}_{+}$. Here $\mathfrak{u}_{\pm}$ is the Lie algebra of the positive/negative unipotent radical of the parabolic corresponding to $L_{G, \lambda}$. Note that as an $L_{G, \lambda}$ representation under conjugation  $\mathfrak{g}$ splits as $\mathfrak{l}$ and $\mathfrak{u}_{-} \oplus \mathfrak{u}_{+}$. The sequence $ \mathfrak{u}_{+} \oplus \mathfrak{u}_{-} \to \mathfrak{l}[1] \to \mathfrak{g}[1] $ is the tangent sequence of the map $\B L_{G,\lambda} \to \B G$.  Taking mapping stacks we get the following tangent sequence for $\Loc_{L_{G,\lambda}}(S^{1}) \to \Loc_{G}(S^{1})$ with $\mathcal{L}_{S^{1}}$ the universal local system as in Example \ref{s1_forms}
\begin{equation}
    \CC^{*}(S^{1} , \ad_{\mathcal{L}_{S^{1}}} \mathfrak{u}_{-} \oplus \mathfrak{u}_{+}) \to \CC^{*}(S^{1} , \ad_{\mathcal{L}_{S^{1}}} \mathfrak{l}[1]) \to \CC^{*}(S^{1} , \ad_{\mathcal{L}_{S^{1}}}\mathfrak{g}[1])
\end{equation}
 To prove that the map is \'etale on the generic locus it is enough to show that  $\CC^{*}(S^{1} , \ad_{\mathcal{L}_{S^{1}}}\mathfrak{u}_{-} \oplus \mathfrak{u}_{+})$ is acyclic. To show that this complex is acyclic, recall the differential in Example \ref{s1_forms}. Since the differential is a map between free modules of the same rank, it is enough to show that it is injective, as it will then be an isomorphism. If the differential is an isomorphism then the complex has no cohomology. The differential is given by $v \mapsto X_{f}^{-1}v X_{f} -v$ with $X_{f}$ satisfying the genericity condition. Therefore, the eigenvalues in the blocks labelled by $\lambda_i$ are distinct. If this map has a kernel, then $X_{f}$ commutes with a matrix $v \in \mathfrak{u}_{+} \oplus \mathfrak{u}_{-}$ having a non-trivial Jordan block, which is a contradiction.
\par
Finally we prove that the map $\Theta^{g} \colon \LocB^{g}_{L_{G, \lambda}} \to \LocB_{G}$ is oriented.
Now we have two volume forms $\omega_{L}$ and $\omega_{G}$, which we want to compare. We have that $\omega_{L} = f\Theta^{*} \omega_{G}$ where $f$ is the function induced on determinants from the isomorphism $\Theta^{*} \mathbb{L}_{\LocB_{G}} \cong \mathbb{L}_{\LocB^{g}_{L_{G, \lambda}}}$. Note that the function $f$ is in degree $0$, therefore to compute it is enough to consider the truncation $\Loc^{g}_{L_{G, \lambda}}$. Furthermore, we only need to compare the volume form on the reduced locus to check that that $\Theta^{g}$ is a map of oriented $d$-critical stacks. This implies that we can compare the values at closed points of $\Loc^{g}_{L_{G, \lambda}}$. To compute $f$, consider the fiber sequence $\Theta^{g*} \mathbb{L}_{\LocB_{G}} \to  \mathbb{L}_{\LocB^{g}_{L_{G, \lambda}}} \to \mathbb{L}_{\LocB^{g}_{L_{G, \lambda}}/ \LocB_{G}}$. Upon taking determinants we get isomorphisms
    \begin{equation}
         \det \mathbb{L}_{\LocB^{g}_{L_{G, \lambda}}} \cong  \Theta^{g*}\det \mathbb{L}_{\LocB_{G}} \otimes \det \mathbb{L}_{\LocB^{g}_{L_{G, \lambda}} / \LocB_{G}} \cong  \Theta^{*g}\det \mathbb{L}_{\LocB_{G}} \otimes \mathcal{O}_{\LocB^{g}_{L_{G, \lambda}}}
    \end{equation}
     Since the map is \'etale $\mathbb{L}_{\LocB^{g}_{L_{G, \lambda}} / \LocB_{G}}$ is acyclic and thus we have $\det \HHf (\mathbb{L}_{\LocB^{g}_{L_{G, \lambda}} / \LocB_{G}}) \cong \mathcal{O}_{\LocB^{g}_{L_{G, \lambda}}}$.  

The function $f$ evaluated at a local system $\mathcal{L} \in \Loc_{L_{G, \lambda}}$ will then be given by the torsion (see Definition \ref{torsion})  
$$\psi(\CC_{*}(T^{3} , \ad_{\mathcal{L}}\mathfrak{u}^{*}_{-} \oplus \mathfrak{u}^{*}_{+}[-1]))$$ since it is induced by the isomorphism
$$\det \mathbb{L}_{\LocB^{g}_{L_{G, \lambda}}, \mathcal{L}} \cong \det \HHf (\mathbb{L}_{\LocB^{g}_{L_{G, \lambda}} / \LocB_{G}})_{\mathcal{L}} \cong \mathbb{C}.$$
Recall the description of the cotangent complex of local systems on $T^3$ in Example \ref{cotangent_t3}.
 Note that the signs $N((\CC_{*}(T^{3} , \ad_{\mathcal{L}}\mathfrak{u}^{*}_{-} \oplus \mathfrak{u}^{*}_{+}[-1]))$ are just $1$ because the sums defining them are both trivial using the description in Example \ref{cotangent_t3}. Furthermore, since we are working with an acyclic complex we can just pick all the $h$ terms to be $1$ and so we do not have to consider them in the calculation.  Let $C =  \mathfrak{u}^{*}_{-} \oplus \mathfrak{u}^{*}_{+}$. Let us denote the three monodromy operators that appear in the definition of the differentials by $t_{1},t_{2},t_{3}$ and the three matrices defining the monodromy of $\mathcal{L}$ by $X_{1},X_{2},X_{3}$. In particular, we view $t_{i}$ as the operator $ C \to  C$ given by $ f \mapsto \ad^{*}_{X_i}f   - f$. Let us write our cochain complex as 
\begin{equation*}
     C_{3}  \xrightarrow{d_{1}} C_{2} \xrightarrow{d_{2}} C_{1} \xrightarrow{d_{3}} C_{0} \\
\end{equation*}
Start by picking a basis $\hat{c}_{3}$ of $C$. Now let us pick a basis for $C_i$ and the elements $b$
\begin{align*}
    C_{3} & = C \text{ with basis } \hat{c}_{3}  \text{ and } b_{3} = \hat{c}_3 \\
    C_{2} & = C^{3} \text{ with basis } \hat c_2 = \{ (\hat c_{3},0,0), (0,\hat c_{3},0),  (0,0,\hat c_{3}) \} \text{ and } b_{2} = \{ (0, \hat c_{3} , 0 ), (0, 0,  \hat c_{3} ) \}  \\
    C_{1} & = C^{3} \text{ with basis } \hat c_1 = \{ (\hat c_{3},0,0), (0,\hat c_{3},0),  (0,0,\hat c_{3}) \} \text{ and }  b_{2} = \{ (0, \hat c_{3} , 0 ) \} \\
    C_{0} & = C \text{ with basis }\hat c_{0} = \hat c_{3}. 
\end{align*}
 Recall that to compute the torsion we have to compute the determinant of certain change of basis matrices defined in terms of the differential. We do this by computing the following block matrices
    $$[d(b_3) b_{2}: \hat{c}_2]  = \begin{bmatrix}
    t_{3} & 0 & 0 \\
     -t_{2} & I & 0 \\
     t_{1} & 0 & I
  \end{bmatrix},$$  
  $$[d(b_2) b_{1}: \hat{c}_1]  = \begin{bmatrix}
    -t_{3} & 0 & I \\
     0 & -t_{3} & 0 \\
     t_{1} & t_{2} & 0
  \end{bmatrix},$$  
  $$ 
  \det  [d(b_3) b_{2}: \hat{c}_2]  = \det t_{3} , \qquad \det  [d(b_2) b_{1}: \hat{c}_1] = \det t_{3} \det t_{1}, \qquad \det [ b_{0}: \hat{c}_1] =  \det t_{1}.$$
Here we have computed using the determinant formula in \cite[Section 4.2]{block_det}. Hence we can see that taking alternating multiplication we get that the torsion is just $1$ and the map is oriented.
\end{proof}
\begin{remark}
    The idea for the above proposition comes from the following observation. Recall that we have for each Levi a $(-1)$-shifted Lagrangian correspondence \eqref{lagrag_og}. The genericity condition above is designed to ensure that restricted to this locus the map $\LocB_{P} \to \LocB_{L_{\lambda}}$ is an isomorphism. More specifically we have the following diagram
\begin{equation}\label{lag_corresp_g}
\begin{tikzcd}
	& {\LocB^{g,1}_{P}(T^3)} \\
	{\LocB^{g,1}_{G}(T^3)} & {\LocB^{1}_{P}(T^3)} & {\LocB^{g,1}_{L_{G, \lambda}}(T^3)} \\
	{\LocB^{1}_{G}(T^3)} && {\LocB^{1}_{L_{G, \lambda}}(T^3)}
	\arrow[from=2-2, to=3-1]
	\arrow[from=2-2, to=3-3]
	\arrow[hook, from=2-3, to=3-3]
	\arrow["\cong", from=1-2, to=2-3]
	\arrow[from=1-2, to=2-2]
	\arrow[from=1-2, to=2-1]
	\arrow[from=2-1, to=3-1]
\end{tikzcd}
\end{equation}
\end{remark}
\section{Exponential map }
In this section we consider an exponential map that connects the two moduli spaces in the previous two sections. We prove that the map preserves the natural closed $2$-forms on the formal neighborhood of the $0$ representation and the trivial local system. We then show the complex analytic version is a map of $d$-critical loci. Denote $\Map(S^{1}, \XB)$ by $L \XB$. Recall that there is always a canonical constant loops map $\XB \to L \XB$.  Recall that $\To[-1] \B G \cong \mathfrak{g} /G$ and $L \B G \cong G/G$ then there is an exponential map if we complete along $\B G$ on both sides.
\subsection{Exponential map and closed forms}
\begin{prop} \label{compact_gp_res}
The exponential map $\widehat{\To}^{\B \GL_n}[-1] \B \GL_n \to \widehat{L}^{\B \GL_n} \B \GL_n$ preserves the respective closed $2$-forms, where the completion is happening with respect to the $0$-section $\B \GL_n \to \To[-1] \B \GL_n$ and constant loops $\B \GL_n \to L \B \GL_n$.
\end{prop}
\begin{proof}
Start by defining explicitly the forms on  $\mathfrak{gl}_{n}/ \GL_n$ and $\GL_n/\GL_n$. Recall the closed two form in Example \ref{s1_forms} $\omega_{0} + \omega_{1}$ on ${\GL_n/\GL_n}$ as well as the exact $2$-form in Example \ref{tangent_bg} $\gamma_{0}$ on $\mathfrak{gl}_{n} / \GL_n$. Firstly it is enough to consider the question for the map $\widehat{\mathfrak{gl}_{n}}/ \GL_n \to \widehat{\GL_n}/\GL_n$. Once we have the claim we can pullback to the completions at a point which are $\widehat{\mathfrak{gl}_{n}}/ \widehat{\GL_n} \to \widehat{\GL_n}/\widehat{\GL_n}$. Here completions are happening at $0 \in \mathfrak{gl}_{n}$ and $1 \in \GL_n$. Now we can use the Cartan model for equivariant de Rham cohomology here to get explicit models for the de Rham complex
\begin{equation}
    \mathbf{DR}(\widehat{\mathfrak{gl}_{n}} / \GL_n)  = (\Omega^{*} \widehat{\mathfrak{gl}_{n}} \otimes S^{*} \mathfrak{gl}_{n}^{*}[-2])^{\GL_n}  
\end{equation}
We wish to prove that the closed two forms $\gamma = (\gamma_{0} , 0 , \dots)$ and $\omega = (\exp^{*} \omega_{0}, \exp^{*} \omega_{1}, 0 ,\dots)$ define the same cohomology class in $\HH^{1}(A^{2,\cl}(\widehat{\mathfrak{gl}_{n}} / \GL_n))$. Note that pulling back forms by $\exp$ makes sense since we are working formally. To prove this we need to find a $2$-form of degree $0$ $\varpi_{1}$ that satisfies the equations
\begin{align} \label{form_equation}
    \gamma_{0} - \exp^* \omega_{0} & = d \varpi \\
    -\exp^{*} \omega_{1} & = d_{\dR}\varpi \nonumber
\end{align}
In \cite[Lemma 3.3]{group_valued_moment}, the forms $\gamma_{0}$ and $\exp^{*} \omega$ are compared in the setting of compact Lie groups over $\mathbb{R}$. We will now transfer these results to our setting. \par
The maximal compact subgroup of $\GL_n$ is the group of unitary matrices $\UU_{n} \subseteq \GL_n$ with Lie algebra $\mathfrak{u}_{n}$. Now we can similarly define graded mixed complexes over $\mathbb{R}$ 
\begin{align*}
    \mathbf{DR}( \mathfrak{u}_{n} / \UU_{n}) & =(\Omega^{*} \mathfrak{u}_{n} \otimes S^{*} \mathfrak{u}_{n}^{*}[-2])^{\UU_{n}}  \\
    \mathbf{DR}(\widehat{\mathfrak{u}}_{n}/ \UU_{n}) &= (\Omega^{*} \widehat{\mathfrak{u}}_{n} \otimes S^{*} \mathfrak{u}_{n}^{*}[-2])^{\UU_{n}}
\end{align*}
here by $\Omega^{*} \widehat{\mathfrak{u}}_{n}$ we mean differential forms on $\mathfrak{u}_{n}$ with coefficients in the ring $\widehat{C}^{\infty}(\mathfrak{u}_{n})$. Here $\widehat{C}^{\infty}(\mathfrak{u}_{n})$ is the completion of $C^{\infty}(\mathfrak{u}_{n})$ at the $0$ matrix. Using these graded mixed complexes we can define the complexes $A^{2,\cl}$ and thus define closed $2$-forms in this setting as well. We also have $\widehat{C}^{\infty}(\mathfrak{u}_{n}) \otimes_{\mathbb{R}} \mathbb{C} \cong \mathcal{O}(\widehat{\mathfrak{gl}_{n}})$ using that $\mathfrak{u}_{n} \otimes_{\mathbb{R}} \mathbb{C} = \mathfrak{gl}_{n}$. We can define the real counterparts of the forms $\gamma_0$, $\omega_{0}$ and $\omega_{1}$ using the real pairing $\operatorname{Tr} \colon \mathfrak{u}_{n} \times \mathfrak{u}_{n} \to \mathbb{R}$. Note that complexifying this pairing on $\mathfrak{u}_{n}$ we get the trace pairing on $\mathfrak{gl}_{n}$.  We have comparison maps
\begin{align*}
\mathbf{DR}(\mathfrak{u}_{n} / \UU_{n})_{\mathbb{C}} & \to  \mathbf{DR}(\widehat{\mathfrak{u}_{n}} / \UU_{n})_{\mathbb{C}} \xleftarrow{\cong} \mathbf{DR}(\widehat{\mathfrak{gl}_{n}} / \GL_n) \\
A^{2,\cl}( \mathfrak{u}_{n}/\UU_{n})_{\mathbb{C}} & \to  A^{2,\cl}( \widehat{\mathfrak{u}_{n}}/\UU_{n})_{\mathbb{C}} \xleftarrow{\cong} A^{2,\cl}( \widehat{\mathfrak{gl}_{n}}/\GL_n)
\end{align*}
Here $(-)_{\mathbb{C}} = (-) \otimes_{\mathbb{R}} \mathbb{C}$. The last maps in both rows are isomorphisms of graded mixed complexes using that invariants over $\UU_{n}$ are the same as invariants over $\GL_n$. From the statement about pairings it then follows that the real versions of the forms $\gamma$ and $\omega$ are mapped to their complex versions under the comparison maps. Now \cite[Lemma 3.3]{group_valued_moment}  says that $\gamma - \omega = 0 \in \HH^{1}(A^{2,\cl}(\mathfrak{u}_{n}/\UU_{n}))$. But the comparison maps then imply that $\gamma - \omega = 0 \in \HH^{1}(A^{2,\cl} (\widehat{\mathfrak{gl}_{n}} / \GL_n))$.

\end{proof}
 The next proposition proves a well-known relationship between mapping stacks and cotangent bundles that the author could not find a reference for. For the proof we will need the following explicit description of transgressed $1$-forms on $\Map(\ZB , \XB)$. Assume that $\ZB$ satisfies all the finiteness assumptions in subsection \ref{aksz_sect} and has a fundamental class of degree $d$. Consider an $n$-shifted $1$-form $\mathcal{O}_{\XB} \to \mathbb{L}_{\XB}[n]$. The transgressed $1$-form is given by
\begin{align*}
A^{1}(\XB) & \to A^{1}(\Map(\ZB,\XB)) \\
(\mathcal{O}_{\XB} \to \mathbb{L}_{\XB}[n]) \mapsto (\mathcal{O}_{\Map(\ZB, \XB)} & \to \pi_{\#} \ev^{*}\mathcal{O}_{\XB}[-d] \to \pi_{\#}\ev^{*}\mathbb{L}_{\XB}[n-d] \cong \mathbb{L_{\Map(\ZB, \XB)}}[n-d]).
\end{align*}
\begin{prop} \label{loops_co}
    Let $\XB$ be a derived Artin stack and let $\ZB$ satisfy all the finiteness assumptions in subsection \ref{aksz_sect} with a fundamental class of degree $d$. Then we have a symplectic isomorphism $\eta \colon \Map(\ZB, \To^{*}[n]\XB) \to \To^{*}[n-d] \Map(\ZB, \XB)$. The $n-d$ shifted symplectic structure on the left is given by AKSZ and on the right by the canonical form on shifted cotangent stacks.
\end{prop}
\begin{proof}
    Let us start by defining the map $\eta$. Based on the functor of points of $\To^{*}[n-d] \Map(\ZB,\XB)$ we need to define a map $f \colon \Map(\ZB, \To^{*}[n]\XB) \to \Map(\ZB,\XB)$ and a section $s \in \Gamma(f^{*} \mathbb{L}_{\Map(\ZB,\XB)}[n-d])$. The map $f$ is given by taking mapping stacks of the projection $\pi_{\XB} \colon \To^{*}[n] \XB \to \XB$. Now we have the canonical section on $\To^{*}[n]\XB$, $l_{\XB} = \mathcal{O}_{\To^{*}[n]\XB} \to \pi^{*}_{\XB} \mathbb{L}_{\XB}[n]$. Pulling back along $\ev \colon \ZB \times \Map(\ZB, \To^{*}[n]\XB) \to \To^{*}[n]\XB$ and pushing forward by $p \colon \ZB \times \Map(\ZB, \To^{*}[n]\XB) \to \Map(\ZB, \To^{*}[n]\XB)$, $s$ is given by the map
    \begin{align} \label{transgr}
        \mathcal{O}_{\Map(\ZB,\To^{*}[n]\XB)} \xrightarrow{[\ZB]} p_{\#}\mathcal{O}_{\ZB \times \Map(\ZB, \To^{*}[n]\XB)}[-d] \xrightarrow{p_{\#} \ev^{*} l_{\XB}} p_{\#} \ev^{*} \pi^{*}_{\XB} \mathbb{L}_{\XB}[n-d] \cong f^{*} \mathbb{L}_{\Map(\ZB,\XB)}[n-d].
    \end{align}
   By considering the functor of points of both spaces we can prove that $\eta$ is an equivalence. In particular, the data of a map from a derived scheme $S$ to $\Map(\ZB, \To^{*}[n]\XB)$ is given by a map $f \colon S \times \ZB \to \XB$ and $s_f \colon \mathcal{O}_{S \times \ZB} \to f^{*} \mathbb{L}_{\XB}[n]$. On the other hand, a map $S$ to $\To^{*}[n-d] \Map(\ZB, \XB)$ is given by $f \colon S \times \ZB \to \XB$ and $\widetilde{s}_{f}: \mathcal{O}_{S} \to \pi_{S\#} f^{*} \mathbb{L}_{\XB}[n-d]$. Here $\pi_S$ is the projection $S \times \ZB \to S$. The map $\eta$ defined above sends $s_{f}$ to a map $\mathcal{O}_{S} \to \pi_{S\#} f^{*} \mathbb{L}_{\XB}[n-d]$ by applying $\pi_{S\#}$ and precomposing with the fundamental class. Now given $\widetilde{s}_{f} \colon \mathcal{O}_{S} \to \pi_{S\#} f^{*} \mathbb{L}_{\XB}[n-d]$ we can use the natural isomorphism $\pi_{S*} \to \pi_{S \#}[-d]$ and the adjunction between $\pi_{S*}$ and $\pi^{*}_{S}$ to define a map $\mathcal{O}_{S \times \ZB} \to f^{*} \mathbb{L}_{\XB}[n]$. The natural isomorphism $\pi_{S*} \to \pi_{S\#}[-d]$ is defined using the one for $p \colon \ZB \to \pt$. This shows that we have an equivalence on points.  \par Now the $n-d$ shifted symplectic structure on $\To^{*}[n-d] \Map(\ZB,\XB)$ is induced from the canonical $1$-form 
   
   $$\lambda_{\Map(\ZB,\XB)} \colon \mathcal{O}_{\XB} \xrightarrow{l_{\Map(\ZB,\XB)}} \pi^{*}_{\Map(\ZB,\XB)} \mathbb{L}_{\Map(\ZB,\XB)}[n-d] \to \mathbb{L}_{\To^{*}[n-d]\Map(\ZB,\XB)}[n-d].$$ 
   The symplectic form on $\Map(\ZB, T^{*}[n]\XB)$ is given by transgressing the canonical one form $\lambda_{\XB}$ on $\To^{*}[n]\XB$ and taking de Rham differential since the AKSZ construction commutes with the de Rham differential. The transgression of the form is given by 
    \begin{align*} 
        \lambda_{\operatorname{aksz}} \colon \mathcal{O}_{\Map(\ZB,\XB)} \xrightarrow{[\ZB]} p_{\#}\mathcal{O}_{\ZB \times \Map(\ZB, \To^{*}[n]\XB)}[-d] & \xrightarrow{p_{\#} \ev^{*} l_{\XB}} \pi_{\#} \ev^{*} \pi^{*}_{\XB} \mathbb{L}_{\XB}[n-d] \to  \\
        & \xrightarrow{ } \pi_{\#} \ev^{*}\mathbb{L}_{\To^{*}[n]\XB}[n-d]=\mathbb{L}_{\Map(\ZB, \To^{*}[n]\XB)}[n-d].
    \end{align*}
    note that the composition of the first two maps is given by the section $s$ in equation \eqref{transgr}. Pullback commutes with de Rham differential so it is enough to compare the pulled back $1$ form $\eta^{*} \lambda_{\To^{*}[n-d]  \Map(\ZB, \XB)}$ to $\lambda_{\operatorname{aksz}}$. By definition of the map $\eta$ we have $\eta^{*} l_{ \Map(\ZB,\XB)} = s$. Therefore, $\eta^{*} \lambda_{\To^{*}[n-d]  \Map(\ZB, \XB)}$ is given by composing $s$ with the top horizontal maps in the diagram below
    \begin{equation*}
\begin{tikzcd}
	{\eta^{*}\pi^{*}_{\Map(\ZB,\XB)} \mathbb{L}_{\Map(\ZB,\XB)}[n-d]} & {\eta^{*} \mathbb{L}_{\To^{*}{n-d} \Map(\ZB,\XB)}[n-d]} & {\mathbb{L}_{\Map(\ZB, \To^{*}[n] \XB)}[n-d]} \\
	& {f^{*} \mathbb{L}_{\Map(\ZB,\XB)}[n-d]}
	\arrow[from=1-1, to=1-2]
	\arrow[from=1-2, to=1-3]
	\arrow["\cong"', from=1-1, to=2-2]
	\arrow[from=2-2, to=1-3]
\end{tikzcd}
\end{equation*}
However, because the diagram commutes $\lambda_{\operatorname{aksz}}$ is identified with $\eta^{*} \lambda_{\Map(\ZB,\XB)}$.
\end{proof}
 In particular, we can use the previous proposition with $\ZB = M$ for $M$ a closed oriented $d$-manifold. Using that the AKSZ construction is functorial we can deduce the following corollary of Proposition \ref{compact_gp_res}.
\begin{cor} \label{exp_symp_final}
    The map $\exp \colon \widehat{\To^{*}}[-1] \LocB_{\GL_n}(T^2) \to \widehat{\LocB}_{\GL_n}(T^{3})$ induced my taking $\Map(T^{2},-)$ of the map  $\widehat{\To}^{\B \GL_n}[-1] \B \GL_n \to \widehat{L}^{\B \GL_n} \B \GL_n$ in Proposition \ref{compact_gp_res} preserves closed $2$-forms. Here we are completing at the $0$ section $\LocB_{\GL_n}(T^{2}) \to \To^{*}[-1] \LocB_{\GL_n}(T^{2})$ and the constant loops $\LocB_{\GL_n}(T^{2}) \to L \LocB_{\GL_n}(T^{2}) = \LocB_{\GL_n}(T^3)$.
\end{cor}
\begin{proof}
     Start by using the pairing on $\mathfrak{gl}_{n}$ to identify $\To^{*}[1] \B \GL_n \cong \To[-1]\B \GL_n$ along with their $1$-shifted symplectic structures. 
    Consider now taking $\Map(T^{2},-)$ of the map $\exp \colon \widehat{\To}^{\B \GL_n}[-1] \B \GL_n \to \widehat{L}^{\B \GL_n} \B \GL_n$. To avoid clutter write $\Map(T^2, \widehat{\To}^{\B \GL_n}[-1]) = (\widehat{\To}^{\B \GL_n}[-1])^{T^2}$. Let us compute $(\widehat{\To}^{\B \GL_n}[-1])^{T^2}$.
    \begin{align*}
        (\To[-1] \B \GL_n \times_{(\To[-1] \B \GL_n)_{\dR}} (\B \GL_n)_{\dR})^{T^2} & \cong (\To^{*}[1] \B \GL_n \times_{(\To^{*}[1] \B \GL_n)_{\dR}} (\B \GL_n)_{\dR})^{T^2} \\
        & \cong (\To^{*}[1] \B \GL_n)^{T^{2}} \times_{((\To^{*}[1] \B \GL_n)_{\dR})^{T^2}} ((\B \GL_n)_{\dR})^{T^{2}} \\
        & \cong (\To^{*}[1] \B \GL_n)^{T^{2}} \times_{((\To^{*}[1] \B \GL_n)^{T^2})_{\dR}} ((\B \GL_n)^{T^2})_{\dR} \\
        & \cong \To^{*}[-1] (\B \GL_n)^{T^2} \times_{(\To^{*}[-1] (\B \GL_n)^{T^2})_{\dR}} ((\B \GL_n)^{T^2})_{\dR} \\
        & \cong \widehat{\To^{*}}^{\LocB_{\GL_n}(T^2)}[-1] \LocB_{\GL_n}(T^2).
    \end{align*}
We have used that mapping stacks preserve pullbacks in the target. Also we have used that the $(-)_{\dR}$ functor commutes with finite limits, which allows us to commute it with mapping stacks. Finally, we have used Proposition \ref{loops_co} in the second to last isomorphism. To summarise, we have shown that $(\widehat{\To}^{\B \GL_n}[-1])^{T^2}$ is isomorphic to the completion at the $0$ section of $\To^{*}[-1] \LocB_{\GL_n}(T^{2})$. The isomorphism also preserves closed $2$-forms. A similar calculation shows that $(\widehat{L}^{\B GL_n} \B \GL_n)^{T^2} \cong \widehat{\LocB}^{\LocB_{\GL_n}(T^{2})}_{\GL_n}(T^{3})$ is  the completion of $L \LocB_{\GL_n}(T^{2})$ at the constant loops. We now have the following diagram
\begin{equation}
\begin{tikzcd}
	{\To^{*}[-1] \LocB_{\GL_n}(T^{2})} && {\LocB_{\GL_n}(T^{3})} \\
	{\Map(T^{2}, \To[-1] \B \GL_n)} && {\Map(T^{2}, \LocB_{G}(S^{1}))} \\
	{\Map(T^{2}, \widehat{\To}[-1] \B \GL_n)} && {\Map(T^{2}, \widehat{\LocB_{\GL_n}}(S^{1}))} \\
	{\widehat{\To^{*}}[-1] \LocB_{\GL_n}(T^{2})} && {\widehat{\LocB}_{\GL_n}(T^{3})}
	\arrow["\cong", from=2-1, to=1-1]
	\arrow["\cong", from=2-3, to=1-3]
	\arrow["{j_1}", from=3-1, to=2-1]
	\arrow["{\Map(T^2,\exp)}", from=3-1, to=3-3]
	\arrow["{j_2}", from=3-3, to=2-3]
	\arrow["\cong", from=4-1, to=3-1]
	\arrow["\exp", from=4-1, to=4-3]
	\arrow["\cong", from=4-3, to=3-3]
\end{tikzcd}
\end{equation}
In this diagram we are completing with respect to $0$ sections or constant loops respectively. Here $j_i$ preserve closed $2$-forms since they are induced from the \'etale maps $\widehat{\To}[-1]\B \GL_n \to \To[-1]\B \GL_n$ and $\widehat{\Map}(S^{1}, \B \GL_n) \to \Map(S^{1}, \B \GL_n)$, which preserve the closed $2$-forms. The right vertical isomorphism preserves closed forms using Lemma \ref{aksz_completions}. The left vertical isomorphism preserves closed forms because of Lemma \ref{aksz_completions} and the fact that the isomorphism in Proposition \ref{loops_co} preserves closed forms. Using functoriality of the AKSZ construction and Proposition \ref{compact_gp_res} we can conclude that $\exp \colon \widehat{\To^{*}}[-1] \LocB_{\GL_n}(T^2) \to \widehat{\LocB}_{\GL_n}(T^{3})$ preserves closed $2$-forms. 
\end{proof}
 Furthermore, the closed $2$-form on $\widehat{\To^{*}}[-1] \LocB_{\GL_n} (T^{2})$ is the one induced by the open inclusion 
    $$\To^{*}[-1] \LocB_{\GL_n} (T^{2}) \to \To^{*}[-1] \To^{*} \mathfrak{gl}_{n}/\GL_n .$$ The latter space is the derived critical locus of $$\operatorname{Tr}(\widetilde{W})_{n} / \GL_n \colon \Rep_{n}(\widetilde{Q}_{\operatorname{Jor}})/ \GL_n \to \mathbb{A}^{1} .$$ This follows since $\To^{*}[-1] \LocB_{\GL_n} (T^{2})$ is open in $\To^{*}[-1] \To^{*} \mathfrak{gl}_{n}/\GL_n$ and the symplectic form is given just by restriction. The fact that $\To^{*}[-1] \To^{*} \mathfrak{gl}_{n}/\GL_n$  is the derived critical locus of $\operatorname{Tr} (\widetilde{W})_{n} / \GL_n $ follows from Proposition \ref{stacky_kinjo}.
\begin{lem}\label{aksz_completions}
Let $N$ be a closed oriented $n$-manifold and $\YB, \ZB$ derived Artin stacks. Fix a closed $m$-form $\omega$ on $\YB$ and call the transgressed $(m-n)$-form on $\Map(N, \YB)$ $\omega_{\operatorname{aksz}}$. There is a closed $(m-n)$-form $\hat \omega = \gamma^{*} \omega_{\operatorname{aksz}}$ induced from the map  $\gamma \colon \widehat{\Map}^{\Map(N, \ZB)}(N, \YB) \to \Map(N, \YB)$ from the $(m-n)$-form $\omega_{\operatorname{aksz}}$ by restriction. We can also consider the closed $(m-n)$-form $\hat \omega_{\operatorname{aksz}}$ on $\Map(N, \widehat{\YB}^{\ZB})$ obtained by AKSZ from the form $\hat \omega_{res} = \eta^{*} \omega$ given by restriction of $\omega$ along the map $\eta \colon \widehat{\YB}^{\ZB} \to \YB$. Under the isomorphism $\widehat{\Map}^{\Map(N, \ZB)}(N, \YB) \cong \Map(N, \widehat{\YB}^{\ZB})$ the forms  $\hat \omega$ and $\hat \omega_{\operatorname{aksz}}$ coincide. 
\end{lem}
\begin{proof}
First let us note that we can use the fact that mapping stacks preserve pullbacks in the target and that $(-)_{\dR}$ commutes with finite limits to deduce the isomorphism $\widehat{\Map}^{\Map(N, \ZB)}(N, \YB) \cong \Map(N, \widehat{\YB}^{\ZB})$. To avoid clutter denote $\Map(N, \YB) $ by $\YB^{N}$. We can identify the correspondence $( \widehat{\YB}^{\ZB})^{N} \xleftarrow{} N \times (\widehat{\YB}^{\ZB})^{N} \xrightarrow{\ev} \widehat{\YB}^{\ZB}$ with the correspondence $\widehat{(\YB^{N})}^{\ZB^{N}} \xleftarrow{} N \times \widehat{(\YB^{N})}^{\ZB^{N}} \xrightarrow{\ev} \widehat{\YB}^{\ZB}$. Where the latter $\ev$ map is constructed from the following commutative diagram using the definition of completions as pullbacks. 
\begin{equation}
\begin{tikzcd}
	{N \times \widehat{(\YB^{N})}^{\ZB^{N}}} & {N \times (\ZB^{N})_{\dR}} \\
	{N \times \YB^{N}} & {N \times (\YB^{N})_{\dR}} & {\ZB_{\dR}} \\
	\YB & {\YB_{\dR}}
	\arrow[from=1-1, to=1-2]
	\arrow[from=1-1, to=2-1]
	\arrow[from=1-2, to=2-2]
	\arrow[from=1-2, to=2-3]
	\arrow[from=2-1, to=2-2]
	\arrow[from=2-1, to=3-1]
	\arrow[from=2-2, to=3-2]
	\arrow[from=2-3, to=3-2]
	\arrow[from=3-1, to=3-2]
\end{tikzcd}
\end{equation}
Using this we obtain the following commutative diagram where the left square is a pullback
\begin{equation}
\begin{tikzcd}
	{\widehat{(\YB^{N})}^{\ZB^{N}}} & {N \times\widehat{(\YB^{N})}^{\ZB^{N}} } & {\widehat{\YB}^{\ZB}} \\
	{\YB^{N}} & {N \times \YB^{N}} & \YB
	\arrow[from=1-1, to=2-1]
	\arrow[from=1-2, to=1-1]
	\arrow["\ev", from=1-2, to=1-3]
	\arrow[from=1-2, to=2-2]
	\arrow[from=1-3, to=2-3]
	\arrow[from=2-2, to=2-1]
	\arrow["\ev"', from=2-2, to=2-3]
\end{tikzcd}
\end{equation}
Using \cite[Remark 3.1.4]{calaque2022aksz}, the compatibility of the pushforward of differential forms with pullbacks then shows that the AKSZ form on $\widehat{(\YB^{N})}^{\ZB^{N}}$ is identified with the one given by restriction from $\YB^{N}$.
\end{proof}
Note that because we have a statement about completions at $\LocB_{\GL_n}(T^{2})$ we can immediately deduce that the exponential map preserves closed $2$-forms also for completions along the closed points $(A,B,0) \in \To^{*}[-1] \LocB_{\GL_n}(T^{2})$ and $(A,B,I) \in \LocB_{\GL_n}(T^{3})$.
\subsection{Exponential map and volume forms}
Here we describe how the exponential map behaves with respect to volume forms as in \cite[Section 5.1]{naef2023torsion}. Firstly, by \cite[Proposition 5.9,5.17]{naef2023torsion} we can view perfect complexes on $\B \widehat{\mathbb{G}}_{a} \times \XB$ and $\B \mathbb{G}_{a} \times \XB$ as perfect complexes on $\XB$ with an endomorphism and a nilpotent endomorphism respectively. There is a $\mathbb{G}_{m}$ action on $\B \widehat{\mathbb{G}}_{a}$ induced by the one on $\widehat{\mathbb{G}}_{a}$. We also have the following proposition.
\begin{prop} \cite[Proposition 5.12]{naef2023torsion}
    Let $\XB$ be a derived prestack. We have a $\mathbb{G}_{m}$-equivariant isomorphism $\To[-1] \mathbf{X} \cong \Map( \B \widehat{\mathbb{G}}_{a} , \mathbf{X})$. 
\end{prop}
Under this isomorphism given a complex $F \in \QCoh(\XB)$ the complex $\ev^{*}F \in \QCoh (\B \widehat{\mathbb{G}}_{a}) \times \To[-1] \XB$ corresponds to a weight $1$ endomorphism on the complex $p^{*} F \in \QCoh(\To[-1] \XB)$ for $p: \To[-1]\XB \to \XB$.
\begin{defn}[Atiyah Class]
Let $F$ be a bounded above quasicoherent complex on $\mathbf{X}$. Then the Atiyah class is a weight $1$ endomorphism $\at_{F} \colon p^* F \to p^{*}F$. 
\end{defn}
Because $F$ is bounded above the endomorphism $p^{*}F \to p^{*}F$ is equivalent to a map $F \to F \otimes \mathbb{L}_{\XB}[-1]$ by \cite[Theorem 2.5]{monier2021notelinearstacks}.
\begin{remark} \label{orientation_remark}
    Let $T^{2}$ be the $2$-torus. Consider the equivalences 
    $$\To^{*}[-1] \LocB_{\GL_n}(T^{2}) \cong \Map(T^{2}, \To[-1] \B \GL_n) \cong \To[-1] \LocB_{\GL_n}(T^{2}).$$ Here the last equivalence is using the mapping stack description of $\To[-1] \LocB_{\GL_n}(T^{2})$. In \cite[Proposition 5.16]{naef2023torsion} a volume form $\omega_{a}$ on $\To[-1] \XB$ is defined using the using the abelian group structure of $\To^{*}[-1] \XB$ relative to $\XB$ to get $\mathbb{L}_{\To[-1]\XB / \XB} \cong p^{*}_{\XB} \mathbb{L}_{\XB}[1]$ with $p_{\XB} \colon \To[-1]\XB \to \XB$. Comparing this volume form to the orientation defined in \cite[Example 2.15]{dim_red_kinjo} we see that they are defined in the same way hence the volume form on $\To^{*}[-1]\XB$ induces this canonical orientation as in \cite{dim_red_kinjo}. Furthermore, we have a description of $\To^{*}[-1]\LocB_{\GL_n}(T^{2})$ as an open of a derived critical locus, this gives a trivial orientation $\mathbb{Z}/ 2 \mathbb{Z}$ local system by \cite[Lemma 2.19]{dim_red_kinjo}. 
\end{remark}
Now consider the correspondence $\B \widehat{\mathbb{G}}_{a} \to \B \mathbb{G}_{a} \xleftarrow{} \B \mathbb{Z}$, which by taking mapping stacks $\Map(-,\mathbf{X})$ induces a correspondence 
\begin{equation} \label{unipotent_loop_corresp}
    \To[-1] \mathbf{X} \xleftarrow{q^{a}} L^{u} \mathbf{X} \xrightarrow{q^{m}} L \mathbf{X}.
\end{equation} 
Here $L^{u}$ is called the unipotent loop space. Now we can restrict the Atiyah class to the unipotent loop space to get a nilpotent endomorphism $p^* F \to p^* F$. Therefore, we can evaluate on any invertible power series $f \in \mathbb{C}[[x]]^{\times}$, to get an automorphism $f(\at_{F})$ and an invertible function $\det f(\at_{F}) \in \mathcal{O}^{\times}_{L^{u} \mathbf{X}}(L^{u} \mathbf{X})$. Then we have the following theorem
\begin{thm}\cite[Theorem 5.23]{naef2023torsion}
Let $\omega_{a}$ be the natural volume form on $\To[-1] \mathbf{X}$ and $\omega_{m}$ the natural volume form on $L \mathbf{X}$. We have an equality of volume forms on $L^{u} \mathbf{X}$
\begin{equation}
    q^{*}_{m} \omega_{m} = q^{*}_{a} \omega_{a} \cdot \det (\frac{\at_{\mathbb{L}_{\mathbf{X}}}}{\exp(\at_{\mathbb{L}_{\mathbf{X}}}) - 1 }).
\end{equation}
\end{thm}
For the $3$ torus we can make the following computation.
\label{at_class}
\begin{prop}
    Let $\mathbf{X} = \LocB_{G} (T^{2})$ with $G$ a reductive group. Then $q^{*}_{m} \omega_{m} = q^{*}_{a} \omega_{a}$.
\end{prop}
\begin{proof}
Note that to compute the function $ \det (\frac{\at_{\mathbb{L}_{\mathbf{X}}}}{\exp(\at_{\mathbb{L}_{\mathbf{X}}}) - 1 })$ it is enough to consider the classical truncation of the correspondence \eqref{unipotent_loop_corresp}. Then for $\XB = \LocB_{\GL_n} (T^{2})$ we get 
\begin{equation}
    \CC_{3}(\GL_n^{2}, \mathfrak{gl}_{n})/\GL_n \xleftarrow{} \CC_{3}(\GL_n^{2}, \widehat{\mathfrak{gl}_{n}}^{\mathcal{N}})/\GL_n \xrightarrow{} \CC_{3}(\GL_n)/\GL_n 
\end{equation}
Here by $\CC_{3}(\GL_n^{2}, \widehat{\mathfrak{gl}_{n}}^{\mathcal{N}})/\GL_n$ we mean the stack of $3$-commuting matrices one of which in the Lie algebra which is completed along the nilpotent matrices in the Lie algebra.
We can compute the restriction of the Atiyah class to the classical truncation of $\To[-1] \LocB_{\GL_n}(T^{2})$ as the map of $\GL_n$-equivariant complexes
\begin{equation}
\begin{tikzcd}
	{\mathcal{O}(\CC_{3}(\GL_n^{2}, \mathfrak{gl}_{n})) \otimes \mathfrak{gl}_{n}^{*}} & {\mathcal{O}(\CC_{3}(\GL_n^{2}, \mathfrak{gl}_{n})) \otimes (\mathfrak{gl}_{n}^{*})^{2}} & {\mathcal{O}(\CC_{3}(\GL_n^{2}, \mathfrak{gl}_{n})) \otimes \mathfrak{gl}_{n}^{*}} \\
	{\mathcal{O}(\CC_{3}(\GL_n^{2}, \mathfrak{gl}_{n})) \otimes \mathfrak{gl}_{n}^{*}} & {\mathcal{O}(\CC_{3}(\GL_n^{2}, \mathfrak{gl}_{n})) \otimes (\mathfrak{gl}_{n}^{*})^{2}} & {\mathcal{O}(\CC_{3}(\GL_n^{2}, \mathfrak{gl}_{n})) \otimes \mathfrak{gl}_{n}^{*}}
	\arrow[from=1-1, to=1-2]
	\arrow["{\ad^{*}_{Z}}"', from=1-1, to=2-1]
	\arrow[from=1-2, to=1-3]
	\arrow["{\ad^{*}_{Z} \oplus \ad^{*}_{Z}}", from=1-2, to=2-2]
	\arrow["{\ad^{*}_{Z}}", from=1-3, to=2-3]
	\arrow[from=2-1, to=2-2]
	\arrow[from=2-2, to=2-3]
\end{tikzcd}
\end{equation}
The form of the cotangent complex can be obtained in the same way as for the $3$-torus in Example \ref{cotangent_t3}. Here $Z$ is the matrix of variables corresponding to $\mathfrak{gl}_{n}$ and $\ad^{*}_{Z}$ is the coadjoint action. Now restricting to $\CC_{3}(\GL_n^{2}, \widehat{\mathfrak{gl}_{n}})$ amounts to taking $Z$ to be nilpotent. Then we can apply the function $f(x) = \frac{x}{ \exp(x) - 1}$ and take the determinant. However, in this case because of the grading we have 
$$\det^{1} (f (\ad_{Z})) \cdot \det^{-2} (f (\ad_{Z})) \cdot \det^{1} (f (\ad_{Z})) = 1.$$
\end{proof}
\subsection{Exponential map and d-critical loci}
\begin{notation} \label{notation_expdcrit}
In this subsection, let $G$ be $\GL_n$ or a Levi subgroup $L \subsetneq \GL_n$ and let $\CC_{3}(G^{2},\mathfrak{g})$ be the space of pairwise commuting triples in $G^{2} \times \mathfrak{g}$. Throughout this subsection call $\CC_{3}(G^{2}, \mathfrak{g})/G = M_{\mathfrak{g}} / G$ and $\Loc_{G}(T^{3}) = M_{G} / G$. \par
We also refer to the natural symplectic forms or $d$-critical structures on $M_{\mathfrak{g}} / G$ as \textbf{additive} and the corresponding ones on $M_{G} / G$ as \textbf{multiplicative}. We will denote their derived enhancements by $\mathbf{M}_{\mathfrak{g}} / G = \To^{*}[-1] \LocB_{G}(T^{2})$  and $\mathbf{M}_{G} / G =  \LocB_{G}(T^{3})$.
\end{notation}
    Note that $M_{\mathfrak{g}} / G$ is the classical truncation of $\To^{*}[-1] \LocB_{G}(T^{2})$ with its induced d-critical structure. We want to prove that the exponential map $\exp \colon M_{\mathfrak{g},\an} / G_{\an} \to M_{G,\an} / G_{\an}$ defines a \' etale cover that preserves the d-critical structures. To show this we need to show that 
\begin{align}\label{exp_section}
\exp^{\star} s_{m} = s_{a}
\end{align}
 for the sections $s_{a} \in \Gamma(M_{\mathfrak{gl}_{n},\an} , S^{0}_{M_{\mathfrak{gl}_{n}},\an})^{\GL_n}$ and $s_{m} \in \Gamma(M_{\GL_n,\an} , S^{0}_{M_{\GL_n, \an} })^{\GL_{n}}$ that control the d-critical structures.  We will show this by checking at the level of stalks of the sheaves $S^{0}_{M_{\mathfrak{gl}_{n}}, \an}$ and $S^{0}_{\GL_{n},\an}$  by using formal geometry and equivariance properties. It is enough to compare the stalks of $s_{a}$ and $\exp^{\star} s_{m}$ at closed orbits in $M_{ \mathfrak{gl}_{n},\an}$. Namely, let $y \in M_{ \mathfrak{gl}_{n},\an}$ be a point in a non-closed orbit. Then there is a closed orbit $C$ with some point $c \in C$ in the closure $\overline{\GL_{n} y}$. We know that $(\exp^{\star}s_{m})_{c} = s_{a,c}$ implies there is some analytic open around $c$ for which $\exp^{\star}s_{m} = s_{a}$. However, under the conjugation action we can always move the point $y$ to be inside this analytic open thus we must have that $(\exp^{\star}s_{m})_{y} = s_{a,y}$ as well. To carry out our strategy we will need the results of Section \ref{s0_formal_sect}  relating the stalks of the $S^{0}_{X}$ sheaf at a point $x$ and the formal completion $S^{0}_{\widehat{X}^{x}}$ for $X$ a d-critical stack.   \par
 \begin{lem}\label{section_equiv}
     Let $s_{a}, (\exp^{\star}s_{m}) \in \Gamma(M_{\mathfrak{gl}_{n}, \an},S^{0}_{M_{\mathfrak{gl}_{n}, \an}})^{\GL_{n, \an}}$ and denote the stalks of the sections at $(A,B,0) \in M_{\mathfrak{gl}_{n}, \an}/\GL_{n}$ by $(\exp^{\star}s_{m})_{(A,B,0)}$  and $s_{a, (A,B,0)}$. \par Then $(\exp^{\star}s_{m})_{(A,B,0)} = 
     s_{a, (A,B,0)}$ implies $(\exp^{\star}s_{m})_{(A,B,\lambda I)} = 
     s_{a, (A,B,\lambda I )}$ for $\lambda \in \mathbb{C}$.
 \end{lem}
 \begin{proof}
     We begin by noting that there is a $\mathbb{G}_{a}$ action on $\CC_{3}(\mathfrak{gl}_{n})/ \GL_n$ given by $\lambda \cdot (A,B,C) = (A,B,C+ \lambda I)$. This action restricts to $M_{\mathfrak{gl}_{n}}/ \GL_n$ and preserves the $d$-critical structure since it also preserves the potential $\operatorname{Tr}(C[A,B])$. Similarly, there is a $\mathbb{G}_{m}$ action on $\Loc_{\GL_n}(T^{3})$ given by $\lambda \cdot (A,B,C) = (A,B, \lambda C)$. This action scales the added loop via the identification $L \Loc_{\GL_n}(T^{2}) \cong \Loc_{\GL_n}(T^{3})$. This action also preserves the $d$-critical structures due to the fact it is the restriction of the action of $\Loc_{\ZZ(\GL_n)}(T^3)$ on $\Loc_{\GL_n}(T^{3})$. The latter action preserves the symplectic structure via Proposition \ref{center_Wact}. Therefore, the sections $s_a$ and $s_m$ are $\mathbb{G}_{a}$ and $\mathbb{G}_{m}$ equivariant respectively. We can also immediately check that the exponential map $\exp \colon M_{\mathfrak{gl}_{n}, \an} / \GL_{n, \an} \to M_{\GL_{n}, \an} / \GL_{n, \an}$ is equivariant with respect to the group homomorphism $\exp \colon \mathbb{C} \to \mathbb{C}^{*}$. We can now conclude using the functoriality of pullback maps on $S^0$ sheaves.
 \end{proof}
Before starting the proof, we will use the following well known lemma to compute the \'etale locus of $\exp \colon M_{\mathfrak{gl}_{n}, \an}/ \GL_{n, \an} \to M_{\GL_{n}, \an}/ \GL_{n, \an}$. 
\begin{lem}[Derivative of the exponential map] 
Let $\exp  \colon  \mathfrak{gl}_{n} \to \GL_{n}$ be the exponential map. The derivative of $\exp$ is 
\begin{equation}\label{der_exp}
    d \exp_{X} Y = \exp(X) \frac{1- \exp(-\ad_{X})}{\ad_{X}}Y 
\end{equation}
then $\exp$ is \'etale when $\ad_{X}$ is invertible. Equivalently the exponential map is \'etale for all $X \in \mathfrak{gl}_{n}$ that satisfy
\begin{equation}\label{etale_cond}
   \lambda_{i} - \lambda_j \neq 2 \pi i k  \text{ for any two eigenvalues of } X \text{ and }  k \in \mathbb{Z} \setminus 0. 
\end{equation}
\end{lem}
\begin{defn}[\'Etale locus] \label{etale_locus_c3gln}
 We write $\mathfrak{gl}^{\et}_n$ for the space of matrices in $\mathfrak{gl}_n$ that satisfy condition \eqref{etale_cond}. This leads us to define the \'etale locus $M^{\et}_{\mathfrak{gl}_{n}, \an}=\CC_{3}(\GL^{2}_{n}, \mathfrak{gl}^{\et}_{n})$ of $\CC_{3}(\GL^{2}_{n}, \mathfrak{gl}_{n})$ as the open subspace of triples of pairwise commuting matrices $(x,y,z) $, where $z$ satisfies condition \eqref{etale_cond}.   
\end{defn}
Before we start with the main theorem of this section let us introduce a version of the genericity conditions \ref{generic_loci} and stratifications \ref{stratifications} on $M_{\mathfrak{g}}/ G$. We can define these analogously as is defined for $M_{\GL_n}$ except we work with the good moduli space $\So^{n}(\mathbb{G}^{2}_{m} \times \mathbb{G}_{a})$ instead of $\So^{n}\mathbb{G}^{3}_{m} $.  Similarly as in the multiplicative case we can prove that the maps $ \Theta_{a} \colon M_{\mathfrak{l}_{\lambda}}/ L_{\lambda} \to M_{\mathfrak{gl}_{n}}/ \GL_n$ are \'etale and preserve $d$-critical structures. These statements follow by considering the derived enhancement $\mathbf{M}_{\mathfrak{gl}_{n}} / \GL_n = \To^{*}[-1] \LocB_{\GL_n}(T^2) $. The same arguments as in the proof of Proposition \ref{lag_corresp_gen} work since we can rewrite $\To^{*}[-1] \LocB_{\GL_n}(T^2) \cong \Map(T^{2}, \To[-1] \B \GL_n)$.
 \begin{thm}\label{exp_dcrit}
The map $\exp \colon \CC_{3} ( \GL^{2}_{n} , \mathfrak{gl}^{\et}_{n})/ \GL_{n} \to \Loc_{\GL_{n}} (T^{3})$ is an \'etale map of oriented complex analytic $d$-critical loci. 
\end{thm}
\begin{proof}
We first show that the map $\CC_{3}(\GL^{2}_{n}, \mathfrak{gl}^{\et}_{n})/ \GL_{n} \to \Loc_{\GL_{n}}(T^{3})$ is \'etale. We know that this map is given by taking $\Map(T^{2}, \mathfrak{gl}^{\et}_{n} /\GL_{n} \to \GL_{n}/ \GL_{n})$. We can make sense of this mapping stack for complex analytic stacks by taking iterated inertia stacks. Then we use that inertia stacks preserve \'etale maps. This finishes the proof and we also have that the map on atlases $\CC_{3}(\GL^{2}_{n}, \mathfrak{gl}^{\et}_{n}) \to \CC_{3}(\GL_{n})$ is \'etale.  \par
To avoid clutter switch to notation \ref{notation_expdcrit}. Let us turn to proving that the map preserves $d$-critical structures. We start by proving that 
\begin{equation} \label{s_sect_equation_thm7}
    (\exp^{\star}s_{m})_{0} = s_{a,0}
\end{equation}
, where by $0$ we mean any element $(A,B,0) \in M^{g}_{\mathfrak{gl}_{n}}$ for $A,B \in \GL_n$. We also denote by $1$ the element $(A,B,I)$ in $M_{\GL_{n}}$. Using Corollary \ref{exp_symp_final} we can see that the exponential preserves the $(-1)$-shifted closed $2$-forms. We will now deduce the statement by using the results of Section \ref{s0_formal_sect} to perform a chase along the following diagram. In particular, we repeatedly use diagrams \eqref{an_to_formal} and \eqref{functoriality}   \par
\begin{equation}
\begin{tikzcd}
	{\mathcal{A}^{2,cl}(\MB_{\mathfrak{gl}_{n}}/\GL_n,-1)} & {\Gamma(M_{\mathfrak{gl}_{n}},S^{0}_{M_{\mathfrak{gl}_{n}}})^{\GL_{n}}} & {\Gamma(M_{\mathfrak{gl}_{n}},S^{0}_{M_{\mathfrak{gl}_{n}}})} & {S^{0}_{M_{\mathfrak{gl}_{n}},0}} \\
	{\mathcal{A}^{2,cl}(\widehat{\MB}^{0}_{\mathfrak{gl}_{n}}/\GL_{n},-1)} & {S^{0}_{\widehat{M}^{0}_{\mathfrak{gl}_{n}}/\GL_{n}}} & {S^{0}_{\widehat{M}^{0}_{\mathfrak{gl}_{n}}}} & {S^{0}_{M_{\mathfrak{gl}_{n}, \an},0}} \\
	{\mathcal{A}^{2,cl}(\widehat{\MB}^{1}_{\GL_{n}}/\GL_{n},-1)} & {S^{0}_{\widehat{M}^{1}_{\GL_{n}}/\GL_{n}}} & {S^{0}_{\widehat{M}^{1}_{\GL_{n}}}} & {S^{0}_{M_{\GL_{n}, \an},1}} \\
	{\mathcal{A}^{2,cl}(\MB_{\GL_{n}}/\GL_{n},-1)} & {\Gamma(M_{\GL_{n}},S^{0}_{M_{\GL_{n}}})^{\GL_{n}}} & {\Gamma(M_{\GL_{n}},S^{0}_{M_{\GL_{n}}})} & {S^{0}_{M_{\GL_{n}},1}}
	\arrow[from=1-1, to=1-2]
	\arrow[from=1-1, to=2-1]
	\arrow[hook, from=1-2, to=1-3]
	\arrow[from=1-2, to=2-2]
	\arrow[from=1-3, to=1-4]
	\arrow[from=1-3, to=2-3]
	\arrow[hook', from=1-4, to=2-4]
	\arrow[from=2-1, to=2-2]
	\arrow[hook, from=2-2, to=2-3]
	\arrow[hook', from=2-4, to=2-3]
	\arrow["\exp", from=3-1, to=2-1]
	\arrow[from=3-1, to=3-2]
	\arrow["\exp", from=3-2, to=2-2]
	\arrow[hook, from=3-2, to=3-3]
	\arrow["\exp", from=3-3, to=2-3]
	\arrow["\exp", from=3-4, to=2-4]
	\arrow[hook', from=3-4, to=3-3]
	\arrow[from=4-1, to=3-1]
	\arrow[from=4-1, to=4-2]
	\arrow[from=4-2, to=3-2]
	\arrow[hook, from=4-2, to=4-3]
	\arrow[from=4-3, to=3-3]
	\arrow[from=4-3, to=4-4]
	\arrow[hook, from=4-4, to=3-4]
\end{tikzcd}
\end{equation}
In more detail: in the first column of the above diagram we work with the derived enhancements of the spaces $M_{\mathfrak{g}}$ and $M_{G}$. We use Corollary \ref{exp_symp_final} to get that the exponential map pulls back the additive closed form to the multiplicative one. Then by commutativity of the first column of squares we can also deduce that the equation \eqref{s_sect_equation_thm7} also holds for the $S^{0}$ sheaves of the formal completions of the classical truncations at $0$ and $1$. The $d$-critical structures $s_a$ and $s_m$ on the analytifications of $M_{\mathfrak{gl}_n}$ and $M_{\GL_{n}}$ are induced from algebraic ones via the vertical maps in the rightmost column. By Lemma \ref{ss_ssheaf_inj} the map 
$$S^{0}_{M_{\mathfrak{gl}_{n}, \an,0}} \to S^{0}_{\widehat{M}^{0}_{\mathfrak{gl}_{n}}}$$ is injective and by Lemma \ref{stack_S_sheafcomm} it is then enough to check that they are the same under the exponential map by first embedding into $S^{0}_{\widehat{M}_{\mathfrak{gl}_{n}}}$. \par
Using the fact that the sections $s_{a}$ and $\exp^{\star}s_{m}$ are $\mathbb{G}_{a}$ invariant by Lemma \ref{section_equiv} we can then also conclude that $(\exp^{\star}s_{m})_{(A,B, \lambda I)} = s_{a,(A,B,\lambda I )}$ for $\lambda \in \mathbb{C}$. \par 
Now let us fix a closed point $x \in M_{\mathfrak{gl}_{n}, \an}/\GL_{n, \an}$ corresponding to a closed orbit. Write $y = \exp(x)$. Then using an additive version of the stratification \ref{stratifications} $x$ is in the image of the \'etale map $\Theta_{a} \colon M^{g}_{\mathfrak{l}_{\lambda}}/ L_{\lambda} \to M_{\mathfrak{gl}_{n}}/ \GL_n$ for some $\lambda$. The map $\Theta_{m} \colon M^{g}_{L_{\lambda}}/ L_{\lambda} \to M_{\GL_n}/ \GL_n$ induces a commutative diagram
\begin{equation}\label{levi_exp_square}
\begin{tikzcd}
	{M^{g}_{L_{\lambda}}} & {M_{GL_{n}}} \\
	{\widehat{M}^{g,x}_{L_{\lambda}}} & {\widehat{M}^{x}_{GL_{n}}}
	\arrow[from=1-1, to=1-2]
	\arrow[from=2-1, to=1-1]
	\arrow["\cong", from=2-1, to=2-2]
	\arrow[from=2-2, to=1-2]
\end{tikzcd}
\end{equation}
We have a similar diagram in the additive case
\begin{equation}\label{levi_exp_square_lie}
\begin{tikzcd}
	{M^{g}_{\mathfrak{l}_{\lambda}}} & {M_{\mathfrak{gl}_{n}}} \\
	{\widehat{M}^{x,g}_{\mathfrak{l}_{\lambda}}} & {\widehat{M}^{x}_{\mathfrak{gl}_{n}}}
	\arrow[from=1-1, to=1-2]
	\arrow[from=2-1, to=1-1]
	\arrow["\cong", from=2-1, to=2-2]
	\arrow[from=2-2, to=1-2]
\end{tikzcd}
\end{equation}
Furthermore, the following diagram commutes
\begin{equation}
\begin{tikzcd}
	{M^{\et,g}_{\mathfrak{l}_{\lambda}, \an}/L_{\lambda, \an}} & {M^{g}_{L_{\lambda}, \an}/L_{\lambda, \an}} \\
	{M_{\mathfrak{gl}_{n}, \an}/\GL_{n, \an}} & {M_{\GL, \an}/\GL_{n, \an}}
	\arrow["{\exp_{L_{\lambda}}}", from=1-1, to=1-2]
	\arrow[from=1-1, to=2-1]
	\arrow[from=1-2, to=2-2]
	\arrow["{\exp_{\GL_n}}", from=2-1, to=2-2]
\end{tikzcd}
\end{equation}
This induces the following diagram on stalks of $S$ sheaves
\begin{equation} \label{exp_inj_s_sheaf}
\begin{tikzcd}
	{S^{0}_{M_{\mathfrak{l}_{\lambda}},\an, x}} & {S^{0}_{M_{L_{\lambda}},\an, y}} \\
	{S^{0}_{M_{\mathfrak{gl}_{n}}, \an ,x}} & {S^{0}_{M_{\GL_{n}}, \an , y}}
	\arrow["{\exp_{L_{\lambda}}}"', from=1-2, to=1-1]
	\arrow[hook, from=2-1, to=1-1]
	\arrow[hook, from=2-2, to=1-2]
	\arrow["{\exp_{GL_{n}}}"'{pos=0.4}, from=2-2, to=2-1]
\end{tikzcd}
\end{equation}
where the vertical maps preserve $d$-critical structures by Proposition \ref{lag_corresp_gen} and are injective.  Injectivity follows by using diagram  \eqref{functoriality} and the fact that we have an isomorphism on formal $S^{0}$ sheaves induced by the multiplicative and additive diagrams \eqref{levi_exp_square} and \eqref{levi_exp_square_lie}. Since the map is injective we can use diagram \eqref{exp_inj_s_sheaf} to  prove $(\exp^{\star}_{\GL_{n}} s_{m})_{x} = s_{a,x}$ by using that 
$$(\exp^{\star}_{L_{\lambda}} s_{m})_{\Theta_{a}(x)} = s_{a,\Theta_{a}(x)}$$ for a Levi subgroup $L_{\lambda}$ of $\GL_{n}$. Recall that  $\mathbf{M}_{L_{\lambda}} / L_{\lambda} = \prod_{i} \mathbf{M}_{\GL_{\lambda_{i}}}/ \GL_{\lambda_{i}}$ and $\mathbf{M}_{\mathfrak{l}_{\lambda}} / L_{\lambda} = \prod_{i} \mathbf{M}_{\mathfrak{gl}_{\lambda_{i}}}/ \GL_{\lambda_{i}}$. Therefore, we can use the description of the d-critical structure of products in Example \ref{dcrit_products} and that the exponential map factors into products. For $n =2$ the only non-trivial Levi is the torus for which the stacks are smooth and the result $(\exp^{\star}_{L_{\lambda}} s_{m})_{x} = s_{a,x}$ is automatic. For $n \geq 3$ we can deduce the claim by induction. In particular, assume we have proved the claim for $n$ then any Levi $L \subsetneq \GL_{n+1}$ will be a product of $\GL_m$ with $m \leq n$. Then we can use diagram \eqref{exp_inj_s_sheaf} given that we already now the claim for all the $\GL_m$ by assumption.   \par
We conclude by checking the orientations using volume forms and orientations on analytic stacks as in Definition \ref{naef_saf_vol}. The \'etale map $\exp$ induces a volume form $\exp^{*} \omega_{m}$ from $\Loc_{\GL_{n}}(T^{3})$ on  $\CC_{3}(\GL^{2}_{n}, \mathfrak{gl}^{\et}_{n}) / \GL_{n}$. There is also a volume form $\omega_a$ on $\CC_{3}(\GL^{2}_{n}, \mathfrak{gl}^{\et}_{n}) / \GL_{n}$ as in Remark \ref{orientation_remark}. We will show that $\omega_{a} = \exp^{*} \omega_{m}$. The difference between the two volume forms is some invertible function $g$ on $\CC_{3}(\GL^{2}_{n}, \mathfrak{gl}^{\et}_{n}) / \GL_{n}$. First we claim that 
\begin{equation*}
    \mathcal{O}(\CC_{3}(\GL^{2}_{n}, \mathfrak{gl}_{n})/ \GL_{n}) \cong  \mathcal{O}(H^{2} \times \mathfrak{h})^{W}.
\end{equation*}
Now $\mathcal{O}(X_{\GL_{n}, \an}) \cong \mathcal{O}((H^{2} \times \mathfrak{h} /\!\!/ W)_{\an}) \cong  \mathcal{O}((H^{2} \times \mathfrak{h})_{\an})^{W} $. This follows due to the proof of  \cite[Theorem 8 page 17]{fossum1989invariant}. This means that $\mathcal{O}((H^{2} \times \mathfrak{h})_{\an})^{W}$ is an integral domain and we can compute the function $g$ by computing it on a formal completion at the trivial local system. However, now we can use Proposition \ref{at_class} to conclude that the function $g = 1$ and hence the map preserves volume forms. Therefore, $\exp$ also preserves orientations. Finally, by Remark \ref{orientation_remark} we see that the orientation induced by $\omega_a$ gives the trivial $\mathbb{Z}/2 \mathbb{Z}$ local system on $\CC_{3}(\GL^{2}_{n}, \mathfrak{gl}^{\et}_{n})/\GL_{n}$.
\end{proof}
Finally  the map $\exp \colon \mathfrak{gl}^{\et}_{n} \to \GL_{n}$ is still surjective and we check that the induced map $\exp \colon M_{\mathfrak{gl}_{n, \an}}/ \GL_{n, \an} \to M_{\GL_{n, \an}}/ \GL_{n, \an}$ is surjective, giving us an \'etale cover.
 \begin{lem} \label{exp_sur_3comm}
     The map $\exp \colon \CC_{3}(\GL^{2}_{n} , \mathfrak{gl}^{\et}_{n}) \to \CC_{3}(\GL_{n})$ is surjective.
 \end{lem}
 \begin{proof}
      Let $(X,Y,Z)$ be $3$-pairwise commuting matrices in $\GL_{n}$, then since $\exp \colon \mathfrak{gl}^{\et}_{n} \to \GL_{n}$ is surjective there exists $z \in \mathfrak{gl}^{\et}_{n}$ with $\exp (z) = Z$. Now as in the proof of the main theorem in \cite{commuting_exp} we can express the operator $\ad_{z}$
     as 
     \begin{equation}
        \ad_{z}(-) =  \frac{\ad_{z}(-)}{\exp{\ad_{z}(-)} - I} \circ (\exp(-z) (-) \exp(z)- I) \text{ where } \frac{\ad_{z}}{\exp{\ad_{z}} - I}
     \end{equation}
     
     is the inverse of the operator appearing in equation \eqref{der_exp}. Because we have restricted to the locus defined by equation \eqref{etale_cond}, this makes sense. Then we have 
     \begin{equation*}
         \ad_{z}(X) = \frac{\ad_{z}(X)}{\exp{\ad_{z}}(X) - X} \circ (\exp(-z) X \exp(z) - X)
     \end{equation*}
     but since $X$ commutes with $\exp(z)$ we have $ \ad_{z}(X) = 0$. The same holds for $Y$ so we are done. 
 \end{proof}

\section{Cohomological integrality for the 3-torus} \label{section_coh_int}
 Recall the set up in Subsection \ref{loc_stack_t3}. In particular, the good moduli spaces $X$ in Definition \ref{good_moduli_loc}. In the remaining sections of the paper $G= \GL_n, \SL_n, \PGL_n$ and $\lambda$ is a partition of $n$ with length $l$. In this section we prove cohomological integrality for the stack of local systems of the $3$-torus. 
 \subsection{Formulating cohomological integrality for \texorpdfstring{$\GL_n$}{TEXT}, \texorpdfstring{$\SL_n$}{TEXT} and \texorpdfstring{$\PGL_n$}{TEXT}}
 Write $\varphi_{m}$ for the DT sheaf on $\coprod_{n \in \mathbb{N}} \Loc_{\GL_{n}}$ and $\pi_{m} \colon \coprod_{n \in \mathbb{N}} \Loc_{\GL_{n}} \to \coprod_{n \in \mathbb{N}} \So^{n} \mathbb{G}^{3}_{m} = \So \mathbb{G}^{3}_{m}$.
\begin{thm}[Cohomological integrality] \label{sym_coh_3tor}
We have an equivalence in $\DD^{+}_{c}(\So \mathbb{G}^{3}_{m})$
$$\JHmp \Jm \cong \Sym_{\boxdot} (\BPS_m \otimes \HHf(\B \mathbb{G}_{m})[-1])$$
with $\BPS_m \coloneqq \pH^{1} \JHmp \Jm  \cong \bigoplus_{n}(\Delta \colon \mathbb{G}^{3}_{m} \to \So^{n} \mathbb{G}^{3}_{m})_{*}  \mathbb{Q}_{\mathbb{G}^{3}_{m}}[3]$.
\end{thm}
We first reformulate cohomological integrality in a form that is applicable to $G = \GL_{n}$, $\SL_n$ or $\PGL_n$.
Just as in the additive case of Proposition \ref{coh_add_prop}, we have the following proposition, which is a reformulation of Theorem \ref{sym_coh_3tor}.
\begin{prop}
     Fix standard Levi subgroups in $\GL_n$ corresponding to a partition $\lambda$ of $n$. Write $\varphi_{\GL_{n}}$ the DT sheaf on $\Loc_{\GL_{n}}$ and $\pi_{\GL_{n}} \colon  \Loc_{\GL_{n}} \to X_{\GL_{n}}= \So^{n} \mathbb{G}^{3}_{m}$. Cohomological integrality is equivalent to the following statement for all $n$.
    \begin{equation}\label{coh_gln}
        \pi_{\GL_{n}*} \varphi_{\GL_{n}} \cong \bigoplus_{L_{\GL_{n},\lambda} \subseteq \GL_{n}} (\theta_{*}\BPS_{L_{\GL_{n},\lambda}} 
 \otimes \HHf(\B \ZZ(L_{\GL_{n},\lambda}))[- \dim \ZZ(L_{\GL_{n},\lambda})])^{W_{L_{\lambda}}}
    \end{equation}
\end{prop}
From this we can write down a cohomological integrality theorem for $\SL_n$ and $\PGL_n$.
\begin{thm}[Cohomological Integrality for $\SL_n$, and $\PGL_n$] \label{coh_proof_pgsln}
Let $G= \SL_n$ or if $n$ is prime, $G=\PGL_n$. Denote by $\Loc^{1}_{G}$ the connected component of the trivial local system in $\Loc_{G}$, $\pi_{G} \colon \Loc^{1}_{G} \to X^{1}_{G}$ the good moduli space and $\varphi^{1}_{G}$ the restriction of the DT sheaf to $\Loc^{1}_{G}$.  Then we have that
    \begin{equation} \label{coh_sln}
        \pi_{G,*} \varphi^{1}_{G} \cong \bigoplus_{L_{G, \lambda} \subseteq G} (\theta_{*} \BPS_{L_{G, \lambda}} 
 \otimes \HHf(\B \ZZ(L_{G, \lambda}))[- \dim \ZZ(L_{G, \lambda})])^{W_{L_{\lambda}}}
    \end{equation}
is an equivalence for $G$.
Here $\BPS_{L_{G, \lambda}} = \Delta_{*} \mathbb{Q}_{\ZZ(L_{G, \lambda})^{3}} [ 3\dim \ZZ(L_{G, \lambda}) ]$ with $\Delta \colon \ZZ^{3}(L_{G, \lambda}) \to X^{1}_{L_{G, \lambda}}$.
\end{thm}
\begin{remark}
    The restriction to the component of the trivial local system is only necessary for the case of $\PGL_n$. In other cases, $\Loc_G$ or $X_{G}$ are connected. Therefore, we will drop the $1$ from the notation if we are working with $\GL_n$ or $\SL_n$. \par
    We expect that the above Theorem is also true for $\PGL_n$ for all $n$. However, in the current proof we need to compare the geometries of $X_{\SL_n}$ and $X_{\PGL_n}$ and we use that $n$ is prime. See Subsection \ref{sln_to_pgln_gms}. 
\end{remark}
For prime $n$, using Lemma \ref{prime_pgln} we can incorporate the contributions of the non-trivial components of $\Loc_{\PGL_n}$ and then use cohomological integrality for $\Loc_{\SL_n}$ and the the trivial component of $\Loc_{\PGL_n}$ to compare $\HHf(\Loc_{\SL_n} , \varphi_{\SL_n})$ and $\HHf(\Loc_{\PGL_n} , \varphi_{\PGL_n})$. Then we can deduce the following corollary proved in subsection \ref{pgln_subsect}.
\begin{cor}[Langlands duality for prime rank]\label{langlands}
Let $n$ be prime. We have an isomorphism of graded vector spaces $\HHf(\Loc_{\SL_{n}} , \varphi_{\SL_{n}}) \cong \HHf(\Loc_{\PGL_{n}}, \varphi_{\PGL_{n}})$.
\end{cor}
Our strategy will consist of the following steps.
\begin{strat} \label{strategy} \,
    \begin{enumerate} 
\item Show that $\pi_{\GL_{n}*} \varphi_{\GL_n}$ is a pure complex of mixed Hodge modules and compute the supports appearing in the Saito decomposition using the exponential map. The supports will be the smooth strata, as in Definition \ref{stratifications}, $X^{\lambda}_{\GL_n} \subseteq X_{\GL_n}$.
\item for $G=\GL_n , \SL_n$ examine the geometry of  the maps $\theta \colon X_{L_{G,\lambda}} \to X_{G}$ induced by the inclusion $L_{G, \lambda} \to G$. Restricted to certain loci $\widetilde{Z}^{3}(L_{G, \lambda}) \subseteq X_{L_{G,\lambda}}$ the maps $ \widetilde{Z}^{3}(L_{G, \lambda}) \xrightarrow{\theta_{\lambda}} X^{\lambda}_{G}$ become $W_{L_{\lambda}}: 1$ covers. We can restrict to perverse pieces with support $X^{\lambda}_{G}$. On $X^{\lambda}_{G}$ these become some local systems $\mathcal{L}^{\lambda}_{i}$. Using the cover $\theta_{\lambda}$ we  then compute the $\mathcal{L}^{\lambda}_{i}$ in terms of perverse pieces of $\pi_{L_{G, \lambda}*} \varphi_{G, \lambda}$, supported on $\widetilde{Z}^{3}(L_{G, \lambda})$ (Proposition \ref{levis_final}). These contributions turn out to be copies of the BPS sheaves $\BPS_{L_{G , \lambda}}$ with some action of $W_{L_{G, \lambda}}$. (Part $(1)$ of Lemma \ref{loc_sys_comp}). See also Example \ref{2_4_coh_int}.
    \item Show that the local systems $\mathcal{L}^{\lambda}_{i}$ appearing in Step $2$ agree with the ones from equation \eqref{coh_sln}. In other words, we compute the right hand side 
    \begin{equation*}
        (\theta_{*}\BPS_{L_{G,\lambda}} 
 \otimes \HHf(\B \ZZ(L_{G,\lambda}))[- \dim \ZZ(L_{G,\lambda})])^{W_{L_{\lambda}}}
    \end{equation*}
    of equation \eqref{coh_sln} in terms of $\IC$ sheaves supported on the strata $X^{\lambda}_{G}$ (Lemma \ref{loc_sys_comp}). We then explicitly compare the local systems that appear to the local systems $\mathcal{L}^{\lambda}_i$ in step $2$. For $\GL_n$ this completes the proof of cohomological integrality.
    \item Deduce purity of $\pi_{\SL_n *}\varphi_{\SL_n}$ for $\SL_n$  from $\GL_n$ (Proposition \ref{purity_sln}) and use Step 2 and Step 3.
    \item Prove integrality for $\PGL_n$ for prime $n$ using integrality for $\SL_n$. The restriction to primes is necessary to compare $X_{\SL_n}$ and $X_{\PGL_n}$ and also to compute the non-trivial components of $\Loc_{\PGL_n}$.
\end{enumerate}
\end{strat}
Let us first introduce a stratification on $X_{G}$. Recall that $X_{G} = H^{3}_{G} / \! \! / W$ so that we can view an element $ x \in X_{G}$ as a triple of diagonal matrices $(D_{1},D_{2},D_{3})$ up to permutation with $D_{i} \in H_{G}$. We say a diagonal $n \times n$ matrix has eigenvalues distinct according to $\lambda$ if up to permutation it can be written as
\begin{equation*}
   D = \diag( \underbrace{x_{1}, \dots , x_{1} , }_{\lambda_{1}\text{ times}}  \underbrace{x_{2}, \dots , x_{2} , }_{\lambda_{2}\text{ times}} \dots \underbrace{x_{l} , \dots , x_{l}}_{\lambda_{l}\text{ times}} )
\end{equation*}
\begin{defn}[Stratifications] \label{stratifications}
Let $\lambda$ be a partition of $n$ of length $l$.
\begin{enumerate}
    \item \textbf{Stratification for} $\GL_n$: \\
Let us define a stratification of $X_{\GL_{n}}$ by setting 
$$X^{\lambda}_{\GL_n} = \So^{n}_{\lambda} \mathbb{G}^{3}_{m} = \{ \sum^{l}_{i =1 } \lambda_{i} x_{i}   \mid \lambda_{i} \in \lambda =(\lambda_{1}, \dots, \lambda_{l}), x_{i} \neq x_{j} \in \mathbb{G}^{3}_{m} \},$$ 
\begin{equation*}
   \text{with }\sum^{l}_{i} \lambda_{i} x_{i} = \{ \underbrace{x_{1}, \dots , x_{1} , }_{\lambda_{1}\text{ times}}  \underbrace{x_{2}, \dots , x_{2} , }_{\lambda_{2}\text{ times}} \dots \underbrace{x_{l} , \dots , x_{l}}_{\lambda_{l}\text{ times}} \} \in \So^{n} \mathbb{G}^{3}_{m} .
\end{equation*}
This defines a locally closed smooth subscheme.
Define an open dense $X^{\lambda,g}_{\GL_n}$ of $X^{\lambda}_{\GL_n}$ by the condition that if  $(D_{1},D_{2},D_{3}) \in X^{\lambda, g}_{\GL_n}$, then there is an $1 \leq i \leq 3$ such that $D_{i}$ has eigenvalues that are distinct according to $\lambda$.
\item \textbf{Stratification for} $\SL_n$: \\
We define a stratification of $X_{\SL_n}$ by setting $X^{\lambda}_{\SL_{n}} = X_{\SL_{n}} \cap X^{\lambda}_{\GL_{n}}  $. Similarly $X^{\lambda , g} _{\SL_n} = X_{\SL_n} \cap X^{\lambda,g}_{\GL_n}$.
\end{enumerate}
\end{defn}
\begin{remark} 
     Given that $X^{\lambda}_{\GL_{n}}$ is a locally closed subvariety of $X_{\GL_{n}}$, $X^{\lambda}_{\SL_{n}}$  is a locally closed subvariety of $X_{\SL_{n}}$. However, the strata of $X_{\SL_{n}}$ are no longer connected in general. \par
\end{remark}
\begin{remark}\label{special_part}
Note that for $x = (D_{1},D_{2},D_{3}) \in X_{\GL_n}$ each element $D_{i}$ has elements corresponding to some partition $\lambda^{i}$. If there is an $1 \leq i \leq 3$ such that $\lambda^{i} = (1, \dots , 1)$, then all the elements $z_{i}$ of $x = (z_{1}, \dots z_{n}) \in \So^{n} \mathbb{G}^{3}_{m}$ with $z_{i} \in \mathbb{G}^{3}_{m}$ will be distinct. This follows since $z_{k} = (d_{1k},d_{2k},d_{3k})$ where $d_{lk}$ is the $k$-th element of $D_{l}$. Therefore, $x \in X^{(1, \dots 1)}_{\GL_n}$. The same argument then shows this is also true for $X_{\SL_n}$.
\end{remark}
\begin{defn}[Centers] \label{tilde_centers_def}
    Recall that $\ZZ^{3}(L_{\GL_{n}, \lambda}) = \prod^{l}_{i = 1} \mathbb{G}^{3}_{m}$. 
    \begin{enumerate}
        \item \textbf{Centers for} $\GL_n$: \\
         We define 
\begin{equation} \label{tilde_z}
    \widetilde{\ZZ}^{3}(L_{\GL_n, \lambda}) = \{  (x_{i}, \dots , x_{l}) \in \ZZ^{3}(L_{\GL_{n},\lambda}) \mid x_{i} \neq x_{j} \text{ for } i \neq j\}.
\end{equation}
We then define the space $\widetilde{\ZZ}^{3,g}(L_{\GL_n , \lambda})$ as the space of ordered triples of matrices $(D_{1},D_{2},D_{3})$ in the centre $\ZZ^{3}(L_{\GL_n , \lambda})$ such that at least one of the $D_{i}$ have distinct eigenvalues according to the blocks indexed by $\lambda$.
\item \textbf{Centers for} $\SL_n$: \\
We define  
\begin{align} \label{tilde_center_sln}
    \widetilde{\ZZ}^{3}(L_{\SL_n, \lambda}) &= \widetilde{\ZZ}^{3}(L_{\GL_n, \lambda}) \cap L^{3}_{\SL_{n}, \lambda} \\
    \widetilde{\ZZ}^{3,g}(L_{\SL_n, \lambda}) & = \widetilde{\ZZ}^{3,g}(L_{\GL_n, \lambda}) \cap L^{3}_{\SL_{n}, \lambda}.
\end{align}
\item \textbf{Centers for} $\PGL_n$: \\
\begin{align*}
    \widetilde{\ZZ}^{3}(L_{\PGL_n, \lambda}) & = \widetilde{\ZZ}^{3}(L_{\SL_n, \lambda})/ \mu^{3}_{n} \\
    \widetilde{\ZZ}^{3,g}(L_{\PGL_n, \lambda}) & = \widetilde{\ZZ}^{3,g}(L_{\SL_n, \lambda})/ \mu^{3}_{n}.
\end{align*}
    \end{enumerate}
\end{defn}
These are all open subvarieties of $\ZZ^{3}(L_{G , \lambda})$ and hence have dimension $\dim \ZZ^{3}(L_{G , \lambda})$.
\begin{ex} \label{z_tilde}
Let us consider $\GL_4$ with $\lambda = (2,2)$ then $L_{\GL_4 , \lambda} = \GL_2 \times \GL_2$ and the space  $\widetilde{\ZZ}^{3}(L_{\GL_{4}, \lambda})$ is
    \begin{align*}
        \biggl \{ \begin{bmatrix}
    \gamma^{x}_{1} & 0 & 0 & 0   \\
     0 &  \gamma^{x}_{1} & 0 & 0 \\
      0 &  0 &  \gamma^{x}_{2} & 0 \\
      0 &  0 &  0 & \gamma^{x}_{2} \\
  \end{bmatrix} ,    \begin{bmatrix}
    \gamma^{y}_{1} & 0 & 0 & 0   \\
     0 &  \gamma^{y}_{1} & 0 & 0 \\
      0 &  0 &  \gamma^{y}_{2} & 0 \\
      0 &  0 &  0 & \gamma^{y}_{2} \\
  \end{bmatrix} ,   \begin{bmatrix}
    \gamma^{z}_{1} & 0 & 0 & 0   \\
     0 &  \gamma^{z}_{1} & 0 & 0 \\
      0 &  0 &  \gamma^{z}_{2} & 0 \\
      0 &  0 &  0 & \gamma^{z}_{2} \\
  \end{bmatrix} \mid  
   (\gamma^{x}_{i} , \gamma^{y}_{i} , \gamma^{z}_{i}) \neq (\gamma^{x}_{j} , \gamma^{y}_{j} , \gamma^{z}_{j}) \text{ if } i \neq j \biggr \}.
\end{align*}
     We can see that the relative Weyl group $W_{L_{\lambda}} = \So_2$ acts freely on $\widetilde{\ZZ}^{3}(L_{\GL_{4}, \lambda})$ because of this condition. In this case $\widetilde{\ZZ}^{3,g}(L_{\GL_4, \lambda}) = \widetilde{\ZZ}^{3}(L_{\GL_4, \lambda})$ but generally these two spaces will be different.  \par To get the space $\widetilde{\ZZ}^{3}(L_{\SL_4 , \lambda})$ we would impose that all the matrices are of determinant $1$. \par
Let $n=3$. We will illustrate the difference between $\widetilde{\ZZ}^{3}(L_{\GL_3 , (1,1,1)})$ and $\widetilde{\ZZ}^{3,g}(L_{\GL_3 , (1,1,1)})$. Note that the matrix
\begin{equation*}
    (\begin{bmatrix}
    \gamma^{x}_{1} & 0 & 0    \\
     0 &  \gamma^{x}_{1} & 0  \\
      0 &  0 &  \gamma^{x}_{1}  \\
  \end{bmatrix} ,    \begin{bmatrix}
    \gamma^{y}_{1} & 0 & 0    \\
     0 &  \gamma^{y}_{1} & 0  \\
      0 &  0 &  \gamma^{y}_{2}  \\
  \end{bmatrix} ,   \begin{bmatrix}
    \gamma^{z}_{1} & 0 & 0    \\
     0 &  \gamma^{z}_{2} & 0  \\
      0 &  0 &  \gamma^{z}_{2}  \\
  \end{bmatrix})
\end{equation*}
with $\gamma^{x}_{1} \neq \gamma^{y}_{2} \neq \gamma^{y}_{1} \neq \gamma^{z}_{1} \neq \gamma^{z}_{2}$ is contained in $\widetilde{\ZZ}^{3}(L_{\GL_3 , (1,1,1)})$ but is \emph{not} contained in $\widetilde{\ZZ}^{3,g}(L_{\GL_3 , (1,1,1)})$ since none of the three matrices have distinct eigenvalues according to the partition $(1,1,1)$.
\end{ex}
\begin{lem} \label{tilde_center_pull_lemma}
Let $G = \GL_n, \SL_n$. Consider the induced map $\ZZ^{3}(L_{G, \lambda}) \subseteq H^{3}_{G} \to X_{G}$. We have a diagram where both squares are pullbacks
 \begin{equation} \label{tilde_center_pullbacks}
\begin{tikzcd}
	{\widetilde{\ZZ}^{3,g}(L_{G, \lambda})} & {\widetilde{\ZZ}^{3}(L_{G, \lambda})} & {\ZZ^{3}(L_{G, \lambda})} \\
	{X^{\lambda,g}_{G}} & {X^{\lambda}_{G}} & {\overline{X}^{\lambda}_{G}}
	\arrow[from=1-1, to=1-2]
	\arrow["{\theta_{\lambda}}"', from=1-1, to=2-1]
	\arrow[from=1-2, to=1-3]
	\arrow["\varpi", from=1-2, to=2-2]
	\arrow["\pi", from=1-3, to=2-3]
	\arrow[from=2-1, to=2-2]
	\arrow[from=2-2, to=2-3]
\end{tikzcd}
    \end{equation}
    with $\varpi$ and $\theta_{\lambda}$ $W_{L_{\lambda}} \colon 1$ covers. \par
This implies that the strata $X^{\lambda}_{G}$ are smooth.
\end{lem}
\begin{proof}
We will first show that the image of the restriction of the quotient map $H^{3}_{G} \to X_{G}$ to $\ZZ(L_{G, \lambda})$ is  given by $\overline{X}^{\lambda}_{G}$. Because the quotient map is continuous, closed and surjective, we can compute the closures of $X^{\lambda}_{G}$ by pulling back via the quotient map $H^{3}_{G} \to X_{G}$ and computing the closure there. In particular, using this we can show that the closure is given by 
\begin{align*}
\overline{X}^{\lambda}_{\GL_{n}} & = \{ x=\sum^{n}_{i=1} \lambda_{i}x_{i} \in \So^{n}_{\lambda} \mathbb{G}^{3}_{m}   \} \\
    \overline{X}^{\lambda}_{\SL_{n}} & = \{ x=\sum^{n}_{i=1} \lambda_{i}x_{i} \in \So^{n}_{\lambda} \mathbb{G}^{3}_{m} \mid  x \in \SL^{3}_{n}   \}. 
\end{align*}
Recall that by $x \in \SL^{3}_{n}$ we mean that considered as a triple of diagonal matrices $x$ is in $\SL^{3}_{n}$. In particular, in the closure there is no condition that the $x_{i}$ are distinct. Note that 
\begin{equation*}
    \ZZ^{3}(L_{\GL_n , \lambda}) = \prod^{l}_{i=1} \mathbb{G}^{3}_{m}, \text{ and } H^{3}_{\GL_n} = \prod^{n}_{i=1} \mathbb{G}^{3}_{m}.
\end{equation*}
Then the map $\ZZ^{3}(L_{\GL_n , \lambda}) \to H^{3}_{\GL_n}$ is the same as the map 
\begin{equation*}
    \prod^{l}_{i=1} \mathbb{G}^{3}_{m} \to \prod^{n}_{i=1} \mathbb{G}^{3}_{m}
\end{equation*}
given by 
\begin{equation*}
    (x_{1}, \dots x_{l}) \mapsto  ( \underbrace{x_{1}, \dots , x_{1} , }_{\lambda_{1}\text{ times}}  \underbrace{x_{2}, \dots , x_{2} , }_{\lambda_{2}\text{ times}} \dots \underbrace{x_{l} , \dots , x_{l}}_{\lambda_{l}\text{ times}} ).
\end{equation*}
We can then directly compute that the image of $\ZZ^{3}(L_{\GL_n, \lambda})$ under the quotient map is $ \overline{X}^{\lambda}_{\GL_n}$. By restricting to matrices with determiant $1$ we will then get the same for $\SL_n$. From the description of the map $\ZZ^{3}(L_{G, \lambda}) \to H^{3}_{G}$ we can then see that restricting to $X^{\lambda}_{G}$ implies that the preimage under the quotient is exactly the space $\widetilde{\ZZ}^{3}(L_{G, \lambda})$. Note that the condition defining $\widetilde{\ZZ}^{3}(L_{G, \lambda})$ also ensures that $W_{L_{\lambda}}$ acts freely on $\widetilde{\ZZ}^{3}(L_{G, \lambda})$ and hence on the fiber of the map $\varpi$. This proves that $\varpi$ is a $W_{L_{\lambda}} \colon 1$ cover. The fact that $\theta_{\lambda}$ is a $W_{L_{\lambda}} \colon 1$ cover then follows immediately by  restriction. \par
    Since $\widetilde{\ZZ}^{3}(L_{G, \lambda})$ is smooth then $X^{\lambda}_{G}$ is smooth since it is covered by a smooth variety.
\end{proof}
   For discussion of the stratification for $\PGL_n$ see Definition \ref{strat_pgl_n}. For the analogue of the previous Lemma see Lemma \ref{tilde_center_pull_lemma_pgln} in Subsection \ref{pgln_subsect}.
\begin{remark}
Note that the map $X_{\SL_n} \to X_{\PGL_n}$ does \emph{not} preserve the stratifications on both sides. See Example \ref{sln_to_pgln_ex}. 
\end{remark}
With preliminaries out of the way we start with step 1 and properties of the exponential map.
\subsection{Step 1: Deducing purity using the exponential}
 In this subsection we will use the exponential map to deduce several strong properties of the DT sheaf $\varphi_{\GL_{n}}$ from the additive version as in subsection \ref{jordan}. Before we start let us state two results we will use repeatedly in our arguments.
 \begin{lem}[Descending \'etale morphisms] \label{descending_etale}
Let $f \colon X \to Y$ be a map of affine schemes that is equivariant with respect to a homorphism of finite groups $\phi \colon G \to H$ such that at any point $x \in X$ we have an isomorphism of stabilisers $G_{x} \to H_{f(x)}$. Then for any point $x \in X$ such that $f$ is \'etale at $x$, the induced morphism $\overline{f} \colon X /\!\!/ G \to Y /\!\!/H$ is \'etale at the image of $x$ under the quotient map $X \to X /\!\!/ G$. Furthermore, this result also holds for $X_{\an}$ and $Y_{\an}$ the associated complex analytic spaces and $g \colon X_{\an} \to Y_{\an}$. Here the map $g$ need not be induced from an algebraic map $g^{'} \colon X \to Y$. 
 \end{lem}
 \begin{proof}
     We follow the proof given in \cite[Remark 4.4.4]{alper_moduli_lect}. There, the proof for affine schemes is already explained. Let $f \colon X_{\an} \to Y_{\an}$ be an \'etale map of complex analytic spaces equivariant with respect to the homomorphism $G \to H$ that preserves stabilisers at $x$. Now since $f$ is \'etale at $x \in X_{\an}$, the induced map $\widehat{\mathcal{O}}_{X,x} \to \widehat{\mathcal{O}}_{Y,f(x)}$ on formal completions is a $G \to H$ equivariant isomorphism, which also gives the isomorphism 
     \begin{equation} \label{invariants_formal_iso_eq}
         \widehat{\mathcal{O}}^{G_{x}}_{X,x} \to \widehat{\mathcal{O}}^{H_{f(x)}}_{Y,f(x)}
     \end{equation}
     Using that formal completions of $X /\!\!/ G $ and $(X /\!\!/ G)_{\an}$ agree we get that the formal completion of the quotient $X_{\an}/\!\!/G$ at $x$ is
\begin{equation}
    \widehat{\mathcal{O}}_{(X /\!\!/G)_{\an} ,x} \cong \widehat{\mathcal{O}}_{X /\!\!/G ,x}  \cong \widehat{\mathcal{O}}^{G_{x}}_{X,x}.
\end{equation}
The last isomorphism follows from \cite[Exercise 4.2.15]{alper_moduli_lect}. Using equation \eqref{invariants_formal_iso_eq} we get that the map $ \widehat{\mathcal{O}}_{(X /\!\!/G)_{\an} ,x} \to  \widehat{\mathcal{O}}_{(Y /\!\!/H)_{\an} ,f(x)}$ is an isomorphism thus giving that the induced map $f \colon X_{\an} /\!\!/G \to Y_{\an} /\!\!/ H$ is \'etale at $x$.
 \end{proof}
We also will need the following result about commuting diagrams of stacks and good moduli spaces.
 \begin{prop} \cite[Proposition 6.8]{alper_local_quotient} \cite[Proposition 6.3.30]{alper_moduli_lect} \label{alper_prop}
     Consider a commutative square of algebraic stacks $\mathcal{X}$, $\mathcal{Y}$ and their respective good moduli spaces $X$ and $Y$
     \begin{equation}
\begin{tikzcd}
	{\mathcal{X}} & {\mathcal{Y}} \\
	X & Y
	\arrow["{f^{'}}", from=1-1, to=1-2]
	\arrow["{g^{'}}"', from=1-1, to=2-1]
	\arrow["f"', from=2-1, to=2-2]
	\arrow["g", from=1-2, to=2-2]
\end{tikzcd}
     \end{equation}
Assume that $f^{'} \colon \mathcal{X} \to \mathcal{Y}$ is a separated, representable morphism of noetherian stacks with affine diagonal. If we have 
\begin{enumerate}
    \item $f^{'}$ is \'etale 
    \item $f^{'}$ maps closed points to closed points
    \item $f^{'}$ induces an isomorphism on stabilisers at all closed points
\end{enumerate}
then $f$ is \'etale and the square is cartesian.
 \end{prop}
 
 Recall the \'etale loci in Definition \ref{etale_locus_c3gln}. Define 
 \begin{equation}
    \So^{n,\et }(\mathbb{G}^{2}_{m} \times \mathbb{G}_{a}) = \{(x_{j},y_{j},z_{j})_{ 1 \leq j \leq n} \in \So^{n}(\mathbb{G}^{2}_{m} \times \mathbb{G}_{a}) \mid z_{j} - z_l \neq 2 \pi i k \text{ for } j\neq l \text{ and }  k \in \mathbb{Z} \setminus \{0 \} \}
 \end{equation}
 In particular, we require that the $z$ part of the element satisfies the condition in equation \eqref{etale_cond}.
 We will now prove three technical lemmas on the behaviour of the exponential map. These lemmas are necessary as we cannot immediately use the above proposition since we are working with analytic stacks.
 \begin{lem} \label{stabiliser_exp_diag}
     The exponential map $\exp \colon \mathfrak{gl}^{\et}_{n} \to \GL_n$ preserves stabilisers of diagonalisable matrices in $\mathfrak{gl}^{\et}_{n}$ under the conjugation action of $\GL_n$.
 \end{lem}
 \begin{proof}
     It is enough to check the claim on any representative of a diagonalizable matrix in the orbit under conjugation since the stabilisers are isomorphic under conjugation. Therefore, we may assume that $D \in \mathfrak{gl}^{\et}$ is diagonal with $D = \diag(\gamma_{1} I_{\lambda_{1}}, \dots , \gamma_{l} I_{\lambda_{l}})$ corresponding to some partition $\lambda$ of $n$ of length $l$ and $\gamma_{i} \neq \gamma_{j}$ for $i \neq j$. Now the stabiliser of $D$ only depends on the partition $\lambda$ and not on the values $\gamma_{i}$ so the only way it can change is if $\exp \gamma_{i} = \exp \gamma_{j}$ but this implies that
     \begin{equation}
         \gamma_{i} - \gamma_{j} = 2 \pi i k 
     \end{equation}
     for $k \neq 0$. This would contradict the condition in equation \eqref{etale_cond} so the exponential preserves stabilisers.
 \end{proof}
 \begin{lem} \label{one_variab_exp_lem}
     The commutative square
     \begin{equation}
\begin{tikzcd}
	{\mathfrak{gl}^{\et}_{n}} & {\GL_n} \\
	{\So^{n} \mathbb{G}^{\et}_{a}} & {\So^{n} \mathbb{G}_{m}}
	\arrow["\exp", from=1-1, to=1-2]
	\arrow[from=1-1, to=2-1]
	\arrow[from=1-2, to=2-2]
	\arrow["\exp"', from=2-1, to=2-2]
\end{tikzcd}
     \end{equation}
     is a pullback of complex analytic spaces. 
 \end{lem}
 \begin{proof}
     Recall that we can view $\So^{n} \mathbb{G}^{\et}_{a}$ as the space of diagonal $n \times n$ matrices up to permutation that satisfy the condition in equation \eqref{etale_cond} and $\So^{n} \mathbb{G}_{m}$ as the space of diagonal matrices up to permutation with non-zero entries. The pullback $P = \So^{n} \mathbb{G}^{\et}_{a} \times_{\So^{n} \mathbb{G}_{m}} \GL_{n}$ has points
     \begin{equation}
         \{ (D, A) \in \So^{n} \mathbb{G}^{\et}_{a} \times \GL_n \mid \exp D = A_{ss} \}
     \end{equation}
     where $A_{ss}$ is the diagonal $n \times n$ matrix that contains the eigenvalues of $A$ up to multiplicity. We can define a map
     \begin{align*}
         f \colon \mathfrak{gl}^{\et}_{n} & \to P \\
         x & \mapsto (x_{ss}, \exp x)
     \end{align*}
     where $x_{ss}$ is the diagonal matrix containing eigenvalues of $x$ up to multiplicity.  The map $f$ is surjective since the exponential map is still surjective once restricted to $\mathfrak{gl}^{\et}_{n}$. We will now show that $f$ is injective. Assume that we have $x, y \in \mathfrak{gl}^{\et}_{n}$ such that
     \begin{equation}
         (x_{ss}, \exp x) = (y_{ss}, \exp y)
     \end{equation}
Now consider Jordan-Chevalley decompositions of $x$ and $y$
\begin{equation}
    x = s_{x} + n_{x} \quad y = s_{y} + n_{y}
\end{equation}
with $s_{-}$ diagonalisable $n_{-}$ nilpotent and 
\begin{equation*}
    s_{x}n_{x} = n_{x} s_{x} \quad s_{y}n_{y} = n_{y}s_{y}.
\end{equation*}
Then we can take the exponential to get the equation
\begin{equation}
    \exp(s_{x}) \exp(n_{x})= \exp x = \exp y = \exp(s_{y}) \exp(n_{y})  
\end{equation}
here $\exp(s_{x})$ is still diagonalisable and $\exp(n_{x})$ is unipotent so we can use the uniqueness of the multiplicative Jordan-Chevalley decompositions of $\exp x = \exp y$  to conclude that 
\begin{equation} \label{ss_part_exp}
    \exp s_{x} = \exp s_{y} \quad \exp n_{x} = \exp n_{y}.
\end{equation}
The exponential map defines a bijection between the nilpotent cone and the unipotent cone which implies that $n_{x} = n_{y}$. Now because $x_{ss} = y_{ss}$ we have that there exist invertible matrices $S_{1}, S_{2}$ such that
\begin{equation*}
    S^{-1}_{1} s_{x} S_{1} = x_{ss} = y_{ss} = S^{-1}_{2} s_{y} S_{2}
\end{equation*}
Then exponentiating, using the fact that the exponential commutes with conjugation and equation \eqref{ss_part_exp} we get
\begin{equation*}
    \exp s_{x} = S_{1} S^{-1}_{2} \exp s_{x} S_{2} S^{-1}_{1}.
\end{equation*}
However, we know from Lemma \ref{stabiliser_exp_diag} that the exponential preserves stabilisers, which implies that 
\begin{align*}
    s_{x} & = S_{1} S^{-1}_{2}  s_{x} S_{2} S^{-1}_{1} \\
    S^{-1}_{1}s_{x}S_{1} & =  S^{-1}_{2}  s_{x} S_{2}  \\
    S^{-1}_{2}  s_{x} S_{2} & = S^{-1}_{2}  s_{y} S_{2}
\end{align*}
hence $s_{x} = s_{y}$ and therefore $x = y$.
\end{proof}
Consider the following diagram 
 \begin{equation} \label{exp_diagram}
\begin{tikzcd}
	{\CC_{3}(\GL^{2}_{n}, \mathfrak{gl}^{\et}_{n}) / \GL_n} & {\CC_{3}(\GL^{2}_{n}, \mathfrak{gl}_{n}) / \GL_n } & {\Loc_{\GL_n}} \\
	{\So^{n,\et }(\mathbb{G}^{2}_{m} \times\mathbb{G}_{a})} & {\So^{n}(\mathbb{G}^{2}_{m} \times\mathbb{G}_{a})} & {\So^{n}(\mathbb{G}^3_{m})}
	\arrow["\exp", from=1-2, to=1-3]
	\arrow["{\pi_{\mathfrak{gl}_{n}}}"', from=1-2, to=2-2]
	\arrow["{\pi_{\GL_{n}}}", from=1-3, to=2-3]
	\arrow["\exp"', from=2-2, to=2-3]
	\arrow[hook, from=1-1, to=1-2]
	\arrow[from=1-1, to=2-1]
	\arrow[hook, from=2-1, to=2-2]
\end{tikzcd}
 \end{equation}
 \begin{lem}\label{pullback_exp_int}
     The outer commutative square in equation \eqref{exp_diagram} is a pullback diagram and the horizontal compositions are \'etale. Furthermore, the exponential map is surjective restricted to the \'etale locus. 
 \end{lem}
 \begin{proof}
Denote $\CC_{3}(\GL^{2}_{n},\mathfrak{gl}^{\et}_{n})$ by $M_{a}$, $\Loc^{\ff}_{\GL_{n}}(T^{3}) = \CC_{3}(\GL_{n})$ by $M_{m}$, the map $M_{a} \to \So^{n}(\mathbb{G}^{2}_{m} \times \mathbb{G}^{\et}_{a})$ by $\pi_{a} \colon M_{a} \to X_{a}$ and  $\So^{n}(\mathbb{G}^{3}_{m})$ by $X_{m}$. Consider the diagram where every square is a pullback

\begin{equation}
\begin{tikzcd}
	{M_{a}} \\
	{P_{2}} & {M_{m}} \\
	{P_{1}} & {M_{m}/\GL_{n}} \\
	{X_{a}} & {X_{m}}
	\arrow[from=3-2, to=4-2]
	\arrow["\exp"', from=4-1, to=4-2]
	\arrow[from=2-2, to=3-2]
	\arrow["{p^{'}}", from=2-1, to=2-2]
	\arrow["p", from=3-1, to=3-2]
	\arrow[from=3-1, to=4-1]
	\arrow[from=2-1, to=3-1]
	\arrow["\exp", from=1-1, to=2-2]
	\arrow["{\pi_{1}}", from=1-1, to=2-1]
\end{tikzcd}
\end{equation}
First using Lemma \ref{descending_etale} we can show that $\exp \colon X_{a} \to X_{m}$ is an \'etale map since on the chosen locus it preserves stabilizers. Then it follows that the maps $p$ and $p^{'}$ in the diagram are also \'etale. Using the universal property of $P_{2}$, we get a map $\pi_{1} \colon M_{a} \to P_{2}$. From the proof of Theorem \ref{exp_dcrit} we know that $\exp \colon M_{a} \to M_{m}$ is \'etale. The map $\pi_{1}$ is then also \'etale by the $2$ out of $3$ property for \' etale maps.  Write $A_{ss}$ for the diagonal matrix containing the eigenvalues of $A$. The map $\pi_{1}$ is defined 
\begin{align*}
    \pi_{1} \colon M_{a} & \to P_{2} \\
    (a_1, a_2, A_3) &\mapsto ((a_{1 ss}, a_{2ss}, A_{3ss}),(a_{1} ,a_{2}, \exp A_{3}) )
\end{align*}
as the identity on the first two matrices and in the same way as the map $f$ in the proof of Lemma \ref{one_variab_exp_lem} on the matrix we are exponentiating. Since we are only exponentiating along one of the $3$ pairwise commuting matrices the fact that $\pi_{1} \colon M_{a} \to P_{2}$ is a bijection on points immediately reduces to Lemma \ref{one_variab_exp_lem}.  Then since $\pi_1$ is an \'etale bijection it must be an isomorphism. We can use the natural maps $M_{a}/\GL_{n} \to M_{m} / \GL_{n}$ and $M_{a}/ \GL_{n} \to X_{a}$ to define a map $M_{a} / \GL_{n} \to P_{1}$. Then we have the diagram
\begin{equation}
\begin{tikzcd}
	{M_{a}} & {M_{a}} & {M_{m}} \\
	{M_{a}/\GL_{n}} & {P_{1}} & {M_{m}/\GL_{n}}
	\arrow["\operatorname{id}", from=1-1, to=1-2]
	\arrow[from=1-1, to=2-1]
	\arrow[from=1-2, to=1-3]
	\arrow[from=1-2, to=2-2]
	\arrow[from=1-3, to=2-3]
	\arrow[from=2-1, to=2-2]
	\arrow[from=2-2, to=2-3]
\end{tikzcd}
\end{equation}
here the big square is a pullback and the righmost square is a pullback, which by the $2$ out of $3$ property for pullbacks implies that the leftmost square is a pullback. Now we can conclude that we have an isomorphism $M_{a}/ \GL_{n} \to P_{1}$ by using that isomorphisms are local under smooth maps. 
 \end{proof}
 Using the fact that we have a pullback square and an \'etale cover we can use the following Proposition and Theorem \ref{exp_dcrit} to transfer information about the additive DT sheaf to the multiplicative DT sheaf.
 \begin{prop}\cite[Proposition 4.5]{ben2015darboux}
Let $f \colon X \to Y$ be a smooth map of oriented $d$-critical loci of relative dimension $n$. Then we have the natural isomorphism $\varphi_X \cong f^{*}[n] \varphi_{Y}$.
 \end{prop}
 The above result is proven for algebraic $d$-critical loci but the result also holds for complex analytic $d$-critical loci.
 \begin{lem} \cite[Lemma 2.2]{davison2023purity} \label{purity_local}
Let $f \colon X \to Y$ be an \'etale or smooth map of complex analytic spaces. Let $\mathcal{F}$ be a mixed Hodge module on $Y$. Then $\mathcal{F}$ is pure if and only if $f^{*} \mathcal{F}$ is pure.
 \end{lem}
 Putting together the last three claims we get:
 \begin{cor} \label{purity_GL}
Write $\varphi_{\mathfrak{gl}_{n}}$ for the DT sheaf on $\CC_{3}(\mathfrak{gl}_{n})/\GL_n$ and $\pi_{\mathfrak{gl}_n} \colon \CC_{3}(\mathfrak{gl}_{n})/\GL_n \to  \So^{n} \mathbb{G}_{a}^{3}$, the good moduli space map. We can deduce the following properties
\begin{enumerate}
    \item $\exp^{*}\varphi_{\GL_{n}} \cong \varphi_{\mathfrak{gl}_{n}}|_{\CC_{3}(\GL^{2}_{n},\mathfrak{gl}^{\et}_{n})/ \GL_n} $
    \item $\exp^* \pi_{\GL_{n}*} \varphi_{\GL_{n}} \cong (\pi_{\mathfrak{gl}_{n}*} \varphi_{\mathfrak{gl}_{n}})|_{\So^{n,\et }(\mathbb{G}^{2}_{m}\times\mathbb{G}_{a})}$
    \item $\pi_{\GL_{n}*} \varphi_{\GL_{n}}$ is a pure complex of mixed Hodge modules 
    \item $\exp^* \BPS_{\GL_{n}} \cong \BPS_{\mathfrak{gl}_{n}}|_{\So^{n,\et }(\mathbb{G}^{2}_{m}\times\mathbb{G}_{a})}$.
\end{enumerate}
Note that the diagonal $\mathbb{G}^{2}_{m} \times \mathbb{G}_{a} \subseteq \So^{n}(\mathbb{G}^{2}_{m} \times \mathbb{G}_{a})$ also lives inside $\So^{n, \et}(\mathbb{G}^{2}_{m} \times \mathbb{G}_{a})$ since the condition in equation \eqref{etale_cond} holds on the diagonal.
 \end{cor}
 \begin{cor}\label{bps_const}
     $\BPS_{\GL_{n}}$ is the sheaf $\Delta_{*} \mathbb{Q}_{\ZZ^{3}(\GL_{n})} [ \dim \ZZ^{3}(\GL_{n}) ] $ with $\Delta \colon \ZZ^{3}(\GL_{n}) \to X_{\GL_n}$.
 \end{cor}
 \begin{proof}
     Using the previous corollary we know that $\BPS_{\GL_{n}}$ is supported on the diagonal $\mathbb{G}^{3}_{m}$ in $\So^{n} \mathbb{G}^{3}_m$ because $\BPS_{\mathfrak{gl}_{n}}$ is supported on the diagonal $\mathbb{G}^{3}_{a}$ in $\So^{n} \mathbb{G}^{3}_a$ via the additive support lemma \ref{support_lemma_quivers}. Therefore, $\BPS_{\GL_{n}}$ must be some local system $\mathcal{L}$ of rank $1$ on $\mathbb{G}^{3}_{m}$.  Furthermore, using Proposition \ref{center_Wact} we get that $\mathbb{G}^{3}_{m}$ acts on $\varphi_{\GL_n}$ and thus also on on each perverse piece of $\pi_{\GL_{n}*} \varphi_{\GL_{n}}$. In particular, $\mathcal{L}$ is $\mathbb{G}^{3}_{m}$-equivariant. $\mathbb{G}^{3}_{m}$ acts transitively on $\mathbb{G}^{3}_{m}$ therefore, we can conclude that $\mathcal{L}$ must be the trivial local system of rank $1$.
 \end{proof}
 \begin{cor}[ Support lemma for $\GL_{n}$] \label{support_locsys}
 
We have a decomposition 
\begin{equation}
    \pi_{\GL_{n}*} \varphi_{\GL_{n}} \cong \bigoplus_{i \geq 1} \bigoplus_{\lambda}\IC_{X^{\lambda}_{\GL_{n}}}(\mathcal{L}^{\lambda}_{i})[-i]
\end{equation}
for some local systems $\mathcal{L}^{\lambda}_{i}$ on $X^{\lambda}_{\GL_n}$. The index $i$ corresponds to the perverse cohomology degree and $\lambda$ is a partition of $n$ giving the corresponding stratum $X^{\lambda}_{\GL_n}$
 \end{cor}
 \begin{proof}
     Firstly, we can use the purity of $\pi_{\GL_n *} \varphi_{\GL_n}$ to obtain a decomposition into $\IC$ sheaves. To get the statement of the corollary we now have to compute the supports. \par 
     The computation of supports is a direct consequence of the support Lemma \ref{additive_computation_lemma} and the decomposition in equation \eqref{add_double_decomp} in the additive case. We can restrict all the perverse sheaves in the decomposition of $\pi_{\mathfrak{gl}_{n}*}\varphi_{\mathfrak{gl}_{n}}$ on $\So^{n}(\mathbb{G}^{3}_{a})$ to the open set $\So^{n,\et }(\mathbb{G}^{2}_{m} \times \mathbb{G}_{a})$. Now consider a summand $\mathcal{F}$ of the $k$-th perverse cohomology in the Saito decomposition of $\pi_{\GL_{n}*} \varphi_{\GL_{n}}$. We will prove that $\mathcal{F}$ is an $\IC$ sheaf supported on $X^{\lambda}_{\GL_n}$ for some $\lambda$ using Lemma \ref{ic_achar}. 
     Pulling back we claim that 
     \begin{equation}
         \exp^{*} \mathcal{F} \cong \IC_{X^{\lambda}_{\mathfrak{gl}_{n}}} (\mathcal{L}).
     \end{equation}
     The above equation follows since $\exp$ is $t$-exact and so $\exp^{*} \mathcal{F}$ is summand of the $k$-th perverse cohomology  of $\pi_{\mathfrak{gl}_{n} *} \varphi_{\mathfrak{gl}_{n}}$. Then as a consequence of  part $(4)$ of Lemma \ref{additive_computation_lemma} and Proposition \ref{coh_add_prop} we can use equation \eqref{add_double_decomp} where we computed that all the summands in the perverse cohomology of $\pi_{\mathfrak{gl}_{n} *} \varphi_{\mathfrak{gl}_{n}}$ are $\IC$ sheaves with supports $X^{\lambda}_{\mathfrak{gl}_n}$. \par  There is a pullback
     \begin{equation}
\begin{tikzcd}
	{\So^{n,\et }_{\lambda}(\mathbb{G}^{2}_{m} \times \mathbb{G}_{a})} & {\So^{n}_{\lambda}\mathbb{G}^{3}_{m}} \\
	{\So^{n,\et }(\mathbb{G}^{2}_{m} \times \mathbb{G}_{a})} & {\So^{n}\mathbb{G}^{3}_{m},}
	\arrow["j", hook, from=1-2, to=2-2]
	\arrow["{\exp^{'}}", from=1-1, to=1-2]
	\arrow["\exp"', from=2-1, to=2-2]
	\arrow["{j^{'}}"', hook, from=1-1, to=2-1]
\end{tikzcd}
     \end{equation}
Considering the pullback of $j^{*} \mathcal{F}$ under $\exp^{'}$ we get
\begin{equation}
    \exp^{'*}j^{*} \mathcal{F} \cong j^{'*}\exp^{*} \mathcal{F} \cong j^{'*} \IC_{X^{\lambda}_{\mathfrak{gl}_{n}}}(\mathcal{L}) \cong \mathcal{L}.
\end{equation}
 Therefore, $j^{*} \mathcal{F}$ must be a local system since it is a local system when pulled back by an \'etale map. Similarly, $\mathcal{F}$ will have no subobjects or quotients on $\overline{\So^{n} \mathbb{G}^{3}_{m}} \setminus \So^{n}_{\lambda} \mathbb{G}^{3}_{m}$ so we have  $\mathcal{F} = \IC_{X^{\lambda}_{\GL_{n}}} (j^{*} \mathcal{F})$.
 \end{proof}
 In light of these properties, the only thing we have to calculate is the local systems appearing in the Saito decomposition of the pure complex of mixed Hodge modules $\pi_{\GL_{n}*}\varphi_{\GL_{n}} $. 
 \subsection{Step 2: Reduction to Levis  }
 Let us start with computing the terms in equation \eqref{coh_sln}. To do this we need to describe the BPS sheaves on Levis. \par The action of $W_{L_{\lambda}}$ on $L_{G, \lambda}$ induces an action on the  equivariant cohomology $\HHf( \B \ZZ(L_{G, \lambda}))$. We start by describing the centres of the Levis $L_{\SL_n , \lambda}$ and $L_{\PGL_n , \lambda}$ more explicitly. We will also more explicitly compute the action of the relative Weyl group $W_{L_{\lambda}}$ on the equivariant cohomology of the centre.
 \begin{lem}[Computation of centres of Levis] \label{levi_centers}
 Let $L_{\SL_{n}, \lambda}$ and $L_{\PGL_{n}, \lambda}$ be the Levis corresponding to the partition $\lambda$ as in  \eqref{levi_partitions}. We have $\ZZ(L_{\PGL_{n}, \lambda}) \cong \mathbb{G}^{l-1}_{m}$ and $\ZZ(L_{\SL_{n}, \lambda}) \cong \mathbb{G}^{l-1}_{m} \times \mu_{\operatorname{gcd}(\lambda_{1}, \dots , \lambda_{l})}$.
 \end{lem}
 \begin{proof}
We have $\ZZ(L_{\GL_{n}, \lambda}) = \ZZ(\prod^{l}_{i} \GL_{\lambda_{i}}) = \prod^{l}_{i} \ZZ(\GL_{\lambda_{i}})$. Let us start with $\PGL_{n}$. Take $X \in \ZZ(L_{\PGL_{n}, \lambda})$ a block matrix with representatives in $\GL_{\lambda_{i}}$ given by $\widetilde{X}_{i}$. $X$ being central means that we have $[\widetilde{X}_{i},Y_{i}] = t_{Y_{i}} I$ for all $Y_{i} \in \GL_{\lambda_{i}}$ and $t_{Y_i}  \in \mathbb{C}^{*}$. This implies that the image of $\widetilde{X}_{i} \in \PGL_{\lambda_{i}}$ is central and so $\widetilde{X}_{i} = \gamma_{i} I_{\lambda_{i}}$ for some $\gamma \in \mathbb{G}_{m}$.
Therefore $\widetilde{X}$ is of the form $\operatorname{diag}(\gamma_{1} I_{\lambda_{1}}, \dots ,\gamma_{l} I_{\lambda_{l}})$. This shows that $\ZZ(L_{\PGL_{n}, \lambda}) \cong \mathbb{G}^{l-1}_{m}$ since in $\PGL_n$ we can quotient out by one of the $\gamma_{i}$. Now consider $\SL_{n}$
 \begin{equation}
    \ZZ(L_{\SL_{n}, \lambda}) = \{  \operatorname{diag}(x_{1} I_{\lambda_{1}, } \dots x_{l}I_{\lambda_{l}})  \mid \prod^{l}_{i} x^{\lambda_{i}}_{i} = 1  \}. 
 \end{equation}
Note that $\ZZ(L_{\SL_{n} , \lambda})= \ker f_{\lambda}$ where $f_{\lambda}  \colon \mathbb{G}^{l}_{m} \to \mathbb{G}_{m}$ $(\gamma_{1}, \dots , \gamma_{l}) \mapsto \prod_{i} \gamma^{\lambda_{i}}_{i}$. Note that the map $f_{\lambda}$ is determined by the map of characters $f^{*}_{\lambda} \colon \XX^{*}(\mathbb{G}_{m}) \to \XX^{*}(\mathbb{G}^{l}_{m})$, which is given by $\mathbb{Z} \to \mathbb{Z}^{l}$ $1 \mapsto [\lambda_{1}, \dots , \lambda_{l}]$. However, by using Smith normal forms this map is equivalent to the map $\widetilde{f}^{*}_{\lambda} \colon 1 \mapsto [\operatorname{gcd}(\lambda_{1}, \dots , \lambda_{l}),0, \dots , 0]$. Therefore the kernel of $f_{\lambda}$ is the same as the cokernel of $\widetilde{f}_{\lambda}$. Thus $\ZZ(L_{\SL_{n}, \lambda}) \cong \mathbb{G}^{l-1}_{m} \times \mu_{\operatorname{gcd}(\lambda_{1}, \dots , \lambda_{l})}$.
 \end{proof}
 Let $\lambda$ be a non-trivial partition of $n$. Note that we can write
 \begin{equation}
     \HHf(\B \ZZ(L_{G, \lambda})) \cong \XX^{*}(\ZZ(L_{G, \lambda})) \otimes_{\mathbb{Z}} \mathbb{Q}.
 \end{equation}
 $W_{L_{\lambda}}$ acts on 
 $$M_{\lambda} = \XX^{*}(\ZZ(L_{\GL_n, \lambda})) \otimes_{\mathbb{Z}} \mathbb{Q} \cong \bigoplus^{l}_{i=1} \mathbb{Q}_{{\lambda_i}}$$ by permuting the elements in the blocks $\mathbb{Q}_{{\lambda_i}}$ which have the same subscript.
 \begin{ex}
     Consider $n = 5$ and $\lambda = (2,2,1)$ then $M_{\lambda} = \mathbb{Q}_{2} \oplus \mathbb{Q}_{2} \oplus \mathbb{Q}_{1}$ and $W_{L_{\lambda}} \cong \So_{2}$ acts by permuting $\mathbb{Q}^{\oplus 2}_{2}$ and acts trivially on $\mathbb{Q}_{1}$.
 \end{ex}
 Pick a basis $x_{i} \in \mathbb{Q}_{\lambda_i} \subseteq M_{\lambda}$ and $1 \leq i \leq l$. $M_{\lambda}$ splits as $M^{'}_{\lambda} \oplus \mathbb{Q}\sum^{l}_{i=1} \lambda_{i}x_{i}$ where $\mathbb{Q}\sum^{l}_{i=1} \lambda_{i}x_{i}$ is the trivial $1$-dimensional $W_{L_{\lambda}}$ representation spanned by the element $\sum^{l}_{i=1} \lambda_{i}x_{i}$.  We can write 
 $$\HHf(\B \ZZ(L_{\GL_n, \lambda})) \cong \Sym( M^{'}_{\lambda}[-2]) \otimes \Sym( \mathbb{Q}[-2])$$ as $W_{L_{\lambda}}$ representations. Using the basis $x_{i}$ we have $\HHf(\B \ZZ(L_{\GL_n, \lambda})) \cong \mathbb{Q}[x_{1}, \dots , x_{l}]$ with $x_{i}$ in degree $2$. Tensoring by $\mathbb{Q}$ identifies $\XX^{*}(\ZZ(L_{\SL_n, \lambda}))$ and $\XX^{*}(\ZZ(L_{\PGL_n, \lambda}))$ so we will also identify $\HHf( \B \ZZ(L_{\SL_n , \lambda}))$ and $\HHf( \B \ZZ(L_{\PGL_n , \lambda}))$.
 \begin{lem} \label{coh_centers}
  We have for $G= \SL_n$ or $\PGL_n$, $\HHf(\B \ZZ(L_{G, \lambda})) \cong \Sym( M^{'}_{\lambda}[-2]) \xhookrightarrow{} \HHf(\B \ZZ(L_{\GL_n, \lambda}))$ and $W_{L_{\lambda}}$ acts by restriction. Therefore $\HHf(\B \ZZ(L_{\SL_n, \lambda}))$ and $\HHf(\B \ZZ(L_{\PGL_n, \lambda}))$ are polynomial algebras in $l-1$ variables.
 \end{lem}
 \begin{proof}
     Consider the exact sequence 
     $$1 \to \mu_{n} \to \ZZ(L_{\SL_{n}, \lambda}) \times \mathbb{G}_{m} \to \ZZ(L_{\GL_{n}, \lambda}) \to 1.$$ 
     This short exact sequence induces the following map on characters
     \begin{equation}
         \XX^{*}(\ZZ(L_{\GL_n, \lambda})) \to \XX^{*}(\ZZ(L_{\SL_n, \lambda})) \oplus \mathbb{Z}
     \end{equation}
which becomes an isomorphism after tensoring by $- \otimes_{\mathbb{Z}} \mathbb{Q}$. 
We can show that we have a diagram
\begin{equation}
\begin{tikzcd}
	{\XX^{*}(\ZZ(L_{\GL_n, \lambda})) \otimes_{\mathbb{Z}} \mathbb{Q}} & {\XX^{*}(\ZZ(L_{\SL_n, \lambda})) \otimes_{\mathbb{Z}} \mathbb{Q} \oplus \mathbb{Q}} \\
	{M_{\lambda}} & {M^{'}_{\lambda} \oplus \mathbb{Q} \sum\lambda_{i}x_{i}}
	\arrow["\cong"', from=1-1, to=1-2]
	\arrow["\cong", from=1-1, to=2-1]
	\arrow["\cong"', from=1-2, to=2-2]
	\arrow["\cong", from=2-1, to=2-2]
\end{tikzcd}
\end{equation}
Thus giving an isomorphism $\XX^{*}(\ZZ(L_{\SL_n, \lambda})) \otimes_{\mathbb{Z}} \mathbb{Q} \cong M^{'}_{\lambda}$ and therefore an isomorphism
\begin{equation}
    \HHf(\B \ZZ(L_{\SL_n, \lambda}))  \cong \Sym(M^{'}_{\lambda}[-2]). 
\end{equation}
 \end{proof}
   Consider the decomposition of $\pi_{G*} \varphi_{G}$ into perverse pieces. We can consider the summands with support $X^{\lambda}_{G}$ by Corollary \ref{support_locsys}. Restricting these perverse sheaves to $X^{\lambda}_{G}$ we get local systems. Therefore, over $X^{\lambda}_{G}$ we will have a direct sum
\begin{equation}\label{loc_sys_dir_sum_strata}
    \bigoplus_{i \in \mathbb{Z}} \mathcal{L}^{\lambda}_{i} [ \dim X^{\lambda}_{G}-i] 
\end{equation}
of shifted local systems with $i$ corresponding to the perverse cohomology degree of the associated $\IC$ sheaf. Note that for some $i$, $\mathcal{L}^{\lambda}_{i}$ will be $0$ depending on which stratum $X^{\lambda}_{G}$ we have chosen. We will now proceed to computing these local systems. 

 Let $G= \GL_n ,\SL_n$. Recall the generic loci in Definition \ref{generic_loci}, the stratification in Definition \ref{stratifications} and spaces in Definition \ref{tilde_centers_def}. Note that for the purpose of computing the $\IC$ sheaves we can restrict the local systems $\mathcal{L}^{\lambda}_{i}$ to an open dense to compute them. We will therefore instead compute the local systems $\mathcal{L}^{\lambda, g}_{i} $, which are the restrictions of $\mathcal{L}^{\lambda}_{i}$ to the open dense $X^{\lambda ,g}_{G}$ of the stratum $X^{\lambda}_{G}$. The space $X^{\lambda,g}_{G}$ fits into the following diagram where both squares are pullbacks. The genericity condition is necessary to ensure that the maps $\Theta^{g}$ and $\theta^{g}$ are \'etale. 
\begin{equation}\label{pull_levi}
\begin{tikzcd}
	{\Loc_{G}} & {\Loc^{g}_{L_{G,\lambda}}} \\
	{X_{G}} & {X^{g}_{L_{G,\lambda}}} \\
	{X^{\lambda, g}_{G}} & {\widetilde{\ZZ}^{3,g}(L_{G, \lambda})}
	\arrow["{\pi_{G}}"', from=1-1, to=2-1]
	\arrow["\theta^{g}"', from=2-2, to=2-1]
	\arrow["\Theta^{g}"', from=1-2, to=1-1]
	\arrow["{\pi_{L}}", from=1-2, to=2-2]
	\arrow["{j_{\lambda,G}}", from=3-1, to=2-1]
	\arrow["{\theta_{\lambda}}", from=3-2, to=3-1]
	\arrow["{j_{\lambda,L}}"', from=3-2, to=2-2]
\end{tikzcd}
\end{equation}

\begin{prop} \label{levis_final_splitoff}
Let $G = \GL_n$ or $\SL_n$.
    \begin{enumerate}
    \item The restriction of the map $\theta^{g}$ to $\widetilde{\ZZ}^{3,g}(L_{G, \lambda})$ is given by the $W_{L_{\lambda}} \colon 1$ cover $\theta_{\lambda}$ in the diagram \eqref{tilde_center_pullbacks}. 
    \item Pulling back the pushforward of the DT sheaf we have the isomorphism of $W_{L_{\lambda}}$-equivariant sheaves: 
\begin{equation}\label{levi_red_equation}
        \theta^{g*}\pi_{G*} \varphi_{G} \cong  \pi_{L_{G,\lambda}*} \varphi_{L_{G, \lambda}}. 
    \end{equation}
    \end{enumerate}
\end{prop}
\begin{proof}
    We first prove that the top square in diagram \eqref{pull_levi} is a pullback. This follows by using Proposition \ref{alper_prop}. In particular, we need to show the map $\Theta^{g}$ is \'etale, separated, representable, sends closed points to closed points and preserves stabilisers. We already know that the map is \'etale from Proposition \ref{lag_corresp_gen}. Since it is a map of quotient stacks induced from an equivariant map with respect to the inclusion of subgroups $L_{G ,\lambda} \subseteq G$, it is representable  and separable. Furthermore, closed orbits are sent to closed orbits since closed orbits correspond to diagonalizable triples of matrices. Finally, let us consider preservation of stabilisers, first in the case $G= \GL_n$. If $x = (A_{1},A_{2},A_{3}) \in \Loc^{g}_{L_{\GL_{n}, \lambda}}$, the stabiliser of $x$ is some block matrix with blocks according to the partition $\lambda$. The stabiliser of $\Theta^{g}(x)$ could be bigger since we are now acting by $\GL_n$ rather than a Levi subgroup. However, because the genericity condition requires that at least one of the $A_{i}$ has distinct eigenvalues in each block, the stabiliser must also split into blocks according to $\lambda$ in $\Loc_{\GL_n}$ so it is the same as in $\Loc_{L_{\GL_{n}, \lambda}}$. Using the description of Levi subgroups of $\SL_n$ we can use the same argument to show that the stabilisers are also preserved by $\Theta^{g}$ in this case.  Using the last part of Proposition \ref{alper_prop} we get that $\theta^{g}$ is \'etale as well.  \par
      
      To prove part $(1)$ consider the preimage $y$ under $\theta^{g}$ of an element $x \in X^{\lambda,g}_{G}$. By the definition of $X^{\lambda,g}_{G}$ this forces $y$ to be a triple of diagonal matrices of the form $\operatorname{diag}(\gamma_{1} I_{\lambda_{1}}, \dots , \gamma_{l} I_{\lambda_{l}})$ with at least one of the three matrices satisfying that the $\gamma_{i}$ are all distinct.
      We therefore get
      $$(\theta^{g})^{-1}(X^{\lambda,g}_{G}) = \widetilde{\ZZ}^{3,g}(L_{G, \lambda}).$$
    Since $X_{L_{G, \lambda}} = H^{3}_{G}/ W_{L_{G,\lambda}}$ where $W_{L_{G, \lambda}}$ is the Weyl group of the Levi, which is a subgroup of the Weyl group of $G$. The map $\theta^{g} \colon X^{g}_{L_{G , \lambda}} \to X_{G}$ is then essentially a further quotient by the full Weyl group. So we see that the maps 
    \begin{align*}
\widetilde{Z}^{3,g}_{L_{G,\lambda}} & \to H^{3}_{G} \to X_{G} \\
\widetilde{Z}^{3,g}_{L_{G,\lambda}} & \to X^{g}_{L_{G,\lambda}} \to X_{G}
\end{align*} 
are the same.
    Therefore, restricting $\theta^{g}$ to $\widetilde{\ZZ}^{3,g}(L_{G, \lambda})$ we get the map $\theta_{\lambda}$ in Lemma \ref{tilde_center_pull_lemma}, which is a $W_{L_{\lambda}} \colon 1$ cover. (See also Examples \ref{z_tilde} and \ref{geom_covers_example}).  \par
     To prove  part $(2)$ we use Proposition \ref{lag_corresp_gen}  to conclude that $\Theta^{*g} \varphi_{G} \cong \varphi_{L_{G, \lambda}}$. Since the first square in \eqref{pull_levi} is a pullback we get 
      $$\pi_{L_{G, \lambda}*} \varphi_{L_{G, \lambda}} \cong \pi_{L_{G, \lambda}*}\Theta^{*g} \varphi_{G} \cong \theta^{*g}\pi_{G*} \varphi_{G}.$$ 
\end{proof} 
\begin{ex} \label{geom_covers_example}
Let us illustrate the geometry of the map $\theta_{\lambda}$ in the example of $n = 6$ and $\lambda = (2,2,1,1)$. We have $L_{\GL_{6} , \lambda} = \GL_{2} \times \GL_2 \times \mathbb{G}_{m} \times \mathbb{G}_{m}$, the centre of this group is given by matrices of the block diagonal form $\diag(\gamma_{1} I_{2}, \gamma_{2} I_{2} , \gamma_{3}, \gamma_{4})$. The relative Weyl group is given by $\So_{2} \times \So_{2}$. The first factor of $\So_2$ permutes the first two blocks $\gamma_{1}I_{2}$ and $\gamma_{2} I_{2}$ while the second factor permutes the elements $\gamma_{3}$ and $\gamma_{4}$. The relative Weyl group $W_{L_{\lambda}}$ acts diagonally on $\ZZ^{3}(L_{\GL_{6}, \lambda})$. \par
The good moduli space of the Levi is
$$X_{L_{\GL_{6}, \lambda}} = \So^{2} \mathbb{G}^{3}_{m} \times \So^{2} \mathbb{G}^{3}_{m} \times  \mathbb{G}^{3}_{m} \times  \mathbb{G}^{3}_{m}.$$   
We can also write $\ZZ^{3}(L_{\GL_{6}, \lambda}) = \mathbb{G}^{3}_{m} \times \mathbb{G}^{3}_{m} \times \mathbb{G}^{3}_{m} \times \mathbb{G}^{3}_{m}$ consisting of elements $x=(x_{1},x_{2},x_{3},x_{4}) \in \prod^{4}_{i=1}\mathbb{G}^{3}_{m}$. The relative Weyl group $W_{L_{\lambda}} = S_{2} \times S_{2}$ then acts by the first factor permuting the $x_{1}$ and $x_{2}$ and the second factor permuting $x_{3}$ and $x_{4}$.
The space $\widetilde{\ZZ}^{3}(L_{\GL_{6}, \lambda})$ is the open subvariety of $\ZZ^{3}(L_{\GL_{6}, \lambda})$ consisting of elements $x=(x_{1},x_{2},x_{3},x_{4}) \in \prod^{4}_{i=1}\mathbb{G}^{3}_{m}$ such that the $x_{i}$ are all distinct. $\widetilde{\ZZ}^{3,g}(L_{\GL_{6}, \lambda})$ is the open subvariety of $\widetilde{\ZZ}^{3}(L_{\GL_{6}, \lambda})$ given by further requiring that the eigenvalues of different blocks are distinct for at least one of the $3$ matrices. Then we have the $\So_{2} \times \So_{2} : 1$ cover 
$$\widetilde{\ZZ}^{3,g}(L_{\GL_{6}, \lambda})  \to X^{(2,2,1,1),g}_{\GL_4} \subseteq \So^{6}_{(2,2,1,1)} \mathbb{G}^{3}_{m} .$$
We get a $\So_{2} \times \So_{2} : 1$ cover since the group $S_{2} \times S_{2}$ acts freely on $\widetilde{Z}^{3}(L_{\GL_6 , \lambda})$ and thus also on $\widetilde{Z}^{3,g}(L_{\GL_6 , \lambda})$.
\end{ex}
We view the next lemma as the multiplicative version of Lemma \ref{additive_computation_lemma}. In particular, we will compute the right hand side of equation \eqref{coh_sln} in terms of $\IC$ sheaves on $X_{G}$. Recall that we can split $\HHf(\B \ZZ(L_{G, \lambda}))[ - \dim \ZZ(L_{G,\lambda})]$ by cohomological degree into subspaces $V_i$ $i \geq \dim \ZZ(L_{G,\lambda})$. The natural $W_{L_{\lambda}}$ action on $\HHf(\B \ZZ(L_{G, \lambda}))$ perserves cohomological degree so each $V_{i}$ is a $W_{L_{\lambda}}$ subrepresentation.
\begin{lem}[BPS sheaves for Levis in $G$]\label{loc_sys_comp}
Let $G = \GL_n, \SL_n$. If $n$ is prime, we also allow $G=\PGL_n$. The following properties hold for BPS sheaves on Levis
 \begin{enumerate}
     \item $\pi_{L_{G,\lambda}*} \varphi_{L_{G,\lambda}}$ has perverse cohomology bounded below. We define $$\BPS_{L_{G,\lambda}} = \pH^{l} \pi_{L_{G,\lambda}*} \varphi_{L_{G,\lambda}}$$  here $l = \dim \ZZ(L_{G,\lambda})$. Furthermore, $\BPS_{L_{G,\lambda}}$ is a constant sheaf of rank $1$ supported on 
 $$\operatorname{supp}(\BPS_{L_{G,\lambda}}) =  \im(\Delta^{\lambda} \colon \ZZ^{3}(L_{G,\lambda}) \xhookrightarrow{} X_{L_{G,\lambda}}) .$$
     Furthermore, the components of the Saito decomposition of $\pi_{L_{G, \lambda}*} \varphi_{L_{G, \lambda}}$ with supports given by $\ZZ^{3}(L_{G, \lambda})$ are $\BPS_{L_{G,\lambda}} \otimes \HHf(\B \ZZ(L_{G, \lambda}))[- \dim \ZZ(L_{G, \lambda})]$.
\item 
Consider the subspace $V_{i}$ of cohomological degree $i$ in $\HHf(\B \ZZ(L_{G, \lambda}))[ - \dim \ZZ(L_{G,\lambda})]$.  The term 
 \begin{equation} \label{equiv_bps}
     \BPS_{L_{G, \lambda}}  \otimes V_{i} 
 \end{equation}
  has a natural action of $W_{L_{\lambda}}$, which corresponds to the finite dimensional representation $V_{i}$. Pushing forward by $\theta \colon X_{L_{G,\lambda}} \to X_{G}$ and taking invariant part we get
  \begin{equation} \label{mult_levi_ic_comp}
      (\theta_{*} \BPS_{L_{G, \lambda}} 
 \otimes V_{i})^{W_{L_{\lambda}}} \cong \IC_{X^{\lambda,g}_{G}}(\mathcal{K}^{\lambda,g}_{i})[-i]
  \end{equation}
  where 
  \begin{equation} \label{llambda_comp_bpsforlevislem}
      \mathcal{K}^{\lambda,g}_{i}[ \dim X^{\lambda,g}_{G} - i] = (\theta_{\lambda *}(\BPS_{L_{G, \lambda}}
 \otimes V_{i})|_{\widetilde{\ZZ}^{3,g}(L_{G, \lambda})}))^{W_{L_{\lambda}}}
  \end{equation}
     with $\theta_{\lambda}$ defined as in diagrams \eqref{tilde_center_pullbacks} for $\GL_n$ and $\SL_n$ and \eqref{tilde_center_pullbacks_pgln} for $\PGL_n$.
 \end{enumerate}
 \end{lem}
 \begin{remark}
     The above lemma computes the terms on the RHS of equation \eqref{coh_sln} as $\IC$ sheaves of local systems supported on the stratification in Definition \ref{stratifications} for $\GL_n$, $\SL_n$ and for $\PGL_n$ supported on the stratification in Definition \ref{strat_pgl_n}.
 \end{remark}
 \begin{proof}
 Part $(1)$ for the trivial Levi $L_{\GL_n , \lambda} = \GL_n$ follows immediately from Corollaries \ref{purity_GL} and \ref{bps_const} using the exponential map to deduce it from the additive case in Lemma \ref{additive_computation_lemma}. Now we can write
 \begin{equation} \label{box_tens_dt}
     \pi_{L_{\GL_n , \lambda} *} \varphi_{L_{\GL_n , \lambda}} \cong \pi_{\GL_{\lambda_{1}}*} \varphi_{\GL_{\lambda_{1}}} \boxtimes \dots \boxtimes \pi_{\GL_{\lambda_{1}}*} \varphi_{\GL_{\lambda_{1}}}
 \end{equation}
 Then repeatedly using Corollary \ref{support_locsys} to write
 \begin{equation}
    \pi_{\GL_{\lambda_{i}}*} \varphi_{\GL_{\lambda_{i}}} \cong \bigoplus_{k \geq 1} \bigoplus_{\gamma^{i}}\IC_{X^{\gamma^{i}}_{\GL_{\lambda_{i}}}}(\mathcal{L}^{\gamma^{i}}_{k})[-k]
\end{equation}
where now $\gamma^{i}$ is some partition of $\lambda_{i}$ and \emph{not} $n$. 
Recall that $\BPS_{\GL_n} \cong \IC_{X^{(n)}_{\GL_{n}}}(\mathbb{Q}_{X^{(n)}_{\GL_n}})$ with $X^{(n)}_{\GL_n} = \ZZ^{3}(\GL_n)$. Using that external tensor product is compatible with $\IC$ sheaves we get 
\begin{align*}
    \IC_{X^{(\lambda_{1})}_{\GL_{\lambda_{1}}}}(\mathbb{Q}_{X^{(\lambda_i)}_{\GL_{\lambda_1}}}) \boxtimes \dots \boxtimes \IC_{X^{(\lambda_{l})}_{\GL_{\lambda_{l}}}}(\mathbb{Q}_{X^{(\lambda_l)}_{\GL_{\lambda_l}}}) & \cong  \IC_{ \prod X^{(\lambda_{i})}_{\GL_{\lambda_{i}}}}(\mathbb{Q}_{\prod X^{(\lambda_i)}_{\GL_{\lambda_i}}})\\
    & \cong  \IC_{ \ZZ^{3}(L_{\GL_n, \lambda})}(\mathbb{Q}_{\ZZ^{3}(L_{\GL_n , \lambda})}).
\end{align*}
So from equation \eqref{box_tens_dt} we can compute all the summands supported on $\ZZ^{3}(L_{\GL_n , \lambda})$ as a tensor product. From this it follows that the summands supported on $\ZZ^{3}(L_{\GL_n , \lambda})$ are
\begin{align*}
     \BPS_{L_{\GL_n , \lambda}} & \otimes (\HHf( \B \ZZ(\GL_{\lambda_{1}})[-\dim \ZZ(\GL_{\lambda_{1}})]) \otimes \dots \otimes \HHf( \B \ZZ(\GL_{\lambda_{l}}))[-\dim \ZZ(\GL_{\lambda_l})])  \cong \\
     & \cong \BPS_{L_{\GL_n , \lambda}} \otimes \HHf( \B \ZZ(L_{\GL_{n}, \lambda }))[-\dim \ZZ(L_{\GL_{n}, \lambda})].
\end{align*}
\par
 For $\GL_n$, using diagram \eqref{tilde_center_pullbacks} in Lemma \ref{tilde_center_pull_lemma}, the computation in Part $(2)$ and in particular the proof of equation \eqref{mult_levi_ic_comp} here is analogous to Lemma \ref{additive_computation_lemma}. However, we further restrict to the generic locus $\widetilde{\ZZ}^{3,g}(L_{G, \lambda})$.
 The proofs for $\SL_n$ and $\PGL_{n}$ are given in Section \ref{loc_sys_comp_sl} and Section  \ref{loc_sys_comp_pgl} respectively. 
 \end{proof}
 In the following lemma we will consider the action on $W_{L_{\lambda}}$ on the terms $(\BPS_{L_{G,\lambda}}  \otimes \HHf(\B \ZZ (L_{G, \lambda}))[ - \dim \ZZ(L_{G,\lambda})])$. Note that a priori there are two different $W_{L_{\lambda}}$ actions. One coming from the natural action on the sheaf $\varphi_{L_{G, \lambda}}$ and one induced by the natural action on $\HHf(\B \ZZ (L_{G, \lambda}))$. We will now show they are the same.
 \begin{lem} \label{action_lemma_split_off}
     Let $G= \GL_n$ or $\SL_n$. The $W_{L_{ \lambda}}$ action on $(\BPS_{L_{G,\lambda}}  \otimes \HHf(\B \ZZ (L_{G, \lambda}))[ - \dim \ZZ(L_{G,\lambda})])$ induced from the $W_{L_{\lambda}}$-action on $\varphi_{L_{G, \lambda}}$ is equivalent to the natural action of $W_{L_{ \lambda}}$ on 
     $$\HHf(\B \ZZ (L_{G, \lambda}))[ - \dim \ZZ(L_{G,\lambda})].$$
 \end{lem}
 \begin{proof}
     Recall Lemma \ref{coh_centers} and consider first the case of $\GL_n$. Here the relative Weyl group acts by permuting the factors of the DT sheaf 
    $$\varphi_{L_{\GL_{n}, \lambda}} \cong  \varphi_{\GL_{\lambda_{1}}} \boxtimes \dots \boxtimes \varphi_{\GL_{\lambda_{l}}} $$ hence on the pushforward it also acts by permuting the factors $$\pi_{\GL_{\lambda_{1}}*} \varphi_{\GL_{\lambda_{1}}} \boxtimes \dots \boxtimes \pi_{\GL_{\lambda_{l}}*} \varphi_{\GL_{\lambda_{l}}}.$$ Restricting to the support of the BPS sheaf $\BPS_{L_{\GL_{n},\lambda}}$ this gives the permutation action on $\BPS_{L_{\GL_{n},\lambda}} \otimes \HHf (\B \ZZ(L_{\GL_{n}, \lambda}))$. \par For $\SL_n$, we will use Lemma \ref{loc_sys_comp}, which is proven in subsection \ref{sln_subsect}. In particular, we will use that the components of the Saito decomposition of $\pi_{L_{\SL_n, \lambda}*} \varphi_{L_{\SL_n, \lambda}}$ with supports given by $\ZZ^{3}(L_{\SL_n, \lambda})$ are $\BPS_{L_{\SL_n,\lambda}} \otimes \HHf(\B \ZZ(L_{\SL_n, \lambda}))[- \dim \ZZ(L_{\SL_n, \lambda})]$. We have the short exact sequence on centers $1 \to \mu_n \to \ZZ(L_{\SL_n , \lambda}) \times \mathbb{G}_{m}  \to \ZZ(L_{\GL_{n}, \lambda}) \to 1$. The short exact sequence induces the $W_{L_{\lambda}}$-equivariant map $$\nu \colon \ZZ^{3}(L_{\SL_n , \lambda}) \times \mathbb{G}^{3}_{m}  \to \ZZ^{3}(L_{\GL_{n}, \lambda})$$ We know from Proposition \ref{diagram_sln} that the map $\eta \colon \LocB_{\SL_{n} \times \mathbb{G}_{m}} \to \LocB_{\GL_{n}}$ is $W_{L_{\lambda}}$-equivariant and also gives an isomorphism $\eta^{*} \varphi_{L_{\GL_n, \lambda}} \cong \varphi_{L_{\SL_n, \lambda}} \boxtimes \varphi_{\mathbb{G}_{m}}$ of $W_{L_{\lambda}}$-equivariant perverse sheaves. We can write
\begin{equation} \label{wl_equiv_eq}
    (\pi_{L_{\SL_{n},\lambda}} \times \pi_{\mathbb{G}_{m}})_{*} (\varphi_{L_{\SL_{n},\lambda}} \boxtimes \varphi_{\mathbb{G}_{m}}) \cong p^{*}_{X}\pi_{L_{\SL_{n},\lambda}*}\varphi_{L_{\SL_{n},\lambda}} \otimes p^{*}_{\mathbb{G}^{3}_{m}}\pi_{\mathbb{G}_{m}*} \varphi_{\mathbb{G}_{m}}.
\end{equation}
From the first part of this proof $W_{L_{\lambda}}$ acts by permutation on 
$$\BPS_{L_{\GL_{n}, \lambda}} \otimes \HHf(\B \ZZ(L_{\GL_{n}, \lambda}))[- \dim \ZZ(L_{\GL_n , \lambda})],$$ then $W_{L_{\lambda}}$ also acts by permutation on $  \nu^{*}\BPS_{L_{\GL_{n}, \lambda}} \otimes \HHf(\B \ZZ(L_{\GL_{n}, \lambda}))[- \dim \ZZ(L_{\GL_n , \lambda})]$ since $\nu$ is $W_{L_{\lambda}}$-equivariant.
    Restricting to components only supported on $\ZZ^{3}(L_{\SL_n, \lambda}) \times \mathbb{G}^{3}_{m}$ the above equation \eqref{wl_equiv_eq} implies an equivalence of $W_{L_{\lambda}}$-equivariant sheaves 
     \begin{align}
          p^{*} & \BPS_{L_{\SL_{n}, \lambda}} \otimes  \HHf(\B \ZZ(L_{\SL_{n}, \lambda})) [- \dim \ZZ(L_{\SL_n , \lambda})] \otimes \mathbb{Q}_{\ZZ^{3}(L_{\SL_n , \lambda}) \times \mathbb{G}^{3}_{m}}[2] \otimes \HHf( \B \mathbb{G}_{m})  \cong  \\
          & \cong \nu^{*} \BPS_{L_{\GL_{n}, \lambda}} \otimes \HHf(\B \ZZ(L_{\GL_{n}, \lambda}))[- \dim \ZZ(L_{\GL_n , \lambda})] \nonumber
      \end{align}
       with $p \colon \ZZ^{3}(L_{\SL_{n}, \lambda}) \times \mathbb{G}^{3}_{m} \to \ZZ^{3}(L_{\SL_{n}, \lambda})$ a $W_{L_{\lambda}}$-equivariant projection. Using Lemma \ref{coh_centers}, $W_{L_{\lambda}}$ acts trivially on the $\HHf(\B \mathbb{G}_{m})$ factor. So the $W_{L_{\lambda}}$ action on $p^{*}  \BPS_{L_{\SL_{n}, \lambda}} \otimes  \HHf(\B \ZZ(L_{\SL_{n}, \lambda}))[- \dim \ZZ(L_{\GL_n , \lambda})] $ is the natural action of $W_{L_{\lambda}}$ on $\HHf(\B \ZZ(L_{\SL_{n}, \lambda}))[- \dim \ZZ(L_{\GL_n , \lambda})]$
       as in Lemma \ref{coh_centers}. Pulling back by the $W_{L_{\lambda}}$-equivariant inclusion
       \begin{equation}
           \iota \colon \ZZ^{3}(L_{\SL_{n}, \lambda}) \to \ZZ^{3}(L_{\SL_{n}, \lambda}) \times \mathbb{G}^{3}_{m}
       \end{equation}
       we get
       \begin{equation}
           \iota^{*} p^{* } \BPS_{L_{\SL_{n}, \lambda}} \otimes  \HHf(\B \ZZ(L_{\SL_{n}, \lambda})) [- \dim \ZZ(L_{\SL_n , \lambda})] \cong \BPS_{L_{\SL_{n}, \lambda}} \otimes  \HHf(\B \ZZ(L_{\SL_{n}, \lambda}))[- \dim \ZZ(L_{\SL_n , \lambda})]
       \end{equation}
        Therefore, $W_{L_{\lambda}}$ acts on $\BPS_{L_{\SL_{n}, \lambda}} \otimes  \HHf(\B \ZZ(L_{\SL_{n}, \lambda}))$ in the natural way as in Lemma \ref{coh_centers}. 
 \end{proof}
 \begin{remark} \label{fin_gp_rem}
    Let $X$ and $Y$ be smooth varieties and $f \colon X \to Y$ be a $K \colon 1$ cover by a finite group $K$. Any local system $\mathcal{L}$ on $Y$ is determined by $f^{*} \mathcal{L}$ and the $K$-equivariant structure on $f^{*} \mathcal{L}$. Knowing the $K$-equivariant structure, we can recover $\mathcal{L}$ as $(f_{*} f^{*} \mathcal{L})^{K}$. Where $(-)^{K}$ is taking the invariants.  
 \end{remark}
  We will use the above remark in the proof of the following proposition, which computes the restrictions of the local systems $\mathcal{L}^{\lambda}_{i}$ in equation \eqref{loc_sys_dir_sum_strata} to $X^{\lambda, g}_{G}$. 
\begin{prop}[Reduction to Levis] \label{levis_final}
Assume that Lemma \ref{loc_sys_comp} has been established for $G= \GL_n$ and $\SL_n$. 
Fix a local system $\mathcal{L}^{\lambda}_{i}$ in equation \eqref{loc_sys_dir_sum_strata}. We have the following isomorphism of $W_{L_{\lambda}}$-equivariant local systems
\begin{equation} \label{bps_loc_strat}
    \theta^{*}_{\lambda}\mathcal{L}^{\lambda,g}_{i}[\dim X^{\lambda,g}_{G} - i] \cong  (\BPS_{L_{G,\lambda}}  \otimes V_{i})|_{\widetilde{\ZZ}^{3,g}(L_{G, \lambda})}
\end{equation}
where $V_{i}$ is the degree $i$ subspace in $\HHf(\B \ZZ(L_{G, \lambda}))$ with its natural $W_{L_{\lambda}}$ action. The degrees $i$ such that $\mathcal{L}^{\lambda,g}_{i}$ are non-zero in equation \eqref{loc_sys_dir_sum_strata} are in one-to-one correspondence with the cohomological degrees of $\HHf(\B \ZZ(L_{G, \lambda}))[- \dim \ZZ(L_{G, \lambda})]$. The local system $\mathcal{L}^{\lambda,g}_{i}$ then satisfies
\begin{equation} \label{mult_loc_sys_red}
    \mathcal{L}^{\lambda,g}_{i}[\dim X^{\lambda,g}_{G} - i] \cong (\theta_{\lambda*}(\BPS_{L_{G,\lambda}}  \otimes V_{i})|_{\widetilde{\ZZ}^{3,g}(L_{G, \lambda})})^{W_{L_{ \lambda}}}.
\end{equation}
 \end{prop}
 \begin{proof}  
     Since $\theta^{g*}$ is \'etale it preserves the perverse $t$-structure. Then using equation \eqref{levi_red_equation} we have that for any $k$
\begin{equation} \label{eq_red_lev}
    \theta^{g*} (\pH^{k}\pi_{G*} \varphi_{G} )\cong \pH^{k}(\theta^{g*}\pi_{G*} \varphi_{G}) \cong  \pH^{k} \pi_{L_{G,\lambda}*} \varphi_{L_{G, \lambda}}.
\end{equation}
 Furthermore, by Propositions \ref{purity_GL}, \ref{purity_sln},  $\pi_{G*} \varphi_{G}$ and $\pi_{L_{G,\lambda}*} \varphi_{L_{G, \lambda}}$ are pure complexes of mixed Hodge modules for $G = \GL_n , \SL_n$. Therefore, both $\pi_{G*} \varphi_{G}$ and $\pi_{L_{G,\lambda}*} \varphi_{L_{G, \lambda}}$ decompose into a direct sum of their shifted perverse pieces.  Fix a perverse degree $k$ and take a summand $\mathcal{F}^{\lambda}_{k}$ of $\pH^{k} \pi_{G*} \varphi_{G}$ that has support on $X^{\lambda}_{G}$.  In particular,
\begin{align*}
    \mathcal{F}^{\lambda}_{k} & \cong \IC_{X^{\lambda,g}_{G}}(\mathcal{L}^{\lambda,g}_{k})  \text{ for } X^{\lambda,g}_{G} \subseteq X_{G}. 
\end{align*}
Equation \eqref{eq_red_lev} then implies that
\begin{equation}
     \theta^{g*}\IC_{X^{\lambda,g}_{G}}(\mathcal{L}^{\lambda,g}_{k}) \cong \IC_{\theta^{g-1}(X^{\lambda,g}_{G})}(\theta^{g*} \mathcal{L}^{\lambda,g}_{k}).  
\end{equation}
is some summand of the $k$-th perverse cohomology of $\pi_{L_{G, \lambda}*} \varphi_{L_{G, \lambda}}$. Using diagram \eqref{pull_levi} and Proposition \ref{levis_final_splitoff} we get that $\theta^{g-1}(X^{\lambda}_{G}) = \widetilde{\ZZ}^{3,g}(L_{G, \lambda}) \subseteq \ZZ^{3}(L_{G, \lambda})$. So pulling back by $\theta^{g}$ we get the summands of the decomposition of $\pi_{L_{G,\lambda}*} \varphi_{L_{G, \lambda}}$ supported on $\ZZ^{3}(L_{G, \lambda})$, which by part $(1)$ of Lemma \ref{loc_sys_comp} (which we assume has been proven ) are given by 
\begin{equation}
    \BPS_{L_{G,\lambda}} \otimes \HHf(\B \ZZ(L_{G, \lambda}))[- \dim \ZZ(L_{G, \lambda})].
\end{equation}
Therefore, using that the lower square in diagram \eqref{pull_levi} commutes and $j^{*}_{\lambda , G}\IC_{X^{\lambda,g}_{G}}(\mathcal{L}^{\lambda,g}_{k}) \cong \mathcal{L}^{\lambda,g}_{k}[\dim X^{\lambda,g}_{G}]$ we get
\begin{align}  \label{proof_mid_levi_eq1}
    \theta^{*}_{\lambda}( \bigoplus_{k} \mathcal{L}^{\lambda,g}_{k} [\dim X^{\lambda,g}_{G}-k] ) & \cong (\BPS_{L_{G,\lambda}}  \otimes \HHf(\B \ZZ( L_{G, \lambda}))[ - \dim \ZZ(L_{G,\lambda})])|_{\widetilde{\ZZ}^{3,g}(L_{G, \lambda})}  \\
    \label{proof_mid_levi_eq2} \theta^{*}_{\lambda}\mathcal{L}^{\lambda,g}_{k}[\dim X^{\lambda,g}_{G} -k] & \cong  (\BPS_{L_{G,\lambda}}  \otimes V_{k})|_{\widetilde{\ZZ}^{3,g}(L_{G, \lambda})}.
\end{align}
Because $\theta_{\lambda}$ is \'etale we have the following equation for perverse pieces supported on $X^{\lambda}_{G}$ for any $k$
\begin{equation}
    \pH^{k}(\pi_{L_{G, \lambda}*} \varphi_{L_{G,\lambda}}) \cong \theta^{g*} \pH^{k}(\pi_{G*} \varphi_{G}).
\end{equation}
Therefore, the terms $\mathcal{L}^{\lambda,g}_{k}$ and $\BPS_{L_{G,\lambda}}  \otimes V_{k}$ are in $1:1$ correspondence. \par
We now explain why equation \eqref{bps_loc_strat} is an equivalence of $W_{L_{\lambda}}$-local systems. Note that all of the maps in diagram \eqref{pull_levi} are $W_{L_{\lambda}}$-equivariant. Because of Propositions \ref{levis_final_splitoff} and \ref{lag_corresp_gen} we know that the map $\Theta^{g}$ is $W_{L_{\lambda}}$-oriented  and $W_{L_{\lambda}}$-equivariant, so we get an equivalence of $W_{L_{\lambda}}$-equivariant sheaves
\begin{equation}
    \Theta^{g*} \varphi_{G} \cong \varphi_{L_{G, \lambda}}.
\end{equation}
Pushing forward to the good moduli space we have induced actions by $W_{L_{\lambda}}$. Since $\theta^{g}$ is $W_{L_{\lambda}}$-equivariant we then get an equivalence of $W_{L_{\lambda}}$-equivariant complexes of sheaves
\begin{equation}
    \theta^{g*}\pi_{G*} \varphi_{G} \cong \pi_{L_{G, \lambda}*}\varphi_{L_{G, \lambda}}.
\end{equation}
The fact that all the maps in diagram \eqref{pull_levi} are $W_{L_{\lambda}}$-equivariant then allows us to upgrade the equivalences in equations \eqref{proof_mid_levi_eq1}, \eqref{proof_mid_levi_eq2} to equivalences of $W_{L_{\lambda}}$-equivariant local systems. Finally, we can use Remark \ref{fin_gp_rem} for the map $\theta_{\lambda}$ to get
\begin{equation} 
    \mathcal{L}^{\lambda,g}_{i}[\dim X^{\lambda,g}_{G} - i] \cong (\theta_{\lambda*}(\BPS_{L_{G,\lambda}}  \otimes V_{i})|_{\widetilde{\ZZ}^{3,g}(L_{G, \lambda})})^{W_{L_{ \lambda}}}.
\end{equation}
as required in equation \eqref{mult_loc_sys_red}.
 \end{proof}
\begin{ex} \label{2_4_coh_int}
Consider the example of $\GL_2$. In this case the only non-trivial Levi is the maximal torus $H_{\GL_2}$ corresponding to the partition $\lambda = (1,1)$. Let us illustrate the computation of the local systems in equation \eqref{loc_sys_dir_sum_strata} in this case. The diagram \eqref{pull_levi} becomes
\begin{equation}
\begin{tikzcd}
	{\CC_{3}(\GL_n)/\GL_n} & { (\mathbb{G}^{3}_{m} \times \mathbb{G}^{3}_{m})^{g}  / \mathbb{G}^{2}_{m}} \\
	{\So^{2}(\mathbb{G}^{3}_{m})} & {(\mathbb{G}^{3}_{m} \times \mathbb{G}^{3}_{m})^{g}} \\
	{\So^{2}(\mathbb{G}^{3}_{m}) \setminus \Delta} & {(\mathbb{G}^{3}_{m} \times \mathbb{G}^{3}_{m}) \setminus \Delta}
	\arrow[from=1-1, to=2-1]
	\arrow["{\Theta^{g}}"', from=1-2, to=1-1]
	\arrow[from=1-2, to=2-2]
	\arrow["{\theta^{g}}", from=2-2, to=2-1]
	\arrow["{j_{\lambda}}", from=3-1, to=2-1]
	\arrow[from=3-2, to=2-2]
	\arrow["{\theta_{\lambda}}", from=3-2, to=3-1]
\end{tikzcd}
\end{equation}
here $(\mathbb{G}^{3}_{m} \times \mathbb{G}^{3}_{m})^{g}$ is the space of triples of diagonal $2 \times 2$ matrices such that atleast one matrix in the triple has distinct eigenvalues, $\So^{2} \mathbb{G}^{3}_{m} \setminus \Delta$ is the complement of the diagonal $\Delta \colon \mathbb{G}^{3}_{m} \to \So^{2} \mathbb{G}^{3}_{m}$ and $(\mathbb{G}^{3}_{m} \times \mathbb{G}^{3}_{m}) \setminus \Delta$ is the complement of the diagonal $\Delta \colon \mathbb{G}^{3}_{m} \to \mathbb{G}^{3}_{m} \times \mathbb{G}^{3}_{m} $ in $\mathbb{G}^{3}_{m} \times \mathbb{G}^{3}_{m}$. In this particular case $(\mathbb{G}^{3}_{m} \times \mathbb{G}^{3}_{m})^{g}$ is the same as $(\mathbb{G}^{3}_{m} \times \mathbb{G}^{3}_{m}) \setminus \Delta$. In particular, for $\GL_2$ we do not need to further restrict to a generic locus. Similarly $X^{(1,1),g}_{\GL_2} = X^{(1,1)}_{\GL_2} = \So^{2} \mathbb{G}^{3}_{m} \setminus \Delta$. However, for higher $n$ the two spaces will not be the same and we will need to restrict to the generic locus. The map $\theta_{\lambda}$ is a $\So_{2} : 1$ cover.
We have
\begin{equation}
    \HHf(\B \ZZ(\mathbb{G}^{2}_m)) \cong \mathbb{Q}[u_{1},u_{2}]
\end{equation}
with $u_{i}$ in degree $2$ and the relative Weyl group $\So_{2}$ acts by permuting the $u_{i}$. Then we see that we can decompose $\mathbb{Q}[u_{1},u_{2}][-2]$ into pieces $V_{k}$ spanned by homogenous polynomials in $u_{i}$ of degree $k-2$, $k\geq 2$. $V_{k}$ is in cohomological degree $2k+2$. Note that since the $\So_{2}$ action preserves cohomological degree we get an action on $V_{k}$. In low degrees we can more explicitly write down the equation  \eqref{bps_loc_strat} 
\begin{align} \label{s2_reps}
    \theta^{*}_{\lambda} \mathcal{L}^{\lambda}_{2}[6-2] \cong \BPS_{H_{\GL_2}} \otimes V_{2} & = \BPS_{H_{\GL_2}} \otimes \mathbb{Q}[-2] \cong \BPS_{H_{\GL_2}} \otimes W_{\operatorname{triv}}[-2] \\
    \theta^{*}_{\lambda} \mathcal{L}^{\lambda}_{4}[6-4] \cong \BPS_{H_{\GL_2}} \otimes V_{4} & =  \BPS_{H_{\GL_2}} \otimes (\mathbb{Q}u_{1} \oplus \mathbb{Q}u_{2})[-2]   \cong  \BPS_{H_{\GL_2}} \otimes W_{\operatorname{perm}}[-2] \nonumber \\
    \theta^{*}_{\lambda} \mathcal{L}^{\lambda}_{6}[6-6] \cong \BPS_{H_{\GL_2}} \otimes V_{6} & = \BPS_{H_{\GL_2}} \otimes( (\mathbb{Q}u^{2}_{1} \oplus \mathbb{Q}u^{2}_{2}) \oplus \mathbb{Q}u_{1}u_{2})[-2]  \nonumber \\
    & \cong \BPS_{H_{\GL_2}} \otimes (W_{\operatorname{perm}} \oplus W_{\operatorname{triv}})[-2] \nonumber
\end{align}
with $\dim X^{(1,1)}_{\GL_2} = \dim \ZZ^{3}(H_{\GL_2}) = 6$. Here $W_{\operatorname{triv}}$ is the trivial one-dimensional representation of $\So_{2}$ and $W_{\operatorname{perm}}$ is the permutation representation of $\So_{2}$ on $\mathbb{Q}^{2}$.
In this case we can see that the only degrees that will appear in equation \eqref{loc_sys_dir_sum_strata} will be $2+2n$ $n \geq 0$. 
In low degrees we will then get the following local systems on $\So^{2} \mathbb{G}^{3}_{m} \setminus \Delta$ using equation \eqref{s2_reps}
\begin{align*}
    \mathcal{L}^{\lambda}_{2} & \cong \mathbb{Q}_{\So^{2} \mathbb{G}^{3}_{m} \setminus \Delta} \\
    \mathcal{L}^{\lambda}_{4} & \cong \mathbb{Q}_{\So^{2} \mathbb{G}^{3}_{m} \setminus \Delta} \oplus \mathcal{K}_{2}
\end{align*}
where $\mathcal{K}_{2}$ is a $\mathbb{Z} / 2 \mathbb{Z}$-local system. \par
To prove cohomological integrality for $\GL_2$ we just have to consider two strata corresponding to the partitions $(2)$ and $(1,1)$. We know the contributions of the stratum corresponding to $(2)$ from Corollary \ref{support_locsys}, which come from copies of the BPS sheaf. Running the above argument for the maximal torus we can finish the proof of cohomological integrality by comparing the local systems $\mathcal{L}^{(1,1)}_{i}$ to the ones appearing in equation \eqref{llambda_comp_bpsforlevislem}. \par
Let us illustrate briefly how the higher $n$ cases work with the example of $n = 4$. We have the following partitions of $4 \colon$ $\lambda^{1} = (4)$, $\lambda^{2} = (3,1)$, $\lambda^{3} = (2,1,1)$, $\lambda^{4} = (2,2)$ and $\lambda^{5} = (1,1,1,1)$. The BPS sheaf is supported on the centre $\ZZ^{3}(\GL_4)$ which corresponds to the partition $\lambda^{1}$, while $\lambda^{5}$ corresponds to the maximal torus of $\GL_{4}$. To prove cohomological integrality we have to use Proposition \ref{levis_final} $4$ times 
 and compare with Lemma \ref{loc_sys_comp} for all the non-trivial partitions. Note that the relative Weyl group is trivial for $\lambda^{2}$ so the cover $\theta_{\lambda^{2}}$ actually becomes an isomorphism.
\end{ex}
\begin{proof}[Proof of cohomological integrality for $\GL_n$]
 We combine the previous lemmas and propositions following the strategy \ref{strategy}. From Corollary \ref{support_locsys} we know that $\pi_{\GL_{n}*} \varphi_{\GL_{n}}$ splits into $\IC$ sheaves with support $X^{\lambda}_{\GL_{n}}$ for some partition $\lambda$ of $n$. On a fixed stratum $X^{\lambda}_{\GL_n}$ we have shifted local systems $\mathcal{L}^{\lambda}_{i}[\dim X^{\lambda}_{\GL_n} -i]$ as in Equation \eqref{loc_sys_dir_sum_strata}. We will compare the restrictions $\mathcal{L}^{\lambda ,g}_{i}$ of these local systems to the generic locus $X^{\lambda,g}_{\GL_n}$ term by term with the local systems that appear on the RHS of Equation \eqref{coh_gln}. The local systems $\mathcal{K}^{\lambda,g}_{i}$ on the RHS of Equation \eqref{coh_gln} are computed in part $(2)$ of  Lemma \ref{loc_sys_comp} and are controlled by the local system corresponding to the $W_{L_{\lambda}}$ representation $V_{i}$. On the other hand we have computed the local systems $\mathcal{L}^{\lambda,g}_{i}$ in Proposition \ref{levis_final}. By comparing equations \eqref{mult_loc_sys_red} and \eqref{llambda_comp_bpsforlevislem} we can see that the local systems $\mathcal{K}^{\lambda, g}_{i}$ and $\mathcal{L}^{\lambda,g}_{i}$ are isomorphic.
\end{proof}

\section{Cohomological integrality for \texorpdfstring{$\SL_n$}{TEXT}, \texorpdfstring{$\PGL_n$}{TEXT} and Langlands duality}
 
\subsection{Proof of integrality for \texorpdfstring{$\SL_n$}{TEXT}} \label{sln_subsect}
We now consider the stack of $\SL_n$-local systems. Recalling Definition \ref{framed_locsys}, we write $\Loc_{G} = \Loc^{\ff}_{G} / G$.
Even though $\ZZ(L_{\SL_n, \lambda})$ can be disconnected we will still say $\dim \ZZ(L_{\SL_n, \lambda}) = l-1$ since every connected component has the same dimension. Using the short exact sequence   
\begin{equation} \label{levi_ses_sln}
  1 \to \mu_n \to L_{\SL_n, \lambda} \times \mathbb{G}_{m} \to L_{\GL_n , \lambda} \to 1  
\end{equation}
We get a diagram 
\begin{equation}\label{ses_m_square}
\begin{tikzcd}
	{\Loc^{\ff}_{L_{\SL_n , \lambda}} \times \mathbb{G}^{3}_{m}\times \B \mu_{n}} & {\Loc^{\ff}_{L_{\SL_n , \lambda}} \times \mathbb{G}^{3}_{m}} \\
	{\Loc^{\ff}_{L_{\SL_n , \lambda}} \times \mathbb{G}^{3}_{m}/(L_{\SL_n , \lambda } \times \mathbb{G}_{m})} & {\Loc^{\ff}_{L_{\SL_n , \lambda}} \times \mathbb{G}^{3}_{m}/L_{\GL_n , \lambda }}
	\arrow["{\eta_{1}}", from=2-1, to=2-2]
	\arrow["p", from=1-2, to=2-2]
	\arrow["{\widetilde{\eta_{1}}}", from=1-1, to=1-2]
	\arrow["{\widetilde{p}}", from=1-1, to=2-1]
\end{tikzcd}
\end{equation}
Note that $L_{\GL_n , \lambda }$ acts on $\Loc^{\ff}_{L_{\SL_n}, \lambda}$ since conjugation by elements in $\GL_n$ preserves determinant $1$ matrices. We take the trivial action of $L_{\GL_n , \lambda} $ on $\mathbb{G}^{3}_{m}$. We also allow the trivial Levis $L_{\SL_n , \lambda} = \SL_n$ and $L_{\GL_n , \lambda} = \GL_n$  which correspond to  $\lambda = (n)$.
\begin{lem}
    The diagram \eqref{ses_m_square} is a pullback.
\end{lem}
\begin{proof}
Let 
\begin{equation}
    1 \to K \to G \to H \to 1
\end{equation}
be a short exact sequence of algebraic groups and assume that $G$ acts on a scheme $X$ with $K$ acting trivially so that the action descends to $H$. We then have
\begin{equation} \label{ses_pull_standard_stacks}
    X/G \times_{X/H} X \cong X \times \B K
\end{equation} 
where the map $X/G \to X/ H$ is induced by the identity $X \to X$. The lemma follows by equation \eqref{ses_pull_standard_stacks} applied to $1 \to \mu_n \to L_{\SL_n , \lambda} \times \mathbb{G}_{m} \to L_{\GL_n , \lambda} \to 1$ and $X = \Loc^{\ff}_{L_{\SL_n , \lambda}} \times \mathbb{G}^{3}_{m}$. 
\end{proof}
The short exact sequences \eqref{levi_ses_sln} induce by taking mapping stacks a commutative diagram of stacks and good moduli spaces 
\begin{equation}\label{sln_square}
\begin{tikzcd}
	{\Loc_{L_{\SL_n , \lambda} \times \mathbb{G}_{m}}} & {(\Loc^{\ff}_{L_{\SL_n , \lambda}} \times \mathbb{G}^{3}_{m}) / L_{\GL_n , \lambda}} & {\Loc_{L_{\GL_n , \lambda }}} \\
	& {X_{L_{\SL_n , \lambda}} \times \mathbb{G}^{3}_{m}} & {X_{L_{\GL_n , \lambda}}}
	\arrow["{\eta_{2}}", from=1-2, to=1-3]
	\arrow["\pi"', from=1-2, to=2-2]
	\arrow["{\widetilde{\eta}_{2}}"', from=2-2, to=2-3]
	\arrow["{\pi_{\mathfrak{gl}_{n}}}", from=1-3, to=2-3]
	\arrow["{\eta_{1}}", from=1-1, to=1-2]
	\arrow["\eta", bend left = 20, from=1-1, to=1-3]
	\arrow["{\pi_{L_{\SL_n , \lambda}} \times \pi_{\mathbb{G}_{m}}}"', from=1-1, to=2-2]
\end{tikzcd}
\end{equation}
\begin{ex} \label{sl2_etale_coverex}
Let us consider the map $\widetilde{\eta}_{2}$ in the example $n =2$. Then we have
\begin{align*}
        X_{\SL_{2}} \times \mathbb{G}^{3}_{m}  & \to X_{\GL_2} \\
        ((D_{1},D_{2},D_{3}), (\gamma_1, \gamma_2 , \gamma_3)) & \mapsto (\gamma_1 D_1 , \gamma_2 D_2 , \gamma_3 D_3) 
    \end{align*}
where we view $D_{i}$ as diagonal matrices in the maximal torus $H_{\SL_2}$ up to permutation by the Weyl group. Now fix matrices $(C_{1}, C_{2},C_{3}) \in X_{\GL_2}$, which we view as diagonal matrices in the maximal torus $H_{\GL_2}$. Now writing $\alpha_{i} = (\det C_{i})^{\frac{1}{2}}$ the preimage under $\widetilde{\eta}_{2}$ consists of 
\begin{equation}
    (( \alpha^{-1}_{1}C_{1}, \,   \alpha^{-1}_{2}C_{2}, \,   \alpha^{-1}_{3}C_{3}), ( \alpha_{1},  \alpha_{2} ,  \alpha_{3})).
\end{equation}
Noticing that for each $\alpha_i$ we can equivalently plug in $- \alpha_i$ we see that there are in total $8$ elements in the preimage, which gives that $\widetilde{\eta}_{2}$ is a $\mu^{3}_{2} \colon 1$ cover.
\end{ex}
\begin{prop} \label{diagram_sln}
We have 
\begin{enumerate}
    \item The map $\eta \colon \LocB_{L_{\SL_{n}, \lambda } \times \mathbb{G}_{m}} \to \LocB_{L_{\GL_{n}, \lambda }}$ induced by $L_{\SL_{n}, \lambda } \times \mathbb{G}_{m} \to L_{\GL_{n}, \lambda }$ is $(-1)$-shifted symplectic, \'etale and oriented.
    \item The square in the diagram \eqref{sln_square} is a pullback and $\widetilde{\eta}_{2}$ is a $\mu^{3}_{n} \colon 1$ cover.
\end{enumerate}
\end{prop}
\begin{proof}
The fact that $\eta$ and $\eta_{1}$ are \'etale follows from the fact that the map $\B (L_{\SL_{n}, \lambda } \times \mathbb{G}_{m}) \to \B L_{\GL_n , \lambda}$ is \'etale and mapping stacks preserve \'etale maps. We can see that the map $\B (L_{\SL_{n}, \lambda } \times \mathbb{G}_{m}) \to \B L_{\GL_n , \lambda}$ is $2$-symplectic using the decomposition of Lie algebras $\mathfrak{l}_{\mathfrak{sl}_n, \lambda} \oplus \mathbb{C} \cong \mathfrak{l}_{\mathfrak{gl}_n, \lambda}$. Since the AKSZ construction preserves symplectic maps we get that $\eta$ is $(-1)$-symplectic.
    We can directly see from the form of the cotangent complexes of $\B L_{\GL_n , \lambda}$ and $\B L_{\SL_n , \lambda} \times \mathbb{G}_m$  that the induced map $\eta^{*} \mathbb{L}_{\LocB_{L_{\SL_{n}, \lambda }}} \to \mathbb{L}_{\LocB_{L_{\SL_{n}, \lambda } \times \mathbb{G}_{m}}}$  is the isomorphism induced from the isomorphism $\mathfrak{l}_{\mathfrak{sl}_n, \lambda} \oplus \mathbb{C} \to \mathfrak{l}_{\mathfrak{gl}_n, \lambda}$. Since the map on cotangent complexes is an honest isomorphism of complexes and not a quasi-isomorphism the map $\eta$ immediately preserves volume forms and thus orientations. By the same reason the orientation is also $W_{L_{\lambda}}$-invariant. \par 
    We prove part $2$ using Proposition \ref{alper_prop}. Firstly the square commutes since it is induced by a $L_{\GL_{n}, \lambda }$ equivariant map $\Loc^{\ff}_{L_{\SL_{n}, \lambda }} \times \mathbb{G}^{3}_{m} \to \Loc^{\ff}_{L_{\GL_{n}, \lambda }}$. This also gives that $\eta_2$ is separated and representable. $\eta_{2}$ is \'etale given that both $\eta$ and $\eta_{1}$ are and using the $2$ out of $3$ property. Closed points are clearly preserved under this map. Finally, because we are quotienting $\Loc^{\ff}_{L_{\SL_{n}, \lambda }} \times \mathbb{G}^{3}_{m}$ by $L_{\GL_{n}, \lambda }$ and not $L_{\SL_{n}, \lambda } \times \mathbb{G}_{m}$ the stabilizers are also preserved. The fact that $\widetilde{\eta}_{2}$ is a $\mu^{3}_{n}$ cover follows by a direct computation similar to Example \ref{sl2_etale_coverex}. Indeed, recall from Definition \ref{good_moduli_loc} that $X_{G} = H^{3}_{G} / \! \! / W$. Note that for us when $G = L_{\GL_n , \lambda}$, $W = \prod \So_{\lambda_i}$. So we can think of the map $\widetilde{\eta}_{2}$ as 
    \begin{align*}
        X_{L_{\SL_n , \lambda}} \times \mathbb{G}^{3}_{m}  & \to X_{L_{\GL_n ,\lambda}} \\
        ((D_{1},D_{2},D_{3}), (\gamma_1, \gamma_2 , \gamma_3)) & \mapsto (\gamma_1 D_1 , \gamma_2 D_2 , \gamma_3 D_3) 
    \end{align*}
    where $D_{i}$ are diagonal matrices in the maximal torus of the Levi $L_{\SL_n , \lambda}$. 
    Now writing $\alpha_{i} = (\det C_{i})^{\frac{1}{n}}$ for one of the $n$-th roots of unity of $\det C_{i}$ the preimage under $\widetilde{\eta}_{2}$ consists of 
\begin{equation}
    (( \alpha^{-1}_{1}C_{1}, \,   \alpha^{-1}_{2}C_{2}, \,   \alpha^{-1}_{3}C_{3}), ( \alpha_{1},  \alpha_{2} ,  \alpha_{3})).
\end{equation}
where for $\alpha_i$ we can plug in any of the $n$-th roots of unity. Therefore there are a total of $|\mu^{n}_3|$ different choices. In other words, $\mu^{3}_n$-acts freely on the preimage since $\mu^{3}_n$ acts freely on $\mathbb{G}^{3}_m$ by multiplication despite the fact that $\mu^{3}_n$ \emph{does not} act freely on $X_{L_{\SL_n, \lambda}}$. See Example \ref{sln_to_pgln_ex} for more on this subtlety.
\end{proof}
\begin{prop}[Purity for $\SL_n$] \label{purity_sln}
Let $\lambda$ be a partition of $n$ of length $l$.
$\pi_{L_{\SL_n , \lambda}*} \varphi_{L_{\SL_n , \lambda}}$ is a pure complex of mixed Hodge modules with perverse cohomology bounded below with lowest non-zero degree $l-1 = \dim \ZZ(L_{\SL_{n}, \lambda})$. 
\end{prop}
\begin{proof}
    Consider the diagram \eqref{sln_square}. We will start by proving that
    \begin{equation}
        (\pi_{L_{\SL_n , \lambda}} \times \pi_{\mathbb{G}_{m}})_{*} (\varphi_{L_{\SL_n , \lambda}} \boxtimes \varphi_{\mathbb{G}_{m}}) \cong \widetilde{\eta}^{*}_{2}\pi_{L_{\GL_n , \lambda}*} \varphi_{L_{\GL_n , \lambda}}
    \end{equation}
    We first prove that $\eta_{1*} (\varphi_{L_{\SL_n , \lambda}} \boxtimes \varphi_{\mathbb{G}_{m}}) \cong \eta^{*}_{2} \varphi_{L_{\GL_n , \lambda}}$. To do this we use the pullback square \eqref{ses_m_square} to deduce that
\begin{align} \label{basechangearg}
     p^{*} \eta_{1*}(\varphi_{L_{\SL_n , \lambda}} \boxtimes \varphi_{\mathbb{G}_{m}}) &\cong p^{*}\eta_{1*}\eta^{*}_{1}\eta^{*}_{2} \varphi_{L_{\GL_n , \lambda}} \qquad \text{(since $\eta^* \cong \eta^{*}_{1}\eta^{*}_2 $ \text{ and } $ \eta^{*}\varphi_{L_{\GL_n , \lambda}} \cong \varphi_{L_{\SL_n , \lambda}} \boxtimes \varphi_{\mathbb{G}_{m}}$ )} \\
    & \cong  \widetilde{\eta}_{1*}\widetilde{p}^{*} \eta^{*}_{1} \eta^{*}_{2} \varphi_{L_{\GL_n , \lambda}} \qquad \text{(base change along \eqref{ses_m_square})} \nonumber \\
    & \cong \widetilde{\eta}_{1*} \nonumber \widetilde{\eta}^{*}_{1}p^{*}\eta^{*}_{2}\varphi_{L_{\GL_n , \lambda}} \qquad ( \widetilde{\eta}^{*}_{1}p^{*} \cong \widetilde{p}^{*} \eta^{*}_{1} ) \nonumber \\
    & \cong p^{*} \eta^{*}_{2} \varphi_{L_{\GL_n , \lambda}}. \nonumber
\end{align}
Pulling back by $\widetilde{\eta}_{1} \colon \Loc^{\ff}_{L_{\SL_n , \lambda}} \times \mathbb{G}^{3}_{m}\times \B \mu_{n} \to \Loc^{\ff}_{L_{\SL_n , \lambda}} \times \mathbb{G}^{3}_{m}$ on sheaves we get the trivial $\mu_{n}$-equivariant structure. Therefore, the last equality above follows since $\widetilde{\eta}^{*}_{1}p^{*}\eta^{*}_{2}\varphi_{L_{\GL_n , \lambda}}$ has the trivial $\mu_{n}$-equivariant structure and pushing forward by $\widetilde{\eta}_{1}$ only picks up the cohomology of $\B \mu_{n}$, which is trivial. The functor $p^{*}$ is conservative since it is the pullback from a quotient stack so is the forgetful functor from $L_{\GL_n , \lambda}$-equivariant sheaves to sheaves. Therefore, we can conclude that  $ \eta_{1*}(\varphi_{L_{\SL_n , \lambda}} \boxtimes \varphi_{\mathbb{G}_{m}}) \cong \eta^{*}_{2} \varphi_{L_{\GL_n , \lambda}}$.
Since the square in diagram \eqref{sln_square} is a pullback we get the equation $\pi_{*}\eta^{*}_{2} \varphi_{L_{\GL_n , \lambda}} \cong \widetilde{\eta}^{*}_{2}\pi_{L_{\GL_n , \lambda}*} \varphi_{L_{\GL_n , \lambda}}$. Hence we have
\begin{equation} \label{final_big_basechange}
    (\pi_{L_{\SL_n , \lambda}} \times \pi_{\mathbb{G}_{m}})_{*} (\varphi_{L_{\SL_n , \lambda}} \boxtimes \varphi_{\mathbb{G}_{m}}) \cong \pi_{*}\eta_{1*} (\varphi_{L_{\SL_n , \lambda}} \boxtimes \varphi_{\mathbb{G}_{m}}) \cong \pi_{*} \eta^{*}_{2} \varphi_{L_{\GL_n , \lambda}} \cong \widetilde{\eta}^{*}_{2}\pi_{L_{\GL_n , \lambda}*} \varphi_{L_{\GL_n , \lambda}}.
\end{equation}
This implies that  $(\pi_{\SL_{n}} \times \pi_{\mathbb{G}_{m}})_{*} (\varphi_{\SL_{n}} \boxtimes \varphi_{\mathbb{G}_{m}})$ is pure as a complex of mixed Hodge modules since $\pi_{L_{\GL_n , \lambda}*} \varphi_{L_{\GL_n , \lambda}}$ is pure and $\widetilde{\eta}_{2}$ preserves purity by Lemma \ref{purity_local} since it is \'etale. We can rewrite 
\begin{equation}
    (\pi_{L_{\SL_{n},\lambda}} \times \pi_{\mathbb{G}_{m}})_{*} (\varphi_{L_{\SL_{n},\lambda}} \boxtimes \varphi_{\mathbb{G}_{m}}) \cong p^{*}_{X}\pi_{L_{\SL_{n},\lambda}*}\varphi_{L_{\SL_{n},\lambda}} \otimes p^{*}_{\mathbb{G}^{3}_{m}}\pi_{\mathbb{G}_{m}*} \varphi_{\mathbb{G}_{m}}
\end{equation}
using the projections 
$$p_X \colon X_{L_{\SL_{n}, \lambda}} \times \mathbb{G}^{3}_{m} \to X_{L_{\SL_{n}, \lambda}} \text{ and } p_{\mathbb{G}^{3}_{m}} \colon X_{L_{\SL_{n}, \lambda}} \times \mathbb{G}^{3}_{m} \to \mathbb{G}^{3}_{m}.$$
\par
Because $\Loc_{\mathbb{G}_{m}}$ is smooth we have 
\begin{equation}
    \pi_{\mathbb{G}_{m}*} \varphi_{\mathbb{G}_{m}} \cong \mathbb{Q}_{\mathbb{G}^{3}_{m}}[2] \otimes \HHf(\B \mathbb{G}_{m}).
\end{equation}
 Therefore, we can take the summand $\mathbb{Q}_{X_{L_{\SL_n , \lambda}} \times \mathbb{G}^{3}_{m}}[2]$
 \begin{equation}
     p^{*}_{X}\pi_{L_{\SL_{n},\lambda}*}\varphi_{L_{\SL_{n},\lambda}} \otimes (\mathbb{Q}_{X_{L_{\SL_n , \lambda}} \times \mathbb{G}^{3}_{m}}[2] \otimes \HHf(\B \mathbb{G}_{m}))
 \end{equation}
 to get $p^{*}_{X} \pi_{L_{\lambda,\SL_{n}*}} \varphi_{L_{\lambda,\SL_{n}}}[2]$. This implies that $p^{*}_{X} \pi_{L_{\lambda,\SL_{n}*}} \varphi_{L_{\lambda,\SL_{n}}}$ is pure and thus that $\pi_{L_{\lambda,\SL_{n}*}} \varphi_{L_{\lambda,\SL_{n}}}$ is pure since it is pure under the smooth projection $p_{X}$. In the last step we used the locality of purity as in Lemma \ref{purity_local}.
\end{proof}
\begin{prop}[Supports for $\SL_n$] \label{supports_sln}
We have a decomposition 
\begin{equation} \label{decomps_sln_pf}
    \pi_{\SL_{n}*} \varphi_{\SL_{n}} \cong \bigoplus_{i \geq 0} \bigoplus_{\lambda}\IC_{X^{\lambda}_{\SL_{n}}}(\mathcal{L}^{\lambda}_{i})[-i]
\end{equation}
for some local systems $\mathcal{L}^{\lambda}_{i}$ on $X^{\lambda}_{\SL_n}$. The index $i$ corresponds to the perverse cohomology degree and $\lambda$ is a partition of $n$ giving the corresponding stratum $X^{\lambda}_{\SL_n}$.
\end{prop}
\begin{proof}
    We start by computing the lowest non-zero perverse degree of $\pi_{L_{\SL_{n}, \lambda}*} \varphi_{L_{\SL_{n}, \lambda}}$. Note that since $p^{*}_{X}$ is smooth of relative dimension $3$ we have 
\begin{equation} \label{proj_equation}
    p^{*}_{X}\pH^{i}\pi_{L_{\SL_{n}, \lambda}*} \varphi_{L_{\SL_{n}, \lambda}} \cong \pH^{i-3} p^{*}_{X} \pi_{L_{\SL_{n}, \lambda}*} \varphi_{L_{\SL_{n}, \lambda}}
\end{equation}

Tensoring by $\pi_{\mathbb{G}_{m}*}\varphi_{\mathbb{G}_{m}}$ and using equation \eqref{final_big_basechange}  we get
\begin{align} \label{sum_equation}
    \pH^{i} \widetilde{\eta}^{*}_{2}\pi_{L_{\GL_{n}, \lambda *}} \varphi_{L_{\GL_n, \lambda}} & \cong  \pH^{i}( p^{*}_{X} \pi_{L_{\SL_{n}, \lambda}*} \varphi_{L_{\SL_{n}, \lambda}} \otimes \mathbb{Q}_{X_{L_{\SL_n , \lambda}} \times \mathbb{G}^{3}_{m}}[2] \otimes \HHf(\B \mathbb{G}_{m})) \\
    & \cong  \bigoplus_{j \geq 0} \pH^{i-2j+2} p^{*}_{X} \pi_{L_{\SL_{n}, \lambda}*} \varphi_{L_{\SL_{n}, \lambda}} \nonumber   \qquad \text{( writing } \HHf(\B \mathbb{G}_{m} ) \cong \oplus_{j \geq 0} \mathbb{Q}[-2j]) \\
    & \cong \bigoplus_{j \geq 0 }p^{*}_{X}\pH^{i-2j-1}\pi_{L_{\SL_{n}, \lambda}*} \varphi_{L_{\SL_{n}, \lambda}} \nonumber \qquad \text{( using equation \eqref{proj_equation})}.
\end{align}
We know from the $\GL_n$ version of  Lemma \ref{loc_sys_comp} that  $\pi_{L_{\GL_n , \lambda}*} \varphi_{L_{\GL_n , \lambda}}$ has lowest degree $l = \dim \ZZ(L_{\GL_n, \lambda})$. Let us first consider the case $L_{\SL_n , \lambda} = \SL_n$. In this case $l = 1$ and  $\pH^{1} \widetilde{\eta}^{*}_{2}\pi_{\GL_{n} *} \varphi_{\GL_n}$ is the lowest perverse degree of $\pi_{\GL_{n} *} \varphi_{\GL_n}$. Therefore, by plugging in $i = -1 , 0 $ into equation \eqref{sum_equation} we get that $p^{*}_{X} \pi_{L_{\SL_{n}, \lambda}*} \varphi_{L_{\SL_{n}, \lambda}}$ has no perverse cohomology in negative degrees and by plugging in $i = 1$ we see that the lowest non-zero perverse cohomology is in degree $0$. The general case follows similarly by plugging in $i \leq l$ into Equation \eqref{sum_equation} we see that $\pi_{L_{\SL_{n}, \lambda}*} \varphi_{L_{\SL_{n}, \lambda}}$ has non-zero perverse cohomology in lowest degree $l-1$. \par
   To prove the decomposition \eqref{decomps_sln_pf} we proceed as in Corollary \ref{support_locsys}. 
Let $\mathcal{F}$ be some summand of the $k$-th perverse cohomology of $\pi_{\SL_{n}*} \varphi_{\SL_{n}}$. We will use Lemma \ref{ic_achar} to show that $\mathcal{F}$ is an $\IC$ sheaf supported on $X^{\lambda}_{\SL_n}$.  Pulling back $\mathcal{F} $ by $p_{X}$ we have to get some summand of the $k$-th perverse cohomology of $\widetilde{\eta}^{*}_{2} \pi_{\GL_n *} \varphi_{\GL_n}$ so we get by Corollary \ref{support_locsys}
\begin{equation} \label{proj_sln_ic}
    p^{*}_{X} \mathcal{F} \cong \widetilde{\eta}^{*} \IC_{X^{\lambda}_{\GL_n}}(\mathcal{L}^{\lambda}_{k}) \cong \IC_{X^{\lambda}_{\SL_n} \times \mathbb{G}^{3}_{m}}(\widetilde{\eta}^{*}_{2}\mathcal{L}^{\lambda}_{k}).
\end{equation}
The last isomorphism follows by noting that the pullback under $\widetilde{\eta}_{2}$ of $X^{\lambda}_{\GL_{n}}$ is $X^{\lambda}_{\SL_{n}} \times \mathbb{G}^{3}_{m}$. So $p^{*}_{X} \mathcal{F}$ must be supported on $\overline{X}^{\lambda}_{\SL_{n}} \times \mathbb{G}^{3}_{m}$ for some $\lambda$ and so $\mathcal{F}$ is supported on $\overline{X}^{\lambda}_{\SL_{n}}$. Then we want to show that $\mathcal{F}$, pulled back to $X^{\lambda}_{\SL_{n}}$, is a local system. We have a commutative diagram
   \begin{equation} \label{supp_sln_diag}
\begin{tikzcd}
	{X^{\lambda}_{\SL_{n}} \times \mathbb{G}^{3}_{m}} & {X_{\SL_{n}} \times \mathbb{G}^{3}_{m}} \\
	{X^{\lambda}_{\SL_{n}}} & {X_{\SL_{n}}}
	\arrow[hook, from=1-1, to=1-2]
	\arrow["{p^{\lambda}_{X}}", from=1-1, to=2-1]
	\arrow["{p_{X}}", from=1-2, to=2-2]
	\arrow[hook, from=2-1, to=2-2]
\end{tikzcd}
   \end{equation}
We can use equation \eqref{proj_sln_ic} to deduce that $p^{*}_{X} \mathcal{F}$ is a local system restricted to $X^{\lambda}_{\SL_{n}} \times \mathbb{G}^{3}_{m}$. This implies that $\mathcal{F}$ is a local system when restricted to $X^{\lambda}_{\SL_{n}}$ since we know it is once pulled back by the projection $p^{\lambda}_{X}$, using the commutative diagram \ref{supp_sln_diag}. The sheaf $\mathcal{F}$ will have no quotients or subobjects supported on $\overline{X}^{\lambda}_{\SL_n}  \setminus X^{\lambda}_{\SL_n}$ since the shifted pullback $p^{*}_{X}[3]$ is exact for the perverse $t$-structure and  $p^{*}_{X}[3] \mathcal{F}$ is an $\IC$ sheaf supported on $X^{\lambda}_{\SL_n} \times \mathbb{G}^{3}_{m}$.
\end{proof}
\begin{ex}[BPS sheaves for $\SL_n$]
    From the above proposition we can calculate the BPS sheaves for $\SL_n$ which are constant sheaves supported on the image of $\ZZ^{3}(\SL_n) = \mu^{3}_{n} \xhookrightarrow{} X_{\SL_{n}}$. So we get skyscrapers supported on $n^{3}$ points and furthermore the BPS sheaf only contributes in degree $0$ of $\pi_{\SL_n *} \varphi_{\SL_n}$.
\end{ex}
\begin{proof}[Proof of Lemma \ref{loc_sys_comp} for $\SL_n$] \label{loc_sys_comp_sl}
    From the previous Proposition \ref{supports_sln} we know that $\pi_{L_{\SL_n , \lambda}*} \varphi_{L_{\SL_n , \lambda}}$ has perverse cohomology bounded below with lowest perverse piece $l-1$, so  we  define 
    \begin{equation}
        \BPS_{L_{\SL_n , \lambda}} = \pH^{l-1}\pi_{L_{\SL_n , \lambda}*} \varphi_{L_{\SL_n , \lambda}}.
    \end{equation}
     Then the same argument for computing supports as in the proof of Proposition \ref{supports_sln} works, using the fact that $\supp \BPS_{L_{\GL_n , \lambda}} = \ZZ^{3}(L_{\GL_n , \lambda})$. We therefore get,
    \begin{equation}
        \supp \BPS_{L_{\SL_n , \lambda}} = \ZZ^{3}(L_{\SL_n , \lambda}).
    \end{equation}
    By the same argument as in Corollary \ref{bps_const}, we can then use the transitive action of $\ZZ^{3}(L_{\SL_{n}, \lambda})$ to prove that the BPS sheaf is constant of rank $1$. The same argument will work to show that any perverse piece supported on $\ZZ^{3}(L_{\SL_n, \lambda})$ is a constant sheaf. To compute the contributions of the Saito decomposition of $\pi_{L_{\SL_{n}, \lambda}*} \varphi_{L_{\SL_{n}}, \lambda}$ supported on $\ZZ^{3}(L_{\SL_n , \lambda})$ we will repeatedly use Part $(1)$ in the $\GL_n$ version of Lemma \ref{loc_sys_comp} and equation \eqref{sum_equation} to split off the extra $\HHf(\B \mathbb{G}_{m})$ factor in 
    $$\HHf(\B\ZZ(L_{\GL_n, \lambda})) \cong  \HHf(\B\ZZ(L_{\SL_n, \lambda})) \otimes \HHf(\B \mathbb{G}_{m}).$$ This is a lengthy computation so we give the proof in its own Lemma \ref{split_off_lemsln}.  \par
    The rest of the calculation in part $2$ is analogous to the $\GL_n$ case of the Lemma \ref{loc_sys_comp}. 
    For the computation we will need that the map $\theta \colon X_{L_{\SL_{n},\lambda}} \to X_{\SL_{n}}$ is finite. This follows from the fact that $X_{L_{\GL_{n}, \lambda}} \to X_{\GL_{n}}$ is finite, $X_{L_{\SL_{n},\lambda}} \to X_{L_{\GL_{n},\lambda}}$ is a closed immersion, composition of finite maps is finite and the map $\theta$ is a composition of these two maps. To finish mimicking the proof of the $\GL_n$ version of the Lemma we can use the pullback diagrams in Lemma \ref{tilde_center_pull_lemma}. The computation is then the same as in the $\GL_n$ case.
\end{proof}
\begin{lem} \label{split_off_lemsln}
    The components of the Saito decomposition of $\pi_{L_{\SL_n, \lambda}*} \varphi_{L_{\SL_n, \lambda}}$ with supports given by $\ZZ^{3}(L_{\SL_n, \lambda})$ are $\BPS_{L_{\SL_n,\lambda}} \otimes \HHf(\B \ZZ(L_{\SL_n, \lambda}))[- \dim \ZZ(L_{\SL_n, \lambda})]$.
\end{lem}
\begin{proof}
    Let us compute all of the perverse pieces of $\pi_{L_{\SL_n , \lambda}*} \varphi_{L_{\SL_n , \lambda}}$ that are supported on $\ZZ^{3}(L_{\SL_n , \lambda})$. Just for this proof call these 
    $$\pH^{m}( \SL_n) $$ to avoid clutter. We know from the $\GL_n$ version of Lemma \ref{loc_sys_comp} that the pieces supported on $\ZZ^{3}(L_{\GL_n , \lambda})$ are $\BPS_{L_{\GL_n , \lambda}} \otimes \HHf(\B \ZZ(L_{\GL_n ,\lambda}))[- \dim L_{\SL_{n}, \lambda}]$. These terms are in cohomological degrees $l + 2m$ for $m \geq 0$ with dimension $\binom{m+l}{l}$, the number of homogeneous polynomials in $l$-variables. Our goal is to prove using equation \eqref{sum_equation}  that on $\ZZ^{3}(L_{\SL_n , \lambda})$ we have terms with dimension $\binom{m+l-1}{l-1}$, the number of homogeneous polynomials in $(l-1)$-variables.  \par
    Again to avoid clutter we write 
    $$\pH^{m}(\BPS_{L_{\GL_n , \lambda}} \otimes \HHf(\B \ZZ(L_{\GL_n ,\lambda}))[- \dim L_{\SL_{n}, \lambda}]) = \pH^{m}( \GL_n)$$
     Now we can restrict equation \eqref{sum_equation} to the stratum $\ZZ^{3}(L_{\SL_n , \lambda}) \times \mathbb{G}^{3}_{m}$  to get
    \begin{equation} \label{pf_mideq}
     \widetilde{\eta}^{*}_{2}\pH^{l+k}(\GL_n) \cong \bigoplus_{j \geq 0 }p^{*}_{X} \pH^{l-1+k-2j} ( \SL_n). 
    \end{equation}
    Recall that $\widetilde{\eta}^{*}_{2}\pH^{l+k}(\GL_n) = 0$ for $k < 0$ so plugging  $k= -1$ and $k= -2$ into equation \eqref{pf_mideq} we can deduce
    \begin{equation*}
        p^{*}_{X} \pH^{l-1-m} ( \SL_n) = 0 
    \end{equation*}
    for $m > 0$. Similarly $\widetilde{\eta}^{*}_{2}\pH^{l+2k+1}(\GL_n) = 0$ for $k \geq 0$ so again from equation \eqref{pf_mideq} we can deduce that
    \begin{equation*}
        p^{*}_{X} \pH^{l-1+(2m+1)} ( \SL_n) = 0 
    \end{equation*}
    for $m \geq 0$.
    Therefore, we can simplify equation \eqref{pf_mideq} to 
     \begin{equation}
     \widetilde{\eta}^{*}_{2}\pH^{l+2k}(\GL_n) \cong \bigoplus_{0 \leq j \leq k }p^{*}_{X} \pH^{l-1+2k-2j} ( \SL_n). 
    \end{equation}
We now prove by induction that $\rank \pH^{l-1 +2k }   ( \SL_n) = \binom{k +l-1}{k}$, the number of homogenenous polynomials in $l-1$ variables of degree $k$. For $k= 0$ this follows from the fact that $\rank \BPS_{L_{\GL_n, \lambda}} = \rank \BPS_{L_{\SL_n , \lambda}} = 1$. Now assume the statement for $k$. We can use the identities
    \begin{equation}
        \binom{k+ l-1 }{k} = \binom{k-1+l}{k-1} + \binom{k+l -2}{k}
    \end{equation}
    for dimensions of homogenous polynomials in $l$ variables related in terms of polynomials in $l-1$ variables repeatedly to see that
    \begin{equation}
          \binom{k+ l-1 }{k} = \sum^{k}_{j \geq 0}  \binom{j+ l-1 }{j}
    \end{equation}
    so the number of homogenous polynomials in $l$ variables of degree $k$ is the sum of the number of homogenous polynomials in $l-1$ variables of degree from $0$ to $k$.
    Then we have
    \begin{align*}
        \rank \pH^{l+2k+2}( \GL_n) & = \sum_{0 \leq j \leq k+1} \rank \pH^{l-1 +2k +2 -2j }   ( \SL_n) \\
        & =  \rank \pH^{l-1 +2k +2 }   ( \SL_n) + \sum_{0 \leq j \leq k} \rank \pH^{l-1 +2k -2j }   ( \SL_n)
    \end{align*}
Now using the induction assumption we can write
\begin{equation*}
    \binom{k+1 + l-1 }{k+1} = \rank \pH^{l-1 +2k +2 }   ( \SL_n) + \sum_{0 \leq j \leq k} \binom{j+l-1}{j}
\end{equation*}

    Therefore, we have $\pH^{*} (\SL_n)$ is non-zero in degrees $l-1 + 2k$ with rank $\binom{l-1 +k}{l-1}$ which is exactly the dimensions and degrees of 
    $$\BPS_{L_{\SL_n, \lambda}}  \otimes \HHf(\B \ZZ (L_{\SL_n, \lambda}))[- \dim \ZZ(L_{\SL_n , \lambda})].$$ Note  $\HHf(\B \ZZ (L_{\SL_n, \lambda}))$ is a polynomial algebra in $l-1$ variables as computed in Lemma \ref{coh_centers}.
\end{proof}
 \begin{proof}[Proof of cohomological integrality for $\SL_n$]
 Now that we have established Lemma \ref{loc_sys_comp} for $\SL_n$ we also have Proposition \ref{levis_final}. Therefore, we can follow strategy \ref{strategy} and use the same argument as in the $\GL_{n}$ case.
 \end{proof}
 \subsection{Comparison of \texorpdfstring{$\SL_n$}{TEXT} and \texorpdfstring{$\PGL_n$}{TEXT} good moduli spaces} \label{sln_to_pgln_gms}
 We have a projection map $\SL_n \to \GL_n \to \PGL_n$, which is an \'etale $\mu_n \colon 1 $ cover. This map induces an \'etale map $\B \SL_n \to \B \PGL_n$, which is also $2$-shifted symplectic. Similarly there is an \'etale map $\B L_{\SL_n , \lambda } \to \B L_{\PGL_n \lambda }$  This induces an \'etale map $\LocB_{ L_{\SL_n , \lambda }} \to \LocB_{ L_{\PGL_n , \lambda }}$. The image of this map only hits the trivial component $\LocB^{1}_{ L_{\PGL_n , \lambda }}$. In this subsection we will consider the geometry of the induced map $\eta \colon X_{\SL_n} \to X_{\PGL_n}$ on good moduli spaces. We start with an example
\begin{ex} \label{sln_to_pgln_ex}
Let us consider the example $n =2$. Then we have
\begin{align*}
        \eta \colon X_{\SL_{2}}   & \to X_{\PGL_2} \\
        (D_{1},D_{2},D_{3}) & \mapsto ( \overline{D}_1 ,  \overline{D}_2 ,  \overline{D}_3) 
    \end{align*}
where we view $D_{i}$ as diagonal matrices in the maximal torus $H_{\SL_2}$ up to permutation by the Weyl group and $\overline{D}_{i}$ as their image in $\PGL_2$ where they land in $H_{\PGL_2}$. Now consider the element 
\begin{equation} \label{badpoint_2}
   x = (\diag(1,1) , \diag(1,1), \diag(i,-i)) \in X_{\PGL_2} 
\end{equation}
 then note that under the map $\eta \colon X_{\SL_2} \to X_{\PGL_2}$ the element $x$ has preimage of size smaller than $8$ since $\diag(i,-i)$ and $-1 \cdot \diag(i,-i)$ define the same element under permutation. However, a generic element in $X_{\PGL_2}$ will be covered by $\mu^{3}_{2}$, which has $8$ elements. Therefore, $\eta$ will not be an \'etale cover but only a finite map. Note however that only $x$ that contain matrices in the centre or the matrix  $\diag(i,-i)$ will have smaller preimage. This follows because $\diag(i,-i)$ is the only matrix that has non trivial $\mu_2$ stabiliser via the action of $\mu_2$ on $H_{\SL_2} / \! \! / W$.
\end{ex}
To work with the map $\eta$ we will need to throw out the bad points such as in equation \eqref{badpoint_2}. We will now examine how the map $\eta$ behaves with respect to the stratification $X^{\lambda}_{\SL_n}$. For the rest of this section assume $n$ is prime and fix $\omega = e^{2 \pi i / n}$ a $n$-th root of unity and $W \cong \So_{n}$ the Weyl group of $\SL_n$. Denote by $A$ the matrix
\begin{align} \label{special_matrices}
A & = \diag(1 , \omega , \dots , \omega^{n-1}) \in H_{\GL_n}. \\ 
A_{S} & =  W \cdot  A \subseteq H_{\SL_n}  \text{ for } n \geq 3.  \nonumber \\
A_{S} & =  W \cdot  iA \subseteq H_{\SL_2}.  \nonumber
\end{align}
Note that $\det A = (-1)^{n+1}$ so since $n$ is prime $A \in \SL_n$ for $n \geq 3$. We see that the image of $A_{S}$ under the quotient map $H_{\SL_n} \to H_{\SL_n} / \! \! / W $ is just a single point. \par
 We will now prove a few technical lemmas we need to check when we can restrict the map $\eta$ to become an \'etale cover. This extra subtlety arises because the map $H_{\SL_n} \to H_{\PGL_n}$ does not preserve stabilisers under the action of $W$. See also Example \ref{sln_to_pgln_ex}. \par
 We now will write down the set of points where the map $\eta \colon X_{\SL_n} \to X_{\PGL_n}$ fails to be \'etale. We start by proving that the map $\eta$ is finite.
\begin{lem} \label{pgln_finite_map}
    Let $\lambda$ be any partition of $n$. 
    \begin{enumerate}
        \item The map $\eta \colon X_{L_{\SL_n, \lambda}} \to X_{L_{\PGL_n , \lambda}}$ is finite.
        \item $\eta$ is \'etale when restricted to an open subset $U \subseteq X_{L_{\SL_n , \lambda}}$ with $\ZZ^{3}(L_{\GL_n}) \subseteq U$.
    \end{enumerate}
\end{lem}
\begin{proof}
    Note that we have the \'etale $\mu^{3}_n$ cover $H_{L_{\SL_n , \lambda}} \to H_{L_{\PGL_n, \lambda}}$. Then we have a commutative diagram 
    \begin{equation}
\begin{tikzcd}
	{H^{3}_{L_{\SL_n , \lambda}}} & {H^{3}_{L_{\PGL_n , \lambda}}} \\
	{X_{L_{\SL_n , \lambda}}} & {X_{L_{\PGL_n , \lambda}}}
	\arrow[from=1-1, to=1-2]
	\arrow[from=1-1, to=2-1]
	\arrow[from=1-2, to=2-2]
	\arrow[from=2-1, to=2-2]
\end{tikzcd}
    \end{equation}
    since the quotient maps are finite the fact that $\eta$ is finite then follows from the $2$ out of $3$ property for finite maps. This proves the first part. \par
    For the second part we can use that stabiliser preservation is an open condition. In particular, since we are working with a DM stack $H^{3}_{\PGL_n} / W$, the map $I_{H^{3}_{\PGL_n} / W} \to H^{3}_{\PGL_n} / W$ will be proper. Here $I_{H^{3}_{\PGL_n} / W}$ is the inertia stack. Therefore, we can use \cite[Proposition 2.5]{alper_local_quotient} or \cite[Proposition 3.5]{Rydh_2013} to conclude that the locus of stabiliser preserving points $U^{'} \subseteq H^{3}_{L_{\SL_n , \lambda}}$ is open. Write $\varpi \colon H^{3}_{L_{\SL_n} , \lambda} \to X_{L_{\SL_n , \lambda}}$ for the quotient map. Then by Proposition \ref{alper_prop} restricted to $U = \varpi(U^{'})$ the map $\eta$ is  \'etale. $U$ is open since $X_{L_{\SL_n , \lambda}}$ has the quotient topology.  Indeed, $U^{'}$ is a $W$-invariant subset so $\varpi^{-1}(\varpi(U^{'})) = U^{'}$ and therefore $U$ is open. Now the Weyl group acts trivially on both the centre $\ZZ^{3}(L_{\SL_n , \lambda})$ and $\ZZ^{3}(L_{\PGL_n , \lambda})$. Therefore, the stabilisers of the centres will be preserved under the map $H^{3}_{L_{\SL_n , \lambda}} \to H^{3}_{L_{\PGL_n , \lambda}}$ so they will be contained in $U^{'}$.
\end{proof}
\begin{defn} \label{bad_locus}
\,
\begin{enumerate}
    \item  Define the set $\bb  \subseteq H^{3}_{\SL_n}$ as the set where $x = (D_{1}, D_{2}, D_{3}) \in \bb $ if \emph{all} the $D_{i}$ are in  $A_{S} \cup \ZZ( \SL_n)$  and \emph{at least} one of the $D_i$ \emph{must} be in $A_{S}$.
    \item Denote by $\bb_{S}$ the image of $\bb $ under $H^{3}_{\SL_n} \to X_{\SL_n}$.
    \item Denote by $\bb^{'}_{P}$ the image of $R$ under the map $H^{3}_{\SL_n} \to H^{3}_{\PGL_n}$ and by $\bb_{P}$ the image of $\bb_{S} $ under $X_{\SL_n} \to X_{\PGL_n}$.
\end{enumerate}
\end{defn}

\begin{notation}
Denote by $\lambda^{H}$ the partition $(1, \cdots , 1)$. 
\end{notation}
We will first consider what happens under the action of $W \times \mu_n$ on $H_{\SL_n}$. We will then use this to compute for the $3$-torus. We will want to use Proposition \ref{alper_prop} so we need to understand the stabilisers of the action of $W \times \mu^{3}_{n}$ on $H^{3}_{\SL_n}$.
\begin{lem} \label{1_loop_etale}
Under the action of $W \times \mu_n$ on $H_{\SL_n}$ the only elements $x$ that can be stabilised by $(\sigma, \zeta)$ with $\zeta \neq 1$  are given by the elements $A_{S}$. Therefore, these are the only elements $x \in H_{\SL_n}$ such that $\operatorname{Stab}(x)$ is not contained in  the subgroup $W \times 1$. 
\end{lem}
\begin{proof}
    For $n =2$ we can directly compute to prove the lemma. Namely as in Example \ref{sln_to_pgln_ex} we see that the only matrix in $H_{\SL_2}$ that is stabilised by $(\sigma, -1)$ is $\diag(i,-i)$.  So assume that $n \geq 3$. Fix an element $x \in H_{\SL_n}$, up to permutation we can write it as 
    \begin{equation}
        x= \diag( \underbrace{x_{1}, \dots , x_{1} , }_{\lambda_{1}\text{ times}}  \underbrace{x_{2}, \dots , x_{2} , }_{\lambda_{2}\text{ times}} \dots \underbrace{x_{l} , \dots , x_{l}}_{\lambda_{l}\text{ times}} ) 
    \end{equation}
    with $x_{i} \neq x_{j}$ if $i \neq j$ for some partition $\lambda$ of $n$.  Assume that $(\sigma , \zeta)$ stabilises $x$ so that we have $(\sigma, \zeta) \cdot x = \zeta \sigma(x) = x$. Assume that $\sigma$ fixes a block corresponding to some  $\lambda_i$. This will imply 
    \begin{equation}
        x_{i} = \zeta \sigma(x_i) = \zeta x_{i} .
    \end{equation}
Hence, $\zeta$ = 1. Therefore, if $\zeta \neq 1$ the permutation $\sigma$ must permute the blocks corresponding to $\lambda_i$. However, because $n$ is prime $\gcd(\lambda_1 , \dots , \lambda_l) = 1$. This means that unless $\lambda = \lambda^{H}$ the permutation $\sigma$ must mix at least one block with another. Assume now that $\lambda \neq \lambda^{H}$. This implies that the following equations will have to hold for some $k  \neq l$ and $m$
\begin{align*}
    \zeta \cdot x_{k} & = x_{m} \\
    \zeta \cdot x_{l} & = x_{m}.
\end{align*}
So after applying the permutation both $x_{k}$ and $x_{l}$ are in the $\lambda_{m}$ block of the original element $x$. The above equation implies that $x_{k} = x_{l}$ but this contradicts our assumption that $x_{k}$ and $x_{l}$ are distinct. So we can conclude that if $(\sigma, \zeta)$ stabilises $x$ we must have $ (\sigma, \zeta) = (\sigma,1)$. Note that if $G_{x} \subseteq W \times 1$ then this will also be true for any $G_{(\sigma, \zeta) \cdot x}$ since the two stabiliser groups are conjugate. \par
By the argument we gave above we can see that the only possible matrices that may have stabiliser not contained in $W \times 1$ are the ones with blocks corresponding to $\lambda^{H}$. In particular, this implies that all the elements $x_{i}$ are distinct. Pick a diagonal matrix $D$  with distinct entries $x_{i}$. If $ \zeta \sigma D = D$, then $\sigma$ does not fix any $x_i$ and so since $n$ is prime must be an $n$-cycle. Indeed, decompose $\sigma$ into disjoint cycles and assume that there is a cycle $\sigma^{'}$ of length $k <n$, containing $i$. Then 
\begin{equation}
    x_{i}= \zeta x_{\sigma^{'}(i)}
\end{equation}
iterating this equation we will get $x_{i} = \zeta^{k} x_{i}$, which implies that $\zeta^{k} = 1$. This contradicts the fact that $\zeta$ is a primitive root of unity of $n$.
Since $\sigma$ is an $n$-cycle we can then write for any $l$
\begin{align*}
    x_{l} & = x_{\sigma^{m_{l}}(1)} \text{ for some $m_{l}$} \\
    x_{\sigma^{m_{l}}(1)} & = \omega^{k_{l}} x_{1} \text{ for some $1 \leq k_{l} \leq n-1$}
\end{align*}
Therefore, $x_{i} = \omega^{k_i} x_{1}$ for any $x_{i}$ and also $\omega^{k_i} \neq \omega^{k_{j}}$ for $i \neq j$. Since $D \in \SL_n$ this implies that $x^{n}_{1} =1$ and so $D \in A_{S}$.
\end{proof}

We now define a stratification for $X_{\PGL_n}$ using Definition \ref{stratifications} and the spaces $\bb$ in Definition \ref{bad_locus}.
\begin{defn}[Stratification for $\PGL_n$] \label{strat_pgl_n}
    Recall the map $\eta \colon X_{\SL_n} \to X_{\PGL_n}$. We define a stratification for $X_{\PGL_n}$ in the following way
    \begin{enumerate}
        \item $X^{\lambda}_{\PGL_n}   = \eta(X^{\lambda}_{\SL_n})$, $\lambda \neq \lambda^{H}$.
        \item $X^{\lambda,g}_{\PGL_n}   = \eta(X^{\lambda,g}_{\SL_n})$,  $\lambda \neq \lambda^{H}$ .
        \item $X^{\lambda^{H}}_{\PGL_n} = \eta(X^{\lambda^{H}}_{\SL_n} \setminus \bb_{S})$.
        \item $X^{\lambda^{H},g}_{\PGL_n} = \eta(X^{\lambda^{H},g}_{\SL_n} \setminus \bb_{S})$.
        \item  $X^{\mu}_{\PGL_n} = \eta(\bb_{S}) = \bb_{P}$.
    \end{enumerate}
\end{defn}
\begin{notation}
    We will refer to the various restrictions of the map $\eta \colon X_{\SL_n} \to X_{\PGL_n}$ by the same symbol to reduce clutter.
\end{notation}
Note that we can write the affine GIT quotient $H^{3}_{\PGL_n} / \! \! / W$ as $H^{3}_{\SL_n} / \! \! / (W \times \mu^{3}_n)$. Then the map $\eta \colon X_{\SL_n} = H^{3}_{\SL_n} / \! \! / W \to  H^{3}_{\PGL_n} / \! \! /  W = X_{\PGL_n}$ is induced by the map
    \begin{equation} \label{rando_eq_2}
        H^{3}_{\SL_n} \xrightarrow{ \id} H^{3}_{\SL_n},
    \end{equation}
     which equivariant with respect to the homomorphism $W \to W \times \mu^{3}_n$ $\sigma \mapsto (\sigma, (1,1,1))$.  Here $\mu^{3}_n$ acts as the centre of $\SL_n$ on each factor. We say that an element $(D_1,D_2,D_3) \in H^{3}_{\SL_n}$  corresponds to a partition $\lambda$ if it is in the stratum $X^{\lambda}_{\SL_n}$ under the quotient map. The following Lemma explains the relation to the $\SL_n$ stratification.
\begin{lem}[Stratifications for $\PGL_n$ and \'etale covers] \label{pgln_strat_etale_lem}
Let $n$ be prime.  
\begin{enumerate}
    \item The map in  equation \eqref{rando_eq_2}
    \begin{enumerate}
        \item  \text{preserves stabilisers restricted to} points corresponding to partitions $\lambda \neq \lambda^{H}$.
        \item  preserves stabilisers restricted to points corresponding to the partition $\lambda^{H}$ except at the points in $\bb$.
    \end{enumerate}
 \item The induced map $\eta \colon X_{\SL_n} \setminus \bb_{S} \to X_{\PGL_n} \setminus \bb_{P}$ is a $\mu^{3}_n \colon 1$ cover, in particular $\mu^{3}_n$ acts freely on $X_{\SL_n} \setminus \bb_{S}$. Restricting, we will also get $\mu^{3}_n \colon 1$ covers: 
 \begin{align} \label{etale_eq_strat}
     X^{\lambda}_{\SL_n} & \to X^{\lambda}_{\PGL_n}  \quad \lambda \neq \lambda^{H} \\
      X^{\lambda^{H}}_{\SL_n} \setminus \bb_{S} & \to X^{\lambda^{H}}_{\PGL_n}. \nonumber \\
       X^{\lambda,g}_{\SL_n} & \to X^{\lambda,g}_{\PGL_n}  \quad \lambda \neq \lambda^{H} \nonumber  \\
      X^{\lambda^{H},g}_{\SL_n} \setminus \bb_{S} & \to X^{\lambda^{H},g}_{\PGL_n}. \nonumber 
 \end{align}
\end{enumerate}
\end{lem}
\begin{proof}
Part $(1)$ is an analogue of the statement of Lemma \ref{1_loop_etale} in the $3$-dimensional case. To start, assume we have $x = (D_{1} , D_{2} , D_{3}) \in H^{3}_{\SL_n}$ and $x$ corresponds to some partition $\lambda^{x}$. Each $D_{i}$ has blocks corresponding to some partitions $\lambda^i$.  If $x$ is stabilised, we must have
\begin{equation}
    (\sigma, (\zeta_{1}, \zeta_{2} , \zeta_{3})) \cdot (D_{1} , D_{2}, D_{3}) = (\zeta_{1} \sigma(D_{1}), \zeta_{2} \sigma (D_{2}) , \zeta_{3} \sigma (D_{3})) = (D_{1},D_{2},D_{3}).
\end{equation}
Assume that $\lambda^{x} \neq \lambda^{H}$, then we claim that the map in equation \eqref{rando_eq_2} preserves the stabiliser of $x$. Note that since we have assumed that $\lambda^{x} \neq \lambda^{H}$, by Remark \ref{special_part} we must also have that $\lambda^{i} \neq \lambda^{H}$ for all $i$. Therefore, we can apply the result of part $1$ of Lemma \ref{1_loop_etale} $3$ times for $D_{i}$ $1 \leq i \leq 3$ to prove part $(1)$ $(a)$. \par
Now assume that $\lambda^{x} = \lambda^{H}$. Assume that $(\sigma, (\zeta_{1}, \zeta_{2}, \zeta_{3}))$ stabilises $x$
for $(\zeta_{1}, \zeta_{2}, \zeta_{3}) \neq (1,1,1)$. This implies that for some $i$ we have $\zeta_{i}\sigma(D_{i}) = D_{i}$ with $\zeta_{i} \neq 1$. Therefore, by Lemma \ref{1_loop_etale} we get that $D_{i} \in A_{S}$ and that $\sigma$ is an $n$-cycle.  Because $\sigma$ is an $n$-cycle it can only stabilise elements in the centre. Therefore, if $\zeta_{i} =1$ $D_{i}$ must be in the centre and in $A_{S}$ if $\zeta_{i} \neq 1$. This implies $x \in R$. \par
Part $2$ follows by part $(1)$ and  Proposition \ref{alper_prop} applied to the map 
\begin{equation} \label{rando_map}
    (H^{3}_{\SL_n} \setminus \bb) / W \to (H^{3}_{\SL_n} \setminus \bb) / (W \times \mu^{3}_n) .
\end{equation}
Where again we quotient by $W$ on the left and $W \times \mu^{3}_{n}$ on the right. We have now established that the map $\eta \colon X_{\SL_n} \setminus \bb_{S} \to X_{\PGL_n} \setminus \bb_{P}$ is \'etale so to show it is a $\mu^{3}_n \colon 1$ cover we just have to consider a preimage. Since it is given by quotenting by $\mu^{3}_n$ it has at most $|\mu_n|^{3}$ elements. Now assume that $\mu^{3}_n$ does not act freely on the preimage of an element $y \in X_{\PGL_n} \setminus \bb_{P}$. Then it must be that some element $z \in X_{\SL_n} \setminus \bb_{S}$ with $\eta(z) = x$ has $(\zeta_{1}, \zeta_{2}, \zeta_{3}) \cdot z = \sigma (z)$ for some permutation $\sigma $. But this would imply that the  stabilisers of the lift of $z$ to $H^{3}_{\SL_n}$ are not preserved under the map in equation \eqref{rando_map}.   
We can now just restrict the map $\eta$ to the various strata to get that the maps in equation \eqref{etale_eq_strat} are also $\mu^{3}_n \colon 1$-covers.
\end{proof}
Lemma \ref{pgln_strat_etale_lem} then gives that the  stratification in Definition \ref{strat_pgl_n} is by locally closed smooth subvarieties in $X_{\PGL_n}$.
The next lemma proves similar results but for the map induced by the inclusion of Levi subgroups $L_{\PGL_n , \lambda} \subseteq \PGL_n$. Recall the spaces in Definition \ref{tilde_centers_def}.
\begin{lem} \label{levi_etale_lemma_pgln}
    We get $W_{L_{\lambda}} \colon 1$ covers
\begin{align} 
    \widetilde{Z}^{3}(L_{\PGL_n , \lambda}) & \to X^{\lambda}_{\PGL_n} \quad \lambda \neq \lambda^{H} \\
    \widetilde{Z}^{3}(L_{\PGL_n , \lambda^{H}}) \setminus \bb^{'}_{P} & \to X^{\lambda^{H}}_{\PGL_n}. \nonumber \\
    \widetilde{Z}^{3,g}(L_{\PGL_n , \lambda}) & \to X^{\lambda,g}_{\PGL_n} \quad \lambda \neq \lambda^{H} \\
    \widetilde{Z}^{3,g}(L_{\PGL_n , \lambda^{H}}) \setminus \bb^{'}_{P} & \to X^{\lambda^{H},g}_{\PGL_n}. \nonumber
\end{align}
\end{lem}
\begin{proof}
    Let $\lambda \neq \lambda^{H}$ and consider the $W_{L_{\lambda}} \colon 1$ cover $\widetilde{\ZZ}^{3}(L_{\SL_n , \lambda}) \to X^{\lambda}_{\SL_n}$. This map is equivariant with respect to $\mu^{3}_n$ and $\mu^{3}_n$ acts freely on both sides by Lemma \ref{pgln_strat_etale_lem}. Taking the quotient we get the map $\widetilde{\ZZ}^{3}(L_{\PGL_n , \lambda}) \to X^{\lambda}_{\PGL_n}$, which must also be a $W_{L_{\lambda}} \colon 1$ cover. Indeed, assume that some point $x \in \widetilde{\ZZ}^{3}(L_{\PGL_n , \lambda})$ has non-trivial $W_{L_{\lambda}}$-stabiliser. Then it must be that 
    $$\sigma(\tilde{x}) = (\zeta_{1}, \zeta_{2}, \zeta_{3}) \cdot \tilde{x}$$ or equivalently 
    $$(\zeta^{-1}_{1}, \zeta^{-1}_{2}, \zeta^{-1}_{3}) \cdot \sigma(\tilde{x}) = \tilde{x}$$
    for some $\sigma \in W_{L_{\lambda}}$, some $\zeta = (\zeta_{1}, \zeta_{2}, \zeta_{3}) \in \mu^{3}_{n}$ and some lift $\tilde{x} \in \widetilde{\ZZ}^{3}(L_{\SL_n , \lambda})$ of $x$. However, by the same arguments as in the proofs of  Lemma \ref{1_loop_etale} this cannot happen. Namely, the permutation $\sigma$ permutes blocks of the same size and so because $n$ is prime and $\gcd (\lambda_1 , \cdots \lambda_l) = 1$, $\sigma$ will have to fix some blocks in the matrices in $\tilde{x}$, which implies that $\zeta = (1,1,1)$. However, $W_{L_{\lambda}}$ acts freely on $\widetilde{\ZZ}^{3}(L_{\SL_n , \lambda})$ so $\sigma(\tilde{x}) = \tilde{x}$ implies $\sigma = \id$. The case of $\lambda = \lambda^{H}$ follows the same way once we have taken out the points in $\bb^{'}_{P}$.
\end{proof}
\subsection{Proof of integrality for \texorpdfstring{$\PGL_n$}{TEXT} and Langlands duality} \label{pgln_subsect}
In this subsection we finish the proof of cohomological integrality for $\PGL_n$ with $n$ prime. After establishing cohomological integrality for $\PGL_n$ we prove Langlands duality of the cohomology of DT sheaves for $\SL_n$ and $\PGL_n$. \par
We have a short exact sequence
\begin{equation} \label{ses_levi_new_pgl}
    1 \to \mu_n \to L_{\SL_n, \lambda} \to L_{\PGL_n , \lambda} \to 1 
\end{equation}
Note that the results of Proposition \ref{center_Wact} still apply to $L_{\PGL_n , \lambda}$  since $\mu_n$ is central in $L_{\SL_n}$ and via the short exact sequence \eqref{ses_levi_new_pgl} we can write $L_{\PGL_n ,\lambda} \cong L_{\SL_n , \lambda} / \mu_n$. Therefore, we can say that $\LocB_{L_{\SL_n, \lambda}}$ is a $\LocB_{\mu_n}$-torsor over $\LocB^{1}_{L_{\PGL_n, \lambda}}$. Note however, that $\mu_n$ is \emph{not} the centre of $L_{\SL_n, \lambda}$. 
Consider the diagram
\begin{equation}
\begin{tikzcd}
	{\Loc_{\SL_n}} & {\Loc^{1}_{\PGL_n}} \\
	{X_{\SL_n}} & {X_{\PGL_n}}
	\arrow["{\widetilde{\eta}}", from=1-1, to=1-2]
	\arrow["{\pi_{\SL_n}}"', from=1-1, to=2-1]
	\arrow["{\pi_{\PGL_n}}", from=1-2, to=2-2]
	\arrow["\eta"', from=2-1, to=2-2]
\end{tikzcd}
\end{equation}
Note that we showed before in Proposition \ref{center_Wact} that the map $\widetilde{\eta}$ is a $\Loc_{\ZZ(\SL_n)}$-torsor and in particular $\varphi_{\SL_n}$ has an equivariant structure with respect to this action. Also the $\Loc_{\ZZ(\SL_n)}$ actions on $\Loc^{1}_{\PGL_n}$ and $X_{\PGL_n}$ are trivial, so it makes sense to take invariants. 
\begin{remark} \label{bmu_n_invariants}
    Note that $\Loc_{\ZZ(\SL_n)} \cong \mu^{3}_{n} \times \B \mu_{n}$ so taking invariants with respect to $\Loc_{\ZZ(\SL_n)}$ is the same as taking invariants with respect to $\mu^{3}_{n}$. This follows because $\B \mu_{n}$ has no cohomology and so an equivariant structure on a sheaf with respect to $\B \mu_n$ is always trivial. In other words the category of $\B \mu_n$ equivariant sheaves is just isomorphic to the original category of sheaves.  Therefore, taking $\B \mu_n$ invariants has no effect. In particular, we will be able to say that taking invariants of sheaves and pushing forward commutes since we are taking invariants with respect to the finite group $\mu^{3}_n$.
\end{remark}
We therefore have
\begin{equation} \label{torsor_eq}
    (\widetilde{\eta}_{*} \varphi_{\SL_n})^{\Loc_{\ZZ(\SL_n)}} \cong \varphi^{1}_{\PGL_n}.
\end{equation}
Here $(-)^{\Loc_{\ZZ(\SL_n)}}$ is invariants with respect to the group $\Loc_{\ZZ(\SL_n)}$.
\begin{align} \label{push_invariants_comp}
    \pi_{\PGL_n *} \varphi^{1}_{\PGL_n} & \cong \pi_{\PGL_n *}(\widetilde{\eta}_{*} \varphi_{\SL_n})^{\Loc_{\ZZ(\SL_n)}}   \\
    & \cong (\pi_{\PGL_n *} \widetilde{\eta}_* \varphi_{\SL_n})^{\Loc_{\ZZ(\SL_n)}} \qquad (\text{since invariants and pushforward commutes}) \nonumber \\
    & \cong (\eta_{*}\pi_{\SL_n *} \varphi_{\SL_n})^{\Loc_{\ZZ(\SL_n)}}. \nonumber 
\end{align}
We also have the  analogous equation for Levis obtained by the same argument 
\begin{equation} \label{invariants_levi_pgln}
    \pi_{L_{\PGL_n , \lambda} *} \varphi^{1}_{L_{\PGL_n , \lambda}} \cong (\eta_{*}\pi_{L_{\SL_n , \lambda} *} \varphi_{L_{\SL_n , \lambda}})^{\Loc_{\mu_n}}.
\end{equation}
Now by Lemma \ref{pgln_finite_map}  the map $\eta$ is finite and we can compute the invariants $(\eta_{*}\pi_{\SL_n *} \varphi_{\SL_n})^{\Loc_{\ZZ(\SL_n)}}$ using the decomposition
in Proposition \ref{supports_sln}. Note that we can write
\begin{equation}
    \IC_{X^{\lambda}_{\SL_n}}(\mathcal{L}) \cong \IC_{U}(\mathcal{L}|_{U})
\end{equation}
for some open dense $U \subseteq X^{\lambda}_{\SL_n}$.
In particular, by part $1$ of Lemma \ref{pgln_strat_etale_lem}, which compares the $\SL_n$ and $\PGL_n$ strata under the map $\eta$, we can replace $X^{\lambda^{H}}_{\SL_n}$ with the open dense $X^{\lambda^{H}}_{\SL_n} \setminus \bb_{S}$.  Then, using the finiteness of $\eta$ we can show that the supports of the pushforward $\eta_{*} \pi_{\SL_n *} \varphi_{\SL_n}$ must be of the form $X^{\lambda}_{\PGL_n}$. 
Using this we can write the Saito decomposition for $\PGL_n$ 
\begin{equation} \label{pgln_saito_decomp}
     \pi_{\PGL_n *} \varphi^{1}_{\PGL_n} = \bigoplus_{i \geq 0} \bigoplus_{\lambda} \IC_{X^{\lambda}_{\PGL_n}}(\mathcal{L}^{\lambda}_{i})[-i].
\end{equation}
To prove integrality we are reduced to computing the local systems $\mathcal{L}^{\lambda}_{i}$. We will then compute the $\mathcal{L}^{\lambda}_{i}$ by using the \'etale covers in part $(2)$ of Lemma \ref{pgln_strat_etale_lem}.
We now prove an analogue of Lemma \ref{tilde_center_pull_lemma} for $\PGL_n$
\begin{lem} \label{tilde_center_pull_lemma_pgln}
 Consider the induced map $\ZZ^{3}(L_{\PGL_n, \lambda}) \subseteq H^{3}_{\PGL_n} \to X_{\PGL_n}$. For $\lambda \neq \lambda^{H}$ we have a diagram where both squares are pullbacks
 \begin{equation} 
 \label{tilde_center_pullbacks_pgln}
\begin{tikzcd}
	{\widetilde{\ZZ}^{3,g}(L_{\PGL_n, \lambda})} & {\widetilde{\ZZ}^{3}(L_{\PGL_n, \lambda})} & {\ZZ^{3}(L_{\PGL_n, \lambda})} \\
	{X^{\lambda,g}_{\PGL_n}} & {X^{\lambda}_{\PGL_n}} & {\overline{X}^{\lambda}_{\PGL_n}}
	\arrow[from=1-1, to=1-2]
	\arrow["{\theta_{\lambda}}"', from=1-1, to=2-1]
	\arrow[from=1-2, to=1-3]
	\arrow["\varpi", from=1-2, to=2-2]
	\arrow["\pi", from=1-3, to=2-3]
	\arrow[from=2-1, to=2-2]
	\arrow[from=2-2, to=2-3]
\end{tikzcd}
    \end{equation}
    with $\varpi$ and $\theta_{\lambda}$ $W_{L_{\lambda}} \colon 1$ covers. For $\lambda = \lambda^{H}$ we have
    a diagram whose squares are pullbacks
    \begin{equation}
\begin{tikzcd}
	{\widetilde{\ZZ}^{3,g}(L_{\PGL_n, \lambda^{H}}) \setminus R^{'}_{P}} & {\widetilde{\ZZ}^{3}(L_{\PGL_n, \lambda}) \setminus R^{'}_{P}} & {\ZZ^{3}(L_{\PGL_n, \lambda^{H}}) } \\
	{X^{\lambda^{H},g}_{\PGL_n}} & {X^{\lambda^{H} }_{\PGL_n}} & {X_{\PGL_n}}
	\arrow[from=1-1, to=1-2]
	\arrow["{\theta_{\lambda}}"', from=1-1, to=2-1]
	\arrow[from=1-2, to=1-3]
	\arrow["\varpi", from=1-2, to=2-2]
	\arrow["\pi", from=1-3, to=2-3]
	\arrow[from=2-1, to=2-2]
	\arrow[from=2-2, to=2-3]
\end{tikzcd}
    \end{equation}
\end{lem}
\begin{proof}
The proof is analogous to the case of $\GL_n$ and $\SL_n$ in Lemma \ref{tilde_center_pull_lemma}. We first need that the map $X_{L_{\PGL_{n} , \lambda}} \to X_{\PGL_{n}}$ is finite. This follows since it is a quotient of the finite map $X_{L_{\SL_{n}, \lambda}} \to X_{\SL_{n}}$. Firstly, we need to show that the map $\ZZ^{3}(L_{\PGL_{n}, \lambda}) \to X_{\PGL_{n}}$ lands in $\overline{X}^{\lambda}_{\PGL_{n}}$. Here as in the $\SL_n$ case, we can compute the closure using the quotient map $H^{3}_{\PGL_n} \to X_{\PGL_n}$. We can consider the closure in the maximal torus $H^{3}_{\SL_n}$ of $\SL_n$ by first pulling back along the quotient by $W$ and then by $\mu^{3}_n$. Here $W$ is the Weyl group of $\PGL_n$, not the relative Weyl group of $L_{\PGL_n , \lambda}$.
The fact that $\varpi$ and $\theta_{\lambda}$ are $W_{L_{\lambda}} \colon 1$ covers follows from Lemma \ref{levi_etale_lemma_pgln}.
\end{proof}
\begin{proof}[Proof of Lemma \ref{loc_sys_comp} for $\PGL_{n}$] \label{loc_sys_comp_pgl}
From the equation \eqref{invariants_levi_pgln} and the fact that $\pi_{L_{\SL_n , \lambda}*} \varphi^{1}_{L_{\SL_n , \lambda}}$ has perverse cohomology bounded below with lowest perverse piece $l-1$ we know that the complex $\pi_{L_{\PGL_n , \lambda}*} \varphi^{1}_{L_{\PGL_n , \lambda}}$ has perverse cohomology bounded below with lowest perverse piece $l-1$. Therefore,  we define 
    \begin{equation}
        \BPS_{L_{\PGL_n , \lambda}} = \pH^{l-1}\pi_{L_{\PGL_n , \lambda}*} \varphi^{1}_{L_{\PGL_n , \lambda}}.
    \end{equation}
We can use that $\operatorname{supp}(\BPS_{L_{\SL_{n}, \lambda}}) = \ZZ^{3}(L_{\SL_{n}, \lambda})$, the fact that $\eta$ is finite and part $2$ of Lemma \ref{pgln_finite_map} to show
\begin{equation}
    \ZZ^{3}(L_{\PGL_{n}, \lambda}) = \eta(\ZZ^{3}(L_{\SL_{n}, \lambda})) = \supp (\BPS_{L_{\PGL_n , \lambda}}). 
\end{equation}
 In particular, we have used part $2$ of Lemma \ref{pgln_finite_map} and the fact that the map $\eta$ restricted to $U$ is \'etale to prove that the lowest perverse piece of $\pi_{L_{\PGL_n , \lambda}*} \varphi^{1}_{L_{\PGL_n , \lambda}}$ pulls back to the lowest perverse piece of $\pi_{L_{\SL_n , \lambda}*} \varphi_{L_{\SL_n , \lambda}}$. This allows us to deduce the second equality.  Knowing the support we can make the same argument as in Corollary \ref{bps_const} and use the transitive action of $\ZZ^{3}(L_{\PGL_{n}, \lambda})$ to show that the BPS sheaf is constant of rank $1$. The same argument will work to show that any perverse piece supported on $\ZZ^{3}(L_{\PGL_n, \lambda})$ is a constant sheaf. Therefore, we can see that the terms in the decomposition of $\pi_{\PGL_n *}\varphi^{1}_{\PGL_n}$ supported on $\ZZ^{3}(L_{\PGL_n , \lambda})$ are 
 \begin{equation}
     \BPS_{L_{\PGL_n , \lambda}} \otimes \HHf(\B \ZZ(L_{\PGL_n , \lambda})) [- \dim \ZZ(L_{\PGL_n , \lambda})].
 \end{equation}
\par 
To prove part $(2)$ we can mimic the proof of the $\GL_n$ or $\SL_n$ version of Lemma \ref{loc_sys_comp} by using Lemma \ref{tilde_center_pull_lemma_pgln}. Then the rest of the computation of the local systems is the same as in the $\SL_n$ or $\GL_n$ case.  
\end{proof}
\begin{remark}
    Note that the stratum $X^{\mu}_{\PGL_n}$ does not appear in the above lemma for $\PGL_n$ since we will not use it for the the proof of integrality.
\end{remark}
\begin{ex}[BPS sheaves for $\PGL_n$]
    From the above Lemma we can calculate the BPS sheaves for $\PGL_n$, which are constant sheaves supported on the image of $\ZZ^{3}(\PGL_n) = \pt \xhookrightarrow{} X_{\PGL_{n}}$. So we get a skyscraper supported on the trivial local system and furthermore the BPS sheaf only contributes in perverse degree $0$ of $\pi_{\PGL_n *} \varphi^{1}_{\PGL_n}$.
\end{ex}
As in the $\SL_n$ and $\GL_n$ case we will compute the restrictions of the local systems $\mathcal{L}^{\lambda}_{i}$ to the generic loci $X^{\lambda,g}_{\PGL_n}$.
 \begin{proof}[Proof of cohomological integrality for $\PGL_n$ ] 
We will compute the local systems in the Saito decomposition \eqref{pgln_saito_decomp} of $\pi_{\PGL_n *} \varphi^{1}_{\PGL_n}$. 
We get the following diagram for $\lambda \neq \lambda^{H}$
 \begin{equation} \label{pgln_mu_ndiag}
\begin{tikzcd}
	{X^{\lambda,g}_{\SL_n}} & {\widetilde{\ZZ}^{3,g}(L_{\SL_n, \lambda})} \\
	{X^{\lambda,g}_{\PGL_n}} & {\widetilde{\ZZ}^{3,g}(L_{\PGL_n, \lambda})}
	\arrow["\eta"', from=1-1, to=2-1]
	\arrow["{\theta_{\lambda2}}", from=1-2, to=1-1]
	\arrow["{\eta^{'}}", from=1-2, to=2-2]
	\arrow["{\theta_{\lambda1}}"', from=2-2, to=2-1]
\end{tikzcd}
 \end{equation}
 Here  by Lemma \ref{levi_etale_lemma_pgln} $\theta_{\lambda1}$ is a $W_{L_{\lambda}} \colon 1$ cover and by Lemma \ref{pgln_strat_etale_lem} $\eta$ is a $\mu^{3}_{n} \colon 1$ cover . Consider a local system $\mathcal{L}^{\lambda}_i$ appearing in equation \eqref{pgln_saito_decomp}, we will first compute the pullback $\theta^{*}_{\lambda1} \mathcal{L}^{\lambda,g}_i$ as a $W_{L_{\lambda}}$-equivariant local system on $\widetilde{\ZZ}^{3,g}(L_{\PGL_n, \lambda})$. Using equation \eqref{torsor_eq} we know that
 \begin{equation} \label{pgl_n_to_sln_locsys}
     \eta^{*} \mathcal{L}^{\lambda,g}_{i} \cong \mathcal{L}^{\lambda,g, \SL_n}_{i}
 \end{equation}
 Where $\mathcal{L}^{\lambda, g, \SL_n}_{i}$ is the corresponding local system in the $\SL_n$ decomposition as in Proposition \ref{supports_sln} and the isomorphism follows by equation \eqref{push_invariants_comp}. We want to show that 
 \begin{equation}
     \theta^{*}_{\lambda 1} \mathcal{L}^{\lambda,g}_i \cong (\BPS_{L_{\PGL_n , \lambda}} \otimes V_{i})|_{\widetilde{\ZZ}^{3,g}(L_{\PGL_n , \lambda})},
 \end{equation}
 where $V_{i}$ is the $i$th cohomological degree piece of $\BPS_{L_{\PGL_n , \lambda}} \otimes \HHf(\B \ZZ(L_{\PGL_n, \lambda}))$ with the natural $W_{L_{\lambda}}$ structure as in Lemma \ref{coh_centers}. This is a version of Proposition \ref{levis_final} for $\PGL_n$. We have
 \begin{align*}
     \theta^{*}_{ \lambda 1} \mathcal{L}^{\lambda,g}_i & \cong (\eta^{'}_{*} \eta^{'*}\theta^{*}_{ \lambda 1} \mathcal{L}^{\lambda,g}_i)^{\mu^{3}_n}  \\
     & \cong (\eta^{'}_{*} \theta^{*}_{ \lambda 2}\eta^{*} \mathcal{L}^{\lambda,g}_i)^{\mu^{3}_n} \qquad \text{(since } \eta^{'*}\theta^{*}_{\lambda 1} \cong \theta^{*}_{\lambda 2}\eta^{*}) \\
     & \cong (\eta^{'}_{*} (\BPS_{L_{\SL_n , \lambda}} \otimes V_{i})|_{\widetilde{\ZZ}^{3,g}(L_{\SL_n , \lambda})})^{\mu^{3}_n} \qquad \text{by Proposition \ref{levis_final} for } \SL_n \text{ and  equation } \ref{pgl_n_to_sln_locsys}  \\ 
     & \cong (\BPS_{L_{\PGL_n , \lambda}} \otimes V_{i})|_{\widetilde{\ZZ}^{3,g}(L_{\PGL_n , \lambda})} \qquad  \text{ by Lemma \ref{loc_sys_comp} for  $\PGL_n$}. 
 \end{align*}
 First, recall that we can ignore any $\B \mu_n$ factors by Remark \ref{bmu_n_invariants}. In more detail, pulling back by $\eta$ we get some local system on $X^{\lambda,g}_{\SL_n}$. By the statement of Proposition \ref{levis_final} for $\SL_n$ we know that we get BPS sheaves of $L_{\SL_n , \lambda}$ when we further pullback by $\theta_{\lambda 2}$. Pushing forward by $\eta^{'}$ and taking invariants we will then get the sheaf $(\BPS_{L_{\PGL_n}, \lambda} \otimes V_i)|_{\widetilde{\ZZ}^{3,g}(L_{\PGL_n , \lambda})}$ as in the proof of Lemma \ref{loc_sys_comp} for $\PGL_n$. By Lemma \ref{action_lemma_split_off} the action on $\BPS_{L_{\SL_n , \lambda}} \otimes V_{i}$ is the natural one of $W_{L_{\lambda}}$ on equivariant cohomology. Since the map $\eta^{'}$ is $W_{L_{\lambda}}$-equivariant and the $\mu^{3}_{n}$ and $W_{L_{\lambda}}$ actions on $\widetilde{\ZZ}^{3,g}(L_{\SL_n , \lambda})$ commute we also get the natural $W_{L_{\lambda}}$-equivariant structure on $\BPS_{L_{\PGL_n}, \lambda} \otimes V_i$ with $W_{L_{\lambda}}$ acting as in Lemma \ref{coh_centers}. Now we can conclude that 
 \begin{equation} \label{compare}
     \mathcal{L}^{\lambda,g}_{i} \cong (\theta_{\lambda 1*} \theta^{*}_{\lambda 1} \mathcal{L}^{\lambda,g}_{i})^{W_{L_{\lambda}}} \cong (\theta_{ \lambda 1*} (\BPS_{L_{\PGL_n, \lambda}} \otimes V_{i})|_{\widetilde{\ZZ}^{3,g}(L_{\PGL_n , \lambda})} )^{W_{L_{\lambda}}}.
 \end{equation}
  By comparing equations \eqref{llambda_comp_bpsforlevislem} and \eqref{compare} we can finish the proof of integrality for  $\PGL_n$ and $\lambda \neq \lambda^{H}$. The argument for $\lambda= \lambda^{H}$ will be similar. The difference is that  we will have to remove the bad locus $\bb$ to get \'etale covers. In particular, we will now have the diagram 
  \begin{equation}
\begin{tikzcd}
	{X^{\lambda_{H},g}_{\SL_n} \setminus \bb_{S}} & {\widetilde{\ZZ}^{3,g}(L_{\SL_n, \lambda_{H}}) \setminus \bb} \\
	{X^{\lambda_{H},g}_{\PGL_n} \setminus \bb_{P}} & {\widetilde{\ZZ}^{3,g}(L_{\PGL_n, \lambda_{H}}) \setminus \bb^{'}}
	\arrow["\eta"', from=1-1, to=2-1]
	\arrow["{\theta_{\lambda^{H}_{2}}}", from=1-2, to=1-1]
	\arrow["{\eta^{'}}", from=1-2, to=2-2]
	\arrow["{\theta_{\lambda^{H}_{1}}}", from=2-2, to=2-1]
\end{tikzcd}
  \end{equation}
we can now repeat the same argument as above.
 \end{proof}
 We can also incorporate the contributions of the other components of $\Loc_{\PGL_n}$ when $n$ is prime. By Lemma \ref{prime_pgln} the twisted stacks in this case are just $\B \mu^{2}_{n}$ so the contribution on the level of good moduli spaces is just a constant sheaf on a point.
 \begin{proof}[Proof of Corollary \ref{langlands}]
    Once we have established Theorem \ref{coh_proof_pgsln}, we can pushforward all the way to the point and compute cohomology of the DT sheaves. More precisely, we can use the same argument as in the proof of \cite[Theorem 1]{mhm_sym} and use the reducibility of representations of $W_{L_{\lambda}}$ and the compatibility of the pusforward functor and composition of the maps $X_{L_{G,\lambda}} \to X_{G} \to \pt $ to get the decomposition 
    \begin{equation}
        \HHf( \Loc^{1}_{G} , \varphi^{1}_{G} )  \cong \bigoplus_{L_{G, \lambda} \subseteq G} ( \BPSo_{L_{G, \lambda}} 
 \otimes \HHf(\B \ZZ(L_{G, \lambda}))[- \dim \ZZ(L_{G, \lambda})])^{W_{L_{\lambda}}}
    \end{equation}
Here $\BPSo_{L_{G, \lambda}} = \HHf(X_{L} , \BPS_{L_{G, \lambda}})$ and $G = \SL_n$ or $\PGL_n$.
    Note that by Lemma \ref{levi_centers} for $n$ prime we have an isomorphism $\ZZ(L_{\SL_{n}, \lambda}) \cong \ZZ(L_{\PGL_{n},\lambda})$ so $\BPSo_{L_{\SL_{n}, \lambda}} \cong \BPSo_{L_{\PGL_{n}},\lambda}$. There is also an isomorphism of relative Weyl groups. Therefore, we have
\begin{align*}
\HHf(\Loc_{\SL_{n}} , \varphi_{\SL_{n}})  & \cong \HHf(\ZZ^{3}(\SL_{n}), \mathbb{Q}_{\ZZ^{3}(\SL_{n})}) \oplus     
\\
& \bigoplus_{\substack{L_{\SL_{n},\lambda} \subseteq \SL_{n} \\ L_{\SL_{n},\lambda} \neq \SL_{n}}} ( \BPSo_{L_{\SL_{n},\lambda}} 
 \otimes \HHf(\B \ZZ(L_{\SL_{n},\lambda}))[- \dim \ZZ(L_{\SL_n,\lambda})])^{W_{L_{\lambda}}}
\end{align*}
and
\begin{align*}
\HHf(\Loc_{\PGL_{n}} , \varphi_{\PGL_{n}})  & \cong \HHf(\ZZ^{3}(\PGL_{n}), \mathbb{Q}_{\ZZ^{3}(\PGL_{n})})  \oplus \HHf(\coprod_{\zeta \neq 1} \Loc^{\zeta}_{\PGL_{n}}, \mathbb{Q}) \oplus &   \\ 
& \bigoplus_{\substack{L_{\PGL_{n},\lambda} \subseteq \PGL_{n} \\ L_{\PGL_{n},\lambda} \neq \PGL_{n}}} ( \BPSo_{L_{\PGL_n,\lambda}}
 \otimes \HHf(\B \ZZ(L_{\PGL_{n},\lambda}))[- \dim \ZZ(L_{\PGL_n,\lambda})] )^{W_{L_{\lambda}}}. 
 \end{align*}
Since the centers and Weyl groups are identified, the contributions of the Levis are also identified. Now the contribution of the center of $\SL_{n}$ is  $\mathbb{Q}^{|\mu_{n}|^3}$. By Lemma \ref{prime_pgln} $\Loc_{\PGL_{n}}$ has $|\mu_{n}|^{3}- 1$ twisted components, which on the level of cohomology are just points so they contribute $\mathbb{Q}^{|\mu_{n}|^{3} - 1}$. Finally, the center of $\PGL_{n}$ is trivial so it contributes $\mathbb{Q}$ and we have shown the claim.
\end{proof}

\printbibliography
\end{document}